\documentclass[12pt,amsfonts]{article}

\usepackage{latexsym}
\usepackage{amsfonts} 
\usepackage{amssymb}
\usepackage{graphicx}

\newtheorem{theorem}{Theorem}[section]
\newtheorem{lemma}[theorem]{Lemma}
\newtheorem{corollary}[theorem]{Corollary}
\newtheorem{proposition}[theorem]{Proposition}

\newtheorem{definition}[theorem]{Definition}
\newcommand{\bd}[1]{\begin{definition}\label{#1}\rm}
\newcommand{\ed}{\end{definition}}
\newcommand{\bt}[1]{\begin{theorem}\label{#1}}
\newcommand{\et}{\end{theorem}}
\newcommand{\bprop}[1]{\begin{proposition}\label{#1}}
\newcommand{\eprop}{\end{proposition}}
\newcommand{\bcor}[1]{\begin{corollary}\label{#1}}
\newcommand{\ecor}{\end{corollary}}

\newcommand{\D}{\displaystyle}
\newcommand{\T}{\textstyle}
\newcommand{\lra}{\longrightarrow}

\newcommand{\stack}[2]{\raisebox{-2pt} 
{\renewcommand{\arraystretch}{.01} 
\begin{tabular}{c} 
$#2$\\$\scriptscriptstyle #1$ 
\end{tabular} 
}}

\newcommand{\vp}{\varphi}
\newcommand{\ve}{\varepsilon}
\newcommand{\nid}{\noindent}
\newcommand{\qed}{\hfill$\Box$} 

\def\1{\, {\rm I}\mskip-10mu 1} 
\,

\renewcommand{\t}[1]{\tilde{#1}} 
\newcommand{\bmu}{\mbox{\boldmath${\mu}$}}
\newcommand{\bnu}{\mbox{\boldmath${\nu}$}} 

\newcommand{\sbnu}{\mbox{\tiny\boldmath${\nu}$}} 
\newcommand{\sbmu}{\mbox{\tiny\boldmath${\mu}$}} 

\begin{document}
\title{{Absolute Continuity under Time Shift 
of Trajectories and Related Stochastic Calculus}}  
\par
\author{J\"org-Uwe L\"obus
\\ Matematiska institutionen \\ 
Link\"opings universitet \\ 
SE-581 83 Link\"oping \\ 
Sverige 
}
\date{}
\maketitle
{\footnotesize
\noindent
\begin{quote}
{\bf Abstract}
The paper is concerned with a class of two-sided stochastic processes 
of the form $X=W+A$. Here $W$ is a two-sided Brownian motion with 
random initial data at time zero and $A\equiv A(W)$ is a function 
of $W$. Elements of the related stochastic calculus are introduced. 
In particular, the calculus is adjusted to the case when $A$ is a jump 
process. Absolute continuity of $(X,P_{\sbnu})$ under time shift of 
trajectories is investigated. For example under various conditions on 
the initial density with respect to the Lebesgue measure, $m=d\bnu/dx$, 
and on $A$ with $A_0=0$ we verify  
\begin{eqnarray*}
\frac{P_{\sbnu}(dX_{\cdot -t})}{P_{\sbnu}(dX_\cdot)}=\frac{m(X_{-t})} 
{m(X_0)}\cdot\prod_i\left|\nabla_{W_0}X_{-t}\right|_i
\end{eqnarray*} 
a.e. where the product is taken over all coordinates. Here $\sum_i\left( 
\nabla_{W_0}X_{-t}\right)_i$ is the divergence of $X_{-t}$ with respect 
to the initial position. Crucial for this is the {\it temporal homogeneity} 
in the sense that $X\left(W_{\cdot +v}+A_v\1\right)=X_{\cdot+v}(W)$, $v\in 
{\Bbb R}$, where $A_v\1$ is the trajectory taking the constant value $A_v 
(W)$. 

By means of such a density, partial integration relative to the generator 
of the process $X$ is established. Relative compactness of sequences of 
such processes is established. 
\noindent 


\end{quote}
}

\section{Introduction, Basic Objects, and Main Result}
\setcounter{equation}{0}

Let us begin with Haar functions over $[0,1]$. Next we do the same on 
every interval $[a,a+1]$ for all integers $a$ and consider these 
functions being defined on the entire real axis by extending with zero 
outside of $[a,a+1]$. We end up with a system of functions $H_i$, $i\in 
{\Bbb N}$, for which we can consider all finite linear combinations 
\begin{eqnarray*} 
\sum_i\xi_i\int_0^t H_i(s)\, ds\, ,\quad t\in {\Bbb R}, 
\end{eqnarray*}
where $\xi_i$, $i\in {\Bbb N}$, is a sequence of independent $N(0,1)$
distributed random variables. The L\'evy-Ciesielsky representation of 
standard Brownian motion tells us now that under proper selection of 
these finite linear combinations we obtain in the limit a two-sided 
Brownian motion $W$ with $W_0=0$. The convergence is uniform on every 
compact subset of ${\Bbb R}$, almost surely. 
\medskip 

{\bf The initial example. } Denote by $\vp$ the density function with 
respect to the Lebesgue measure of an $N(0,1)$ random variable. Let us 
consider two individual trajectories $W^1$ and $W^2$ obtained by the 
outcomes $\xi_i=x_i^1$ and, respectively, $\xi_i=x_i^2$, $i\in {\Bbb N}$. 
The weight of $W^2$ relative to $W^1$ is 
$ 
\prod_{i=1}^\infty\frac{\vp(x_i^2)}{\vp(x_i^1)}
$
%
provided that this infinite product converges properly. If we assume that 
$W_0$ is a random variable independent of $\xi_i$, $i\in {\Bbb N}$, with 
positive density $m$ then the weight of $W^2$ relative to $W^1$ is 
\begin{eqnarray*} 
\frac{m(W^2_0)}{m(W^1_0)}\cdot\prod_{i=1}^\infty\frac{\vp(x_i^2)} 
{\vp(x_i^1)}\, . 
\end{eqnarray*}
Now let $W^2=W^1_{\cdot -1}$. We obtain by reordering of indices that 
the weight of $W^2=W^1_{\cdot -1}$ relative to $W^1$ is 
\begin{eqnarray*} 
\frac{m(W^2_0)}{m(W^1_0)}\cdot\prod_{i=1}^\infty\frac{\vp(x_i^2)}{\vp( 
x_i^1)}=\frac{m(W^1_{-1})}{m(W^1_0)}\, . 
\end{eqnarray*}
Since the above ``reordering" followed by canceling is just a heuristic 
argument, all what follows from this relation here in the motivating 
part, is of heuristic nature as well. 
\medskip

The question the paper addresses is what happens if we have a process 
of type $X=W+A(W)$ for which we would be interested in such ``relative 
weights" under time shift. We suppose that $X$ is an injective function 
of $W$. 
Let us take a look at the very particular case when (i) $A_0(W)=0$ for 
all trajectories $W$, (ii) the process $A(W)$ is independent of $W_0$, 
and (iii) we have $X_{\cdot -1}(W)\equiv W_{\cdot -1}+A_{\cdot -1}(W)=X 
(W_{\cdot -1}+A_{-1}\1)$ where $A_{-1}\1$ is the function taking the 
value $A_{-1}$, constant in time. Then it is still plausible that the 
relative weight of $X_{\cdot -1}$ with respect to $X$ is determined by 
the weight of $W_{\cdot -1}+A_{-1}\1$ relative to $W$, which is by the 
initial example and (ii)  
\begin{eqnarray*} 
\frac{m(X_{-1})}{m(X_0)}=\frac{m(W_{-1}+A_{-1})}{m(W_0)}\, . 
\end{eqnarray*}

Indeed, given $W_{\cdot -1}+A_{-1}(W)\1$, according to (iii), we add 
$A_{\cdot-1}(W)-A_{-1}(W)\1=A(W_{\cdot -1}+A_{-1}\1)$ to obtain $X_{\cdot 
-1}(W)=X(W_{\cdot -1}+A_{-1}\1)$. On the other hand, given $W$ we add $A 
\equiv A(W)$ to obtain $X\equiv X(W)$. 
%
We may interpret (iii) as some sort of temporal homogeneity of the in 
general non-adapted process $X$. 
\medskip 

Below we are using a more general concept of {\it temporal homogeneity} 
of $X$. Also, conditions (i) and (ii) will be dropped. However in order 
to develop a framework that finally results in a Radon-Nikodym density 
formula for the measure of the trajectories of $X$ under time shift we 
need to formulate a set of conditions. On the one hand, these conditions 
must allow to apply methods from stochastic / Mallianvin calculus and, 
on the other hand, cover the class of stochastic processes and particle 
systems we want to discover. Next, we will give verbal answers to the 
following two questions. What sort of stochastic processes and particle 
systems are we interested in? What is the major technical challenge we 
have to overcome to prove the density formula. 
\medskip 

{\bf Switching media. } Let us continue the concept of the relative weight 
of $X_{\cdot -1}$ with respect to $X$ on an intuitive level. We are now 
interested in a $d$-dimensional process of the form $X=W+A$ for which (i) 
$A_0(W)=0$ for all trajectories $W$, (ii') $A_s(W)\equiv\left(A^i_s(W_1, 
\ldots ,W_d)\right)_{i=1,\ldots d}$ is for all $s\in {\Bbb R}$ continuously 
differentiable with respect to the arguments $W_1,\ldots ,W_d$, and (iii) 
we have $X_{\cdot -1}(W)\equiv W_{\cdot -1}+A_{\cdot -1}(W)=X(W_{\cdot -1} 
+A_{-1}\1)$. Let $\bnu$ denote the measure whose density with respect to 
the $d$-dimensional Lebesgue measure is $m$, the density of $W_0$. The $d 
$-dimensional volume element is in this sense compatibly denoted by $dW_0 
$. Since we do no longer suppose that the process $A(W)$ is independent of 
$W_0$, it is reasonable that we turn to the Radon-Nikodym derivative ${d 
\bnu(W_{-1}+A_{-1})}/{d\bnu(W_0)}$ rather than to look at some ratio of 
densities. We obtain formally for the relative weight of $W_{\cdot -1}+ 
A_{-1}\1$ with respect to $W$ 
\begin{eqnarray*} 
\frac{d\bnu(W_{-1}+A_{-1})}{d\bnu(W_0)} 
&&\hspace{-.5cm}=\left.m(W_{-1}+A_{-1})\frac{d\left(W_{-1}+A_{-1}(W) 
\vphantom{l^1}\right)}{dW_0}\right/m(W_0) \\ 
&&\hspace{-.5cm}=\frac{m(W_{-1}+A_{-1})}{m(W_0)}\prod_{i=1}^d\left|1 
+\nabla_{W_0}A_{-1}(W)\vphantom{l^1}\right|_i \\ 
&&\hspace{-.5cm}=\frac{m(X_{-1})}{m(X_0)}\prod_{i=1}^d\left|\nabla_{ 
W_0}X_{-1}(W)\vphantom{l^1}\right|_i\, . 
\end{eqnarray*}
Condition (iii) gives rise to consider this term also as the 
relative weight of $X_{\cdot -1}$ with respect to $X$. A 
motivating example for (iii) is the following. Let $i_{1,+}\equiv 
\chi_{{\Bbb R}^d_{1,+}}$ denote the indicator function for the set 
${\Bbb R}^d_{1,+}:=\left\{(x_1,\ldots ,x_d):x_1\ge 0\, ,\ x_2, 
\ldots ,x_n\in {\Bbb R}\right\}$, let $i^\ve_{1,+}$ be some smoothly 
mollified modification, and let $a\in {\Bbb R}$. Let $X$ be the 
solution to the It\^o SDE 
\begin{eqnarray*} 
X_s=W_s+a\int_0^si^\ve_{1,+}(X_r)\, d(W_{1,r},0,\ldots ,0)=:W_s+A_s 
\, ,\quad s\in {\Bbb R}. 
\end{eqnarray*}
We have (i) by definition and (ii'), because of the above 
mollification, by a classical theorem of Blagoveshchenskij and 
Freidlin, cf. \cite{PS98}, Theorem 1.23. For $s\in {\Bbb R}$ it 
follows that with $W^s=W_{\cdot +s}+A_s\1$ 
\begin{eqnarray*} 
X_\cdot\left(W^s\right)&&\hspace{-.5cm}=W^s_\cdot+a\int_0^\cdot 
i^\ve_{1,+}\left(X_r\left(W^s\right)\right)\, d\left(W^s_{1,r},0, 
\ldots ,0\right) \\ 
&&\hspace{-.5cm}=W^s_\cdot+a\cdot\int_0^\cdot i^\ve_{1,+}\left(X_r 
\left(W^s\right)\right)\, d(W_{1,s+r},0,\ldots ,0) 
\end{eqnarray*}
and 
\begin{eqnarray*} 
X_{\cdot +s}(W)&&\hspace{-.5cm}=W_{\cdot +s}+a\int_0^{\cdot +s} 
i^\ve_{1,+}\left(X_r(W)\right)\, d(W_{1,r},0,\ldots ,0) \\ 
&&\hspace{-.5cm}=W_{\cdot +s}+A_s+a\int_s^{\cdot +s}i^\ve_{1,+} 
\left(X_r(W)\right)\, d(W_{1,r},0,\ldots ,0) \\
&&\hspace{-.5cm}=W^s_\cdot+a\int_0^\cdot i^\ve_{1,+}(X_{r+s}(W)) 
\, d(W_{1,r+s},0,\ldots ,0) 
\end{eqnarray*}
which implies $X_\cdot\left(W^s\right)=X_{\cdot +s}(W)$, $s\in { 
\Bbb R}$. In other words, we have (iii). The choice $a\in (-1,0)$ 
(respectively $a\in (0,\infty)$) lets $X$ move slower (respectively 
faster) than $W$ on $\left\{(x_1,\ldots ,x_d):\right.$ $\left.x_1 
\ge\delta\, ,\ x_2,\ldots ,x_n\in {\Bbb R}\right\}$. Here $\delta 
>0$ is some small parameter depending on the above mollification. 
\medskip 

{\bf The particle collision type jumps. } The stochastic processes 
and particle systems we are focusing on may involve jumps. Definition  
\ref{Definition1.7} formulates the mathematical conditions on the 
jumps. In this paragraph we give a physical motivation for Definition 
\ref{Definition1.7}. 

Let us consider a system of $n$ marked particles $\{X_{o1},\ldots , 
X_{on}\}$ in some $d_o$-dimensional domain $D_o$. If we model the 
particles just as points in $D_o$ and consequently the collection of 
all $n$ particles as a point in ${\Bbb R}^{n\cdot d_o}$ then the set 
of all potential collisions is the set $Z_o$ defined as follows. For 
$z_1,\ldots ,z_n\in D_o$, $z=(z_1,\ldots ,z_n)$, and $i\in\{1,\ldots 
,n\}$ introduce $z^{ij}:=(z_1,\ldots ,z_{i-1},z_j,z_{i+1},\ldots ,z_n 
)$, $i,j\in \{1,\ldots ,n\}$, $i<j$. Now, let 
\begin{eqnarray*}
Z_o(i,j):=\{z^{ij}:\, z\in D_o^n\}\, ,\quad i<j\, ,\quad Z_o:= 
\bigcup_{i<j}Z_o(i,j)\, . 
\end{eqnarray*}
If we look at particles as small $d_o$-dimensional objects then a 
counterpart is $Z_o^\ve(i,j):=\left\{z\in {\Bbb R}^{n\cdot d_o}:| 
z-z'|\le\ve,\, z'\in Z_o(i,j)\right\}$, $i<j$, and $Z^\ve_o:= 
\bigcup_{i<j}Z^\ve_o(i,j)$. In this case however, we choose from 
the {\it physical particle} some canonical inner point, the {\it 
mathematical particle}; for example the midpoint if the physical 
particle is a small $d_o$-dimensional ball. 

Now we extend the notion of a particle to an object having a $d_o 
$-dimensional geometric location and a $d_v$-dimensional velocity, 
a vector belonging to the difference $D_v$ of the $d_v$-dimensional 
open ball with center 0 and some radius $v_{\rm max}<\infty$ and 
the closed ball with center 0 and some radius $v_{\rm min}<v_{\rm 
max}$. Usually we have $d_o=d_v$. The system of mathematical 
particles we are now interested in is $\{X_1,\ldots ,X_n\}$ where 
$X_i:=(X_{oi},X_{vi})$ and $X_{vi}$ is the velocity of $X_{oi}$. 
We set $d:=d_o+d_v$ and $D:=D_o\times D_v$. In this sense, the set 
of all potential collisions of particle $i$ with particle $j$, 
$i<j$, is now $Z(i,j):=Z_o(i,j)\times D_v$ or, if the physical 
location of the particle is a small $d_o$-dimensional set, $Z^\ve 
(i,j):=Z_o^\ve (i,j)\times D_v$. Respectively, we set $Z:=\bigcup_{ 
i<j}Z(i,j)$ and $Z^\ve:=\bigcup_{i<j}Z^\ve(i,j)$. 

Let $\xi(s)$, $s\in {\Bbb R}$, be an ${\Bbb R}^{n\cdot d}$-valued 
trajectory describing the path of a system of $n$ mathematical 
particles. Furthermore, for $x\in {\Bbb R}^{n\cdot d}$, write 
$x\1\equiv x\1(s)$ for the trajectory which takes the value $x$ 
for all $s\in {\Bbb R}$. By the concept physical and mathematical 
particles, $\xi_s$, $s\in {\Bbb R}$, never reaches the boundary 
$\partial D^n$. Near the boundary $\partial D^n$ we assume that 
the trajectory $\xi$ models some reflection / adhesion / diffusion 
mechanism of the particle(s) that is (are) closest to $\partial D$. 

For small $z\in {\Bbb R}^{n\cdot d}$, let us also consider the 
{\it parallel trajectory} $\xi_z(s)$, $s\in {\Bbb R}$, whenever it 
meaningfully exists. That is $\xi_z=\xi+z\1$, as long as $\xi+z\1$ 
has not got {\it close} to $\partial D^n$ and is again $\xi+z\1$ as 
soon as it is not close to $\partial D^n$ anymore. The latter could 
also mean ``as soon as two individual mathematical particles with 
the same location perform a collision rather than following some 
intrinsic mechanism near the boundary". 

Until the end of this paragraph we would like to keep the 
interpretation that for $z=(z_1,\ldots ,z_n)$ and $z_i=(z_{oi},z_{ 
vi})$, the component $z_{oi}\in {\Bbb R}^{d_o}$ stands for location 
and the component $z_{vi}\in {\Bbb R}^{d_v}$ stands for velocity of 
$z_i$. We make the following observation. It is reasonable to assume 
that for the trajectory $\xi$ and the parallel trajectory $\xi+z\1$ 
we have $z_{vi}=0$, $i\in\{1,\ldots ,n\}$. However we would like to 
have $\xi_x$ be well-defined for $x$ in some neighborhood of $x=0$ 
and set therefore artificially $\xi_x:=\xi_{(x_o,0)}$ where $x_o:= 
(x_{1o},\ldots ,x_{no})\in {\Bbb R}^{n\cdot d_o}$ is the vector of 
the locations of the $x_i\in {\Bbb R}^d$, $i\in\{1,\ldots ,n\}$.  

Let $\tau(x)$ denote the first time after 0 for $\xi_x$ two particles 
collide, that is, the first time after 0 the trajectory $\xi_x$ hits 
$Z^\ve$. The function $x\to\tau(x)$ is then a piecewise infinitely 
boundedly differentiable function in some neighborhood of $x=0$. The 
term {\it piecewise} refers to the different possibilities for the 
first collision to be caused by the particles $i$ and $j$, $i<j$. On 
the pieces, $\nabla_x\tau(x)\in{\Bbb R}^{n\cdot d}$ is by the above 
convention well-defined and we have $\langle\nabla_x\tau(0),e 
\rangle_{{\Bbb R}^{n\cdot d}}=0$ for any vector $e=((e_{1o},e_{1v}), 
\ldots ,(e_{no},e_{nv}))$ with $e_{1o}=\ldots =e_{no}=0$ where $e_{ 
io}$, $i\in\{1,\ldots ,n\}$, are the locations of the $e_i\in {\Bbb 
R}^d$, $i\in\{1,\ldots ,n\}$. 

By the trajectory $\xi$ we follow $n$ particles. For those we may 
assume that immediately after the first collision of two, the vector 
of the velocities has changed by a jump while the vector of the 
locations is unchanged. Thus, for $\Delta\xi_\tau:=\xi_{\tau +}- 
\xi_{\tau-}$ we get 
\begin{eqnarray*} 
\langle\nabla_x\tau(0),\Delta\xi_\tau\rangle_{{\Bbb R}^{n\cdot d}}=0 
\, . 
\end{eqnarray*}
This is motivation for Definition \ref{Definition1.7} (jv). 
Condition (jv) of Definition \ref{Definition1.7} is at the same 
time a technical hypothesis. 
\medskip 

{\bf The process $Y$. } Instead of (i) and (iii) formulated here 
in the introductory section, for the actual analysis in the paper 
we will use condition of (3) of Subsection 1.2. This condition 
does no longer require $A_0=0$. It includes a process $Y$ with 
$A_0=Y_0$ such that with $W^s:=W_{\cdot +s}+\left(A_s(W)-Y_s(W) 
\vphantom{l^1}\right)\1$ we have the relation $X_\cdot\left(W^s 
\right)=X_{\cdot +s}(W)$, $s\in {\Bbb R}$. Below we give a reason 
why (i) and (iii) is not enough. 

Let $X$ be a $n\cdot n$-dimensional stochastic particle process 
of the form $X=W+A$ where $A\equiv A(W)$ is a right-continuous 
pure jump process. That is, $A$ is constant in time until it 
jumps and thereafter constant again until it jumps again. The 
jump times and jumps are measurable functions of $W$. Assume that 
the jumps are organized as in the previous paragraph. This means 
in particular that whenever $X$ hits $Z^\ve$ at some jump time it 
performs a jump within $Z^\ve$ and keeps moving afterwards as it 
was a $n\cdot n$-dimensional Brownian motion. For proper modeling 
we must here assume that jump times do not accumulate. However 
the following scenario is possible. The trajectory $W$ passes 
through $Z^\ve$ in a subset of $Z^\ve$ no jump is allowed to 
reach. This happens at some negative time and no further jump 
occurs until time zero. Consequently, $A_0=0$ which is equivalent 
to $X_s=W_s$ in some non-positive time interval $s\in (\sigma,0]$, 
is not possible. In order to construct a injection $W\to X$ we 
must allow $A_0\neq 0$. In that case $X_\cdot\left(W^s\right)=X_{ 
\cdot +s}(W)$, $s\in {\Bbb R}$, for $W^s:=W_{\cdot +s}+A_s\1$ 
may fail.  
\medskip  

{\bf The integral $\int_{s\in I} F(s,W)\, d\dot{W}_s$ } for some 
interval $I$. This type of integral appears unavoidably in the 
beginning of the proof of our main result, namely in the second 
last line of (\ref{3.25}). We handle this integral by projecting 
$W$ via the L\'evy-Ciesielsky representation to a piecewise linear 
path. The critical operations on $\int_{s\in I} F(s,W)\, d\dot{W 
}_s$ will be carried out under the projection. Then we apply 
several extension techniques but mainly the approximation theorem 
Theorem \ref{Theorem2.1} by Gihman and Skorohod, cf. \cite{GS66}. 
The calculus on the projections is developed in Subsections 2.1, 
2.3, 3.1, and 3.2. 
\bigskip

{\bf Organization of the paper. } After this motivating part, we 
collect all symbols and definitions we frequently need in 
the paper. The concept of two-sided Brownian motion is taken from 
Imkeller \cite{Im96}. Symbols and definitions just for momentary 
use will be introduced when they are needed. In Subsection 1.2 we 
formulate a list of conditions we need in order to state the main 
result, the change of measure formula under time shift of 
trajectories. The main results are presented in Subsection 1.3. 
Immediate corollaries are formulated and proved in this subsection. 

Section 2 contains the stochastic calculus on the above mentioned 
projections of trajectories. The results of Subsections 2.1 and 
2.3 are important for the proof of the change of measure formula. 
Subsection 2.2 formulates and proves a related stochastic calculus 
where the time shift is  replaced by a certain class of trajectory 
valued flows. Here we refer to Mayer-Wolf and Zakai \cite{M-WZ05} 
and Smolyanov and Weizs\"acker \cite{SW99} for related ideas. 

If the trajectory flow is a translation along the time axis, we 
shift piecewise linear trajectories, obtained by projections on 
the L\'evy-Ciesielsky representation. This operation requires 
particular attention to objects like gradient, stochastic integral, 
and generator of the flow. The related assertions are proved in 
Subsections 3.1 and 3.2. 

The proof of the change of measure formula involves two 
approximations. The first one corresponds to the projection via the 
L\'evy-Ciesielsky representation. The second one is the 
approximation of the trajectory valued flow $s\to W^s=W_{\cdot +s}+ 
(A_s-Y_s)\1$ if $A-Y$ has jumps by a sequence of continuous flows. 
The change of measure formula is proved in Subsection 3.4. We would 
like to refer to Bogachev and Mayer-Wolf \cite{BM-W99}, Buckdahn 
\cite{Bu94}, Kulik and Pilipenko \cite{KP00}. 

Sections 4 and 5 contain two applications. In Section 4 we prove 
integration by parts formulas for operators which are for, in 
general, non-Markov processes counterparts to generators of Markov 
processes. In Section 5 we arrive at the application the paper 
was initiated by, relative compactness of a class of abstract 
particle systems. This application may explain the setup, conditions, 
and notations in the paper. It has been motivated by concrete 
physical modeling as in  Caprino, Pulvirenti, and Wagner \cite{CPW98} 
as well as a more abstract mathematical and physical treatment as in 
Graham and M\'el\'eard \cite{GM97}. 
%
\bigskip 

Before starting this preliminary part of the paper, let us emphasize 
that most of the terms related to the Malliavin calculus and stochastic 
integration are explained in Section 6. Section 6 is independent and 
should at least be browsed prior to reading this introduction. 

\subsection{Analytical setting}

{\bf The symbol $F$. } Let $n,d\in {\Bbb N}$. In order to simplify the 
notation, we will use the letter $F$ in order to denote both, the space 
of all real $n\cdot d$-dimensional vectors $(z_1,\ldots ,z_n)$ with $z_1, 
\ldots ,z_n\in {\Bbb R}^d$ and the space of all real $n\cdot d\times n 
\cdot d$-matrices. Correspondingly, we will use the symbol $\langle\, 
\cdot\, ,\, \cdot\, \rangle_F$ in order to denote the square of a vector 
norm or of a matrix norm on $F$ which we assume to be compatible but do 
not specify. We also assume $\langle\, \cdot\, ,\, \cdot\, \rangle_F^{1 
/2}$ to be submultiplicative if $F={\Bbb R}^{n\cdot d}\otimes {\Bbb R}^{
n\cdot d}$. 

Let $\lambda_F$ denote the Lebesgue measure on $F$ and let $(e_j)_{j 
=1, \ldots ,n\cdot d}$ be a standard basis in ${\Bbb R}^{n\cdot d}$. 
Moreover, let ${\sf e}=\sum_{i=1}^{n\cdot d}e_j$ if $F={\Bbb R}^{n\cdot 
d}$. Let ${\sf e}$ be the $n\cdot d\times n\cdot d$-dimensional unit 
matrix if $F={\Bbb R}^{n\cdot d}\otimes {\Bbb R}^{n\cdot d}$.  

Furthermore, depending on the format of the entries, let us define 
the product $\langle\, \cdot\, ,\, \cdot\, \rangle_{F\to F}$. If both 
entries belong to $F={\Bbb R}^{n\cdot d}$, we set 
\begin{eqnarray*} 
\langle (x_1,\ldots ,x_{n\cdot d}),(y_1,\ldots ,y_{n\cdot d})\rangle_{F 
\to F}:=\sum_{i=1}^{n\cdot d}x_iy_i\cdot e_i\, , \quad x_1,\ldots ,x_{n 
\cdot d},y_1,\ldots ,y_{n\cdot d}\in {\Bbb R}. 
\end{eqnarray*} 
If at least one of the entries belongs to $F={\Bbb R}^{n\cdot d}\otimes 
{\Bbb R}^{n\cdot d}$ and not otherwise stated, let $\langle\, \cdot\, , 
\, \cdot\, \rangle_{F\to F}$ denote the usual matrix-vector, vector-matrix, 
or matrix-matrix multiplication. 

Set $\1(s):={\sf e}$, $s\in {\Bbb R}$, and for $B\in {\cal B}({\Bbb 
R})$ set $\1_B(s)={\sf e}$ if $s\in B$ and $\1_B(s)=0$ if $s\not\in B$. 
If $x\in F$ then for simplicity, we will use the short notation 
\begin{eqnarray*} 
x\1\equiv\langle x,\1(\, \cdot\, )\rangle_{F\to F}\, .  
\end{eqnarray*} 

It will be always clear from the context which of the spaces, ${\Bbb 
R}^{n\cdot d}$ or ${\Bbb R}^{n\cdot d}\otimes {\Bbb R}^{n\cdot d}$, 
we talk about when we use the symbol $F$ and which form of 
multiplication $\langle\, \cdot\, ,\, \cdot\, \rangle_{F\to F}$ we use. 

\bigskip

{\bf Function spaces and orthogonal projections. } 
Among others, which will be introduced upon usage, we will work with 
the spaces $L^2({\Bbb R};F)$ and $L^2_{\rm loc}({\Bbb R};F)$ of all 
quadratically integrable and, respectively, local quadratically 
integrable $F$-valued functions on ${\Bbb R}$. For simplicity, the 
latter space will be abbreviated just by $L^2$ if no ambiguity is 
possible. For $\int\langle a(s),b(s)\rangle_{F\to F}\, ds$, $a,b 
\in L^2({\Bbb R};F)$, we shall use the symbol $\langle a,b\rangle_{ 
L^2\to F}$. Moreover, for an integral of type $\int_{-\infty}^\infty 
\langle a,db\rangle_F$ or $\int_{-\infty}^\infty\langle a,db\rangle_{F 
\to F}$ we shall use the short notation $\langle a,db\rangle_{L^2}$ or 
$\langle a,db\rangle_{L^2\to F}$, again if no ambiguity is possible. 

In addition, let $C^k({\Bbb R};F)$ be the space of all $k$ times 
differentiable $F$-valued functions on ${\Bbb R}$. When adding a 
subscript zero, $C_0^k({\Bbb R};F)$, we restrict ourself to compactly 
supported elements. Furthermore, let $C_a({\Bbb R};F)$ denote the 
space of all absolutely continuous $F$-valued functions such that, 
for $f\in C_a({\Bbb R};F)$, its Radon-Nikodym derivative admits a 
cadlag version. In fact, the cadlag version is at the same time the 
right derivative of $f\in C_a({\Bbb R};F)$ and will be denoted by 
$f'$. When restricting to functions defined just on some subset 
${\cal S}\subset {\Bbb R}$, we simply replace ${\Bbb R}$ in these 
symbols by ${\cal S}$. 

Let $B_{b;{\rm loc}}({\Bbb R};F)$ denote the space of all Lebesgue 
measurable $F$-valued functions on ${\Bbb R}$ possessing a version 
that is bounded on every finite subinterval of ${\Bbb R}$. 

Finally, let $C^k(E)$ be the space of all $k$ times continuously 
differentiable real functions on $E\subseteq F$. When adding here a 
subscript $b$, $C_b^k(E)$, we just count in the bounded functions of 
$C^k(E)$ with bounded derivatives up to order $k$. 
\medskip 

\medskip 

Let $t>0$. Moreover, let 
\begin{eqnarray*}
D_1(s)&&\hspace{-.5cm}=t^{-\frac12}\, , \quad s\in [0,t), \vphantom{\sum_1} 
\quad D_1(s)=0\, , \quad s\in {\Bbb R}\setminus [0,t), \nonumber \\
D_{2^m+k}(s)&&\hspace{-.5cm}= 
\left\{\begin{array}{rl}
t^{-\frac12}\cdot 2^{\frac{m}{2}} & \mbox{if }s\in [\frac{k-1}{2^m}\cdot t, 
\frac{2k-1}{2^{m+1}}\cdot t)\vphantom{\displaystyle\sum_1} \\
-t^{-\frac12}\cdot 2^{\frac{m}{2}} & \mbox{if }s\in [\frac{2k-1}{2^{m+1}} 
\cdot t,\frac{k}{2^m}\cdot t)\vphantom{\displaystyle\sum_1} \nonumber \\ 
0 & \mbox{if }s\in {\Bbb R}\setminus [\frac{k-1}{2^m}\cdot t,\frac{k}{2^m} 
\cdot t)
\end{array}\right.\, , \nonumber \\ &&\hspace{5cm}\vphantom{\sum^1} k\in \{ 
1,\ldots ,2^m\}\, , \ m\in {\Bbb Z}_+. 
\end{eqnarray*}
Furthermore, 
let $f_j$, $j\in\{1,\ldots ,n\cdot d\}$, be a 
standard basis in ${\Bbb R}^{n\cdot d}$, and define 
\begin{eqnarray*}
E_{n\cdot d\cdot(q-1)+j}:=D_q\cdot f_j\, , \quad q\in {\Bbb N}\, , \ j\in 
\{1,\ldots ,n\cdot d\}, 
\end{eqnarray*}
as well as 
\begin{eqnarray*}
E_k^{(w)}(s):=E_k(s+(w-1)t)\, , \quad s\in {\Bbb R},\ k\in {\Bbb N},\ 
w\in {\Bbb Z}, 
\end{eqnarray*}
and rename the elements of $\{E_k^{(w)}:k\in {\Bbb N},\ w\in {\Bbb Z}\}$ 
as $\{H_1,H_2,\ldots\, \}$. 
\medskip

Let $I(m,r)$ be the set of all indices such that $\{H_i:i\in I(m,r)\}$ is 
the set of all $H_j$ with $H_j(s)=0$, $s\not\in [-rt,(r-1)t]$, and supp 
$H_j$ is a closed interval of the form $[\frac{k-1}{2^{m'}}\cdot t,\frac 
{k}{2^{m'}}\cdot t]$, $0\le m'\le m-1$, $m\in {\Bbb Z}_+$, $r\in {\Bbb N}$. 
Furthermore, set $J(m):=\bigcup_{r\in {\Bbb N}}I(m,r)$, $m\in {\Bbb Z}_+$ 
and $I(r):=\bigcup_{m\in {\Bbb Z}_+}I(m,r)$, $r\in {\Bbb N}$. 
\medskip

Let $p_{m,r}$, resp. $p_{m}$, denote the orthogonal projection from $L^2_{ 
\rm loc}({\Bbb R};F)$ to the linear subspace spanned by $\{H_i:i\in I(m,r) 
\}$, resp. $\{H_i:i\in J(m)\}$. 
\medskip

Let ${\cal R}\equiv {\cal R}(r):=[-rt,(r-1)t]$. For calculations under 
the measure $Q^{(m,r)}_{\sbnu}$ we will use the set 
\begin{eqnarray*}
H\equiv H^{(r)}:=\left\{(g,y):g\in L^2({\cal R};F),\ y\in F\right\}\, . 
\end{eqnarray*}
Let $l_{\cal R}$ denote that real measure on ${\Bbb R}$ which is the 
Lebesgue measure on ${\cal R}$ and zero on ${\Bbb R}\setminus {\cal R}$. 
In addition, identify $H$ with the space $\left\{(g,y):g\in L^2({\Bbb 
R},l_{{\cal R}};F),\ y\in F\right\}$. The inner product in $H$ is 
$\langle h_1,h_2\rangle_H=\langle f_1,f_2\rangle_{L^2}+\langle x_1,x_2 
\rangle_F$ but let us also use the notation 
\begin{eqnarray*}
\, \ \langle h_1,h_2\rangle_{H\to F}&&\hspace{-.5cm}:=\langle f_1,f_2 
\rangle_{L^2\to F}+\langle x_1, x_2\rangle_{F\to F}\, , \quad h_i=(f_i, 
x_i)\in H,\ i\in \{1,2\}. 
\end{eqnarray*} 
%
 
Again, we embed $H$ into $C({\Bbb R};F)$ by $j(f,x):=(\int_0^\cdot f(s)\, 
ds,x)\equiv x\1+\int_0^\cdot f(s)\, ds$. The latter gives also rise to set 
\begin{eqnarray*}
j^{-1}g:=(g',g(0))\, , \quad g\in C_a({\Bbb R};F)\ \mbox{\rm with }g'\in 
L^2({\Bbb R},l_{{\cal R}};F). 
\end{eqnarray*} 
Finally, let $H^{(m,r)}$ resp. $H^{(m)}$ be the linear space spanned by 
$\{(H_i,x):i\in I(m,r),\, x\in F\}$ resp. $\{(H_i,x):i\in J(m),\, x\in F\}$. 
\medskip 

\subsection{Probabilistic setting}

{\bf Two-sided Brownian motion with random initial data. } 
Let $\bnu$ be a probability distribution on $(F,{\cal B}(F))$, here $F= 
{\Bbb R}^{n\cdot d}$. Let the probability space $(\Omega,{\cal F},Q_{\sbnu} 
)$ be given by the following. 
\begin{itemize} 
\item[{(i)}] $\Omega =C({\Bbb R};F)$, the space of all continuous functions 
$W$ from ${\Bbb R}$ to $F$. Identify $C({\Bbb R};F)\equiv \{\omega\in C({ 
\Bbb R};F):\omega(0)=0\}\times F$, here also $F={\Bbb R}^{n\cdot d}$. 
\item[{}] \vspace{-1.1cm}
\item[{(ii)}] ${\cal F}$, the $\sigma$-algebra of Borel sets with respect 
to uniform convergence on compact subsets of ${\Bbb R}$.  
\item[{}] \vspace{-1.1cm}
\item[{(iii)}] $Q_{\sbnu}$, the probability measure on $(\Omega ,{\cal F})$ 
for which, when given $W_0$, both $(W_s-W_0)_{s\ge 0}$ as well as $(W_{-s}-
W_0)_{s\ge 0}$ are independent standard Brownian motions with state space 
$(F,{\cal B}(F))$. Assume $W_0$ to be independent of $(W_s-W_0)_{s\in {\Bbb 
R}}$ and distributed according to $\bnu$. 
\end{itemize} 
In addition, we will assume that the natural filtration $\{{\cal F}_u^v= 
\sigma(W_\alpha-W_\beta :u\le\alpha,\beta\le v)\times \sigma(W_0):-\infty 
<u<v<\infty\}$ is completed by the $Q_{\sbnu}$-completion of ${\cal F}$. 
\medskip 

For the measure $\bnu$ we shall assume the following throughout the paper.  
\begin{itemize}
\item[{(i)}] $\bnu$ admits a density $m$ with respect to $\lambda_F$ and 
$m$ is symmetric with respect to the $n$ $d$-dimensional components of $F$.  
\item[{}] \vspace{-1.1cm}
\item[{(ii)}]  $m$ is supported by a set $\overline{D^n}$ where $D$ is a 
bounded $d$-dimensional domain such that 
\begin{eqnarray*} 
0<m\in C^1(D^n)\quad\mbox{\rm and}\quad\lim_{D^n\ni x\to\partial D^n} 
m(x)=0\, . 
\end{eqnarray*}
\item[{}] \vspace{-1.1cm}
\item[{(iii)}] Furthermore, let us assume that there exist $Q\in (1,\infty)$ 
such that, for all $y\in F$,  
\begin{eqnarray*} 
\left.\frac{d} 
{d\lambda}\right|_{\lambda =0}\frac{m(\cdot +\lambda y)}{m}=\left\langle 
\frac{\nabla m}{m},y\right\rangle_F\quad \mbox{\rm exists in }L^Q(F,\bnu)\, . 
\hspace{-.7cm}
\end{eqnarray*} 
\end{itemize} 
\medskip

{\bf Measure, gradient, and stochastic integral under projections. } 
Recalling the projections on subspaces spanned by sets of Haar functions, 
for $(f,x)\in L^2\times F$, we set $q_{m,r}(f,x):=(p_{m,r}f,x)$ and $q_{m} 
(f,x):=(p_{m}f,x)$. For $\Omega\ni W=y+\sum_{i=1}^\infty x_i\cdot\int_0^\cdot 
H_i(u)\, du$ such that the sum converges uniformly on finite intervals we 
will use the notation 
\begin{eqnarray*} 
q_{m,r}j^{-1}(W):=\left(\sum_{i\in I(m,r)} x_i\cdot H_i,y\right)\quad\mbox{\rm 
and}\quad q_{m}j^{-1}(W):=\left(\sum_{i\in J(m)} x_i\cdot H_i,y\right) 
\end{eqnarray*} 
which is compatible with the operator $j$ defined in the end of Subsection 
1.1 and the just defined projections $q_{m,r}$ and $q_{m}$. The sum $\sum_{i 
\in J(m)} x_i\cdot H_i$ is to be understood in the sense $L^2_{\rm loc}({\Bbb 
R};F)$. We are now able to introduce the projections 
\begin{eqnarray*} 
\pi_{m,r}(W):=jq_{m,r}j^{-1}(W)\quad\mbox{\rm and}\quad\pi_m(W):=jq_mj^{-1} 
(W)\, , \quad W\in\Omega .  
\end{eqnarray*} 
Furthermore, we introduce the measure $Q^{(m,r)}_{\sbnu}$ on $(\Omega,{\cal 
B}(\Omega))$ by 
\begin{eqnarray*}
Q^{(m,r)}_{\sbnu}\left(\displaystyle\Omega\setminus\left\{\pi_{m,r}(W):W\in 
\Omega\right\}\right)=0 
\end{eqnarray*} 
and, on $\left\{\pi_{m,r}(W):W\in\Omega\right\}$, by
\begin{eqnarray*} 
Q^{(m,r)}_{\sbnu}=Q_{\sbnu}\circ \pi_{m,r}^{-1}\, ,\quad m\in {\Bbb Z}_+, 
\ r\in {\Bbb N}.  
\end{eqnarray*} 
Similarly, we define the measure $Q^{(m)}_{\sbnu}$ on $(\Omega,{\cal B} 
(\Omega))$ by 
\begin{eqnarray*}
Q^{(m)}_{\sbnu}\left(\displaystyle\Omega\setminus\left\{\pi_m(W):W\in\Omega 
\right\}\right)=0 
\end{eqnarray*} 
and, on $\left\{\pi_{m}(W):W\in\Omega\right\}$, by 
\begin{eqnarray*} 
Q^{(m)}_{\sbnu}=Q_{\sbnu}\circ\pi_{m}^{-1}\, ,\quad m\in {\Bbb Z}_+ . 
\end{eqnarray*} 
Likewise, but without indicating this in the notation, let us also define 
the measure $Q^{(m,r)}_{\sbnu}$ on Borel subsets $\Omega'\in {\cal B}(\Omega)$ 
with $\left\{\pi_{m,r}(W):W\in\Omega\right\}\subseteq \Omega'$ and the measure 
$Q^{(m)}_{\sbnu}$ on $\Omega'$ whenever $\left\{\pi_m(W):W\in\Omega\right\} 
\subseteq\Omega'$. 
\medskip
\begin{definition}\label{Definition1.1}
{\rm (a) Let $f\in C_b^1({\Bbb R^{\#I(m,r)+n\cdot d}})$, $\langle {\rm e} 
,W\rangle$ be the vector of all $\langle e_j,W_0\rangle_F$, $j\in \{1, 
\ldots,$ $n\cdot d\}$, and $\langle {\rm H},W\rangle$ be the vector of all 
$\langle H_i,dW\rangle_{L^2}$, $i\in I(m,r)$, $W\in C({\Bbb R};F)$. Let 
$\vp(W):=f\left(\langle {\rm H},W\rangle,\langle {\rm e},W\rangle\right)$. 
Define the {\em gradient of a cylindrical function} with respect to the 
measure $Q^{(m,r)}_{\sbnu}$ by 
\begin{eqnarray*}
{\Bbb D}\vp(W)(j(g,x))&&\hspace{-.5cm}:=\sum_{i\in I(m,r)}f_i\left(\langle 
{\rm H},W\rangle ,\langle {\rm e},W\rangle\right)\cdot\langle H_i,g 
\rangle_{L^2} \\ 
&&\hspace{-.0cm}+\sum_{j=1}^{n\cdot d}f_j\left(\langle {\rm H},W\rangle , 
\langle {\rm e},W\rangle\right)\cdot\langle e_j,x \rangle_F
\end{eqnarray*} 
where $f_i$ denotes the first order derivative of $f$ relative to $i\in 
I(m,r)$ and $f_j$ is the first order derivative of $f$ relative to $j\in 
\{1,\ldots ,n\cdot d\}$. 
}
\end{definition}

We remark that Definition \ref{Definition1.1} (a) is in compliance with 
(\ref{6.4}) and (\ref{6.5}) by choosing $l_i:=(-dH_i,0)$, $i\in I(m,r)$, 
and $l_j:=(0,e_j)$, $j\in\{1,\ldots ,n\cdot d\}$. Let ${\cal C}\equiv 
{\cal C}(\Omega)$ denote the set of all cylindrical functions of the above 
form. Let $Q$ be the number appearing in the definition of the density $m$ 
and $1/p+1/Q\le 1$. ${\cal C}$ is dense in $L^p(\Omega,Q^{(m,r)}_{\sbnu})$. 
As an adaption of Proposition \ref{Proposition6.3}, the operator $(D ,
{\cal C}):=(j^\ast\circ {\Bbb D},{\cal C})$ is closable on $L^p(\Omega,
Q^{(m,r)}_{\sbnu};H)$. 
\medskip

\nid
{\bf Definition 1.1 continued.}
(b) Let $(D,D_{p,1})\equiv (D,D_{p,1}(Q^{(m,r)}_{\sbnu}))$ denote 
the closure of $(D ,{\cal C}):=(j^\ast\circ {\Bbb D},{\cal C})$ on 
$L^p(\Omega,Q^{(m,r)}_{\sbnu};H)$. \\ 
(c) If no ambiguity is possible we will also use the symbol $(D, 
D_{p,1})$ to denote the (vectors of) gradients of functions of type 
$\Omega\to F$. \\ 
(d) Let $\vp\in D_{p,1}$ be of type $\Omega\to {\Bbb R}$. If $\vp$ 
is for $W\in\Omega$ of the form $\vp(W-W_0\1,W_0)$ and the $n\cdot 
d$-dimensional vector of the directional derivatives of $\vp$ in the 
directions of the components of $W_0$ exists in the sense of 
differentiation in $L^p(\Omega,Q^{(m,r)}_{\sbnu};F)$ then let it be 
denoted by $\nabla_{W_0}\vp\equiv\nabla_{W_0}\vp(W)$. When emphasizing 
that $D\vp(W)\in H$ for $Q^{(m,r)}_{\sbnu}$-a.e. $W\in\Omega$ or 
stressing the particular form of the gradient we will use the notation 
\begin{eqnarray*}
D\vp(W)&&\hspace{-.5cm}\equiv\left((D\vp)_1(W),(D\vp)_2(W)\vphantom{ 
\dot{f}}\right) 
\end{eqnarray*} 
or 
\begin{eqnarray*}
D\vp(W)&&\hspace{-.5cm}\equiv\left((D\vp)_1(W),\nabla_{W_0}\vp(W) 
\vphantom{\dot{f}}\right)\, . 
\end{eqnarray*} 
(e) Let $f\in\bigcup_{k=1}^\infty C_b^1({\Bbb R^{k+n\cdot d}})$ and $W 
\in C({\Bbb R};F)$. Let $\langle {\rm e},W\rangle$ be the vector of all 
$\langle e_j,W_0\rangle_F$, $j\in\{1,\ldots,n\cdot d\}$. Furthermore, 
let $k\in {\Bbb N}$, $i_1,\ldots,i_k\in J(m)$, and $\langle {\rm H},W 
\rangle$ be the vector of all $\langle H_{i_l},dW\rangle_{L^2}$, $l\in  
\{1,\ldots ,k\}$. Let $\vp(W):=f\left(\langle {\rm H},W\rangle,\langle 
{\rm e},W\rangle\right)$. Define here the {\em gradient of a 
cylindrical function} with respect to the measure $Q^{(m)}_{\sbnu}$ by 
\begin{eqnarray*}
{\Bbb D}\vp(W)(j(g,x))&&\hspace{-.5cm}:=\sum_{l\in \{1,\ldots ,k\}}f_l 
\left(\langle{\rm H},W\rangle ,\langle {\rm e},W\rangle\right)\cdot 
\langle H_{i_l},g\rangle_{L^2} \\ 
&&\hspace{-.0cm}+\sum_{j=1}^{n\cdot d}f_j\left(\langle {\rm H},W\rangle 
,\langle {\rm e},W\rangle\right)\cdot\langle e_j,x \rangle_F
\end{eqnarray*} 
where $f_l$ denotes the first order derivative of $f$ relative to the 
index $i_l$, $l\in \{1,\ldots ,k\}$, and $f_j$ is the first order 
derivative of $f$ relative to $j\in\{1,\ldots ,n\cdot d\}$. 
Parts (b)-(d) can now be modified to the definition of $(D,D_{p,1}(Q^{ 
(m)}_{\sbnu}))$ by replacing $I(m,r)$ with $J(m)$, $Q^{(m,r)}_{\sbnu}$ 
with $Q^{(m)}_{\sbnu}$, and $H$ with $\left\{(g,y):g\in L^2({\Bbb R};F) 
,\ y\in F\right\}$ (with inner product $\langle h_1,h_2\rangle=\langle 
f_1,f_2\rangle_{L^2}+\langle x_1,x_2\rangle_F$). 
\medskip 

Referring to part (d), it is important to note that, for $\vp\in 
D_{p,1}$, $D\vp$ is $Q^{(m,r)}_{\sbnu}$-a.e a linear combination of 
$\langle H_i,\cdot\rangle_{L^2}$, $i\in I(m,r)$, and $\langle e_j, 
\cdot\rangle_F$, $j\in\{1,\ldots ,n\cdot d\}$, with random 
coefficients. For this, see also Lemma \ref{Lemma2.2} (a) below. 
\medskip

\nid
{\bf Definition 1.1 continued.}
(f) Let $1/q'+1/Q\le 1$, $1/{p'}+1/{q'}=1$, and $\xi\in L^{p'}\left( 
\Omega,Q^{(m,r)}_{\sbnu};H\right)$. We say that $\xi\in {\rm Dom}_{p'} 
(\delta)$ if there exists $c_{p'}(\xi)>0$ such that 
\begin{eqnarray*} 
\int\langle D\vp,\xi\rangle_H\, dQ^{(m,r)}_{\sbnu}\le c_{p'}(\xi)\cdot 
\|\vp\|_{L^{q'} (\Omega,Q^{(m,r)}_{\sbnu})}\, , \quad\vp\in D_{q',1} 
(Q^{(m,r)}_{\sbnu}). 
\end{eqnarray*} 
In this case, we define the stochastic integral $\delta(\xi)\equiv 
\delta^{(m,r)}(\xi)$ by 
\begin{eqnarray*} 
\int\delta(\xi)\cdot\vp\, dQ^{(m,r)}_{\sbnu}=\int\langle\xi,D\vp 
\rangle_H\, dQ^{(m,r)}_{\sbnu}\, , \quad\vp\in D_{q',1}(Q^{(m,r)}_{ 
\sbnu})\, . 
\end{eqnarray*} 
The representations of Theorem \ref{Theorem6.7} apply. Furthermore, 
$\delta (\xi)\in L^{p'}\left(\Omega,Q^{(m,r)}_{\sbnu}\right)$ by the 
Hahn-Banach theorem. Furthermore, in order to motivate this definition, 
we note that $D_{q',1}(Q^{(m,r)}_{\sbnu})$ is defined for $1/q'+1/Q\le 
1$, recall also Proposition \ref{Proposition6.3}. 
\medskip 

{\bf Gradients with respect to an individual trajectory } Let $\xi 
\equiv\xi(W)$ be an $\hat{F}$-valued random variable with $\hat{F}$ 
being a metric space. 
\begin{definition}\label{Definition1.2}
{\rm (a) For fixed $W\in\Omega$ and $x\in F$, let 
\begin{eqnarray*} 
\nabla_x\xi(W+x\1):=\sum_{j=1}^{n\cdot d}\frac{\partial}{\partial x_j} 
\xi(W+x\1)\cdot (0,e_j) 
\end{eqnarray*} 
whenever the partial derivatives $\frac{\partial}{\partial x_j}\xi(W+x\1) 
=\lim_{h\to 0}\frac1h\left(\vphantom{\displaystyle i^1}\xi(W+x\1+he_j\1)- 
\xi(W+\right.$ $\left.x\1)\vphantom{\displaystyle i^1}\right)$ exist in 
$\hat{F}$ for all $j\in\{1,\ldots ,n\cdot d\}$. \\ 
(b) Set 
\begin{eqnarray*} 
\nabla_{W_0}\xi(W):=\left.\nabla_x\right|_{x=W_0}\xi(W-W_0\1+x\1)\, . 
\end{eqnarray*} 
(c) For $W\in\Omega$ and $\kappa_i:=\int_0^\cdot H_i(v)\, dv$ let 
\begin{eqnarray*} 
\nabla_G\xi(W):=\sum_{i\in I(r)}\frac{\partial}{\partial\kappa_i}\xi(W) 
\cdot (H_i,0) 
\end{eqnarray*} 
whenever the partial derivatives $\frac{\partial}{\partial\kappa_i}\xi(W) 
=\lim_{h\to 0}\frac1h\left(\vphantom{\displaystyle i^1}\xi(W+h\kappa_i)- 
\xi(W)\vphantom{\displaystyle i^1}\right)$ exist in $\hat{F}$ for all $i 
\in I(r)$ and the sum converges in $H$. \\ 
(d) For $W\in\Omega$, set 
\begin{eqnarray*} 
\nabla_H\xi(W):=\nabla_G\xi(W)+\nabla_{W_0}\xi(W)\, . 
\end{eqnarray*} 
} 
\end{definition}
\medskip 

{\bf Finite dimensional divergence } Let $\xi\equiv\xi(W)$ be an $F\equiv 
{\Bbb R}^{n\cdot d}$-valued random variable. Let $\nabla_{W_0}$ be either 
as in part (d) of Definition \ref{Definition1.1} or as in part (d) of 
Definition \ref{Definition1.2}. Set 
\begin{eqnarray*} 
\nabla_{d,W_0}\xi(W):=\left(\left(\nabla_{W_0}\left(\xi(W)\vphantom{l^1} 
\right)_1\right)_1,\ldots ,\left(\nabla_{W_0}\left(\xi(W)\vphantom{l^1} 
\right)_{n\cdot d}\right)_{n\cdot d}\right)^T\, . 
\end{eqnarray*} 
\begin{definition}\label{Definition1.3}
{\rm The {\em divergence} of $\xi$ {\em with respect to the coordinates of} 
$W_0$ is 
\begin{eqnarray*} 
\left\langle {\sf e},\nabla_{d,W_0}\xi\vphantom{l^1} \right\rangle_F 
\equiv\sum_{j=1}^{n\cdot d}\left\langle e_j,\nabla_{W_0}\left\langle\xi, 
e_j\vphantom{l^1} \right\rangle_F\vphantom{\dot{f}}\right\rangle_F . 
\end{eqnarray*} 
} 
\end{definition}
\bigskip 

{\bf The spaces $I^p$, $E_{q,1}$, and $K_{q,1}$; flows belonging to 
${\cal F}_{p,1}(W)$. } 
\begin{definition}\label{Definition1.4}
{\rm Let $1<p<\infty$ and $I^p\equiv I^p(Q^{(m,r)}_{\sbnu})$ denote 
the set of all $F$-valued processes $Y\equiv Y(W)$ such that
\begin{itemize}
\item[{(i)}] $Y(W)=Y_a(W)+Y_j(W)$ where 
\item[{(ii)}] $Y_a(W)$ is an absolutely continuous process with $j^{-1} 
Y_a\in L^p(\Omega,Q^{(m,r)}_{\sbnu};H)$ and 
\item[{(iii)}] $Y_j(W)$ is a cadlag pure jump process with $\frac12 
\left((Y_j)_{0-}+(Y_j)_0\right)=0$ that is of locally finite variation 
on ${\cal R}$, $V_{-rt}^{(r-1)t}(Y_j)=\sum_{w\in {\cal R}}\sum_{i=1}^{n 
\cdot d}\left|((Y_j)_w-(Y_j)_{w-})_i\right|\in L^p (\Omega,Q^{(m,r)}_{ 
\sbnu})$. 
\end{itemize} 
}
\end{definition}
\begin{definition}\label{Definition1.5}
{\rm (a) Let ${\cal F}_{p,1}(W)\equiv {\cal F}_{p,1}(W,Q^{(m,r)}_{ 
\sbnu})$, $1<p<\infty$, denote the set of all $\Omega$-valued random 
flows $W^\rho\equiv W+g(\rho,W)$, $\rho\in {\Bbb R}$, on $\left(\Omega 
,Q^{(m,r)}_{\sbnu}\right)$, with 
\begin{itemize}
\item[(i)] $\dot{g}(\rho,\cdot)\in jH:=\{jh:h\in H\}$ given by the {\it 
mixed derivative} 
\begin{eqnarray*}
\dot{g}(\rho ,\cdot)(s)&&\hspace{-.5cm}=\frac{d^\pm}{d\rho}\, g(\rho , 
\cdot)(s) \\ 
&&\hspace{-.5cm}:=\frac12\left(\frac{d^-}{d\rho}\, g(\rho ,\cdot)(s)+ 
\frac{d^+}{d\rho}\, g(\rho ,\cdot)(s)\right)\, ,\quad s\in {\cal R},\ 
\rho\in {\Bbb R}\quad Q^{(m,r)}_{\sbnu}\mbox{\rm -a.e.} 
\end{eqnarray*}
\item[(ii)] $j^{-1}\dot{g}(\rho,\cdot)\in L^p(\Omega,Q^{(m,r)}_{\sbnu};H 
)$, $\rho\in {\Bbb R}$. 
\end{itemize} 
(b) Let $Q$ be the number appearing in the definition of the density $m$ 
and $1/q+1/Q<1$. Furthermore, let $1<p<\infty$ with $1/p+1/q=1$. Introduce 
\begin{eqnarray*}
E_{q,1}&&\hspace{-.5cm}\equiv E_{q,1}\left(Q^{(m,r)}_{\sbnu}\right) \\ 
&&\hspace{-.5cm}:=\left\{\vp\in D_{q,1}:\left.\frac{d^\pm}{d\sigma}\right 
|_{\sigma =0}\vp(W^\sigma)=\left\langle D\vp(W),j^{-1}\dot{g}(0,W)\right 
\rangle_H \right. \\ 
&&\hspace{4cm}\left.\vphantom{\left.\frac{d^\pm}{d\sigma}\right|_{ 
\sigma =0}}Q^{(m,r)}_{\sbnu}\mbox{-a.e.}\ \mbox{\rm for all}\ (W^\rho)_{ 
\rho\in {\Bbb R}}\in {\cal F}_{p,1}(W)\right\}\, . 
\end{eqnarray*}
}
\end{definition}
\begin{definition}\label{Definition1.6} 
{\rm Let $Q$ be the number appearing in the definition of the density 
$m$ and $1/q+1/Q<1$. Let $K_{q,1}\equiv K_{q,1}\left(Q^{(m,r)}_{\sbnu} 
\right)$ be the set of all processes $Y$ for which 
\begin{itemize}
\item[(i)] $Y_s\in E_{q,1}$, $s\in {\cal R}$, 
\item[(ii)] $Q^{(m,r)}_{\sbnu}$-a.e. we have $D_tY_\cdot\in jH$ and $t\to 
j^{-1}D_tY_\cdot$ has a right continuous version for which the linear 
span of $\{j^{-1}D_tY_\cdot:t\in {\cal R}\}$ is dense in $H$. Furthermore, 
$j^{-1}DY_\cdot\in L^q(\Omega,Q^{(m,r)}_{\sbnu};H\otimes H)$ and
\item[(iii)] we have $\langle j^{-1}DY_\cdot,h\rangle_{H\to F}\in L^p 
(\Omega,Q^{(m,r)}_{\sbnu};H)$, and $\langle DY_\cdot,h\rangle_{H\to F}= 
j\langle j^{-1}DY_\cdot ,h\rangle_{H\to F}$ for all $h\in L^p (\Omega,Q^{ 
(m,r)}_{\sbnu};H)$, $1/p+1/q=1$.  
\end{itemize} 
}
\end{definition} 
Let $h\in H$. It follows from (ii) that $Q^{(m,r)}_{\sbnu}$-a.e. we have 
\begin{eqnarray*}
\left\langle\left|j^{-1}D_\cdot Y\right|,h\right\rangle_{H\to F}\in H\, , 
\quad\left\langle\left|j^{-1}DY_\cdot\right|,h\right\rangle_{H\to F}\in 
H\, ,\quad\left\langle j^{-1}DY_\cdot,h\right\rangle_{H\to F}\in H\, .  
\end{eqnarray*} 

{\bf The process $X$. } 
Let $A\equiv A(W)$ be a stochastic process with trajectories that 
belong to $B_{b;{\rm loc}}({\Bbb R};F)$ and are cadlag. Introduce 
\begin{eqnarray*} 
X:=W+A\, .  
\end{eqnarray*} 
This induces path wise a measurable map $X=u(W):\Omega\equiv C({\Bbb 
R};F)\to B_{b;{\rm loc}}({\Bbb R};F)$ which we assume to be injective. 
The measure $P_{\sbmu}:=Q_{\sbnu}\circ u^{-1}$ is the law of $X$. 
Since $Q_{\sbnu}(W_0\in D^n)=1$ it is reasonable to set $A(W)\equiv 0$ 
whenever $W_0\not\in D^n$. Otherwise we will assume that $X_s\in D^n$, 
$s\in {\Bbb R}$. 

The subscript $\bmu$ indicates the distribution of $X_0=W_0+A_0$.  
We will use the symbol $X$ when we primarily want to refer to $X 
\equiv X(W)$ as a stochastic process, and we will use the symbol 
$u\equiv u(W)$ in order to refer to the map $u:C({\Bbb R};F)\to 
B_{b;{\rm loc}}({\Bbb R};F)$. Let us define $\Omega^u:=\{u(W):W 
\in C({\Bbb R};F)\}$ and 
\begin{eqnarray*}
P^{(m,r)}_{\sbmu}:=Q^{(m,r)}_{\sbnu}\circ u^{-1}\quad\mbox{\rm and} 
\quad P^{(m)}_{\sbmu}:=Q^{(m)}_{\sbnu}\circ u^{-1}\, .  
\end{eqnarray*} 

{\bf Jumps of $X$. } Let $\Delta X_w:=X_w-X_{w-}$. For $W\in\Omega$ 
with $W_0=0$ and $X=u(W)$ let $0<\tau_1\equiv\tau_1(X)<\tau_2\equiv 
\tau_2(X)\ldots\, $ denote the jump times of $X$ on $(0,\infty)$, 
i.e., the times $w$ for which $\Delta X_w\neq 0$. Accordingly, let 
$\ldots\, \tau_{-2}\equiv\tau_{-2}(X)<\tau_{-1}\equiv\tau_{-1}(X) 
\le 0$ denote the jump times on $(-\infty ,0]$. 
\medskip

Based on this, let us introduce the numbering of the jumps of $X=u 
(W)$ for general $W_0\in F$. 
\begin{definition}\label{Definition1.7}
{\rm Suppose that for each $k\in {\Bbb Z}\setminus\{0\}$ there is a 
map $\tau_k:\Omega\equiv\{W\in\Omega:W_0=0\}\times F\to {\Bbb R}$ 
satisfying the following. 
\begin{itemize} 
\item[(j)] For each fixed $W\in\Omega$ with $W_0=0$ and $k\in {\Bbb 
Z}\setminus\{0\}$ the map $F\ni x\to\tau_k(u\circ (W+x\1))$ is {\it 
piecewise continuously differentiable} on $F$. That is, there exist 
finitely many mutually exclusive open sets $F_1\equiv F_1(k),F_2 
\equiv F_2(k),\ldots\, $ with piecewise $C^1$-boundary and $\overline{ 
F}=\bigcup_i\overline{F_i}$ such that 
\begin{eqnarray*}
\left.\tau_k(u\circ (W+\, \cdot\, \1)){\vphantom{\dot{f}}}\right|_{ 
F_i}\in C^2_b(F_i) 
\, , \quad i,k\in {\Bbb N}. 
\end{eqnarray*} 
In addition, if $x_0\in\bigcup_i\partial F_i$ then there exist $i_0 
\equiv i_0(k)$ such that $\tau_k(u\circ (W+\, \cdot\, \1))$ as well 
as $\nabla_x\tau_k(u\circ (W+\, \cdot\, \1))$ can continuously be 
extended from $F_{i_0}$ to $x_0$. 
\item[(jj)] For each $W\in\Omega$ and $X=u(W)$ the sequence $\ldots 
\, \tau_{-2}\equiv\tau_{-2}(X)<\tau_{-1}\equiv\tau_{-1}(X)<\tau_1 
\equiv\tau_1(X)<\tau_2\equiv\tau_2(X)\ldots\, $ is the sequence of 
all {\em jump times of} $X$. 
\end{itemize} 

For each $W\in\Omega$ and $X=u(W)$ let $\tau_0\equiv\tau_0(X):=\sup\{ 
\tau_k:\tau_k\le 0,\ k\in {\Bbb Z}\setminus\{0\}\}$. 
}
\end{definition}
However keep in mind that $\tau_{-1}\le\tau_0<\tau_1$ is not necessarily 
true. For every finite subinterval $S$ of ${\Bbb R}$ and all $W\in\Omega$ 
with $W_0=0$, let us introduce 
\begin{eqnarray*}
G(S;W):=\bigcup_{k,i}\overline{\left\{x\in\partial F_{i}(k):\tau_k(u 
\circ(W+x\1))\in S\vphantom{l^1}\right\}}\, .   
\end{eqnarray*} 
Let us furthermore define $G(W):=\bigcup_SG(S;W)$ where the union is 
taken over all finite subinterval $S$ of ${\Bbb R}$ and stress that 
we have $\bnu (G(W))=0$. Moreover, to simplify the notation we will 
use $X_{\tau_k}(W)$, $A_{\tau_k}(W)$, etc. rather than $X_{\tau_k\circ 
u}(W)$, $A_{\tau_k\circ u}(W)$ etc. 
\medskip 

\nid
{\bf Remark. (1)} For given $W\in\Omega$ with $W_0=0$, let $\{G_i:i 
\in {\Bbb N}\}$ denote the set of all components of $F\setminus G(W 
)$. It is a consequence of Definition \ref{Definition1.7} that, for 
given $W\in\Omega$ with $W_0=0$ and $x_0\in D^n$ there exists a {\it 
time re-parametrization $\sigma(s,x)\equiv\sigma(s,x;W,x_0)$, $s\in 
{\Bbb R}$, $x\in F$, with respect to $W$ and $x_0$} such that we 
have the following. 
\begin{itemize} 
\item[(i)]
\begin{eqnarray*}
\left.\sigma {\vphantom{\dot{f}}}\right|_{{\Bbb R}\times G_i}\in C^2_b 
({\Bbb R}\times G_i)\, , \quad i\in {\Bbb N}.  
\end{eqnarray*} 

\item[(ii)] For $x_1\in\bigcup_i\partial G_i$ there exist $i_1$ 
independent of $x_0$ such that $\sigma(s,\cdot)$ as well as $\nabla_x 
\sigma(s,\cdot)$ can continuously be extended from $G_{i_1}$ to $x_1$, 
$s\in {\Bbb R}$. 

\item[(iii)] $\sigma(s,x_0)=s$ for all $s\in {\Bbb R}$.  
\item[(iv)] For every $x\in F$, $\sigma(\cdot,x)\in C^1({\Bbb R};{\Bbb R} 
)$ and, coordinate wise, 
\begin{eqnarray*}
\frac{\partial}{\partial s}\sigma (s,x)>0\, ,\quad s\in S.  
\end{eqnarray*} 
\item[(v)] $\D\sigma\left(\tau_k\circ u(W+x_0\1),x\vphantom{l^1}\right) 
=\tau_k\circ u\left(W+x\1\vphantom{l^1}\right)$, $x\in F$, $k\in {\Bbb Z} 
\setminus\{0\}$.   
\end{itemize} 
Keeping property (v) of Remark (1) in mind, the next condition (jjj) 
says that jumps of parallel trajectories in $\Omega$ have to be 
compatible with each other in a certain sense. This condition is 
trivially satisfied if, for example, parallel trajectories in $\Omega$ 
generate identical jump times for $X$. In fact, the following condition 
(jjj) is crucial if parallel trajectories in $\Omega$ do not necessarily 
generate identical jump times for $X$. Condition (jv) has been discussed 
in the motivating part of this section.
\medskip

\nid 
{\bf Definition 1.7 continued.} For $W\in\Omega$ with $W_0=0$ 
and $x_0\in D^n$ we assume that there exists a time re-parametrization 
$\sigma(s,x)\equiv\sigma(s,x;W,x_0)$, $s\in {\Bbb R}$, $x\in D^n$, with 
respect to $W$ and $x_0$ such that (i)-(v) of Remark (1) and we have the 
following. 
\begin{itemize}
\item[(jjj)] For every finite subinterval $S$ of ${\Bbb R}$ there is 
an $\ve\equiv\ve(S)>0$ such that for all $k\in {\Bbb Z}\setminus\{0\}$ 
with $\tau_k\circ u(W+x_0\1)\in S$ we have 
\begin{eqnarray*}
&&\hspace{-.5cm}\sigma\left(\tau_k\circ u(W+x\1)+\delta,x\vphantom 
{l^1}\right)-\sigma\left(\tau_k\circ u(W+x_0\1),x\vphantom{l^1}\right) 
 \\ 
&&\hspace{.5cm}=\tau_k\circ u\left(W+x\1\vphantom{l^1}\right)+\delta- 
\tau_k\circ u\left(W+x_0\1\vphantom{l^1}\right)
\end{eqnarray*} 
for all $|\delta|<\ve$ and $x\in D^n$. 
\item[(jv)] $\nabla_{W_0}\tau_k\circ u(W)$ is orthogonal in $F$ to 
$\Delta\left(A_{\tau_k}(W)-Y_{\tau_k}(W)\right)$, $k\in {\Bbb Z} 
\setminus\{0\}$. 
\end{itemize} 
\medskip

{\bf Assumptions on $X$. } Let us post a list of assumptions on $X$ 
which we will make use of in  different stages of the paper. Again, 
let $Q$ be the number appearing in the definition of the density $m$. 
\begin{itemize} 
\item[(1)] 
$A$ {\it has on $\bigcup_m\{\pi_mW:W\in\Omega\}$ a local spatial gradient}. 
That is, 
\begin{itemize} 
\item[(i)] for all $W\in\{\pi_{m}V:V\in\Omega\}$ with $W_0=0$, $x\in 
F$, $s\in {\Bbb R}$, and $k\equiv k(W,x)=\max\{l\in{\Bbb Z}\setminus 
\{0\}:\tau_l(u(W+x\1))\le s\}$ the gradient $\nabla_x A_s(W+x\1)$ exists 
and is bounded and continuous on some neighborhood $U_{s,x}\subset F_i 
\times {\Bbb R}\equiv F_i(k)\times {\Bbb R}$ of $(s,x)$ whenever $x\in 
F_i$ and $s>\tau_k(u(W+x\1))$. 

\item[(ii)] If $x\in\bigcup_{i'}\partial F_{i'}(k)$ and $s>\tau_k(u(W+x\1) 
)$ then, with $i_0$ given by (j), there exists a neighborhood $U_s\times 
V_x$ of $(s,x)$ such that the gradient $\nabla_x A_s(W+x\1)$ is 
well-defined, bounded, and continuous on $(U_s\times V_x)\cap ({\Bbb R} 
\times F_{i_0})\equiv (U_s\times V_x)\cap ({\Bbb R}\times F_{i_0}(k))$ 
and can continuously be extended from $U_s\times F_{i_0}$ to $U_s\times 
\{x\}$.

\item[(iii)] If $s=\tau_k(u(W+x\1))$ and $x\in F_i(k)$ or $x\in\partial 
F_{i_0}(k)\equiv \partial F_i(k)$ in the sense of the previous paragraph, 
then there exists a neighborhood $U_{s,x}$ of $(s,x)$ such that the 
gradient $\nabla_x A_s(W+x\1)$ is well-defined, bounded, and continuous 
on $U_{s,x}\cap\{(v,y)\in {\Bbb R}\times F_i:\tau_{k(W,x)}(u(W+y\1))\le 
v<\tau_{k(W,x)+1}(u(W+y\1))\}$ and can continuously be extended from this 
set to $(s,x)$. 

\item[(iv)] The bound on the gradients $\nabla_x A_s(W+x\1)$ is uniform 
for all $W\in\bigcup_m\{\pi_{m}V:V\in\Omega,V_0=0\}$ and $(s,x)\in S\times 
F$ where $S$ is any finite subinterval of ${\Bbb R}$. 

\end{itemize} 
\end{itemize} 
In this way, $\nabla_{W_0}A_s(W)$ is well-defined for all $W\in 
\bigcup_m\{\pi_{m}V:V\in\Omega\}$ and $s\in {\Bbb R}$. Condition (1) is 
rather sophisticated. In simplified situations we may use the following 
condition. 

\begin{itemize} 
\item[(1')] $u_s\in D_{q,1}\equiv D_{q,1}(Q^{(m,r)}_{\sbnu})$, $s\in 
{\Bbb R}$, for some $q$ with $1/q+1/Q<1$. 
\end{itemize} 

We remark that for the well-definiteness of $u_s\in D_{q,1}$ we only 
need $1/q+1/Q\le 1$, cf. Definition \ref{Definition1.1}. However, for 
example, in Lemma \ref{Lemma3.2} (f) below we need $1/q+1/Q<1$. 
For the next condition we recall the fact that the right continuous 
version of a Radon-Nikodym derivative with respect to the time is, 
whenever it exists, nothing but the right derivative. We will use the 
symbol ``$\, '\, "$ for this. 
\begin{itemize} 
\item[(2)] 
\begin{itemize} 
\item[(i)] For all $W\in\Omega$ and $A\equiv A(W)$, $A=A^1+A^2$. $A^1$ 
is a continuous process and $A^2$ is a pure jump process with $A^2_{ 
\tau_{-1}}=A^2_{\tau_1-}=0$ such that, for fixed time $s\in {\Bbb R}$, 
$A^2_s$ jumps only on $Q_{\sbnu}^{(m)}$- and $Q_{\sbnu}$-zero sets. 

Furthermore, on $\bigcup_m\left\{(\pi_m W)_{\cdot +w}:W\in\Omega\, ,\ w 
\in\left[0,\frac{1}{2^m}\cdot t\right)\right\}$, the jumps of $A^2$ do not 
accumulate at finite time and $A^1$ possesses a Radon-Nikodym derivative 
with cadlag version ${A^1}'=A'$. 
\item[(ii)] Moreover, 
\begin{eqnarray*} 
\nabla_{W_0}A'_s(W):=\left.\nabla_x\right|_{x=W_0}A'_s(W-W_0\1+x\1)\, . 
\end{eqnarray*} 
exists and is bounded on $\bigcup_m\left\{(\pi_m W)_{\cdot +w}:W\in 
\Omega\, ,\ w\in\left[0,\frac{1}{2^m}\cdot t\right)\right\}$ uniformly 
for all $s$ belonging to any finite subinterval of ${\Bbb R}$. The same 
holds for the right continuous version of $A_{-\, \cdot}\ $. 
\item[(iii)] 
For $W\in\bigcup_m\left\{(\pi_m V)_{\cdot +w}:V\in\Omega\, , \ w\in\left 
[0,\frac{1}{2^m}\cdot t\right)\right\}$ with $W_0=0$ and $k\in {\Bbb Z} 
\setminus\{0\}$ the map $F\ni x\to\Delta A_{\tau_k}(W+x\1)$ is piecewise 
continuously differentiable on $F$. That is $\Delta A_{\tau_k}(W+\, \cdot 
\, \1)$ satisfies (j) of the definition of the jump times with $\Delta 
A_{\tau_k}(W+\, \cdot\, \1)$ instead of ${\tau_k}(W+\, \cdot\, \1)$. 
\end{itemize} 
\end{itemize} 
\begin{itemize}
\item[(3)] 
\begin{itemize}
\item[(i)] There exists $Y\equiv Y(W)$ such that $X\equiv X(W)$ is a 
{\it temporally homogeneous} function of $W\in\Omega$ in the sense that 
for $W^u:=W_{\cdot +u}+\left(A_u(W)-Y_u(W)\vphantom{l^1}\right)\1$ we 
have 
\begin{eqnarray*} 
W^0=W\, , \quad A_0\left(\left(W^u\vphantom{l^1}\right)^v\right)=A_0 
(W^{u+v})\, , 
\end{eqnarray*} 
and 
\begin{eqnarray*} 
X_{\cdot+v}(W)=X(W^v)\, , \quad u,v\in {\Bbb R}. 
\end{eqnarray*} 
\end{itemize} 
\end{itemize} 
In the subsequent part (ii) we use the fact that, according to part 
(i) of this condition, the jump times of $u(W^v)$ are the jump times 
of $X=u(W)$ minus $v$. For a subset $S\subset {\Bbb R}$, we define $
S\ominus v=\{s\in {\Bbb R}:s+v\in S\}$. 
\begin{itemize}
\item[(3)] 
\begin{itemize}
\item[(ii)] For every finite subinterval $S$ of ${\Bbb R}$ and $W\in\{ 
\pi_mV:V\in\Omega\}$ with $W_0\not\in G(S;W-W_0\1)$ there exists $\ve>0 
$ such that for all $v\in S$ the Euclidean distance of $W^v_0$ and 
$G(S\ominus v;W^v-W^v_0\1)$ is at least $\ve$. 
\item[(iii)] 
Furthermore, $\sup_{s,W}|Y_s(W)|<\infty$ and $Y(W)\equiv 0$ if $W_0\not 
\in D^n$ and $Y$ has on $\bigcup_m\{\pi_mW:W\in\Omega\}$ a local spatial 
gradient (cf. (1)). In addition, $Y$ has all the properties of $A$ 
mentioned in (2). In particular, the set of the jump times of $Y$ is a 
subset of the set of the jump times of $A$. 
\end{itemize} 
\end{itemize} 

The following condition describes the notion of continuity on $\Omega$ 
we are going to use. 
\begin{itemize} 
\item[(4)] 
\begin{itemize}
\item[(i)] $\Omega\ni W_n\stack{n\to\infty}{\lra}W\in\Omega$ relative 
to the topology of coordinate wise uniform convergence on compact subsets 
of ${\Bbb R}$ and $W_0\not\in G(S;W-W_0\1)$ imply 
\begin{eqnarray*}
A_v(W_n)\stack{n\to\infty}{\lra}A_v(W)\quad\mbox{\rm and}\quad Y_v(W_n) 
\stack{n\to\infty}{\lra}Y_v(W) 
\end{eqnarray*} 
and 
\begin{eqnarray*}
\nabla_{W_0}A_v(W_n)-\nabla_{W_0}Y_v(W_n)\stack{n\to\infty}{\lra} 
\nabla_{W_0}A_v(W)-\nabla_{W_0}Y_v(W)
\end{eqnarray*} 
for all $v\in S\setminus\{\tau_k(X):k\in {\Bbb Z}\setminus\{0\}\}$ and 
any finite subinterval $S\subset{\Bbb R}$. Furthermore, 
\begin{eqnarray*}
\Delta A_{\tau_k}(W_n)-\Delta Y_{\tau_k}(W_n)\stack{n\to\infty}{\lra} 
\Delta A_{\tau_k}(W)-\Delta Y_{\tau_k}(W)\, , \quad k\in {\Bbb Z} 
\setminus\{0\}. 
\end{eqnarray*} 
\item[(ii)] 
\begin{eqnarray*}
A(W_n+x\1)-Y(W_n+x\1)-\left(A(W+x\1)-Y(W+x\1)\vphantom{l^1}\right) 
\stack{n\to\infty}{\lra}0 
\end{eqnarray*} 
in $L^1(S;F)$ uniformly in $x\in F$ for all finite subintervals $S\subset 
{\Bbb R}$. 
\end{itemize} 
\end{itemize} 

Instead of condition (4) the following will be used in particular in 
Section 2 below. 
\begin{itemize} 
\item[(4')] $\Omega\ni W_n\stack{n\to\infty}{\lra}W\in\Omega$ relative 
to the topology of coordinate wise uniform convergence on compact subsets 
of ${\Bbb R}$ implies $X_s(W_n)\equiv u_s(W_n)\stack{n\to\infty}{\lra} 
X_s\equiv u_s(W)$ for all $s\in {\Bbb R}\setminus\{\tau_k(X):k\in {\Bbb 
Z}\setminus\{0\}\}$.
\end{itemize} 

\nid
{\bf Remark. (2)} According to conditions (1) and (3) $A$ and $Y$ have on 
$\{\pi_mV:V\in \Omega\}$ a local spatial gradient. Since by conditions 
(2) and (3), $A^2_{\tau_k}(W+\, \cdot\, \1)$ and $Y^2_{\tau_k}(W+\, \cdot 
\, \1)$ are for $W\in\{\pi_mV:V\in\Omega\}$ with $W_0=0$ and $k\in{\Bbb Z} 
\setminus\{0\}$ piecewise continuously differentiable on $F$, $A^1$ and 
$Y^1$ have on $\{\pi_mV:V\in\Omega\}$ also a local spatial gradient. 
\bigskip

To ease the formulation of Theorem \ref{Theorem1.11} below one could also 
require the following. 
\begin{itemize} 
\item[(4'')] For all $\bigcup_{m}\{\pi_mV:V\in\Omega\}\ni W_n\stack{n 
\to\infty}{\lra}W\in\Omega$ relative to the topology of coordinate wise 
uniform convergence on compact subsets of ${\Bbb R}$ with $W_0\not\in 
G(S;W-W_0\1)$, and all $s\in S\setminus\{\tau_k(X):k\in {\Bbb Z} 
\setminus\{0\}\}$ where $X=u(W)$ there exist the limit  
\begin{eqnarray*}
\lim_{n\to\infty}\nabla_{W_0}A^1_s(W_n)-\nabla_{W_0}Y^1_s(W_n)=:\nabla_{ 
W_0}A^1_s(W)-\nabla_{W_0}Y^1_s(W) 
\end{eqnarray*} 
where $S$ is any finite subinterval of ${\Bbb R}$. 
\end{itemize} 
Assuming (4''), the piecewise continuous differentiability of $\Delta A$ 
and $\Delta Y$, cf. conditions (2) and (3), together with condition (4) 
(i) define $\nabla_{W_0}\Delta A_{\tau_k}-\nabla_{W_0}\Delta Y_{\tau_k}$, 
$k\in {\Bbb Z}\setminus\{0\}$, $Q_{\sbnu}$-a.e. 
\medskip

In order to formulate the next condition, introduce
\begin{eqnarray*} 
f(s):=\left\{ 
\begin{array}{cl}
\D\exp\left\{-1/\left(1-s^2\vphantom{l^1}\right)\vphantom{\dot{f}} 
\right\}\, , & s\in\left(-1,1\right), \\ 
0\, , & s\in {\Bbb R}\setminus\left(-1,1\right),\vphantom{\D\dot{f}} 
\end{array} 
\right. 
\end{eqnarray*} 
and $c_n:=\left(\int_{\Bbb R} f(ns)\, ds\right)^{-1}$ as well as 
$c_{F,n}:=\left(\int_F f(n|x|)\, dx\right)^{-1}$. Define the 
mollifier functions $g_n(s):=c_n\cdot f(ns)\, {\sf e}=nc_1\cdot f(ns) 
\, {\sf e}$, $s\in {\Bbb R}$, and 
\begin{eqnarray*}
\gamma_n(x):=c_{F,n^3}\cdot f\left(n^3|x|\right){\sf e}\, , \quad x\in 
F,\ n\in {\Bbb N}. 
\end{eqnarray*} 
For the motivation of the order $n^3$ in the definition of $\gamma_n$
we refer to the text above (\ref{2}) in Subsection 3.3. In addition, 
set 
\begin{eqnarray*} 
A_s(\, \cdot\, ,\gamma_n)(W):=\int_F\left\langle A_s(W+x\1)\, ,\, 
\gamma_n(x)\vphantom{l^1}\right\rangle_{F\to F}\, dx 
\end{eqnarray*} 
$s\in {\Bbb R}$, $n\in {\Bbb N}$. Likewise define $Y_s(\, \cdot\, , 
\gamma_n)$. Furthermore, let $\kappa_i:=\int_0^\cdot H_i(u)\, du$, 
$i\in {\Bbb N}$. 
\begin{itemize} 
\item[(5)] 
\begin{itemize}
\item[(i)] For all $r\in {\Bbb N}$ and all $s\in {\cal R}\equiv 
{\cal R}(r)$ there exists $h\equiv h_s(r)\in H\equiv H^{(r)}$ such 
that for all $W,w\in jH$ with $w_0=0$, and all $\ve\in (-1,1)$, 
\begin{eqnarray*}
&&\hspace{-.5cm}\left|A_s(W+\ve w)-A_s(W)\vphantom{l^1}\right|\le 
\left|\left\langle w',h_s\right\rangle_{L^2}\right|\cdot\max_j\left| 
A_s(W+\ve e_j\1)-A_s(W)\vphantom{l^1}\right| 
\end{eqnarray*} 
as well as 
\begin{eqnarray*}
&&\hspace{-.5cm}\left|Y_s(W+\ve w)-Y_s(W)\vphantom{l^1}\right|\le 
\left|\left\langle w',h_s\right\rangle_{L^2}\right|\cdot\max_j\left| 
Y_s(W+\ve e_j\1)-Y_s(W)\vphantom{l^1}\right| 
\end{eqnarray*} 
Furthermore $\sup_{s\in {\cal R}}\|h_s\|_{L^2}^2<\infty$.
\item[(ii)] Furthermore, for all $r\in {\Bbb N}$, all $n\in {\Bbb N}$, 
and all $s\in {\cal R}$, $A_s(\, \cdot\, ,\gamma_n)(W)$ and $Y_s(\,  
\cdot\, ,\gamma_n)(W)$ are two times Fr\' echet differentiable on $W 
\in jH\equiv jH^{(r)}$. 
\end{itemize} 
\end{itemize} 
For an $F$-valued random variable $\xi$ we will use the symbol $D_F\xi 
(W)$ to denote the Fr\' echet derivative of $\xi$ at $W\in jH$ if it 
exists. In this case, 
\begin{eqnarray*}
D_F\xi(W)=\nabla_H\xi(W)\equiv\sum_{i\in I(r)}\frac{\partial\xi(W)} 
{\partial\kappa_i}(H_i,0)+\sum_{j'}\frac{\partial\xi(W)}{\partial 
\lambda_{j'}}\left(0,e_{j'}\right)\in H\, . 
\end{eqnarray*} 
Furthermore, we will use the notation 
\begin{eqnarray*}
(D_{F,1})_s\xi(W):=\sum_{i\in I(r)}\frac{\partial\xi(W)}{\partial\kappa_i 
}\cdot H_i(s)\, ,\quad s\in {\cal R}, 
\end{eqnarray*} 
if $\xi$ is on $W\in jH$ Fr\' echet differentiable. For the second order 
Fr\' echet derivative we will use the symbol $D^2_F\xi(W)$. 
\begin{itemize} 
\item[(5)] 
\begin{itemize}
\item[(iii)] For all $W\in jH$, the gradients $\nabla_GA_s(\, \cdot\, , 
\gamma_n)(W)\equiv\left(\nabla_G\right)_rA_s(\, \cdot\, , \gamma_n)(W)$ 
as well as $\nabla_GY_s(\, \cdot\, , \gamma_n)(W)\equiv\left(\nabla_G 
\right)_rY_s(\, \cdot\, , \gamma_n)(W)$ are continuously differentiable 
with respect to the variable $r\in {\cal R}$. 
\end{itemize} 
\end{itemize} 
%
%
\bigskip

\nid
{\bf Remarks.} (on conditions (1)-(4)) {\bf (3)} It follows from condition 
(3) that 
\begin{eqnarray*}
X_v=X_0(W^v)=W_v+(A_v(W)-Y_v(W))+A_0(W^v) 
\end{eqnarray*} 
which implies 
\begin{eqnarray*}
Y_v(W)=A_0(W^v)\, , \quad v\in {\Bbb R}.  
\end{eqnarray*} 
Furthermore, $Y_v(W^u)=A_0\left((W^u)^v\vphantom{l^1}\right)=A_0 
(W^{u+v})=Y_{u+v}(W)$ which means that we have in addition to condition (3) 
also  
\begin{eqnarray*} 
Y_{\cdot+v}(W)=Y(W^v)\, , \quad v\in {\Bbb R}. 
\end{eqnarray*} 
Moreover, 
\begin{eqnarray*}
A_v(W^u)-A_0(W^u)&&\hspace{-.5cm}=X_v(W^u)-X_0(W^u)-\left((W^u)_v-(W^u 
)_0\vphantom{l^1}\right)\vphantom{\left(\dot{f}\right)} \\ 
&&\hspace{-.5cm}=X_v(W^u)-X_0(W^u)-\left(W_{u+v}-W_u\vphantom{l^1}\right) 
\vphantom{\left(\dot{f}\right)} \\ 
&&\hspace{-.5cm}=X_{u+v}(W)-X_u(W)-\left(W_{u+v}-W_u\vphantom{l^1}\right) 
\vphantom{\left(\dot{f}\right)} \\ 
&&\hspace{-.5cm}=A_{u+v}(W)-A_u(W)\, , \quad u,v\in {\Bbb R}. \vphantom{ 
\left(\dot{f}\right)} 
\end{eqnarray*} 

It follows now from condition (3) that $W(\, \cdot\, )$ possesses the {\it 
flow property}, 
\begin{eqnarray*} 
(W^u)^v&&\hspace{-.5cm}=(W^u)_{\cdot +v}+\left(A_v(W^u)-Y_v(W^u)\vphantom 
{l^1}\right)\1\vphantom{\left(\dot{f}\right)} \\  
&&\hspace{-.5cm}=W_{\cdot +u+v}+\left(A_u(W)-Y_u(W)\right)\1+\left(A_v(W^u) 
-Y_{u+v}(W)\vphantom{l^1}\right)\1\vphantom{\left(\dot{f}\right)} \\ 
&&\hspace{-.5cm}=W_{\cdot +u+v}+\left(A_u(W)-A_0(W^u)+A_v(W^u)\vphantom 
{l^1}\right)\1-Y_{u+v}(W)\1\vphantom{\left(\dot{f}\right)} \\ 
&&\hspace{-.5cm}=W_{\cdot +u+v}+\left(A_{u+v}(W)-Y_{u+v}(W)\vphantom{l^1} 
\right)\1 \vphantom{\left(\dot{f}\right)} \\ 
&&\hspace{-.5cm}=W^{u+v}\, , \quad u,v\in {\Bbb R}.\vphantom{\left(\dot{f} 
\right)}
\end{eqnarray*} 

\nid
{\bf (4)} 
If $Y_v(W)=Y_0(W)=A_0(W)$, $v\in {\Bbb R}$, $W\in\Omega$, then condition 
(3) reads with $W^u:=W_{\cdot +u}+\left(A_u(W)-A_0(W)\vphantom{l^1}\right) 
\1$, $u\in {\Bbb R}$, as 
\begin{eqnarray*} 
A_0(W^v)=A_0(W)\, , \quad X_{\cdot+v}(W)=X(W^v)\, , \quad v\in {\Bbb R}. 
\end{eqnarray*} 
If $A_0(W)=0$, $W\in\Omega$, then condition (3) (i) reads as simple as 
\begin{eqnarray*} 
X_{\cdot +v}=X(W_{\cdot +v}+A_v\1)\, , \quad v\in {\Bbb R}. 
\end{eqnarray*} 

\nid
{\bf (5)} By condition (2), for all $m\in {\Bbb N}$ there exists $Q^{(m) 
}_{\sbnu}$-a.e. a representation of $A$ of the form 
\begin{eqnarray*}
A_s(W)-A_0(W)=\int_0^s b(w,X(W))\, dw+\sum_{w\in (0,s]}\Delta X_w(W) 
\, , \quad s\in {\Bbb R}. 
\end{eqnarray*} 
In case of $A_0(W)=0$, $W\in\Omega$, and $G(W-W_0\1)=\emptyset$ temporal 
homogeneity is under condition (2) and the first line of condition (4) 
(i) equivalent to the following. 
\begin{itemize} 
\item[(3')] We have 
\begin{eqnarray*} 
b(s+v,X(W))=b\left(s,X(W^v)\right) 
\end{eqnarray*} 
and 
\begin{eqnarray*} 
\Delta X_{s+v}(W)=\Delta X_s(W^v)\, , \quad s,v\in {\Bbb R},\quad Q^{ 
(m)}_{\sbnu}\mbox{\rm -a.e.},\ m\in {\Bbb N}. 
\end{eqnarray*} 
\end{itemize} 

\nid 
{\bf (6)} Let $W\in\bigcup_m\{\pi_mV:V\in\Omega\}$. Because of $A(W+ 
x\1)\equiv 0$ whenever $W_0+x\not\in D^n$ and $A_s(W+x\1)=X_s(W+x\1)- 
W_s-x$ and $X_s\in D^n$ where $D$ is a bounded domain, $A_s(W+x\1)$ is 
bounded on $(s,x)\in S\times F$ for any finite subinterval $S\subset 
{\Bbb R}$. 

It follows from condition (1) (iv) that both 
\begin{eqnarray*} 
\quad\int\left\langle\vphantom{\dot{f}} A_{s}\left(W+y\, \1\vphantom 
{l^1}\right)\, , \gamma_n(y-x)\right\rangle_{F\to F}\, dy\quad\mbox{ 
\rm and}\quad\nabla_x\int\left\langle\vphantom{\dot{f}}A_{s}\left(W+ 
y\,\1\vphantom{l^1}\right)\, , \gamma_n(y-x)\right\rangle_{F\to F}\, 
dy\hspace{-.5cm}
\end{eqnarray*} 
are bounded on $(s,x)\in S\times F$ for any finite subinterval $S 
\subset {\Bbb R}$. According to condition (4) (i) these relations 
hold even for all $W\in\Omega$. 

By condition (3), the same holds for $Y$. Because of condition 
(3) we have for $s=v+w$ the identity 
\begin{eqnarray*} 
\left(W+y\, \1\vphantom{l^1}\right)^w=W_{\cdot +w}+\left(y+A_w\left(W+ 
y\, \1\vphantom {l^1}\right)-Y_w\left(W+y\, \1\vphantom{l^1}\right) 
\vphantom{\dot{f}}\right)\1  
\end{eqnarray*} 
and furthermore 
\begin{eqnarray*} 
A_{v}\left(\left(W+y\, \1\vphantom{l^1}\right)^w\right)&&\hspace 
{-.5cm}=X_{v}\left(\left(W+y\, \1\vphantom{l^1}\right)^w\right)- 
\left(W+y\, \1\vphantom {l^1}\right)^w_v \\ 
&&\hspace{-.5cm}=X_{s}\left(W+y\, \1\vphantom{l^1}\right)-W_s- 
\left(y+A_w\left(W+y\, \1\vphantom {l^1}\right)-Y_w\left(W+y\, 
\1\vphantom{l^1}\right)\vphantom{\dot{f}}\right) \\ 
&&\hspace{-.5cm}=A_{s}\left(W+y\, \1\vphantom{l^1}\right)-A_w 
\left(W+y\, \1\vphantom {l^1}\right)+Y_w\left(W+y\, \1\vphantom 
{l^1}\right)\vphantom{\dot{f}}\, . 
\end{eqnarray*} 
Thus, both 
\begin{eqnarray*}
\int\left\langle\vphantom{\dot{f}}A_v\left(W_{\cdot +w}+y\, 
\1\vphantom{l^1}\right)\, , \gamma_n(y-x)\right\rangle_{F\to 
F}\, dy 
\end{eqnarray*} 
and 
\begin{eqnarray*}
\nabla_x\int\left\langle\vphantom{\dot{f}}A_v\left(W_{\cdot 
+w}+y\, \1\vphantom{l^1}\right)\, ,\gamma_n(y-x)\right 
\rangle_{F\to F}\, dy 
\end{eqnarray*} 
are bounded on $(v,w,x)\in S\times S\times F$ for any finite 
subinterval $S\subset {\Bbb R}$. Similarly, corresponding relations 
for $Y$ can be verified. 
\medskip 

\nid 
{\bf (7)} Recall the existence of 
$\nabla_{W_0}A'$ as well as $\nabla_{W_0}Y'$ in the sense of condition 
(2) and the piecewise continuous differentiability of $\Delta A$ and 
$\Delta Y$, cf. conditions (2) and (3). These properties yield the 
existence and the right continuity of $s\to\nabla_{W_0}A^1_s(W)+\sum_{ 
\{k\ge 1:\tau_k\circ u\le s\}}\nabla_{W_0}\Delta A_{\tau_k}$ as well as 
$s\to\nabla_{W_0}Y^1_s(W)+\sum_{\{k\ge 1:\tau_k\circ u\le s\}}\nabla_{ 
W_0}\Delta Y_{\tau_k}$ and implies 
\begin{eqnarray*} 
\nabla_{W_0}A_s=\nabla_{W_0}A^1_s(W)+\sum_{\{k\ge 1:\tau_k\circ u\le s\}} 
\nabla_{W_0}\Delta A_{\tau_k}
\end{eqnarray*} 
as well as 
\begin{eqnarray*} 
\nabla_{W_0}Y_s=\nabla_{W_0}Y^1_s(W)+\sum_{\{k\ge 1:\tau_k\circ u\le s\}} 
\nabla_{W_0}\Delta Y_{\tau_k}
\end{eqnarray*} 
for all, without loss of generality, $s>\tau_{-1}\circ u$ and $W\in\left 
\{(\pi_mV)_{\cdot +w}:V\in\Omega\, , \ w\in\left[0,\frac{1}{2^m}\cdot t 
\right)\right\}$ such that $W_0\not\in G(S;W-W_0\1)$. For the consequence 
of this remark to the formulation of Theorem \ref{Theorem1.11} below, 
recall condition (4''). 
\bigskip 

{\bf The spaces $D^u_{q,1}$, $E^u_{q,1}$, and $K^u_{q,1}$; flows belonging 
to ${\cal F}_{p,1}(X)$. } For the rest of this subsection, let $X=u(W)$ be 
a bijection $C({\Bbb R};F)\to C({\Bbb R};F)$. Although this implies $\Omega 
=\Omega^u$ we will use the symbol $\Omega^u$ for a better indication. 
\begin{definition}\label{Definition1.8}
{\rm Let ${\cal C}(\Omega^u)$ denote the set of all cylindrical functions 
of Definition \ref{Definition1.1} with $W$ replaced by $X$. Let $Q$ be the 
number from the definition of the density $m$ and let $1/q+1/Q\le 1$. Let 
${\Bbb D}^u$ be that version of ${\Bbb D}$ that acts on functions $\psi$ 
with arguments $X =u(W)$ rather than $W$. In particular, suppose that the 
operator $(D^u,{\cal C}):=(j^\ast\circ {\Bbb D}^u,{\cal C}(\Omega^u))$ is 
closable on $L^q(\Omega^u,P^{(m,r)}_{\sbmu};H)$. \\ 
(a) Following the lines of Definition \ref{Definition1.1} and replacing 
there $W$ by $X$, we introduce the space $D^u_{q,1}\equiv D^u_{q,1}(P^{ 
(m,r)}_{\sbmu})$ as the closure of $(D^u,{\cal C}(\Omega^u)):=(j^\ast\circ 
{\Bbb D}^u,{\cal C}(\Omega^u))$ on $L^q(\Omega^u,P^{(m,r)}_{\sbmu};H)$. \\ 
(b) If no ambiguity is possible we will also use the symbol $(D^u,D^u_{q, 
1})$ to denote the (vectors of) gradients of functions of type $\Omega^u 
\to F$. \\ 
(c) Let $\psi\in D^u_{q,1}$ be of type $\Omega^u\to {\Bbb R}$. If $\psi$ 
is for $X\in\Omega^u$ of the form $\psi(X-X_0\1,X_0)$ and the $n\cdot d 
$-dimensional vector of the directional derivatives of $\psi$ in the 
directions of the components of $X_0$ exists in the sense of differentiation 
in $L^q(\Omega^u,P^{(m,r)}_{\sbmu};F)$ then let it be denoted by $\nabla_{ 
X_0}\psi\equiv\nabla_{X_0}\psi(X)$. \\ 
(d) Let $\xi\equiv\xi(X)$ be an $F\equiv{\Bbb R}^{n\cdot d}$-valued random 
variable. Set 
\begin{eqnarray*} 
\nabla_{d,X_0}\xi(X):=\left(\left(\nabla_{X_0}\left(\xi(X)\vphantom{l^1} 
\right)_1\right)_1,\ldots ,\left(\nabla_{X_0}\left(\xi(X)\vphantom{l^1} 
\right)_{n\cdot d}\right)_{n\cdot d}\right)^T\, . 
\end{eqnarray*} 
The {\em divergence} of $\xi$ {\em with respect to the coordinates of} 
$X_0$ is 
\begin{eqnarray*} 
\left\langle {\sf e},\nabla_{d,X_0}\xi\vphantom{l^1} \right\rangle_F 
\equiv\sum_{j=1}^{n\cdot d}\left\langle e_j,\nabla_{X_0}\left\langle\xi, 
e_j\vphantom{l^1} \right\rangle_F\vphantom{\dot{f}}\right\rangle_F . 
\end{eqnarray*} 
}
\end{definition} 
\begin{definition}\label{Definition1.9}
{\rm (a) Let $1<p<\infty$ and ${\cal F}_{p,1}(X)\equiv {\cal F}_{p,1}(X, 
P^{(m,r)}_{\sbmu})$ denote the set of all $\Omega^u$-valued random flows 
$X^\rho\equiv X+g(\rho,X)$, $\rho\in {\Bbb R}$, on $\left(\Omega^u,P^{(m 
,r)}_{\sbmu}\right)$ with 
\begin{itemize}
\item[(i)] $\dot{g}(\rho,\cdot)\in jH:=\{jh:h\in H\}$ given by 
\begin{eqnarray*}
\dot{g}(\rho ,\cdot)(s)=\frac{d^\pm}{d\rho}\, g(\rho ,\cdot)(s)\, ,\quad 
s\in {\cal R},\ \rho\in {\Bbb R}\quad P^{(m,r)}_{\sbmu}\mbox{\rm -a.e.} 
\end{eqnarray*}
\item[(ii)] $j^{-1}\dot{g}(\rho ,\cdot)\in L^p(\Omega^u,P^{(m,r)}_{\sbmu} 
;H)$, $\rho\in {\Bbb R}$. 
\end{itemize} 
(b) Let $Q$ be the number from the definition of the density $m$ and let 
$1/q+1/Q<1$ as well as $1/p+1/q=1$. Furthermore, let $(X^\rho)_{\rho\in 
{\Bbb R}}\in {\cal F}_{p,1}(X)$. Introduce 
\begin{eqnarray*}
E^u_{q,1}&&\hspace{-.5cm}\equiv E^u_{q,1}\left(P^{(m,r)}_{\sbmu}\right) \\ 
&&\hspace{-.5cm}:=\left\{\vp\in D^u_{q,1}:\left.\frac{d^\pm}{d\sigma}\right 
|_{\sigma =0}\vp(X^\sigma)=\left\langle D^u\vp(X),j^{-1}\dot{g}(0,X)\right 
\rangle_H\right.\\ 
&&\hspace{4cm}\left.\vphantom{\left.\frac{d^\pm}{d\sigma}\right|_{\sigma =0} 
}P^{(m,r)}_{\sbmu}\mbox{-a.e.}\ \mbox{\rm for all}\ (X^\rho)_{\rho\in {\Bbb 
R}}\in {\cal F}_{p,1}(X)\right\}\, . 
\end{eqnarray*}
}
\end{definition}
\begin{definition}\label{Definition1.10}
{\rm Let $Q$ be the number from the definition of the density $m$ and 
$1/q+1/Q<1$. Let $K^u_{q,1}\equiv K^u_{q,1}\left(P^{(m,r)}_{\sbmu}\right)$ 
be the set of all processes $Y$ for which 
\begin{itemize} 
\item[(i)] $Y_s\in E^u_{q,1}$, $s\in {\cal R}$, 
\item[(ii)] $P^{(m,r)}_{\sbmu}$-a.e. we have $D^u_tY_\cdot\in jH$, $t\in 
{\cal R}$. 
Furthermore, $j^{-1}D^uY_\cdot\in L^q(\Omega^u,P^{(m,r)}_{\sbmu};H\otimes 
H)$,
\item[(iii)] 
$\langle j^{-1}D^uY_\cdot,h\rangle_{H\to F}\in L^p (\Omega^u,P^{(m,r)}_{ 
\sbmu};H)$, and $\langle D^uY_\cdot,h\rangle_{H\to F}=j\langle j^{-1} D^u 
Y_\cdot,h\rangle_{H\to F}$ for all $h\in L^p(\Omega^u,P^{(m,r)}_{\sbmu};H)$, 
$1/p+1/q=1$. 
\end{itemize} 
}
\end{definition} 
{\bf Remarks. (8)} (on Definitions \ref{Definition1.8}-\ref{Definition1.10}) 
When working with one of the spaces $D^u_{q,1}$, $E^u_{q,1}$, or $K^u_{q,1}$ 
we will always suppose that the operator $(D^u,{\cal C}):=(j^\ast\circ {\Bbb 
D}^u,{\cal C}(\Omega^u))$ is closable on $L^q(\Omega^u,P^{(m,r)}_{\sbmu};H)$. 
\medskip  

\nid
{\bf (9)} (on Definitions \ref{Definition1.5} and \ref{Definition1.9}) Let 
$\vp\in {\cal C}(\Omega)$ be a cylindrical function in the sense of Definition 
\ref{Definition1.1} (a), i. e., $\vp(W):=f\left(\langle{\rm H},W\rangle , 
\langle {\rm e},W\rangle\right)$ with $f\in C_b^1({\Bbb R^{\#I(m,r)+n\cdot 
d}})$, $\langle {\rm e},W\rangle$ denoting the vector of all $\langle e_j, 
W_0\rangle_F$, $j\in \{1,\ldots,$ $n\cdot d\}$, and $\langle {\rm H},W 
\rangle$ being the vector of all $\langle H_i,dW\rangle_{L^2}$, $i\in I(m, 
r)$, $W\in C({\Bbb R};F)$. By Definition \ref{Definition1.5} (a) 
\begin{eqnarray*} 
\left.\frac{d^\pm}{d\sigma}\right|_{\sigma =0}\left\langle H_i,dW^\sigma 
\right\rangle_{L^2}=\left\langle H_i,d\dot{g}(0,W)\right\rangle_H 
\end{eqnarray*}
because of the particular structure of $H_i$, $W^\sigma=W+g(\sigma,W)$, and 
$(W^\rho)_{\rho\in {\Bbb R}}\in {\cal F}_{p,1}(W)$. Furthermore,  
\begin{eqnarray*} 
\left.\frac{d^\pm}{d\sigma}\right|_{\sigma =0}\left\langle e_j,W^\sigma_0 
\right\rangle_F=\left\langle e_j,\dot{g}(0,W)(0)\right\rangle_F\, .  
\end{eqnarray*}

Let $f_i$ denote the first order derivative of $f$ relative to $i\in I(m,r 
)$ and let $f_j$ be the first order derivative of $f$ relative to $j\in\{1, 
\ldots ,n\cdot d\}$. The last two relations imply $\vp\in E_{1,q}$ since 
\begin{eqnarray*}
&&\hspace{-.5cm}\left.\frac{d^\pm}{d\sigma}\right|_{\sigma =0}\vp(W^\sigma)= 
\left.\frac{d^\pm}{d\sigma}\right|_{\sigma =0}f\left(\langle{\rm H},W^\sigma 
\rangle,\langle {\rm e},W^\sigma\rangle\right) \\ 
&&\hspace{.5cm}=\sum_{i\in I(m,r)}f_i\left(\langle {\rm H},W\rangle ,\langle 
{\rm e},W\rangle\right)\cdot\left.\frac{d^\pm}{d\sigma}\right|_{\sigma =0} 
\left\langle H_i,dW^\sigma\right\rangle_{L^2} \\ 
&&\hspace{1.0cm}+\sum_{j=1}^{n\cdot d}f_j\left(\langle {\rm H},W\rangle , 
\langle {\rm e},W\rangle\right)\cdot\left.\frac{d^\pm}{d\sigma}\right|_{ 
\sigma =0}\left\langle e_j,W^\sigma_0\right\rangle_F \\ 
&&\hspace{.5cm}=\left\langle D\vp(W),j^{-1}\dot{g}(0,W)\right\rangle_H\quad 
Q^{(m,r)}_{\sbnu}\mbox{-a.e.}\ \mbox{\rm for all}\ (W^\rho)_{\rho\in {\Bbb 
R}}\in {\cal F}_{p,1}(W)\, . 
\end{eqnarray*}
Similarly, $\psi\in E^u_{1,q}$ if $\psi\in {\cal C}(\Omega^u)$. 

\subsection{Main result}

The main purpose of the paper is to specify a class of processes $X\equiv X 
(W)$ and a class of distributions $\bnu$ for which the Radon-Nikodym densities  
\begin{eqnarray*} 
\omega_{-t}(W):=\frac{Q_{\sbnu}(dW^{-t})}{Q_{\sbnu}(dW)}\quad\mbox{\rm and} 
\quad\rho_{-t}(X):=\frac{P_{\sbmu}(dX_{\cdot -t})}{P_{\sbmu}(dX)}\, , \quad 
t\in {\Bbb R}, 
\end{eqnarray*} 
exist and to give a representation of these densities. Let $E_{\sbnu}$ denote 
the expectation with respect to the measure $Q_{\sbnu}$ and let $|\cdot |_i$ 
be the absolute value of the $i$-th coordinate. 
\begin{theorem}\label{Theorem1.11}
Assume (1)-(5). (a) For $t\in {\Bbb R}$ the density $\omega_{-t}$ exists 
and we have $Q_{\sbnu}$-a.e. 
\begin{eqnarray}\label{1.1} 
\omega_{-t}(W)&&\hspace{-.5cm}=\frac{m\left(X_{-t}-Y_{-t}\vphantom{l^1} 
\right)}{m(W_0)}\cdot\prod_{i=1}^{n\cdot d}\left|{\sf e}+\nabla_{d,W_0} 
A_{-t}-\nabla_{d,W_0}Y_{-t}\vphantom{\dot{f}}\right|_i\, .  
\end{eqnarray} 
If, in addition, (4'') is satisfied then with $\D I_{-t}:=\left((-t) 
\wedge\tau_{-1},(-t)\vee\tau_{-1}\vphantom{l^1}\right]$, $b=1$ if 
$\tau_{-1}<-t$, and $b=-1$ if $\tau_{-1}\ge -t$ we have the 
representation 
\begin{eqnarray*} 
\omega_{-t}(W)&&\hspace{-.5cm}=\frac{m(X_{-t}-Y_{-t})}{m(W_0)}\times 
\nonumber \\ 
&&\hspace{0cm}\times\prod_{i=1}^{n\cdot d}\left|{\sf e}+\nabla_{d,W_0} 
\left(A^1_{-t}-Y^1_{-t}\right)+b\cdot\sum_{\{k:\tau_k\circ u\in I_{-t} 
\}}\nabla_{d,W_0}\Delta\left(A_{\tau_k}-Y_{\tau_k}\vphantom{l^1}\right) 
\right|_i\, . 
\end{eqnarray*} 
(b) For $t\in {\Bbb R}$ the density $\rho_{-t}$ exists and we have 
$P_{\sbmu}$-a.e. 
\begin{eqnarray}\label{1.2}
\rho_{-t}(X)&&\hspace{-.5cm}=\frac{\D m\left(X_{-t}-Y_{-t}\circ u^{ 
-1}\vphantom{l^1}\right)}{\D m\left(u^{-1}_0\vphantom{l^1}\right)} 
\cdot\vphantom{\left\{\sum_{\tau_k}\right\}}\prod_{i=1}^{n\cdot d} 
\left|\nabla_{d,W_0}X_{-t}-\nabla_{d,W_0}Y_{-t}\circ u^{-1}\vphantom{ 
\dot{f}}\right|_i\, . 
\end{eqnarray} 
\end{theorem} 

If $A^2=0$, i. e., if $X$, $A$, and $Y$ are continuous processes, 
then the conditions under which we have (\ref{1.1}) and (\ref{1.2}) 
simplify as follows. 
\begin{theorem}\label{Theorem1.12}
Assume the following. 
\begin{itemize} 
\item[(i)] $A$ {\it has on $\bigcup_m\{\pi_mW:W\in\Omega\}$ a local 
spatial gradient}. That is, for all $W\in\{\pi_{m}V:V\in\Omega\}$ with 
$W_0=0$, $x\in F$, and $s\in {\Bbb R}$, the gradient $\nabla_x A_s(W+x 
\1)$ exists such that we have (1) (iv). Furthermore, suppose (1'). 
\item[(ii)] On $\{\pi_mW:W\in\Omega\}$, $A$ possesses a Radon-Nikodym 
derivative with cadlag version $A'$. Furthermore, we have (2) (ii). 
\item[(iii)] We have (3) (i) and (iii) (without jump times) as well 
as $Y_s\in D_{q,1}(Q^{(m,r)}_{\sbnu})$, $s\in {\Bbb R}$, for $q$ 
defined in (1'). 
\item[(iv)] $\bigcup_{m}\{\pi_mV:V\in\Omega\}\ni W_n\stack{n\to\infty} 
{\lra}W\in\Omega$ relative to the topology of coordinate wise uniform 
convergence on compact subsets of ${\Bbb R}$ implies  
\begin{eqnarray*}
A_v(W_n)\stack{n\to\infty}{\lra}A_v(W)\quad\mbox{\rm and}\quad Y_v(W_n) 
\stack{n\to\infty}{\lra}Y_v(W) 
\end{eqnarray*} 
and 
\begin{eqnarray*}
\nabla_{W_0}A_v(W_n)-\nabla_{W_0}Y_v(W_n)\stack{n\to\infty}{\lra} 
\nabla_{W_0}A_v(W)-\nabla_{W_0}Y_v(W)
\end{eqnarray*} 
for all $v\in S$ and any finite subinterval $S\subset{\Bbb R}$. Also 
suppose (4) (ii). 
\item[(v)] For all $r\in {\Bbb N}$, all $n\in {\Bbb N}$, and all $s 
\in {\cal R}$, $A_s(W)$ and $Y_s(W)$ are two times Fr\' echet 
differentiable on $W\in jH$. Furthermore, for all $W\in jH$, the 
gradients $\nabla_GA_s(W)\equiv\left(\nabla_G\right)_rA_s(W)$ as well 
as $\nabla_GY_s(W)\equiv\left(\nabla_G\right)_rY_s(W)$ are continuous 
with respect to the variable $r\in {\cal R}$. 
\end{itemize} 
Then we have (\ref{1.1}) and (\ref{1.2}). 
\end{theorem} 

We stress that the technical conditions on the density $m$ formulated in 
the beginning of Subsection 1.2 are crucial in the proofs of Theorem 
\ref{Theorem1.11} and Theorem \ref{Theorem1.12}. A significant relaxation 
of the conditions on $m$ is possible once we have determined the 
representations (\ref{1.1}) and (\ref{1.2}) for $\omega_{-t}(W)$ and 
$\rho_{-t}(X)$ under the conditions of Subsection 1.2. This two-step 
procedure is the reason for presenting the relaxation in a separate 
corollary. 
\begin{corollary}\label{Corollary1.13} 
Suppose (\ref{1.1}) for some $\bar{m}$ satisfying the conditions on the 
Lebesgue density of $W_0$ formulated in the beginning of Subsection 1.2. 
Let $m$ be a strictly positive, everywhere on $D^n$ defined real function 
which is on $D^n$ integrable respect to $\lambda_F$. Suppose also $\int_{ 
D^n}m\, dx=1$ and let $Q_{\sbnu}$ be defined accordingly. Then, for $Q_{ 
\sbnu}$ and $P_{\sbmu}:=Q_{\sbnu}\circ u^{-1}$, relations (\ref{1.1}) and 
(\ref{1.2}) hold true.
\end{corollary} 
Proof. Let $F$ be a bounded measurable function on $\Omega$. Set $\bar 
{\bnu}:=\bar{m}\, dx$ and let $Q_{\bar{\sbnu}}$ be defined accordingly. 
Furthermore, let $\bar{\omega}$ be the expression defined in (\ref{1.1}) 
with $m$ replaced by $\bar{m}$. We have 
\begin{eqnarray*} 
\int F\, dQ_{\sbnu}&&\hspace{-.5cm}=\int F(W)\cdot\frac{m(W_0)}{\bar{m} 
(W_0)}\, dQ_{\bar{\sbnu}} \\ 
&&\hspace{-.5cm}=\int F\left(W^{-t}\right)\cdot\frac{m(W^{-t}_0)}{\bar{m 
}(W^{-t}_0)}\cdot\bar{\omega}_{-t}(W)\, dQ_{\bar{\sbnu}} \\ 
&&\hspace{-.5cm}=\int F\left(W^{-t}\right)\cdot\frac{m\left(X_{-t}-Y_{-t 
}\right)}{\bar{m}\left(X_{-t}-Y_{-t}\right)}\cdot\bar{\omega}_{-t}(W)\cdot 
\frac{\bar{m}(W_0)}{m(W_0)}\, dQ_{\sbnu} \\ 
&&\hspace{-.5cm}=\int F\left(W^{-t}\right)\cdot\omega_{-t}(W)\, dQ_{\sbnu} 
\, .  
\end{eqnarray*}
\qed 
\medskip

We continue with two corollaries of Theorem \ref{Theorem1.11} and 
Theorem \ref{Theorem1.12} for which we pay particular attention to 
condition (1) (iv) and (4) (i) and, respectively, (i) and (iv) of 
Theorem \ref{Theorem1.12}. The first one is the following disintegration 
formula for the measure $P_{\sbmu}$. For this let $\xi_X:=\{X_{\cdot -t} 
:t\in {\Bbb R}\}$, $u^{-1}(X)\in\Omega$, and ${\cal X}:=\{\xi_X:u^{-1}(X) 
\in\Omega\}$. For $A\in {\cal F}$ set ${\cal A}:=\{\xi_X:X=u(W),\, W\in 
A\}$. Endow ${\cal X}$ with the sub-$\sigma$ algebra ${\cal F}_{\cal X}$ 
of $u\circ {\cal F}$ generated by all such sets ${\cal A}$. Furthermore, 
define the measure $\Gamma$ on $\left({\cal X},{\cal F}_{\cal X}\right)$ 
by $\Gamma({\cal A}):=P_{\sbmu}\left(\{X:\mbox{\rm there exists $\xi\in 
{\cal A}$ such that }X\in\xi\}\right)$. 
\begin{corollary}\label{Corollary1.14} 
Suppose the hypotheses of Theorem \ref{Theorem1.11} or Theorem 
\ref{Theorem1.12} and let $m$ be as in Corollary \ref{Corollary1.13}. 
For $\Gamma$-a.e. $\xi\in {\cal X}$ there exists a measure $\gamma_\xi$ 
on $\xi$ endowed with the trace-$\sigma$ algebra ${\cal F}_{\xi}$ of $u 
\circ {\cal F}$ to $\xi$ such that for $A\in u \circ {\cal F}$ and $A_\xi: 
=\xi\cap A$ 
\begin{eqnarray*} 
P_{\sbmu}(A)=\int_{\xi\in {\cal X}}\gamma_\xi(A_\xi)\, \Gamma(d\xi)\, . 
\end{eqnarray*}
For $\Gamma$-a.e. $\xi\in {\cal X}$ and $X\in\xi$, we have 
\begin{eqnarray*} 
\gamma_\xi\left(\left\{X_{\cdot +\tau_k(X)}:k\in {\Bbb Z}\setminus\{0\} 
\right\}\right)=0\, . 
\end{eqnarray*}
Furthermore, for $\Gamma$-a.e. $\xi\in {\cal X}$, $X\in\xi$, $-t\in {\Bbb 
R}\setminus\{\tau_k(X):k\in {\Bbb Z}\setminus\{0\}\}$, and $X^1=X_{\cdot 
-t}$ we have 
\begin{eqnarray*} 
\frac{\gamma_\xi(dX^1)}{\gamma_\xi(dX)}=\rho_{-t}(X) 
\end{eqnarray*}
where $\rho$ is given by (\ref{1.2}). 
\end{corollary} 
An immediate but useful consequence of Corollary \ref{Corollary1.14} as 
well as Theorem \ref{Theorem1.11} or Theorem \ref{Theorem1.12} is the 
following one. 
\begin{corollary}\label{Corollary1.15} 
Suppose the hypotheses of Theorem \ref{Theorem1.11} or Theorem 
\ref{Theorem1.12} and let $m$ be as in Corollary \ref{Corollary1.13}. Let 
${\cal A}$ be an index set. Assume that 
\begin{itemize}
\item[(i)] for every $\alpha\in {\cal A}$, there is a set $\Omega_\alpha 
\in {\cal F}$ such that $\bigcup_{\alpha\in {\cal A}}\Omega_\alpha=\Omega$ 
and $\alpha_1\neq\alpha_2$ implies $\Omega_{\alpha_1}\cap\Omega_{\alpha_2} 
=\emptyset$ and 
\item[(ii)] $u^{-1}\circ X\in\Omega_\alpha$ implies $u^{-1}\circ X_{\cdot 
-t}\in\Omega_\alpha$ for all $t\in {\Bbb R}$, $\alpha\in {\cal A}$. 
\end{itemize} 
Let $t(\alpha):{\cal A} \to {\Bbb R}$ be a bounded map and define the random 
time $\tau^\alpha\equiv\tau^\alpha(X)$ on $(\Omega,{\cal F})$ by $\tau^\alpha 
\circ u:=t(\alpha)$ on $\Omega_\alpha$, $\alpha\in {\cal A}$. Then we have 

\begin{eqnarray*}
P_\mu(-\tau^\alpha\ \mbox{\rm is a jump time for }X)=0 
\end{eqnarray*}
and 
\begin{eqnarray*} 
\frac{P_{\sbmu}(dX_{\cdot -\tau^\alpha})}{P_{\sbmu}(dX)}=\rho_{-\tau^\alpha 
}(X) \quad P_{\sbmu}\mbox{\rm -a.e. } 
\end{eqnarray*}
where $\rho$ is given by (\ref{1.2}). 
\end{corollary} 

The restriction to bounded domains $D$ is due to techniques used in 
the proofs. A possibility to establish absolute continuity under time 
shift of trajectories of processes with unbounded state space comes 
with the following proposition. For this let $X\equiv X(W)=W+A$ be a 
stochastic process with $A_0=0$ and state space $\check{D}^n$ where 
$\check{D}\subseteq {\Bbb R}^d$ is a possibly unbounded domain. Let 
$D\subset {\Bbb R}^d$ be a bounded domain in the above setting and 
$g:\check{D}^n\to D^n$ be a bijection. Define 
\begin{eqnarray*} 
\hat{W}:=g(W_0)\1+W-W_0\1\, , \quad W\in\Omega\, ,\ W_0\in\check{D}^n, 
\end{eqnarray*}
and 
\begin{eqnarray*} 
\check{W}:=g^{-1}(W_0)\1+W-W_0\1\, ,\quad W\in\Omega\, ,\ W_0\in D^n.
\end{eqnarray*}
Furthermore, let $\hat{A}(W):=g\circ X_\cdot(\check{W})-W$ for $W\in 
\Omega$ with $W_0\in D^n$. This yields 
\begin{eqnarray*} 
\hat{X}(W):=g\circ X_\cdot(\check{W})=W+\hat{A}(W)\, , \quad W\in 
\Omega\, ,\ W_0\in D^n. 
\end{eqnarray*}
The formulation and the proof of the following proposition does not 
take into consideration the distribution of $W_0$. 
\begin{proposition}\label{Proposition1.16}
Assume the following. 
\begin{itemize} 
\item[(i)] $X=W+A$ where $W$ is a two-sided Brownian motion with the 
law of $W_0$ supported by $\overline{\check{D}^n}$ and $X:\Omega\to 
B_{b;{\rm loc}} ({\Bbb R};F)$ is a measurable map. 
\item[(ii)] $A_0=0$. 
\item[(iii)] With $W^v:=W_{\cdot +v}+A_v(W)\1$ we have $X(W^v)=X_{ 
\cdot +v}(W)$, $v\in {\Bbb R}$. 
\end{itemize} 
Then $\hat{A}_0(W)=0$ for all $W\in\Omega$ such that $W_0\in D^n$ and 
with $W\, \hat{}\, ^v:=W_{\cdot +v}+\hat{A}_v(W)\1$ we have 
\begin{eqnarray*} 
\hat{X}(W\, \hat{}\, ^v)=\hat{X}_{\cdot +v}(W)\, ,\quad v\in {\Bbb R}. 
\end{eqnarray*}
\end{proposition} 
Proof. Let $W\in\Omega$ with $W_0\in D^n$. By (i) and (ii), we have 
$\hat{A}_0(W)=g(X_0(\check{W}))-W_0=g(\check{W}_0+A_0(\check{W}))-W_0 
=g(\check{W}_0)-W_0=0$. In addition, 
\begin{eqnarray*} 
\raisebox{4mm}{\rotatebox{180}{\raisebox{-3mm}{\hspace{-2.5mm}$ 
\widehat{\hspace{5mm}}$}}}\hspace{-3.5mm}W\, \hat{}\, ^v_{\hspace
{-2mm}0}=g^{-1}\left(W\, \hat{}\, ^v_{\hspace{-2mm}0}\right)=g^{-1} 
\left(W_v+\hat{A}_v(W)\right)=g^{-1}\left(\hat{X}_v(W)\right)=X_v
(\check{W})=\check{W}_0^v
\end{eqnarray*}
which implies
\begin{eqnarray*} 
\raisebox{4mm}{\rotatebox{180}{\raisebox{-3mm}{\hspace{-2.5mm}$ 
\widehat{\hspace{5mm}}$}}}\hspace{-3.5mm}W\, \hat{}\, ^v=\left( 
\check{W}\right)^v\equiv\check{W}^v\, ,\quad v\in {\Bbb R}. 
\end{eqnarray*}
With (iii) we obtain 
\begin{eqnarray*} 
\hat{X}(W\, \hat{}\, ^v)=g\circ X_\cdot\left(\vphantom{\dot{f}} 
\right.\raisebox{4mm}{\rotatebox{180}{\raisebox{-3mm}{\hspace{- 
2.5mm}$\widehat{\hspace{5mm}}$}}}\hspace{-3.5mm}W\, \hat{}\, ^v 
\left.\vphantom{\dot{f}}\right)=g\circ X_\cdot\left(\check{W}^v 
\right)=g\circ X_{\cdot +v}\left(\check{W}\right)=\hat{X}_{\cdot 
+v}(W)\, ,\quad v\in {\Bbb R}. 
\end{eqnarray*}
\qed
\medskip 

In fact, based on (i) and (iii) of Proposition \ref{Proposition1.16}, 
we have verified condition (3) (i) for $\hat{X}$. In words, this says 
that for $A=0$ temporal homogeneity is invariant under bijective space 
transformation. For the absolute continuity under time shift of 
trajectories for $\hat{X}$ check now the remaining conditions of 
Theorem \ref{Theorem1.11} or Theorem \ref{Theorem1.12}, respective 
Corollary \ref{Corollary1.13} relative to $\hat{X}$. An immediate 
consequence is absolute continuity under time shift of trajectories for 
$X$. 
 
\section{Flows and Logarithmic Derivative Relative to $X$ under Orthogonal 
Projection} 
\setcounter{equation}{0}

In this section we are primarily interested in the analysis of the 
process $X=u(W)$ where we just focus on conditions (1') and 
(4') of Subsection 1.2. In particular, we do not refer to $X=W+A$. 
In Subsection 2.1, we will introduce elements of the calculus on 
orthogonal projections of $W$. Subsections 2.2 and 2.3 will be 
dedicated to the analysis of related flows and logarithmic derivatives. 
\medskip 

In order to prove Theorem \ref{Theorem1.11} we will apply the following 
approximation theorem. 
\begin{theorem}\label{Theorem2.1} (Specification of Theorem 3.1 in 
\cite{GS66}, Theorem 3 in \cite{KP00})
Suppose that $M$ is a separable metric space, $\mu$ is a probability 
measure defined on the Borel $\sigma$-algebra over $M$, and $f_n:M\to 
M$, $n\in {\Bbb N}$, is a sequence of measurable maps. Assume that the 
following conditions are satisfied: 
\begin{itemize} 
\item[(i)] For every $n\in {\Bbb N}$ the measure $\mu\circ f_n^{-1}$ is 
absolutely continuous with respect to the measure $\mu$. 
\item[(ii)] The sequence of densities $d\mu\circ f_n^{-1}/d\mu$, $n\in 
{\Bbb N}$, is uniformly integrable. 
\item[(iii)] $f_n\stack{n\to\infty}{\lra}f$ in the measure $\mu$ for some 
$f:M\to M$. 
\end{itemize}
Then $\mu\circ f^{-1}$ is absolutely continuous with respect to the measure 
$\mu$. If $d\mu\circ f_n^{-1}/d\mu\stack{n\to\infty}{\lra}p$ in the measure 
$\mu$ then $p=d\mu\circ f^{-1}/d\mu$. 
\end{theorem}

\subsection{Elements of the analysis on orthogonal projections of Brownian 
paths}

Starting point is the L\'evy-Ciesielsky construction of $W$. Accordingly, 
there are independent $N(0,1)$-distributed random variables $\xi_1 ,\xi_2 , 
\ldots $ such that $W_s=W_0+\sum_{i=1}^\infty\xi_i\cdot\int_0^s H_i(u)\, 
du$, $s\in {\Bbb R}$, where the sum converges uniformly in $s$ on finite 
subintervals of ${\Bbb R}$, $Q_{\sbnu}$-a.e. It holds that $\xi_i=\langle 
H_i,dW\rangle_{L^2}$ so that we get the L\'evy-Ciesielsky representation 
\begin{eqnarray*}
W_s=W_0+\sum_{i=1}^\infty\langle H_i,dW\rangle_{L^2}\cdot\int_0^s H_i(u) 
\, du\, , \quad s\in {\Bbb R}. 
\end{eqnarray*}
Under the measure $Q^{(m,r)}_{\sbnu}$ we shall consider  
\begin{eqnarray}\label{2.1}
W_s=W_0+\sum_{i\in I(m,r)}\langle H_i,dW\rangle_{L^2}\cdot\int_0^sH_i 
(u)\, du\, , \quad s\in {\Bbb R}, 
\end{eqnarray}
and $j^{-1}W=q_{m,r}j^{-1}W=\sum_{i\in I(m,r)}\langle H_i,dW,W_0\rangle_{ 
L^2}\cdot H_i \in H$. In the proof of Theorem \ref{Theorem1.11} in Section 
3 we will also consider similar projections under the measure $Q^{(m)}_{ 
\sbnu}$. 
\medskip 

{\bf The integral $\int\langle V, d\dot{W}\rangle_F$ under the measure 
$Q^{(m,r)}_{\sbnu}$. } For $-rt\le u<v\le (r-1)t$ and $V:{\cal R}\to F$ 
cadlag define 
\begin{eqnarray*} 
&&\hspace{-.5cm}\int_{w=u}^v\langle V,dH_i\rangle_F \\ 
&&\hspace{.5cm}:=\sum_{w\in (u,v)}\frac12\langle V_{w-}+V_w,\Delta H_i(w) 
\rangle_F+\frac12\langle V_u,\Delta H_i(u)\rangle_F +\frac12\langle V_{v-} 
,\Delta H_i (v)\rangle_F\, , 
\end{eqnarray*}
$i\in I(m,r)$, and similarly the integral $\int_{w=u}^v\langle V, dH_i  
\rangle_{F\to F}$. Furthermore, introduce 
\begin{eqnarray}\label{2.2} 
\int_{w=u}^v \langle V, d\dot{W}\rangle_F:=\sum_{i\in I(m,r)}\langle H_i 
,dW\rangle_{L^2}\cdot\int_{w=u}^v\langle V, dH_i \rangle_F\quad Q^{(m,r) 
}_{\sbnu}\mbox{\rm -a.e.}
\end{eqnarray}
For $i,j\in I(m,r)$, we have 
\begin{eqnarray*} 
\langle H_i,dH_j \rangle_{L^2}+\langle H_j,dH_i\rangle_{L^2}=0\, . 
\end{eqnarray*}
It follows that 
\begin{eqnarray}\label{2.3} 
\sum_{i\in I(m,r)}\langle H_i,dW\rangle_{L^2}\cdot\langle H_i,d\dot{W} 
\rangle_{L^2}&&\hspace{-.5cm}=\sum_{i,j\in I(m,r)}\langle H_i,dW 
\rangle_{L^2}\langle H_j,dW\rangle_{L^2}\cdot\langle H_i,dH_j\rangle_{L^2} 
\nonumber \\ 
&&\hspace{-.5cm}=0\quad Q^{(m,r)}_{\sbnu}\mbox{\rm -a.e.}\vphantom{ 
\sum^0_0} 
\end{eqnarray}

{\bf The process $X$ and the map $u$ under the measure $Q^{(m,r)}_{\sbnu}$. 
} Let 
\begin{eqnarray*} 
\kappa_i&&\hspace{-.5cm}:=\int_0^\cdot H_i (v)\, dv
\, , \quad i\in I(m,r), \\ 
\lambda_j&&\hspace{-.5cm}:=e_j\1(\cdot) 
\, , \quad j\in \{1,\ldots ,n\cdot d\}. 
\end{eqnarray*} 
\begin{lemma}\label{Lemma2.2} 
(a) Let $Q$ be the number used in the definition of the density $m$ 
and $1/q+1/Q\le 1$. Let $\vp\in D_{q,1}$ with respect to the measure 
$Q^{(m,r)}_{\sbnu}$ and $\vp_{n'}\in {\cal C}$, $n'\in {\Bbb N}$, be 
a sequence with $\vp_{n'}\stack{n'\to\infty}{\lra}\vp$ in $D_{q,1}$. 
Choosing a subsequence if necessary $D\vp$ is $Q^{(m,r)}_{\sbnu}$-a.e. 
a finite dimensional object for which 
\begin{eqnarray*}
D\vp(W)&&\hspace{-.5cm}=\lim_{n'\to\infty}D\vp_{n'}(W) \\
&&\hspace{-.5cm}=\lim_{n'\to\infty}\sum_{i\in I(m,r)}\frac{\partial 
\vp_{n'}(W)}{\partial\kappa_i}\cdot (H_i,0)+\lim_{n'\to\infty}\sum_{j 
=1}^{n\cdot d}\frac{\partial\vp_{n'}(W)}{\partial \lambda_j}\cdot (0, 
e_j) \\ 
&&\hspace{-.5cm}=\sum_{i\in I(m,r)}\left\langle D\vp,(H_i,0)\right 
\rangle_H\cdot (H_i,0)+\sum_{j=1}^{n\cdot d}\left\langle D\vp,(0,e_j) 
\right\rangle_H\cdot (0,e_j)\, . 
\end{eqnarray*} 
In particular, $\vp=\vp\circ\pi_{m,r}$ and $D\vp=D\vp\circ\pi_{m,r}$ 
$Q^{(m,r)}_{\sbnu}$-a.e. 
\medskip 

\nid 
(b)  Assume condition (4') of Subsection 1.2. We have 
\begin{eqnarray*}
\lim_{m\to\infty}X_t(\pi_mW)=X_t\, , \quad t\in {\Bbb R}\setminus\{ 
\tau_i(X):i=\pm 1,\pm 2,\ldots \}\quad Q_{\sbnu}\mbox{\rm -a.e.}
\end{eqnarray*} 
(c) Let $\xi\equiv\xi(W)$ be an $\hat{F}$-valued random variable with 
$\hat{F}$ being a metric space. Assume that $W_n\stack{n\to\infty}{\lra} 
W$ relative to the topology of coordinate wise uniform convergence on 
compact subsets of ${\Bbb R}$ implies $\xi(W_n)\stack{n\to\infty}{\lra} 
\xi(W)$. Let $\psi\in C_b(\hat{F})$ or $\psi\in C(\hat{F})$ and $\sup_{ 
r\in {\Bbb N}}|\psi(\xi(\pi_{m,r}W))|\in L^1(\Omega,Q^{(m)}_{\sbnu})$. 
Then 
\begin{eqnarray*}
\lim_{r\to\infty}\int\psi(\xi(W))\, Q^{(m,r)}_{\sbnu}(dW)=\int\psi(\xi 
(W))\, Q^{(m)}_{\sbnu}(dW)\, . 
\end{eqnarray*} 
\end{lemma} 
Proof. The representation of $D\vp$ in part (a) follows from Definition 
\ref{Definition1.1}. By the L\'evy-Ciesielsky construction of $W$ and the 
above definition of $D\vp$, $\vp$ and $D\vp$ are independent of $\langle 
H_i,dW\rangle_{L^2}$, $i\not\in I(m,r)$. In other words, $\vp=\vp\circ 
\pi_{m,r}$ and $D\vp=D\vp\circ\pi_{m,r}$ $Q^{(m,r)}_{\sbnu}$-a.e. 

Part (b) is an immediate consequence of (4'), cf. Subsection 1.2. For 
part (c) it is sufficient to recall that $Q^{(m,r)}_{\sbnu}=Q_{\sbnu} 
\circ\pi_{m,r}^{-1}$, relation (\ref{2.1}), the corresponding 
representation of $W$ under the measure $Q^{(m)}_{\sbnu}=Q_{\sbnu}\circ 
\pi_{m}^{-1}$, and condition (4'). In fact, it follows now that 
\begin{eqnarray*}
\lim_{r\to\infty}\int\psi(\xi(W))\, Q^{(m,r)}_{\sbnu}(dW)&&\hspace{-.5cm} 
=\lim_{r\to\infty}\int\psi(\xi(\pi_{m,r}W))\, Q^{(m)}_{\sbnu}(dW) \\ 
&&\hspace{-.5cm}=\int\psi(\xi(W))\, Q^{(m)}_{\sbnu}(dW)\, . 
\end{eqnarray*}  
\qed 
\bigskip 

For $s\in {\cal R}=[-rt,(r-1)t]$, we have $Du_s\in H$ by assumption 
(1') of Section 1 with $Du_s=((Du_s)_1,(Du_s)_2)$ and $(Du_s)_1\in 
L^2_{\rm loc}({\Bbb R};F)$ and $(Du_s)_2\in F$ $Q_{\sbnu}^{(m,r)}$-a.e. 
However, $Q_{\sbnu}^{(m,r)}$-a.e., $(Du_s)_1$ is, according to Lemma 
\ref{Lemma2.2}, supported by a subset of ${\cal R}$. So we will 
assume $Q_{\sbnu}^{(m,r)}$-a.e., $(Du_s)_1\in L^2({\cal R};F)$ and 
$(Du_s)_2\in F$. 
\bigskip

Let ${\cal M}^{f,s}({\cal R};F)$ denote the set of all $F$-valued finite 
signed measures on $({\cal R},{\cal B}({\cal R}))$. 
\begin{lemma}\label{Lemma2.3} 
Let $X=u(W)$ be a map $C({\Bbb R};F)\to B_{\rm b;loc}({\Bbb R};F)$. Assume 
(1') of Section 1 and let $\gamma:C({\cal R};F)\to {\cal M}^{f,s}({\cal R}; 
F)$. Assume furthermore that 
\begin{itemize} 
\item[(i)] for every $i\in I(m,r)$, $\langle H_i,d\gamma\rangle_{L^2}$ 
is differentiable in direction of $\kappa_i$, that is 
\begin{eqnarray*}
\langle H_i,d\gamma\rangle_{L^2}\in L^p(\Omega ,Q_{\sbnu}^{(m,r)}) 
\quad\mbox{\rm and}\quad\frac{\partial\langle H_i,d\gamma\rangle_{ 
L^2}}{\partial\kappa_i}\ \mbox{\rm exists in } L^p(\Omega ,Q_{\sbnu}^{ 
(m,r)}) 
\end{eqnarray*} 
where $1/p+1/q=1$ and $q$ is the constant of assumption (1') of Section 1 
and 
\item[(ii)] for every $i'\in \{1,\ldots ,n\cdot d\}$, $\langle e_{i'},d 
\gamma (0)\rangle_F$ is differentiable in direction of $e_{i'}$, that is 
\begin{eqnarray*}
\langle e_{i'},d\gamma(0)\rangle_F\in L^p(\Omega ,Q_{\sbnu}^{(m,r)}) 
\quad\mbox{\rm and}\quad\frac{\partial\langle e_{i'},d\gamma(0)\rangle_F} 
{\partial e_{i'}}\ \mbox{\rm exists in } L^p(\Omega ,Q_{\sbnu}^{(m,r)})\, . 
\end{eqnarray*} 
\item[(iii)]  The linear span of $\{X_s:s\in {\cal R}\}$ is dense in $L^q 
(\Omega,Q^{(m,r)}_{\sbnu};F)$.
\end{itemize} 
Then $\langle (Du_s)_1,d\gamma\rangle_{L^2\to F}+\langle (Du_s)_2,d\gamma 
(0)\rangle_{F\to F}=0$ for all $s\in {\cal R}$ $Q_{\sbnu}^{(m,r)}$-a.e. 
implies that 
\begin{eqnarray*}
\langle H_i,d\gamma\rangle_{L^2}=\langle e_{i'},d\gamma (0) 
\rangle_F=0\, , \quad i\in I(m,r)\, , \ i'\in \{1,\ldots ,n\cdot d\}. 
\end{eqnarray*} 
\end{lemma} 
\medskip 

\nid 
Proof. For $i\in I(m,r)$, let $\beta^{(i,0)}$ denote the logarithmic 
derivative relative to the measure $\langle H_i,d\gamma 
\rangle_{L^2}\, dQ^{(m,r)}_{\sbnu}$. Moreover, for $i'\in \{1,\ldots 
,n\cdot d\}$, let $\beta^{(0,i')}$ denote the logarithmic derivative 
relative to the measure $\langle e_{i'},d\gamma (0)\rangle_F\, dQ^{ 
(m,r)}_{\sbnu}$. Both logarithmic derivatives exist by hypothesis. 
By Lemma \ref{Lemma2.2} (a) we may write 
\begin{eqnarray*} 
&&\hspace{-1cm}\int\langle\left((Du_s)_1,d\gamma\rangle_{L^2\to F} 
+\langle (Du_s)_2,d\gamma (0)\rangle_{F\to F}\right)\, d Q_{\sbnu}^{ 
(m,r)} \\ 
&&\hspace{.5cm}=\sum_{i\in I(m,r)}\int\langle Du_s,(H_i ,0)\rangle_{ 
H\to F}\cdot\langle H_i,d\gamma\rangle_{L^2}\, dQ_{\sbnu}^{(m,r)} 
 \\ 
&&\hspace{1cm}+\sum_{i'\in \{1,\ldots ,n\cdot d\}}\int\langle Du_s, 
(0,e_{i'})\rangle_{H\to F}\cdot\langle e_{i'},d\gamma (0)\rangle_F\, 
dQ_{\sbnu}^{(m,r)} \\ 
&&\hspace{.5cm}=-\sum_{i\in I(m,r)\atop i'\in \{1,\ldots ,n\cdot d\}} 
\int u_s\cdot\left(\beta^{(i)}_{j(H_i ,0)}\cdot\langle H_i, 
d\gamma\rangle_{L^2}+\beta^{(i')}_{j(0,e_{i'})}\cdot\langle e_{i'}, 
d\gamma (0)\rangle_F\right)\, dQ_{\sbnu}^{(m,r)}\, . 
\end{eqnarray*} 
By $\langle (Du_s)_1,d\gamma\rangle_{L^2\to F}+\langle (Du_s)_2,d 
\gamma (0)\rangle_{F´\to F}=0$ for all $s\in {\cal R}$ $Q_{\sbnu}^{ 
(m,r)}$-a.e. and assumption (iii), we have 
\begin{eqnarray}\label{2.4} 
\sum_{i\in I(m,r)\atop i'\in \{1,\ldots ,n\cdot d\}}\left(\beta^{(i) 
}_{j(H_i ,0)}\cdot\langle H_i,d\gamma\rangle_{L^2}+\beta^{ 
(i')}_{j(0,e_{i'})}\cdot\langle e_{i'},d\gamma (0)\rangle_F\right)= 
0\quad Q_{\sbnu}^{(m,r)}\mbox{\rm -a.e.}	 
\end{eqnarray} 
In addition, it holds that 
\begin{eqnarray*} 
&&\hspace{-.5cm}\int\left\langle\vphantom{\left(\dot{f}\right)}\left( 
\vphantom{I^I_I}\langle (Du_s)_1,d\gamma\rangle_{L^2\to F}+\langle 
(Du_s)_2,d\gamma (0)\rangle_{F\to F}\right),W_v\right\rangle_{F\to F} 
\, dQ_{\sbnu}^{(m,r)}\vphantom{\left(\int_T\right)}\nonumber \\
&&\hspace{.5cm}+\sum_{i\in I(m,r)\atop i'\in \{1,\ldots ,n\cdot d\}} 
\left(\int\left\langle\int_0^v H_i (w)\, dw,u_s\right\rangle_{F\to F} 
\cdot\langle H_i,d\gamma\rangle_{L^2}\, dQ_{\sbnu}^{(m,r)} 
\right. \\ 
&&\hspace{4.9cm}+\left.\int\left\langle e_{i'},u_s\right\rangle_{F\to 
F}\cdot\langle e_{i'},d\gamma(0)\rangle_F\vphantom{\sum_T}\, dQ_{\sbnu 
}^{(m,r)}\right) \\ 
&&\hspace{.0cm}=\sum_{i\in I(m,r)\atop i'\in \{1,\ldots ,n\cdot d\}} 
\left(\int\left\langle\langle Du_s,(H_i ,0)\rangle_{H\to F},W_v\right 
\rangle_{F\to F}\cdot\langle H_i,d\gamma\rangle_{L^2}\, dQ_{\sbnu}^{ 
(m,r)}\vphantom{\sum_{T}}\right. \nonumber \\ 
&&\hspace{2cm}+\left.\int\left\langle\langle Du_s,(0,e_{i'})\rangle_{
H\to F},W_v\right\rangle_{F\to F}\cdot\langle e_{i'},d\gamma (0) 
\rangle_F\, dQ_{\sbnu}^{(m,r)}\vphantom{\sum_{T\atop T}}\right. 
\nonumber \\ 
&&\hspace{2.0cm}\left.+\int\left\langle\langle DW_v,(H_i ,0)\rangle_{ 
H\to F},u_s\right\rangle_{F\to F}\cdot\langle H_i,d\gamma\rangle_{L^2 
}\, dQ_{\sbnu}^{(m,r)}\vphantom{\sum_{T}}\right. \\ 
&&\hspace{2.0cm}\left.+\int\left\langle\langle DW_v,(0,e_{i'}) 
\rangle_{H\to F},u_s\right\rangle_{F\to F}\cdot\langle e_{i'},d\gamma 
(0)\rangle_F\, dQ_{\sbnu}^{(m,r)}\vphantom{\sum_{T}}\right) \\ 
&&\hspace{.0cm}=-\int\langle u_s,W_v\rangle_{F\to F}\cdot\sum_{i\in 
I(m,r)\atop i'\in \{1,\ldots ,n\cdot d\}}\left(\beta^{(i)}_{j(H_i ,0) 
}\cdot\langle H_i,d\gamma\rangle_{L^2}+\beta^{(i')}_{j(0,e_{i'} 
)}\cdot\langle e_{i'},d\gamma (0)\rangle_F\right)\, dQ_{\sbnu}^{(m,r)} 
\, .    
\end{eqnarray*} 
By $\langle (Du_s)_1,d\gamma\rangle_{L^2\to F}+\langle (Du_s)_2,d 
\gamma(0)\rangle_{F\to F}=0$ for all $s\in {\cal R}$ $Q_{\sbnu}^{(m, 
r)}$-a.e. and (\ref{2.4}) we conclude now 
\begin{eqnarray*}
&&\hspace{-.5cm}\sum_{i\in I(m,r)\atop i'\in \{1,\ldots ,n\cdot 
d\}}\int\left(\left\langle\int_0^vH_i (w)\, dw,u_s\right\rangle_{F 
\to F}\cdot\langle H_i,d\gamma\rangle_{L^2}\right. \\ 
&&\hspace{5.5cm}\left.\vphantom{\int_m^n}+\left\langle e_{i'},u_s 
\right\rangle_{F\to F}\langle e_{i'},d\gamma (0)\rangle_F\right)\, 
dQ_{\sbnu}^{(m,r)}=0\, ,    
\end{eqnarray*} 
$s,v\in {\cal R}$, which by assumption (iii) proves the lemma. 
\qed 

\subsection{Two flows associated with $X=u(W)$, substitution law, and 
stochastic integral relative to $X$} 

We may write $X\equiv X(W):=u(W)$. Suppose we are given $h$ with 
$h\circ u:C({\Bbb R};F)\to B_{\rm b;loc}({\Bbb R};F)$ and $f(\rho,W) 
(s):{\Bbb R}\times C({\Bbb R};F)\to C({\Bbb R};F)$ such that 
\begin{eqnarray}\label{2.5}
\int_{\sigma=0}^\rho h\circ u(W+f(\sigma,W ))(\cdot )\, d\sigma 
&&\hspace{-.5cm}=u(W+f(\rho,W))-u(W)\vphantom{\int^0_0} \nonumber 
 \\ 
&&\hspace{-.5cm}=X (W+f(\rho,W))-X (W)\, , \quad\, \rho 
\in {\Bbb R},\vphantom{\int^0_0} \\ 
&&\hspace{-2.1cm}f(0,W)\equiv 0\, . \vphantom{\int^0_0} \nonumber 
\end{eqnarray} 
We notice that since $X=u(W)$ is an injection $\Omega\to\Omega^u$ and 
$W^\rho$, $\rho\in {\Bbb R}$, is a random flow on $\Omega$, because 
of 
\begin{eqnarray*}
X^{\rho+\sigma}&&\hspace{-.5cm}=u(W^{\rho+\sigma})u^{-1}=u(W^\rho\circ 
W^\sigma)u^{-1}\vphantom{\dot{f}} \\ 
&&\hspace{-.5cm}=u(W^\rho)u^{-1}\circ u(W^\sigma)u^{-1}=X^\sigma\circ 
X^\rho\, , \quad \rho ,\sigma\in {\Bbb R},  
\end{eqnarray*} 
$X^\rho$, $\rho\in {\Bbb R}$, is a random flow on $\Omega^u$. In this 
situation, there are two random flows in this equation, $W^\rho:=W+f 
(\rho,W)$ and $X^\rho\equiv X^\rho(W):=u(W+f(\rho,W))$, $\rho\in {\Bbb 
R}$. In Section 3 we will be interested in specifying the second one 
to the formal relation $X^\rho:=X_{\cdot +\rho}$, $\rho\in {\Bbb R}$. 
In order to prepare this we are going to develop elements of the 
related stochastic differential and integral calculus in this 
section. 
\medskip

{\bf Assumptions for the remainder of Section 2. } In this section 
we shall assume that for $Q^{(m,r)}_{\sbnu}$-a.e. $W\in\Omega\equiv 
C({\Bbb R};F)$
\begin{eqnarray*}
f(\rho,W)=f^0(\rho,W)+f(\rho,W)(0)\cdot\1\, , \quad \rho 
\in {\Bbb R},
\end{eqnarray*} 
where $f^0(\rho,W)\in C_a({\Bbb R};F)$ with Radon-Nikodym derivative 
in $L^2_{\rm loc}({\Bbb R};F)$ and $f^0(\rho,W)(0)=0$. Furthermore, 
$f(\rho,\cdot )(0)\equiv f(\rho,W)(0)$, $\rho\in {\Bbb R}$, is a 
stochastic process that we abbreviate by 
\begin{eqnarray*}
f(\rho,W)(0)=B_\rho(W)\, .  
\end{eqnarray*} 
In particular, we shall suppose that $B\equiv B_\rho(W)$ as well as 
$B_{-\cdot\, }\equiv B_{-\rho}(W)$ possess $Q^{(m,r)}_{\sbnu}$-a.e. 
a Radon-Nikodym derivative with respect to $\rho\in {\Bbb R}$ with 
cadlag versions $B'$ and $B_{\, -\cdot\, }'\, $. We recall that $B'$ 
and $B_{\, -\cdot\, }'$ are at the same time $Q^{(m,r)}_{\sbnu}$-a.e. 
the right derivatives of $B$ and $B_{-\cdot\, }\, $, respectively. 
Thus the {\it mixed derivative} 
\begin{eqnarray*}
\dot{B}_\rho:=\frac12\left(\frac{d^-}{d\rho}B_\rho+\frac{d^+}{d\rho} 
B_\rho\right)\, , \quad \rho\in {\Bbb R}, 
\end{eqnarray*} 
is well-defined. 
\medskip

Let $1/p+1/q=1$ and $q$ be given by condition (1') of Section 1. Let 
``$\ \dot{}\ $'' stand further on for differentiation with respect to 
$\rho$ and assume in this and the next subsection the following 
(\ref{2.6}) and (\ref{2.7}). We shall suppose that, under the measure 
$Q^{(m,r)}_{\sbnu}$, $W^\rho$, $\rho\in {\Bbb R}$, is {\it transient} 
in the sense that 
\begin{eqnarray}\label{2.6}
W^\sigma\neq W^\rho\ \mbox{\rm for }\sigma\neq\rho\, .  
\end{eqnarray} 

Let ``$\ {}'\ $'' indicate the right continuous version of the 
Radon-Nikodym derivative of $f^0(\rho,\cdot)(s)$ with respect to the 
Lebesgue measure $ds$. We also shall assume that $\left(f^0(\rho,\cdot) 
\right)'$ is $Q^{(m,r)}_{\sbnu}$-a.e. {\it weakly mixed differentiable} 
in $\rho\in {\Bbb R}$ in the sense that for some $\dot{f}(\rho ,)\in 
I^p$, $\rho\in {\Bbb R}$, and all test elements $k=(g,x)\in H$ with $g 
\in C({\cal R};F)$, $x\in F$, we have 
\begin{eqnarray}\label{2.7}
&&\hspace{-.5cm}\left\langle g,d\dot{f}(\rho ,W)\right\rangle_{L^2} 
+\left\langle x,\dot{f}(\rho ,W)(0)\right\rangle_F \nonumber \\ 
&&\hspace{.5cm}=\frac12\left(\frac{d^-}{d\rho}\left\langle g,df^0 
(\rho,W)\right\rangle_{L^2}+\frac{d^+}{d\rho}\left\langle g,df^0 
(\rho,W)\right\rangle_{L^2}\right)+\left\langle x,\dot{B}_\rho 
\right\rangle_F\, ,\quad \rho\in {\Bbb R}.
\end{eqnarray}

{\bf Hypotheses for Subsection 2.2. } Let $X=u(W)$ be a bijection 
$C({\Bbb R};F)\to C({\Bbb R};F)$ and $q$ given by condition (1') of 
Section 1. Assume that the operator $(D^u,{\cal C}):=(j^\ast\circ 
{\Bbb D}^u,{\cal C}(\Omega^u))$ is closable on $L^q(\Omega^u,P^{(m, 
r)}_{\sbmu};H)$. Suppose also (\ref{2.6}). Let $1/p+1/q=1$ and 
$(W^\rho)_{\rho\in {\Bbb R}}\in {\cal F}_{p,1}(W)$ which implies 
\begin{eqnarray}\label{2.8}
\dot{f}(\rho,W)(s)=\dot{f^0}(\rho,W)(s)+\dot{B}_\rho\cdot\1(s)\, , 
\quad s\in {\cal R},\ \rho\in {\Bbb R}. 
\end{eqnarray}
We mention that $(W^\rho)_{\rho\in {\Bbb R}}\in {\cal F}_{p,1}(W)$ 
(cf. Definition \ref{Definition1.5} (a)) is because of (\ref{2.8}) 
stronger than (\ref{2.7}) together with $\dot{f}(\rho ,)\in I^p$, 
$\rho\in {\Bbb R}$. 

Furthermore, we suppose $u\circ W^\rho\in K_{q,1}$ and $u^{-1}\circ 
X^\rho\in K^u_{q,1}$, $\rho\in {\Bbb R}$, where $q$ is given by 
condition (1') of Section 1. 
\medskip

{\bf Relationships among the flows $(W^\rho)_{\rho\in {\Bbb R}}$ 
and $(X^\rho)_{\rho\in {\Bbb R}}$ and between $u$ and $u^{-1}$. } Let 
us abbreviate $Du_s\equiv D(u(s))$, $s\in {\Bbb R}$. We observe that 
the formal relation (\ref{2.5}) gets now, because of $(W^\rho)_{\rho 
\in {\Bbb R}}\in {\cal F}_{p,1}(W)$ and $u\circ W^\rho\in K_{q,1}$, 
$\rho\in {\Bbb R}$, (especially Definition \ref{Definition1.6} (i)), 
a sense by 
\begin{eqnarray}\label{2.9}
&&\hspace{-.5cm}h\circ\left(u\circ W^\rho\right)(W)=\left\langle D 
\left(u\circ W^\rho\right)(W),j^{-1}\dot{f}(0,W)\right\rangle_{H\to 
F}\nonumber \\ 
&&\hspace{.5cm}=\left.\frac{d^\pm}{d\sigma}\right|_{\sigma=0}\left(u 
\circ W^\rho\right)(W^{\sigma})=\left.\frac{d^\pm}{d\sigma}\right|_{ 
\sigma =0}u(W^{\rho +\sigma})(W)\nonumber \\ 
&&\hspace{.5cm}=\frac{d^\pm}{d\rho}u(W^\rho)(W)\, ,\quad\rho\in {\Bbb 
R},\quad\mbox{\rm for }\, Q^{(m,r)}_{\sbnu}\mbox{\rm -a.e. }W\in\Omega
\end{eqnarray} 
and, in particular,  
\begin{eqnarray}\label{2.10}
h\circ u(W)(\, \cdot\, )\vphantom{\sum}=\left\langle Du_\cdot (W), 
j^{-1}\dot{f}(0,W)\right\rangle_{H\to F}=\left.\frac{d^\pm}{d\sigma} 
\right|_{\sigma=0}u_\cdot(W^\sigma)\, . 
\end{eqnarray} 

Let $\nabla_{W_0}u_t$ be the vector of the gradients of the 
components of $X_t\equiv u_t(W)\equiv u_t(W-W_0\1,W_0)$ relative 
to the components of $W_0$ and let $\nabla_{X_0}u^{-1}_t$ denote 
the vector of the gradients of the components of $W_t\equiv u^{ 
-1}_t(X)\equiv u^{-1}_t(X-X_0\1,X_0)$ relative to the components 
of $X_0$, $t\in {\cal R}$. 

\begin{proposition}\label{Proposition2.4} 
Let $X=u(W)$ be a bijection $C({\Bbb R};F)\to C({\Bbb R};F)$. Suppose 
condition (1') of Subsection 1.2. Let $q$ be the number given by in 
this condition (1'). Suppose that the operator $(D^u,{\cal C}):=( 
j^\ast\circ {\Bbb D}^u,{\cal C}(\Omega^u))$ is closable on $L^q( 
\Omega^u,P^{(m,r)}_{\sbmu};H)$. Assume that the following holds.
\begin{itemize}   
\item[(i)] $(W^\rho)_{\rho\in {\Bbb R}}\in {\cal F}_{p,1}(W)$, $1/p+ 
1/q=1$. The flow  $(W^\rho)_{\rho\in {\Bbb R}}$ satisfies (\ref{2.6}). 
\item[(ii)] $u\circ W^\rho\in K_{q,1}$ and $u^{-1}\circ X^\rho\in 
K^u_{q,1}$, $\rho\in {\Bbb R}$. 
\end{itemize} 
(a) We have $\left(u(W^\rho)\right)_{\rho\in {\Bbb R}}=(X^\rho)_{ 
\rho\in {\Bbb R}}\in {\cal F}_{p,1}(X)$. In addition to (\ref{2.9}), 
\begin{eqnarray*}
\dot{f}(\rho ,u^{-1}(X))(\, \cdot\, )=\left\langle D^u(u^{-1}_\cdot 
\circ X^\rho)(X),j^{-1}h(X)\right\rangle_{H\to F}\, ,\quad \rho\in 
{\Bbb R}, 
\end{eqnarray*} 
for $P^{(m,r)}_{\sbmu}$-a.e. $X\in\Omega^u$. \\ 
(b) We have, $Q^{(m,r)}_{\sbnu}$-a.e., 
\begin{eqnarray*}
\1_{[0,s]}(\, \cdot\, )\cdot {\sf e}=\left\langle j^{-1}(D_\cdot u)_1, 
(D^uu^{-1}_s)\circ u\right\rangle_{H\to F}\ \mbox{\rm and}\ {\sf e}= 
\left\langle j^{-1}\nabla_{W_0}u,(D^u u^{-1}_s)\circ u\right\rangle_{H 
\to F}   
\end{eqnarray*} 
and, $P^{(m,r)}_{\sbmu}$-a.e., 
\begin{eqnarray*}
\ \1_{[0,s]}(\, \cdot\, )\cdot {\sf e}=\left\langle j^{-1}(D^u_\cdot u^{- 
1})_1,(Du_s)\circ u^{-1}\right\rangle_{H\to F}\ \mbox{\rm and}\ {\sf e} 
=\left\langle j^{-1}\nabla_{X_0}u^{-1},(Du_s)\circ u^{-1}\right\rangle_{H 
\to F}\, . \hspace{-3mm} 
\end{eqnarray*} 
(c) For $\psi\in E^u_{q,1}$ it holds that 
\begin{eqnarray*}
D_\cdot\psi\circ u=\left\langle j^{-1}D_\cdot u,(D^u\psi)\circ u\right 
\rangle_{H\to F}\quad Q^{(m,r)}_{\sbnu}\mbox{\rm -a.e.} 
\end{eqnarray*} 
Furthermore, $\langle D\psi\circ u,k\rangle_H\in L^1(\Omega, Q^{(m,r)}_{ 
\sbnu})$ if $k\in L^p(\Omega, Q^{(m,r)}_{\sbnu};H)$. 
\end{proposition}
Proof. (a) By (i), the first part of (ii) (especially Definition 
\ref{Definition1.6} (i)), and by (\ref{2.9}) we get, $P^{(m,r)}_{\sbmu 
}$-a.e., 
\begin{eqnarray*}
\frac{d^\pm}{d\rho}X^\rho=\frac{d^\pm}{d\rho}u\left(W^\rho\right)=\left 
\langle D(u\circ W^\rho)(W),j^{-1}\dot{f}(W,0)\right\rangle_{H\to F}=h 
(X^\rho)\, ,\quad\rho\in {\Bbb R}.  
\end{eqnarray*} 
Now (i) and the first part of (ii) (especially Definition 
\ref{Definition1.6} (iii)) imply $h\circ u(W^\rho)\in jH$ $Q^{(m,r)}_{ 
\sbnu}$-a.e. and $j^{-1}h\circ u(W^\rho)\in L^p(\Omega,Q^{(m,r)}_{\sbnu} 
;H)$, $\rho\in {\Bbb R}$. 

In other words, we have $h(X^\rho)\in jH$ $P^{(m,r)}_{\sbmu}$-a.e. and 
$j^{-1}h(X^\rho)\in L^p(\Omega^u,P^{(m,r)}_{\sbmu};H)$, $\rho\in {\Bbb R} 
$. Furthermore, we have by hypothesis $X^\rho=u(W^\rho)\in C({\Bbb R};F) 
$, $\rho\in {\Bbb R}$. We conclude $(X^\rho)_{\rho\in {\Bbb R}}\in {\cal 
F}_{p,1}(X)$. Now, by the second part of (ii), $P^{(m,r)}_{\sbmu}$-a.e., 
\begin{eqnarray*}
\dot{f}(\rho,u^{-1}(X))(\, \cdot\, )&&\hspace{-.5cm}=\frac{d^\pm}{d\rho} 
W^\rho_\cdot =\frac{d^\pm}{d\rho}u^{-1}_\cdot\left(u\left(W^\rho\right) 
\right) \\
&&\hspace{-.5cm}=\left\langle D^u(u^{-1}_\cdot\circ X^\rho)(X),j^{-1}h 
(X)\right\rangle_{H\to F}\, ,\quad \rho\in {\Bbb R}. 
\end{eqnarray*} 
(b) For the set of all flows $W^\cdot\equiv (W^\rho)_{\rho\in {\Bbb R}} 
\in {\cal F}_{p,1}(W)$, with (i) of Definition \ref{Definition1.5} it 
follows that, $Q^{(m,r)}_{\sbnu}$-a.e., 
\begin{eqnarray}\label{2.11}
H\subseteq\left\{j^{-1}\dot{f}(0,\cdot ):W^\cdot\in {\cal F}_{p,1}(W) 
\right\}
\end{eqnarray} 
by just considering the specific flows of type $W^\rho:=W+\rho\cdot jk$, 
$\rho\in {\Bbb R}$, where $k\in H$. In contrast to the formulation of the 
proposition, in the following, the ``." indicates the variables to be 
jointly integrated over. Also, one should read $j^{-1}Du_\cdot\equiv\left 
((Du_\cdot)',Du_0\right)\in H$. Keeping (\ref{2.11}) in mind, the first 
part of (b) follows from 
\begin{eqnarray*}
\dot{f}(0,W)(s)&&\hspace{-.5cm}=\left\langle (D^u_\cdot u^{-1}_s)\circ u,
j^{-1}\left\langle Du_\cdot ,j^{-1}\dot{f}(0,W)\right\rangle_{H\to F}\right 
\rangle_{H\to F} \\ 
&&\hspace{-1.8cm}=\left\langle (D^u_\cdot u^{-1}_s)\circ u,\left\langle j^{ 
-1}Du_\cdot ,j^{-1}\dot{f}(0,W)\right\rangle_{H\to F}\right\rangle_{H\to F} 
 \\ 
&&\hspace{-1.8cm}=\left\langle\left(\left\langle j^{-1}(Du_\cdot)_1,
(D^u_\cdot u^{-1}_s)\circ u\right\rangle_{H\to F},\langle j^{-1}\nabla_{W_0 
}u ,(D^uu^{-1}_s)\circ u\rangle_{H\to F}\vphantom{\dot{f}}\right),j^{-1}\dot 
{f}(0,W)\vphantom{\left(\dot{f}\right)}\right\rangle_{H\to F}\, .    
\end{eqnarray*} 
For the first equality sign, recall (a) and (\ref{2.10}). For the second 
one, we refer to Definition \ref{Definition1.6} (iii). For the third one 
use Fubini's theorem together with Definition \ref{Definition1.6} (ii) 
as well as the remarks after Definition \ref{Definition1.6}. 
\medskip 

Assuming that there was an element $k\in H$, $k\neq 0$, such that $0 
=\langle k,j^{-1}h(X)\rangle_H$ for all $h(X)=\left.\frac{d^\pm}{d 
\sigma}X^\sigma\right|_{\sigma=0}$ on a set of positive $P^{(m,r)}_{ 
\sbmu}$-measure where $X^\rho=u(W^\rho)$, $\rho\in {\Bbb R}$, and 
$W^\cdot\in {\cal F}_{p,1}(W)$, from (\ref{2.9}) it would follow that 
\begin{eqnarray*}
0&&\hspace{-.5cm}=\left\langle k(\, \cdot\, ),j^{-1}\left\langle 
Du_\cdot,j^{-1}\dot{f}(0,W)\right\rangle_{H\to F}\right\rangle_H \\ 
&&\hspace{-.5cm}=\left\langle\left(\left\langle j^{-1}(Du_\cdot)_1, 
k(\, \cdot\, )\right\rangle_{H\to F},\langle j^{-1}\nabla_{W_0}u,k  
\rangle_{H\to F}\vphantom{\dot{I}}\right),j^{-1}\dot{f}(0,W)\right 
\rangle_H\, .  
\end{eqnarray*} 
By Definition \ref{Definition1.6} (ii) and (\ref{2.11}) the linear span 
of 
\begin{eqnarray*} 
\left\{j^{-1}h(X)=j^{-1}\left.\frac{d^\pm}{d\sigma}X^\sigma\right|_{ 
\sigma =0}:\ X^\cdot=u(W^\cdot),\, W^\cdot\in {\cal F}_{p,1}(W)\right\}
\end{eqnarray*} 
is therefore $P^{(m,r)}_{\sbmu}$-a.e. dense in $H$. The second part of 
(b) is now by (a) and (\ref{2.10}) a consequence of  
\begin{eqnarray*}
h(X)(s)&&\hspace{-.5cm}=\left\langle (D_\cdot u_s)\circ u^{-1},j^{-1} 
\left\langle D^u u^{-1}_\cdot ,j^{-1}h(X)\vphantom{\dot{f}}\right 
\rangle_{H\to F}\right\rangle_{H\to F} \\ 
&&\hspace{-.5cm}=\left\langle\left(\vphantom{\dot{f}}\left\langle j^{- 
1}(D^uu^{-1}_\cdot)_1,(D_\cdot u_s)\circ u^{-1}\right\rangle_{H\to F}, 
\right.\right.\\ 
&&\hspace{.5cm}\left.\left.\left\langle j^{-1}\nabla_{X_0}u^{-1},(Du_s) 
\circ u^{-1}\right\rangle_{H\to F}\vphantom{\dot{f}}\right),j^{-1}h(X) 
\right\rangle_{H\to F}\, .  
\end{eqnarray*} 
(c) Let $k\in H$ and consider the specific flow $W^\rho:=W+\rho\cdot jk$, 
$\rho\in {\Bbb R}$. We have, $Q^{(m,r)}_{\sbnu 
}$-a.e.  
\begin{eqnarray*} 
\langle D\psi\circ u,k\rangle_{H\to F}&&\hspace{-.5cm}=\left.\frac{d^\pm} 
{d\rho}\right|_{\rho=0}\psi\circ u(W^\rho) \\ 
&&\hspace{-1.8cm}=\left\langle (D^u\psi)\circ u,j^{-1}\left.\frac{d^\pm} 
{d\rho}\right|_{\rho=0}u(W^\rho)\right\rangle_{H\to F} \\ 
&&\hspace{-1.8cm}=\left\langle (D^u_\cdot\psi)\circ u,j^{-1}\langle 
Du_\cdot,k\rangle_{H\to F}\vphantom{\left(\dot{f}\right)}\right\rangle_{ 
H\to F}\vphantom{\left.\frac{d}{d\rho}\right|_{\rho =0}} \\ 
&&\hspace{-1.8cm}=\left\langle\left(\left\langle j^{-1}(Du_\cdot)_1, 
(D^u_\cdot\psi)\circ u\right\rangle_{H\to F}\, ,\, \left\langle j^{-1} 
\nabla_{W_0} u,(D^u\psi)\circ u\right\rangle_{H\to F}\right),k\vphantom 
{\left(\dot{f}\right)}\right\rangle_{H\to F} \, .  
\end{eqnarray*} 
Furthermore, if $k\in L^p(\Omega, Q^{(m,r)}_{\sbnu};H)$ then  
\begin{eqnarray*} 
\left\|\langle D\psi\circ u,k\rangle_H\right\|_{L^1(\Omega, Q^{(m,r)}_{ 
\sbnu})}&&\hspace{-.5cm}=\left\|\left\langle (D^u_\cdot\psi)\circ u,j^{- 
1}\langle Du_\cdot,k\rangle_{H\to F}\vphantom{\left(\dot{f}\right)}\right 
\rangle_H\right\|_{L^1(\Omega, Q^{(m,r)}_{\sbnu})}\nonumber \\ 
&&\hspace{-.5cm}\le \left\|(D^u\psi)\circ u\right\|_{L^q(\Omega, Q^{(m, 
r)}_{\sbnu};H)}\cdot\left\|j^{-1}\langle Du_\cdot,k\rangle_{H\to F}\right 
\|_{L^p(\Omega, Q^{(m,r)}_{\sbnu};H)}\vphantom{\left(\frac{d}{d\rho} 
\right)} \nonumber \\ 
&&\hspace{-.5cm}<\infty\, , \vphantom{\dot{f}} 
\end{eqnarray*} 
cf. Definition \ref{Definition1.6} (iii). 
\qed
\bigskip

{\bf Stochastic integral relative to $W$. } If $u$ is, just for the next 
equality, the identity and $\dot{f}(0,\cdot)(s)\in D_{p,1}$, $s\in {\cal 
R}$, and $p$ is fixed in $(W^\rho)_{\rho\in {\Bbb R}}\in {\cal F}_{p,1} 
(W)$ as a hypothesis of the present subsection, then $Q^{(m,r)}_{\sbnu} 
$-a.e., 
\begin{eqnarray}\label{2.12}
\delta\left(j^{-1}\dot{f}(0,\cdot )\right)(W)&&\hspace{-.5cm}=\sum_{i\in 
I(m,r)}\left\langle H_i ,d\dot{f}(0,W)\right\rangle_{L^2}\cdot\langle H_i 
,dW\rangle_{L^2}\nonumber \\ 
&&\hspace{-0cm}-\sum_{i\in I(m,r)}\left\langle (H_i,0),D\left\langle H_i, 
d\dot{f}(0,W)\right\rangle_{L^2}\right\rangle_H \nonumber \\ 
&&\hspace{-0cm}-\left\langle\dot{f}(0,W)(0),\frac{\nabla m(W_0)}{m(W_0)} 
\right\rangle_F-\left\langle {\sf e},\nabla_{d,W_0}\dot{f}(0,W)(0)\right 
\rangle_F\, . 
\end{eqnarray} 
Here $\nabla_{W_0}\left\langle e_j,\dot{f}(0,W)(0)\right\rangle_F 
\equiv\nabla_{W_0}\left\langle e_j,\dot{f}(0,(\, \cdot\, ,W_0))(0) 
\right\rangle_F$ is a gradient in the sense of Definition 
\ref{Definition1.1} (d) identifying in the argument $W\equiv (W 
-W_0\1,W_0)$. To be precise, representation (\ref{2.12}) follows, 
on the one hand, from $(W^\rho)_{\rho\in {\Bbb R}}\in {\cal F}_{ 
p,1}(W)$ which says $j^{-1}\dot{f}(0,\cdot )\in L^p(\Omega,Q^{ 
(m,r)}_{\sbnu};H)$. On the other hand, (\ref{2.12}) is a 
consequence of an adaption of Theorem \ref{Theorem6.7} (b) to the 
present situation. In particular recalling the meaning of $Q$ 
relative to $m$ in Subsection 1.2, the stochastic integral 
(\ref{2.12}) with integrator $\dot{f}(0,\cdot)$ is well-defined 
if $p>P$ where $1/P+1/Q=1$. Furthermore, 
\begin{eqnarray*} 
\left\langle H_i,d\dot{f}(0,W)\right\rangle_{L^2}\in D_{p,1}\equiv 
D_{p,1}(Q^{(m,r)}_{\sbnu})
\end{eqnarray*} 
by $\dot{f}(0,\, \cdot\, )(s)\in D_{p,1}$, $s\in {\cal R}$, and the 
particular form of $H_i$. 
\bigskip

{\bf Substitution law and logarithmic derivative relative to $X$. } 
Being now interested in $\beta_{h(X)}(X)$ where $h$ is given by 
(\ref{2.10}), the formula 
\begin{eqnarray}\label{2.13}
\beta_{h}\circ u=-\delta\left(j^{-1}\dot{f}(0,\cdot)\right)\quad 
Q^{(m,r)}_{\sbnu}\mbox{\rm -a.e.}  
\end{eqnarray} 
looks similar to the substitution law as, for example, demonstrated in 
\cite{SW99}, Theorem 2. Of cause, we have to prove it. The proof is of 
instructive character. It provides some insight of how the Wiener space 
representation of the integral (\ref{2.13}) is embedded in the general 
non-Gaussian stochastic calculus based on elementary integrals of the 
form $-\beta_{j(H_i ,0)}(X)$. In fact, the representation of objects 
related to non-Gaussian stochastic calculus by means of Wiener space 
anslysis terms is one of the major techniques in the paper. 
\medskip 

Let $Q$ be the number appearing in the definition of $m$ in Subsection 
1.2 and let $q$ be given by condition (1') of Section 1. Recall 
Definitions \ref{Definition1.8} and \ref{Definition1.1} (b),(f), 
as well as Proposition \ref{Proposition6.3}. 

Suppose that the operator $(D^u,{\cal C}):=(j^\ast\circ {\Bbb D}^u, 
{\cal C}(\Omega^u))$ is closable on $L^q(\Omega^u,P^{(m,r)}_{\sbmu}; 
H)$. Assume $1/p+1/q=1$ and let $p'\ge q$ as well as $1/p'+1/q'=1$ 
and note that this implies $Q>p\ge q'$. The following observation 
is important. If, for $f\in H$, there exists $c_f>0$ such that for all 
$\psi\in D^u_{p',1}$ we have 
\begin{eqnarray}\label{2.14}
\int\langle D^u\psi ,f\rangle_H\, dP^{(m,r)}_{\sbmu}\le c_f\|\psi 
\|_{L^{p'}(\Omega^u ,P^{(m,r)}_{\sbmu})} 
\end{eqnarray} 
then there is a real element $\beta_{jf}\in L^{q'}(\Omega^u,P_{\sbmu 
}^{(m,r)})$, such that 
\begin{eqnarray}\label{2.15}
-\int\psi\beta_{jf}\, dP^{(m,r)}_{\sbmu}=\int\langle D^u\psi,f 
\rangle_H\, dP^{(m,r)}_{\sbmu}\, .  
\end{eqnarray} 
\begin{proposition}\label{Proposition2.5} 
Let $X=u(W)$ be a bijection $C({\Bbb R};F)\to C({\Bbb R};F)$. Assume 
that the conditions of Proposition \ref{Proposition2.4} are satisfied. 
Suppose that the operator $(D^u,{\cal C}):=(j^\ast\circ {\Bbb D}^u, 
{\cal C}(\Omega^u))$ is closable on $L^q(\Omega^u,P^{(m,r)}_{\sbmu}; 
H)$. Here $q$ is given by condition (1') of Section 1. 
\medskip

Let $m'\equiv m'(X-X_0\1,\, \cdot\, )$ denote the conditional density 
of $X_0$ with respect to $\lambda_F$. Assume that for given $X-X_0\1$, 
the density $m'$ satisfies the conditions on $m$ in Subsection 1.2. 
\medskip


Let $Q$ be the number appearing in the definition of $m$ in Subsection 
1.2 and let $1/p+1/q=1$. Suppose $1/p+1/Q<1$. Let moreover $p'\in (1, 
\infty)$ and $1/p'+1/q'=1$ such that $p'<p$ and $q'<Q$. 

Assume the following holds. 
\begin{itemize} 
\item[(i)] We have (\ref{2.14}) for $f=(H_i,0)$, $i\in I(r)$. 
\item[(ii)] $\dot{f}(0,\cdot)(s)\in D_{p,1}$, $h(\cdot)(s)\in E^u_{p,1}$, 
$s\in {\cal R}$. 
\item[(iii)] $\sum_{i\in I(r)}\left(\|\beta_{j(H_i,0)}\|_{L^{q'}(\Omega^u, 
P^{(m,r)}_{\sbmu})}\vee 1\right)\cdot\|\langle H_i,dh\rangle_{L^2}\|_{L^p 
(\Omega^u,P^{(m,r)}_{\sbmu})}<\infty$. 
\end{itemize} 
In order to apply Lemma \ref{Lemma2.3} we also shall suppose that 
\begin{itemize} 
\item[(iv)] the linear span of $\{X_s:s\in {\cal R}\}$ is dense in $L^q 
(\Omega,Q^{(m,r)}_{\sbnu};F)$.
\end{itemize} 
Then we have (\ref{2.13}), 
\begin{eqnarray*}
\beta_{h}\circ u=-\delta (j^{-1}\dot{f}(0,\cdot))\quad Q^{(m,r)}_{\sbnu} 
\mbox{\rm -a.e.}\, , 
\end{eqnarray*} 
and  
\begin{eqnarray}\label{2.16}
\beta_{h(X)}(X)&&\hspace{-.5cm}=\sum_{i\in I(r)}\left(\langle H_i , dh 
(X)\rangle_{L^2}\cdot\beta_{j(H_i ,0)}(X)+\left\langle (H_i,0),D^u 
\langle H_i ,dh(X)\rangle_{L^2}\right\rangle_H\vphantom{\beta_{j(H_i , 
0)}(X)}\right)\nonumber \\ 
&&\hspace{-0cm}+\left\langle h(X)(0),\frac{\nabla_{X_0} m'(X-X_0\1,X_0)} 
{m'(X-X_0\1,X_0)}\right\rangle_F+\left\langle {\sf e},\nabla_{d,X_0}h(X) 
(0)\right\rangle_F 
\end{eqnarray} 
$P^{(m,r)}_{\sbmu}$-a.e. where, with $1/v=1/p+1/q'$, the infinite sum 
$\sum_{i\in I(r)}$ converges in $L^v(\Omega^u,P^{(m,r)}_{\sbmu})$. The 
term in the second line belongs to $L^w(\Omega^u,P^{(m,r)}_{\sbmu})$ 
where $1/w=1/p+1/Q$. 
\end{proposition}
{\bf Remark. (1)} Summarizing all conditions on the exponents, we get 
first $Q>2$ from $1/p+1/q=1$, $1/p+1/Q<1$, and (1'). With $p'\ge q$ 
and $Q>p\ge q'$ which has been noted for the definition in (\ref{2.14}), 
(\ref{2.15}) and $p'<p$ from the formulation of Proposition 
\ref{Proposition2.5} we get also $p\ge q'>q$, $p>p'\ge q$ and therefore 
$Q>p>2>q>1/(1-1/Q)$. 
\medskip

\nid
Proof. {\it Step 1 } First, we shall confirm well-definiteness of 
(\ref{2.16}). Indeed by (i) and (\ref{2.15}), $\beta_{j(H_i ,0)} 
\in L^{q'}(\Omega^u ,P^{(m,r)}_{\sbmu})$ is well-defined. 
Furthermore, by (ii) and the definition of $H_i$,$i\in I(r)$, we 
have $\langle H_i,dh\rangle_{L^2}\in E^u_{p,1}\subset D^u_{p,1}$. Now, 
the  convergence of the infinite sum 
\begin{eqnarray*}
\sum_{i\in I(r)}\langle H_i,dh(X)\rangle_{L^2}\cdot\beta_{j(H_i,0)}(X) 
\end{eqnarray*} 
in (\ref{2.16}) in the norm of $L^v(\Omega^u ,P^{(m,r)}_{\sbmu})$ is a 
consequence of (iii). Furthermore, we mention that according to (ii) 
of the present proposition and Definition \ref{Definition1.10}, (ii),  
\begin{eqnarray*}
\sum_{i\in I(r)}D^u\langle H_i ,dh(X)\rangle_{L^2}\cdot (H_i,0) 
\end{eqnarray*} 
is the coordinate representation of the projection of $D^uj^{-1}h(X) 
\in L^p(\Omega^u,P^{(m,r)}_{\sbmu};H\otimes H)$ to the linear span of 
$\{(H_i,0):i\in I(r)\}$ and therefore an element of $L^p(\Omega^u,P^{ 
(m,r)}_{\sbmu};H)\subseteq L^v(\Omega^u,P^{(m,r)}_{\sbmu};H)$. 
\medskip 

Next, we shall verify that the right-hand side of (\ref{2.16}) is indeed 
$\beta_h$, the 
logarithmic derivative of $P^{(m,r)}_{\sbmu}$ in direction of $h$. Let 
$\psi$ be a cylindrical function on $\Omega^u$ of the form $\psi (X)=g_0 
(X_0)\cdot g_1(X_{t_1}-X_0,\ldots ,X_{t_k}-X_0)$ where $g_0\in C_b^1(F)$, 
$g_1\in C_b^1(F^k)$, and $t_l\in {\cal R}\cap\{z\cdot t/2^m:z\in {\Bbb 
Z}\}$, $l\in \{1,\ldots ,k\}$. Similar to (\ref{6.5}) and (\ref{6.6}), 
\begin{eqnarray}\label{2.17}
\left\langle (D^u\psi)(X_{\cdot}),j^{-1}\rho\right\rangle_H 
&&\hspace{-.5cm}=g_0(X_0)\cdot\sum_{i=1}^k\langle\nabla_i g_1(X_{t_1}-X_0, 
\ldots ,X_{t_k}-X_0),\rho_{t_i}-\rho_0\rangle_F \nonumber \\ 
&&\hspace{-.0cm}+\langle\nabla g_0(X_0),\rho_0\rangle_F\cdot g_1(X_{t_1} 
-X_0,\ldots ,X_{t_k}-X_0)\vphantom{\displaystyle\sum} 
\end{eqnarray} 
with $\rho\equiv (\rho -\rho_0\cdot\1,\rho_0)\in\{j(f,x):(f,x)\in H\}$ 
and $\nabla_i$ denoting the gradient with respect to the $i$-th entry. 
By Proposition \ref{Proposition2.4} (a), $h$ is $P^{(m,r)}_{\sbmu}$-a.e. 
absolutely continuous and $j^{-1}h\in L^p(\Omega^u,P^{(m,r)}_{\sbmu} 
;H)$. We have 
\begin{eqnarray*} 
\left\langle D^u\psi,j^{-1}h\right\rangle_H=\sum_{i\in I(r)}\left\langle 
D^u\psi,(H_i ,0)\right\rangle_H\left\langle H_i ,dh\right\rangle_{L^2}+ 
\left\langle\nabla_{X_0}\psi,h(0)\right\rangle_F 
\end{eqnarray*} 
where the convergence of the right-hand side in $L^p(\Omega^u,P^{(m,r) 
}_{\sbmu})$ can be deduced from (iii) and (\ref{2.17}). Therefore, 
\begin{eqnarray*}
&&\hspace{-.5cm}\int\left\langle D^u\psi,j^{-1}h\right\rangle_H\, dP^{ 
(m,r)}_{\sbmu} \\ 
&&\hspace{.5cm}=\sum_{i\in I(r)}\int\left\langle D^u\psi,(H_i ,0)\right 
\rangle_H\left\langle H_i ,dh\right\rangle_{L^2}\, dP^{(m,r)}_{\sbmu}
+\int\left\langle\nabla_{X_0}\psi,h(0)\right\rangle_F\, dP^{(m,r)}_{ 
\sbmu} \\ 
&&\hspace{.5cm}=-\sum_{i\in I(r)}\int\psi\cdot\left\langle (H_i,0),D^u 
\left\langle H_i ,dh\right\rangle_{L^2}\right\rangle_H\, dP^{(m,r)}_{ 
\sbmu}-\int\psi\cdot\left\langle {\sf e},\nabla_{d,X_0}h(0)\right 
\rangle_F\, dP^{(m,r)}_{\sbmu}  \\ 
&&\hspace{1cm}-\sum_{i\in I(r)}\int\psi\beta_{j(H_i ,0)}\cdot\left 
\langle H_i ,dh\right\rangle_{L^2}\, dP^{(m,r)}_{\sbmu}-\int\psi\cdot 
\left\langle h(0),\frac{\nabla_{X_0} m'}{m'}\right\rangle_F\, dP^{(m,r) 
}_{\sbmu} 
\end{eqnarray*} 
which first yields the existence of the logarithmic derivative of $P^{ 
(m,r)}_{\sbmu}$ in direction of $h$, $\beta_h\in L^v(\Omega^u,P^{(m,r) 
}_{\sbmu})$ (since $v<w$), and then (\ref{2.16}). In addition, we have 
(\ref{2.14}) for $f=(0,e_i)$, $i\in \{1,\ldots ,n\cdot d\}$ with   
\begin{eqnarray*} 
\beta_{j(0,e_i)}(X)=\frac{(\nabla_{X_0} m')_i(X-X_0\1,X_0)}{m'(X-X_0\1, 
X_0)}\quad P^{(m,r)}_{\sbmu}\mbox{\rm -a.e.}  
\end{eqnarray*} 

It is the objective of the remaining proof to point out that 
relation (\ref{2.13}) is a consequence of (\ref{2.12}) and 
(\ref{2.16}). Strictly speaking, we have to demonstrate that, 
$Q^{(m,r)}_{\sbnu}$-a.e., 
\begin{eqnarray}\label{2.18}
&&\hspace{-1cm}-\sum_{i\in I(r)}\left(\langle H_i,dh\circ u\rangle_{ 
L^2}\cdot\beta_{j(H_i,0)}\circ u+\left\langle (H_i,0),\left(D^u\langle 
H_i,dh\rangle_{L^2}\right)\circ u\right\rangle_H\vphantom{\sum}\right) 
\nonumber \\ 
&&\hspace{.5cm}-\sum_i(h\circ u(0))_i\cdot\beta_{j(0,e_i)}\circ u- 
\left\langle{\sf e},\nabla_{d,X_0}h\circ u(0)\right\rangle_H\nonumber 
 \\ 
&&\hspace{-.0cm}=\sum_{j\in I(m,r)}\langle H_j ,d\dot{f}(0,\cdot) 
\rangle_{L^2}\cdot\langle H_j ,dW\rangle_{L^2}-\sum_{j\in I(m,r)} 
\left\langle (H_j,0),D\langle H_j ,d\dot{f}(0,\cdot)\rangle_{L^2} 
\right\rangle_H \nonumber \\ 
&&\hspace{.5cm}-\left\langle\dot{f}(0,\cdot)(0)\, ,\, \frac{\nabla 
m(W_0)}{m(W_0)}\right\rangle_F-\left\langle{\sf e},\nabla_{d,W_0} 
\dot{f}(0,\cdot)(0)\right\rangle_F\, .  
\end{eqnarray} 
We also recall (\ref{2.9}) which leads to 
\begin{eqnarray}\label{2.19}
\left\langle Du_0 ,j^{-1}\dot{f}(0,\, \cdot\, )\right\rangle_{H\to 
F}=h\circ u(0)\, .   
\end{eqnarray} 

\nid 
{\it Step 2 } In this step, we are going to analyze the terms 
$\langle H_i,dh\circ u\rangle_{L^2}$ and $\left(D^u\langle H_i ,dh 
\rangle_{L^2}\right)\circ u$ from the left-hand side of (\ref{2.18}). 
We stress that everything further in this step will hold almost 
everywhere under the measure $Q^{(m,r)}_{\sbnu}$.
\medskip

Recalling, $u^{-1}\in K^u_{q,1}$, (cf. (ii) of Proposition 
\ref{Proposition2.4}) we may define 
\begin{eqnarray}\label{2.20} 
\begin{array}{rl} 
{\bf H}_i&\hspace{-.2cm}:=q_{m,r}j^{-1}\langle D^uu^{-1}_\cdot\circ u, 
(H_i,0)\rangle_{H\to F}\, , \quad i\in I(r), \\ & \\ 
{\bf e}_j&\hspace{-.2cm}:=q_{m,r}j^{-1}\langle D^uu^{-1}_\cdot\circ u,
(0,e_j)\rangle_{H\to F}\, , \quad j\in\{1,\ldots ,n\cdot d\}.   
\end{array} 
\end{eqnarray} 
In fact, we have 
\begin{eqnarray*}
{\bf H}_i,{\bf e}_j\in H^{(m,r)}\quad Q^{(m,r)}_{\sbnu}\mbox{\rm a.e. 
} 
\end{eqnarray*} 
and with Definition \ref{Definition1.10} (iii) 
\begin{eqnarray*}
{\bf H}_i,{\bf e}_j\in L^p\left(\Omega^u,P^{(m,r)}_{\sbmu};H^{(m,r)} 
\right),\ i\in I(r)\, , \quad j\in\{1,\ldots ,n\cdot d\}.   
\end{eqnarray*} 
Moreover, Lemma \ref{Lemma2.2} implies 
\begin{eqnarray*}
\quad\left\langle D_\cdot u,q_{m,r}j^{-1}\langle D^uu^{-1}_\cdot\circ u, 
(H_i ,0)\rangle_{H\to F}\right\rangle_{H\to F}=\left\langle D_\cdot u, 
j^{-1}\langle D^uu^{-1}_\cdot\circ u,(H_i ,0)\rangle_{H\to F}\right 
\rangle_{H\to F} 
\end{eqnarray*} 
where, as in the proof of Proposition \ref{Proposition2.4}, the dots 
indicate the variable to be jointly integrated over. Thus, Proposition 
\ref{Proposition2.4} (b) as well as the arguments of its proof together 
with (\ref{2.20}) yield 
\begin{eqnarray}\label{2.21}
&&\hspace{-1.0cm}(j(H_i,0))_s\vphantom{\frac12} 
\nonumber \\   
&&\hspace{-.5cm}=\left\langle\left(\left\langle j^{-1}(D^uu^{-1}_\cdot)_1 
\circ u,D_\cdot u_s\right\rangle_{H\to F},\left\langle j^{-1}\nabla_{X_0}
u^{-1}\circ u,Du_s\right\rangle_{H\to F}\vphantom{\dot{f}}\right),(H_i,0) 
\vphantom{\left(\dot{f}\right)}\right\rangle_{H\to F}\nonumber \\  
&&\hspace{-.5cm}=\left\langle D_\cdot u_s,\left\langle j^{-1}D^uu^{-1}_\cdot 
\circ u,(H_i,0)\right\rangle_{H\to F}\right\rangle_{H\to F}\vphantom{\frac12} 
\nonumber \\  
&&\hspace{-.5cm}=\left\langle D_\cdot u_s,j^{-1}\left\langle D^uu^{-1}_\cdot 
\circ u,(H_i,0)\right\rangle_{H\to F}\right\rangle_{H\to F}\vphantom{\frac12} 
\nonumber \\ 
&&\hspace{-.5cm}=\langle Du_s,{\bf H}_i\rangle_{H\to F}\, ,\quad s\in {\cal 
R}\, ,\ i\in I(r). \vphantom{\frac12} 
\end{eqnarray} 
Similarly, 
\begin{eqnarray}\label{2.22}
&&\hspace{-1.0cm}(0,e_j)=\langle Du_s,{\bf e}_j\rangle_{H\to F}\, , \quad 
s\in {\cal R}\, , \ j\in\{1,\ldots ,n\cdot d\}. 
\end{eqnarray} 

As a consequence of Proposition \ref{Proposition2.4} (a), $(X^\rho)_{\rho 
\in {\Bbb R}}\in {\cal F}_{p,1}(X)$, we have $h\circ u\in jH$ $Q^{(m,r)}_{ 
\sbnu}$-a.e. Also, we have $Q^{(m,r)}_{\sbnu}$-a.e. for all $s\in {\cal R}$ 
\begin{eqnarray*} 
\left\langle j^{-1}D^uu^{-1}_s\circ u,j^{-1}h\circ u\right\rangle_H 
&&\hspace{-.5cm}=\sum_{i\in I(r)}\left\langle j^{-1}D^uu^{-1}_s\circ u, 
(H_i,0)\right\rangle_H\cdot\left\langle (H_i,0),j^{-1}h\circ u\right 
\rangle_H \\ 
&&+\left\langle (D^uu^{-1}_s\circ u)_2,h\circ u(0)\right\rangle_F\, . 
\end{eqnarray*} 
Thus, the sum 
\begin{eqnarray*} 
&&\hspace{-.5cm}\sum_{i\in I(r)}{\bf H}_i\cdot \left\langle (H_i,0), 
j^{-1}h\circ u\right\rangle_H\quad \mbox{\rm converges }\ Q^{(m,r)}_{ 
\sbnu} \mbox{\rm -a.e. in }H. 
\end{eqnarray*} 
By (\ref{2.19})-(\ref{2.22}) we get therefore 
\begin{eqnarray*} 
&&\hspace{-.5cm}\left\langle D_\cdot u,\sum_{i\in I(r)}{\bf H}_i(\, 
\cdot\, )\cdot\langle H_i ,dh\circ u\rangle_{L^2}+\sum_{j}{\bf e}_j(\, 
\cdot\, )\cdot\langle e_j,h\circ u(0)\rangle_F\right\rangle_{H\to F} \\ 
&&\hspace{.5cm}=\sum_{i\in I(r)}\left\langle D_\cdot u ,{\bf H}_i(\, 
\cdot\, )\right\rangle_{H\to F}\cdot \langle H_i ,dh\circ u\rangle_{L^2} 
+\sum_j\left\langle D_\cdot u ,{\bf e}_j(\, \cdot\, )\right\rangle_{H\to 
F}\cdot\langle e_j,h\circ u(0)\rangle_F \\  
&&\hspace{.5cm}=\sum_{i\in I(r)}j(H_i,0)\langle (H_i,0),j^{-1}h\circ u 
\rangle_H+h\circ u(0)\cdot\1 \\  
&&\hspace{.5cm}=h\circ u=\left\langle Du ,j^{-1}\dot{f}(0,\, \cdot\, ) 
\right\rangle_{H\to F}\, , 
\end{eqnarray*} 
for the second last equality sign we consider the integral induced by the 
embedding $j$ the inner product in $H$ with $\1\in H$. Using Lemma 
\ref{Lemma2.3} it follows that 
\begin{eqnarray}\label{2.23}
\begin{array}{rl}
\D\sum_{i\in I(r)}\langle (H_j,0),{\bf H}_i\rangle_H\cdot\langle H_i,dh 
\circ u\rangle_{L^2}+\D\sum_i\langle (H_j,0),{\bf e}_i\rangle_H\cdot 
(h\circ u(0))_i &\hspace{-.2cm}=\langle H_j,d\dot{f}(0,\, \cdot\, ) 
\rangle_{L^2}\, ,\hspace{-.5cm}\\ & \\ 
\D\sum_{j\in I(r)}\langle (0,e_i),{\bf H}_j\rangle_H\cdot\langle H_j,dh 
\circ u\rangle_{L^2}+\D\sum_{j}\langle (0,e_i),{\bf e}_j\rangle_H\cdot 
(h\circ u(0))_j &\hspace{-.2cm}=\langle e_i,\dot{f}(0,\, \cdot\, )(0) 
\rangle_F\, , \hspace{-.5cm}
\end{array}
\end{eqnarray} 
$j\in I(m,r)$, $i\in\{1,\ldots ,n\cdot d\}$. Furthermore, we mention that 
$\langle H_i,dh\rangle_{L^2}\in E^u_{p,1}$ by (ii) and the special shape 
of $H_i$, $i\in I(r)$. Now, by Proposition \ref{Proposition2.4} (c), the 
arguments of the proof of \ref{Proposition2.4} (b), and (\ref{2.21}) for 
the last equality sign give 
\begin{eqnarray}\label{2.24}
&&\hspace{-.5cm}\sum_{j\in I(m,r)}\langle {\bf H}_i,(H_j ,0)\rangle_H 
\left\langle (H_j ,0),D\langle H_i ,dh\circ u\rangle_{L^2}\right\rangle_H 
\nonumber \\ 
&&\hspace{1.0cm}+\sum_j\langle {\bf H}_i,(0,e_j)\rangle_H\left\langle 0, 
e_j),D\langle H_i,dh\circ u\rangle_{L^2}\right\rangle_H\nonumber \\ 
&&\hspace{.5cm}=\left\langle {\bf H}_i,D\langle H_i ,dh\circ u\rangle_{ 
L^2}\right\rangle_H\vphantom{\int} \nonumber \\ 
&&\hspace{.5cm}=\left\langle{\bf H}_i,\left\langle j^{-1}Du_\cdot,\left( 
D^u_\cdot\langle H_i,dh\rangle_{L^2}\right)\circ u\right\rangle_{H\to F} 
\right\rangle_H\vphantom{\int} \nonumber \\ 
&&\hspace{.5cm}=\left\langle\langle j^{-1}Du_\cdot ,{\bf H}_i\rangle_{H 
\to F},\left(D^u_\cdot\langle H_i,dh\rangle_{L^2}\right)\circ u\right 
\rangle_H\vphantom{\int} \nonumber \\ 
&&\hspace{.5cm}=\left\langle j^{-1}\langle Du_\cdot ,{\bf H}_i\rangle_{H 
\to F},\left(D^u_\cdot\langle H_i,dh\rangle_{L^2}\right)\circ u\right 
\rangle_H\vphantom{\int} \nonumber \\ 
&&\hspace{.5cm}=\left\langle(H_i,0),\left(D^u\langle H_i,dh\rangle_{L^2} 
\right)\circ u\right\rangle_H\, , \quad i\in I(r).  \vphantom{\int}
\end{eqnarray} 
Similarly, with (\ref{2.22}) for the last equality sign we obtain 
\begin{eqnarray}\label{2.25}
&&\hspace{-.5cm}\sum_i\sum_{j\in I(m,r)}\langle {\bf e}_i,(H_j ,0) 
\rangle_H\left\langle (H_j ,0),D\langle e_i,h\circ u(0)\rangle_F 
\right\rangle_H\nonumber \\ 
&&\hspace{1.0cm}+\sum_i\sum_j\langle {\bf e}_i,(0,e_j)\rangle_H 
\left\langle 0,e_j),D\langle e_i,h\circ u(0)\rangle_F\right\rangle_H 
\nonumber \\ 
&&\hspace{.5cm}=\sum_i\left\langle {\bf e}_i,D\langle e_i,h\circ 
u(0)\rangle_F\right\rangle_H\vphantom{\int} \nonumber \\ 
&&\hspace{.5cm}=\sum_i\left\langle{\bf e}_i,\left\langle j^{-1} 
Du_\cdot,D^u_\cdot\langle e_i,h\circ u(0)\rangle_F\right\rangle_{H\to 
F}\right\rangle_H\vphantom{\int} \nonumber \\ 
&&\hspace{.5cm}=\sum_i\left\langle\langle j^{-1}Du_\cdot ,{\bf e}_i 
\rangle_{H\to F},D^u_\cdot\langle e_i,h\circ u(0)\rangle_F\right 
\rangle_H\vphantom{\int} \nonumber \\ 
&&\hspace{.5cm}=\sum_i\left\langle j^{-1}\langle Du_\cdot ,{\bf e}_i 
\rangle_{H\to F},D^u_\cdot\langle e_i,h\circ u(0)\rangle_F\right 
\rangle_H\vphantom{\int} \nonumber \\ 
&&\hspace{.5cm}=\left\langle{\sf e},\nabla_{d,X_0}h\circ u(0)\right 
\rangle_H.\vphantom{\int}
\end{eqnarray} 
{\it Step 3 } We are interested in the logarithmic derivatives $\beta_{ 
j(H_i,0)}$ and $\beta_{j(0,e_{j'})}$. We recall that by (\ref{2.20}) and 
$u^{-1}\in K^u_{q,1}$ we have ${\bf H}_i\in L^p(\Omega,Q^{(m,r)}_{\sbnu} 
;H)$, $i\in I(r)$, cf. Definition \ref{Definition1.10}, (iii). Let $\psi 
\in {\cal C}(\Omega^u)$ be a cylindric cylindrical function in the sense 
of Definitions \ref{Definition1.1} (a) and \ref{Definition1.8}. By Remark 
(9) of Section 1 it holds that $\psi\in E^u_{q,1}$. 

According to assumption (i) of this proposition, (\ref{2.21}), (\ref{2.22}), 
$u\in K_{q,1}$, and Proposition \ref{Proposition2.4} (c) we have 
\begin{eqnarray*}
-\int\psi\cdot\beta_{j(H_i ,0)}\, dP^{(m,r)}_{\sbmu}&&\hspace{-.5cm}=\int 
\langle D^u\psi ,(H_i ,0)\rangle_H\, dP^{(m,r)}_{\sbmu} \\ 
&&\hspace{-.5cm}=\int\left\langle (D^u_\cdot\psi)\circ u,j^{-1}\left\langle 
Du_\cdot,{\bf H}_i\right\rangle_{H\to F}\right\rangle_H\, dQ^{(m,r)}_{\sbnu} 
\\ 
&&\hspace{-.5cm}=\int\left\langle (D^u_\cdot\psi)\circ u,\left\langle j^{-1} 
Du_\cdot,{\bf H}_i\right\rangle_{H\to F}\right\rangle_H\, dQ^{(m,r)}_{\sbnu} 
\\ 
&&\hspace{-.5cm}=\int\left\langle\left\langle (D^u_\cdot\psi)\circ u,j^{-1} 
Du_\cdot\right\rangle_{H\to F},{\bf H}_i\right\rangle_H\, dQ^{(m,r)}_{\sbnu} 
\\ 
&&\hspace{-.5cm}=\int\langle D\psi\circ u,{\bf H}_i\rangle_H\, dQ^{(m,r) 
}_{\sbnu} \\ 
&&\hspace{-.5cm}=\int\psi\circ u\cdot\delta({\bf H}_i)\, dQ^{(m,r)}_{\sbnu}
\end{eqnarray*} 
and similarly 
\begin{eqnarray*}
-\int\psi\cdot\beta_{j(0,e_{j'})}\, dP^{(m,r)}_{\sbmu}&&\hspace{-.5cm}=\int 
\psi\circ u\cdot\delta({\bf e}_{j'})\, dQ^{(m,r)}_{\sbnu}\, . 
\end{eqnarray*} 
This shows 
\begin{eqnarray*}
\beta_{j(H_i ,0)}\circ u=-\delta ({\bf H}_i)\quad Q^{(m,r)}_{\sbnu}\mbox{ 
\rm -a.e.}\, , \quad i\in I(r),  
\end{eqnarray*} 
and 
\begin{eqnarray*}
\beta_{j(0,e_{j'})}\circ u=-\delta ({\bf e}_{j'})\quad Q^{(m,r)}_{ 
\sbnu}\mbox{\rm -a.e.}\, , \quad j'\in \{1,\ldots ,n\cdot d\}.   
\end{eqnarray*} 
Taking into consideration ${\bf H}_i\in H^{(m,r)}$, $i\in I(r)$, as 
well as ${\bf e}_j\in H^{(m,r)}$, $j\in\{1,\ldots ,n\cdot d\}$ 
$Q^{(m,r)}_{\sbnu}$-a.e., we get ${\bf H}_i=\sum_{j\in I(m,r)}\langle 
{\bf H}_i,(H_j,0)\rangle_H\cdot (H_j,0)+\sum_{j'}\langle {\bf H}_i, 
(0,e_{j'})\rangle_H\cdot (0,e_{j'})$ and thus from Theorem 
\ref{Theorem6.7} (b) the representation 
\begin{eqnarray}\label{2.26} 
\beta_{j(H_i ,0)}\circ u&&\hspace{-.5cm}=-\delta({\bf H}_i)\vphantom{ 
\sum} \nonumber \\ 
&&\hspace{-.5cm}=-\sum_{j\in I(m,r)}\langle {\bf H}_i,(H_j,0)\rangle_H 
\cdot\langle H_j,dW\rangle_{L^2} \nonumber \\ 
&&\hspace{-.0cm}+\sum_{j\in I(m,r)}\left\langle (H_j,0),D\langle {\bf 
H}_i,(H_j,0)\rangle_H\right\rangle_H \nonumber \\ 
&&\hspace{-.0cm}+\left\langle ({\bf H}_i)_2\, ,\, \frac{\nabla m(W_0)} 
{m(W_0)}\right\rangle_F+\left\langle{\sf e}\, ,\, \nabla_{d,W_0}({\bf H 
}_i)_2\right\rangle_F\, , \quad i\in I(r).  
\end{eqnarray} 
Similarly, 
\begin{eqnarray}\label{2.27} 
\beta_{j(0,e_i)}\circ u&&\hspace{-.5cm}=-\delta ({\bf e}_i) 
\vphantom{\sum} \nonumber \\ 
&&\hspace{-.5cm}=-\sum_{j\in I(m,r)}\langle {\bf e}_i,(H_j,0)\rangle_H 
\cdot\langle H_j,dW\rangle_{L^2} \nonumber \\ 
&&\hspace{-.0cm}+\sum_{j\in I(m,r)}\left\langle (H_j,0),D\langle {\bf 
e}_i,(H_j,0)\rangle_H\right\rangle_H \nonumber \\ 
&&\hspace{-.0cm}+\left\langle ({\bf e}_i)_2\, ,\, \frac{\nabla m(W_0)} 
{m(W_0)}\right\rangle_F+\left\langle{\sf e}\, ,\, \nabla_{d,W_0}({\bf 
e}_i)_2\right\rangle_F\, , \quad i\in \{1,\ldots ,n\cdot d\}.  
\end{eqnarray} 
{\it Step 4 } We will finally show that (\ref{2.13}) is a consequence of 
(\ref{2.16}).  For this we verify (\ref{2.18}) by using the results 
of the Steps 2 and 3. Especially (\ref{2.24})-(\ref{2.27}) imply that 
we have $Q^{(m,r)}_{\sbnu}$-a.e. 
\begin{eqnarray*} 
&&\hspace{-1cm}-\sum_{i\in I(r)}\left(\langle H_i,dh\circ u\rangle_{ 
L^2}\cdot\beta_{j(H_i,0)}\circ u+\left\langle (H_i,0),\left(D^u\langle 
H_i,dh\rangle_{L^2}\right)\circ u\right\rangle_H\vphantom{\sum}\right) 
\nonumber \\ 
&&\hspace{.5cm}-\sum_i(h\circ u(0))_i\cdot\beta_{j(0,e_i)}\circ u- 
\left\langle{\sf e},\nabla_{d,X_0}h\circ u(0)\right\rangle_H\nonumber 
 \\ 
&&\hspace{-.0cm}=\sum_{i\in I(r)}\langle H_i,dh\circ u\rangle_{L^2} 
\cdot\sum_{j\in I(m,r)}\langle {\bf H}_i,(H_j,0)\rangle_H\cdot\langle 
H_j,dW\rangle_{L^2} \nonumber \\ 
&&\hspace{.5cm}-\sum_{i\in I(r)}\langle H_i,dh\circ u\rangle_{L^2} 
\cdot\sum_{j\in I(m,r)}\left\langle (H_j,0),D\langle {\bf H}_i,(H_j 
,0)\rangle_H\right\rangle_H \nonumber \\ 
&&\hspace{.5cm}-\sum_{i\in I(r)}\langle H_i,dh\circ u\rangle_{L^2} 
\cdot\left(\left\langle ({\bf H}_i)_2\, ,\, \frac{\nabla m(W_0)} 
{m(W_0)}\right\rangle_F+\left\langle{\sf e}\, ,\, \nabla_{d,W_0}({\bf 
H}_i)_2\right\rangle_F\right)\nonumber \\ 
&&\hspace{.5cm}-\sum_{i\in I(r)}\sum_{j\in I(m,r)}\langle {\bf H}_i, 
(H_j,0)\rangle_H\left\langle (H_j,0),D\langle H_i,dh\circ u\rangle_{ 
L^2}\right\rangle_H \nonumber \\ 
&&\hspace{.5cm}-\sum_{i\in I(r)}\ \ \sum_{j'}\ \ \langle {\bf H}_i, 
(0,e_{j'})\rangle_H\left\langle 0,e_{j'}),D\langle H_i,dh\circ u 
\rangle_{L^2}\right\rangle_H\nonumber \\ 
&&\hspace{.5cm}+\sum_i(h\circ u(0))_i\cdot\sum_{j\in I(m,r)}\langle 
{\bf e}_i,(H_j,0)\rangle_H\cdot\langle H_j,dW\rangle_{L^2}\nonumber 
 \\ 
&&\hspace{.5cm}-\sum_i(h\circ u(0))_i\cdot\sum_{j\in I(m,r)}\left 
\langle (H_j,0),D\langle {\bf e}_i,(H_j,0)\rangle_H\right\rangle_H 
\nonumber \\ 
&&\hspace{.5cm}-\sum_i(h\circ u(0))_i\cdot\left(\left\langle ({\bf 
e}_i)_2\, ,\, \frac{\nabla m(W_0)}{m(W_0)}\right\rangle_F+\left 
\langle{\sf e}\, ,\, \nabla_{d,W_0}({\bf e}_i)_2\right\rangle_F 
\right)\nonumber \\ 
&&\hspace{.5cm}-\sum_i\sum_{j\in I(m,r)}\langle {\bf e}_i,(H_j,0)
\rangle_H\left\langle (H_j,0),D(h\circ u(0))_i\right\rangle_H 
\nonumber \\ 
&&\hspace{.5cm}-\sum_i\ \ \sum_{j'}\ \ \langle {\bf e}_i,(0,e_{j'}) 
\rangle_H\left\langle (0,e_{j'}),D(h\circ u(0))_i\right\rangle_H\, .  
\end{eqnarray*} 
Applying the usual product rule, this is equal to 
\begin{eqnarray*} 
&&\hspace{-.5cm}\sum_{i\in I(r)}\langle H_i,dh\circ u\rangle_{L^2} 
\cdot\sum_{j\in I(m,r)}\langle {\bf H}_i,(H_j,0)\rangle_H\cdot\langle 
H_j,dW\rangle_{L^2} \nonumber \\ 
&&\hspace{.5cm}+\sum_i(h\circ u(0))_i\cdot\sum_{j\in I(m,r)}\langle 
{\bf e}_i,(H_j,0)\rangle_H\cdot\langle H_j,dW\rangle_{L^2}\nonumber 
 \\ 
&&\hspace{.5cm}-\sum_{i\in I(r)}\sum_{j\in I(m,r)}\left\langle (H_j, 
0),D\left(\langle {\bf H}_i,(H_j,0)\rangle_H\cdot\langle H_i,dh\circ 
u\rangle_{L^2}\vphantom{\dot{f}}\right)\vphantom{\left(\dot{f}\right)} 
\right\rangle_H \nonumber \\ 
&&\hspace{.5cm}-\sum_i\ \ \sum_{j\in I(m,r)}\left\langle (H_j, 
0),D\left(\langle {\bf e}_i,(H_j,0)\rangle_H\cdot(h\circ 
u(0))_i\vphantom{\dot{f}}\right)\vphantom{\left(\dot{f}\right)} 
\right\rangle_H \nonumber \\ 
&&\hspace{.5cm}-\sum_{i\in I(r)}\ \ \sum_{j'}\ \ \langle {\bf H}_i, 
(0,e_{j'})\rangle_H\left\langle (0,e_{j'}),D\langle H_i,dh\circ u 
\rangle_{L^2}\right\rangle_H\nonumber \\ 
&&\hspace{.5cm}-\sum_i\ \ \sum_{j'}\ \ \langle {\bf e}_i,(0,e_{j'}) 
\rangle_H\left\langle (0,e_{j'}),D(h\circ u(0))_i\right\rangle_H 
\nonumber \\ 
&&\hspace{.5cm}-\sum_{i\in I(r)}\sum_{j'}\ \langle H_i,dh\circ u 
\rangle_{L^2}\cdot\left(\left\langle \langle {\bf H}_i,(0,e_{j'}) 
\rangle_H\cdot e_{j'}\, ,\, \frac{\nabla m(W_0)}{m(W_0)}\right 
\rangle_F\right.\nonumber \\ 
&&\hspace{8.5cm}\left.+\left\langle D\langle {\bf H}_i,(0,e_{j'}) 
\rangle_H\, ,\, (0,e_{j'})\right\rangle_H\vphantom{\frac{\nabla m 
(W_0)}{m(W_0)}}\right)\nonumber \\ 
&&\hspace{.5cm}-\sum_i\ \ \sum_{j'}(h\circ u(0))_i\cdot\left(\left 
\langle\ \langle {\bf e}_i,(0,e_{j'})\rangle_H\cdot e_{j'}\, ,\, 
\frac{\nabla m(W_0)}{m(W_0)}\right\rangle_F\right.\nonumber \\ 
&&\hspace{8.5cm}\left.+\left\langle D\langle {\bf e}_i,(0,e_{j'}) 
\rangle_H\, ,\, (0,e_{j'})\right\rangle_H\vphantom{\frac{\nabla m 
(W_0)}{m(W_0)}}\right)\, .  
\end{eqnarray*} 
Now (\ref{2.23}) yields 
\begin{eqnarray*} 
&&\hspace{-1cm}-\sum_{i\in I(r)}\left(\langle H_i,dh\circ u\rangle_{ 
L^2}\cdot\beta_{j(H_i,0)}\circ u+\left\langle (H_i,0),\left(D^u\langle 
H_i,dh\rangle_{L^2}\right)\circ u\right\rangle_H\vphantom{\sum}\right) 
 \\ 
&&\hspace{.5cm}-\sum_i(h\circ u(0))_i\cdot\beta_{j(0,e_i)}\circ u- 
\left\langle{\sf e},\nabla_{d,X_0}h\circ u(0)\right\rangle_H \\ 
&&\hspace{-.0cm}=\sum_{j\in I(m,r)}\langle H_j ,d\dot{f}(0,\cdot) 
\rangle_{L^2}\cdot\langle H_j ,dW\rangle_{L^2}-\sum_{j\in I(m,r)} 
\left\langle (H_j,0),D\langle H_j ,d\dot{f}(0,\cdot)\rangle_{L^2} 
\right\rangle_H \\ 
&&\hspace{.5cm}-\left\langle\dot{f}(0,\cdot)(0)\, ,\, \frac{\nabla 
m(W_0)}{m(W_0)}\right\rangle_F-\left\langle{\sf e},\nabla_{d,W_0} 
\dot{f}(0,\cdot)(0)\right\rangle_F
\end{eqnarray*} 
which is (\ref{2.18}). 
\qed

\subsection{A chain rule for absolutely continuous $f(\rho,W)(\cdot)$ 
but jumps in $\dot{f}(\rho,W)(\cdot)$} 

We proceed with a definition whose appropriateness will become obvious  
by the subsequent Lemma \ref{Lemma3.1}. For $m\in {\Bbb Z}_+$ and $r\in 
{\Bbb N}$ introduce the following. 
\begin{definition}\label{Definition2.6}
{\rm Let $1<p<\infty$ and let $Y\in I^p$ such that $Y_s$ is constant 
outside of $s\in {\cal R}$. Let $f=f_c+f_j$ where $f_c\in L^q(\Omega,Q^{ 
(m,r)}_{\sbnu};L^2)$ is a $Q^{(m,r)}_{\sbnu}$-a.e. continuous process 
and $f_j$ is a cadlag pure jump process of finite jump variation on ${ 
\cal R}$, $V_{-rt}^{(r-1)t}(f)\in L^q(\Omega ,Q^{(m,r)}_{\sbnu})$, $1/p 
+1/q=1$. Let $x\in L^q(\Omega ,Q^{(m,r)}_{\sbnu};F)$. 
\medskip 

\nid
(a) Let $\vp=(f,x)$ and define 
\begin{eqnarray*} 
&&\hspace{-.5cm}\int_u^v\left\langle f,dY\right\rangle_F+\left\langle x,
{\textstyle\frac12}\left(Y_{0-}+Y_0\right)\right\rangle_F:=\left\langle 
(f\, \1_{[u,v]},x),j^{-1}Y_a\right\rangle_H \\ 
&&\hspace{.5cm}+\sum_{w\in (u,v)}\frac12\langle f_{w-}+f_w,\Delta Y_j(w) 
\rangle_F+\frac12\langle f_u,\Delta Y_j(u)\rangle_F +\frac12\langle f_{v 
-},\Delta Y_j(v)\rangle_F\quad Q^{(m,r)}_{\sbnu}\mbox{\rm -a.e.} 
\end{eqnarray*} 
noting that $\left\langle (f\, \1_{[u,v]},x),j^{-1}Y_a\right\rangle_H 
=\int_u^v\langle f,dY_a\rangle_F+\left\langle x,{\textstyle\frac12} 
\left(Y_{0-}+Y_0\right)\right\rangle_F$ in the sense of 
Riemann-Stieltjes integration. 
\medskip 

\nid 
(b) Let $\vp=(f,x)$ and define $j^{-1}Y$ applied to the test element 
$\vp$ by 
\begin{eqnarray}\label{2.28}
j^{-1}Y(\vp):=\int_{\cal R}\left\langle f,dY\right\rangle_F +\left 
\langle x,{\textstyle\frac12}\left(Y_{0-}+Y_0\right)\right\rangle_F\quad 
Q^{(m,r)}_{\sbnu}\mbox{\rm -a.e.} 
\end{eqnarray} 
(c) For $\vp$ and $Y$ as above, we introduce the inner product $\langle 
\, \cdot\, , \, j^{-1}(\, \cdot\, )\rangle_H$ by 
\begin{eqnarray*} 
\langle\vp ,j^{-1}Y\rangle_H:=j^{-1}Y(\vp)\, . 
\end{eqnarray*} 
Replacing in (\ref{2.28}) $\langle\, \cdot\, , \, \cdot\, \rangle_F$ 
with $\langle\, \cdot\, , \, \cdot\, \rangle_{F\to F}$, we define 
$\langle\vp,j^{-1}Y\rangle_{H\to F}$. 
}
\end{definition}
\medskip 

{\bf Relations} (\ref{2.9}) {\bf and} (\ref{2.10}) {\bf in a 
distributional sense. } Suppose we are given $W^\rho:=W+f(\rho,W)$. 
Furthermore, allow  
\begin{eqnarray}\label{2.29} 
\dot{f}(0,\cdot)\in I^p 
\end{eqnarray} 
where $1/p+1/q=1$ and $q$ is fixed by assumption (1') of Section 1 and 
let (\ref{2.6}) and (\ref{2.7}) be in force. In particular, we recall 
that the symbol $\dot{f}(\rho ,W)$ is used in the strict sense of a {\it 
weak mixed derivative} which is $Q^{(m,r)}_{\sbnu}$-a.e. given for all 
test elements $k=(g,x)\in H$ with cadlag $g:{\cal R}\to F$ of finite 
jump variation, $x\in F$, by (\ref{2.7}), 
\begin{eqnarray*}
\dot{f}(\rho ,W)(k)&&\hspace{-.5cm}=\frac{d^\pm}{d\rho}\left\langle g, 
df^0(\rho,W)\right\rangle_{L^2}+\left\langle x,\dot{B}_\rho\right 
\rangle_F \\ 
&&\hspace{-.5cm}=\frac12\left(\frac{d^-}{d\rho}\left\langle g,df^0(\rho, 
W)\right\rangle_{L^2}+\frac{d^+}{d\rho}\left\langle g,df^0(\rho,W)\right 
\rangle_{L^2}\right)+\left\langle x,\dot{B}_\rho\right\rangle_F\, .
\end{eqnarray*}

One major difference to the analysis in Subsection 2.2 is that we 
assume there $(W^\rho)_{\rho\in {\Bbb R}}\in {\cal F}_{p,1}$ which 
is replaced here by (\ref{2.29}) and (\ref{2.7}) for $\rho=0$. Even 
more restrictive in Subsection 2.2 is the hypothesis $u_s\in E_{q,1}$, 
$s\in {\cal R}$, as part of $u\in K_{q,1}$ in order to establish the 
generator of the flow $(X^\rho)_{\rho\in {\Bbb R}}$ namely relation 
(\ref{2.10}). Note also that for the following proposition we do not 
require that $W^\rho=W+f(\rho ,W)$, $\rho\in {\Bbb R}$, is a flow.  

As a preparation of the analysis in Section 3, we are interested in a 
counterpart of (\ref{2.10}) at $\rho=0$ under the conditions of the 
present subsection. 
\medskip

For this fix $r\in {\Bbb Z}_+$ and define for all $W\in jH$ the terms 
$\xi_i:=\left\langle (H_i,0),j^{-1}W\right\rangle_H$, $i\in I(r)$, $x_j 
:=(W_0)_j$, $j\in \{1,\ldots ,n\cdot d\}$, and 
\begin{eqnarray*} 
Z_s(X;W)\equiv Z_s\left(X;\, \ldots ,\, \xi_i\, , \ldots\, ;\, \ldots, 
\, x_j\, , \ldots\, \right):=X_s\left(\sum_{i\in I(r)}\xi_i\cdot j(H_i, 
0)+x\1\right)\, . 
\end{eqnarray*} 
\begin{proposition}\label{Proposition2.7} 
Assume (1') of Section 1. Let $1<q<\infty$ be the number defined 
in (1') and let $1/p+1/q=1$. Assume the following. 
\begin{itemize} 
\item[(i)] We have (\ref{2.29}), $\dot{f}(0,\cdot)\in I^p$.  

\item[(ii)] For $W^\rho=W+f(\rho,W)$, $\rho\in {\Bbb R}$, we have 
(\ref{2.7}) at $\rho=0$, that is $Q^{(m,r)}_{\sbnu}$-a.e. there 
exists the weak mixed derivative 
\begin{eqnarray*}
\dot{f}(0,W)(k)=\left.\frac{d^\pm}{d\rho}\right|_{\rho=0}\left 
\langle g,df^0(\rho,W)\right\rangle_{L^2}+\left\langle x,\dot{B}_0 
\right\rangle_F\, . 
\end{eqnarray*} 

Furthermore, for all test elements $k=(g,x)\in H$ with $g=g_c+g_j$ 
where $g_c\in C({\cal R} ;F)$ and $g_j:{\cal R}\to F$ being a cadlag 
pure jump function of finite jump variation and $x\in F$ we have the 
representation 
\begin{eqnarray*}
\dot{f}(0,\cdot)(k)= \left\langle k,j^{-1}\dot{f}(0,\cdot)\right 
\rangle_H\, , 
\end{eqnarray*} 
the right-hand side in the sense of Definition \ref{Definition2.6} 
(c).
\item[(iii)] $Q^{(m,r)}_{\sbnu}$-a.e., $Z_s\left(X;\, \ldots ,\, 
\xi_i\, ,\ldots\, ;\, \ldots, \, x_j\, , \ldots\, \right):l^2\to F$ 
is two times Fr\' echet differentiable for all $s\in {\cal R}$. This 
implies for $Q^{(m,r)}_{\sbnu}$-a.e. $W\in\Omega$ the existence of 
the Fr\' echet derivative $D_FX_s(W)$ in the sense of Subsection 1.2 
(after condition (5)) for which we assume that 
\begin{eqnarray}\label{2.30}
D_{F,1}X_s(W)=\left(D_{F,1}X_s(W)\right)_c+\left(D_{F,1}X_s(W)\right 
)_j  
\end{eqnarray} 
where $\left(D_{F,1}X_s(W)\right)_c\in C({\cal R} ;F)$ and $\left(D_{ 
F,1}X_s(W)\right)_j:{\cal R}\to F$ is a cadlag pure jump function of 
finite jump variation. 
\item[(iv)] For $Q^{(m,r)}_{\sbnu}$-a.e. $W\in\Omega$ there is some 
$0<c\equiv c(W)<\infty$ such that 
\begin{eqnarray*}
\left\langle h\, ,\, j^{-1}(W^\rho-W)\right\rangle_H=\left\langle 
h\, ,\, j^{-1}f(\rho,W)\right\rangle_H\le c\|h\|_H\cdot |\rho|\, , 
\quad h\in H, 
\end{eqnarray*} 
as well as 
\begin{eqnarray*}
\left\|j^{-1}(W^\rho-W)\vphantom{l^1}\right\|_H^2=\left\|j^{-1}f(\rho 
,W)\vphantom{l^1}\right\|_H^2\le c\cdot|\rho|\, ,\quad |\rho|<1. 
\end{eqnarray*} 
\end{itemize} 
Then $\D\left.\frac{d^\pm}{d\rho}\right|_{\rho=0}X_s(W^\rho)$ exists 
for $Q^{(m,r)}_{\sbnu}$-a.e. $W\in\Omega$ and all $s\in {\cal R}$ and 
we have 
\begin{eqnarray*}
&&\hspace{-.5cm}h\circ u(W)(s)=\left.\frac{d^\pm}{d\rho}\right|_{\rho= 
0}X_s(W^\rho)\nonumber \\ 
&&\hspace{.5cm}=\left.\frac{d^\pm}{d\rho}\right|_{\rho =0}\left\langle 
D_FX_s(W),j^{-1}f(\rho,W)\right\rangle_H=\left\langle D_FX_s(W),j^{-1} 
\dot{f}(0,W)\right\rangle_{H\to F}\, .  
\end{eqnarray*} 
Furthermore, this point wise calculated expression belongs to $L^1 
(\Omega ,Q^{(m,r)}_{\sbnu};F)$ whenever $\nabla_HX_s\in L^q(\Omega ,
Q^{(m,r)}_{\sbnu};H)$. 
\end{proposition}
{\bf Remark. (2)} The existence of the second order partial directional 
derivatives of $X_s$ with respect to $\kappa_i$, $i\in I(r)$, and all 
$\lambda_j$ such that 
\begin{eqnarray*}
\sum_{i,i'\in I(r)}\left(\frac{\partial^2 X_s}{\partial\kappa_i\partial 
\kappa_{i'}}\right)^2+\sum_{i\in I(r), j'}\left(\frac{\partial^2 X_s} 
{\partial\kappa_i\partial\lambda_{j'}}\right)^2<\infty
\end{eqnarray*} 
together with the continuity of 
\begin{eqnarray*}
D^2_FX_s&&\hspace{-.5cm}\equiv\sum_{i,i'\in I(r)}\frac{\partial^2 
X_s}{\partial\kappa_i\partial\kappa_{i'}}(H_i,0)\times (H_{i'},0) 
+\sum_{j,j'}\frac{\partial^2 X_s}{\partial\lambda_{j}\partial 
\lambda_{j'}}\left(0,e_{j'}\right)\times\left(0,e_{j'}\right) \\ 
&&\hspace{-.0cm}+2\sum_{i\in I(r), j'}\frac{\partial^2 X_s} 
{\partial\kappa_i\partial\lambda_{j'}}(H_i,0)\times (0,e_{j'})\in 
H\otimes H 
\end{eqnarray*} 
in a neighborhood of $W$, and this for $Q^{(m,r)}_{\sbnu}$-a.e. 
$W\in\Omega$ and all $s\in {\cal R}$, implies the first part of 
(iii). 
\medskip 

\nid
Proof. By condition (iii) of this proposition it follows now that 
for every $\ve >0$ there is $\rho_0\equiv\rho_0(W,\ve)$ such that 
\begin{eqnarray}\label{2.31}
&&\hspace{-.5cm}\left|\frac{X_s(W^\rho)-X_s(W)}{\rho}-\left\langle 
D_FX_s(W)\, ,\, \frac{j^{-1}(W^\rho-W)}{\rho}\right\rangle_H\right. 
\nonumber \\ 
&&\hspace{1cm}-\frac12\left.\left\langle\left\langle D^2_FX_s(W)\, , 
\, \frac{j^{-1}(W^\rho-W)}{\rho}\right\rangle_H\, ,\, j^{-1}(W^\rho 
-W)\right\rangle_H\right|\vphantom{\left(\int^\int_\int\right)} 
\nonumber \\ 
&&\hspace{.5cm}\le\ve\, \frac{\left\|j^{-1}(W^\rho-W)\right\|_H^2} 
{|\rho|}\, ,\quad|\rho|\le\rho_0\, .\vphantom{\sum}
\end{eqnarray} 
For the notation, cf. also Subsection 1.1, in particular the comment 
following condition (5). Using condition (iv) of this proposition it 
follows from (\ref{2.31}) that 
\begin{eqnarray*}  
&&\hspace{-.5cm}\left.\frac{d^\pm}{d\rho}\right|_{\rho =0}X_s(W^\rho) 
=\left.\frac{d^\pm}{d\rho}\right|_{\rho =0}\left\langle D_FX_s(W)\, , 
\, j^{-1}W^\rho\right\rangle_H\, .
\end{eqnarray*} 
Taking into consideration conditions (i) and (ii) of the present 
proposition we get for $Q^{(m,r)}_{\sbnu}$-a.e. $W\in\Omega$
\begin{eqnarray}\label{2.32} 
\left.\frac{d^\pm}{d\rho}\right|_{\rho =0}X_s(W^\rho)&&\hspace{-.5cm} 
=\left.\frac{d^\pm}{d\rho}\right|_{\rho =0}\left\langle D_FX_s(W),j^{ 
-1}f(\rho,W)\right\rangle_H\nonumber \\ 
&&\hspace{-.5cm}=\left\langle D_FX_s(W),j^{-1}\dot{f}(0,W)\right 
\rangle_{H\to F}\, .  
\end{eqnarray} 
In addition, we note that $D_FX_s=\nabla_HX_s$ by (iii) of the present 
proposition. Accordingly, the right-hand side of (\ref{2.32}) is by 
(i) of the present proposition an element of $L^1(\Omega,Q^{(m,r)}_{ 
\sbnu};F)$ whenever $\nabla_HX_s\in L^q(\Omega,Q^{(m,r)}_{\sbnu};H)$. 
\qed 
\bigskip 

{\bf Stochastic integral in a distributional sense. } Let us assume 
$\dot{f}(0,\cdot)(s)\in D_{p,1}$, $s\in {\cal R}$. We have 
\begin{eqnarray*}
\left\langle H_i,d\dot{f}(0,\cdot )\right\rangle_{L^2}\in D_{p,1}\, , 
\quad i\in I(m,r),   
\end{eqnarray*} 
by the special shape of the $H_i$. Consequently 
\begin{eqnarray*} 
D\left\langle H_i,d\dot{f}(0,W)\right\rangle_{L^2}\in H\, , \quad i\in 
I(m,r), 
\end{eqnarray*} 
for $Q^{(m,r)}_{\sbnu}$-a.e. $W\in\Omega$. With this we may redefine 
(\ref{2.12}) $Q^{(m,r)}_{\sbnu}$-a.e. by 
\begin{eqnarray}\label{2.33}
\hat{\delta}(j^{-1}\dot{f}(0,\cdot ))(W)\vphantom{\sum_{i\in I(m,r)}} 
&&\hspace{-.5cm}:=\sum_{i\in I(m,r)}\left\langle H_i,d\dot{f}(0,W) 
\right\rangle_{L^2}\cdot\langle H_i,dW\rangle_{L^2}\nonumber \\ 
&&\hspace{-0cm}-\sum_{i\in I(m,r)}\left\langle (H_i,0),D\left\langle 
H_i,d\dot{f}(0,W)\right\rangle_{L^2}\right\rangle_H\nonumber \\ 
&&\hspace{-0cm}-\left\langle\dot{f}(0,W)(0),\frac{\nabla m(W_0)}{m 
(W_0)}\right\rangle_F-\left\langle {\sf e},\nabla_{d,W_0}\dot{f}(0,W) 
(0)\right\rangle_F 
\end{eqnarray} 
where we note that we just have shown the well-definiteness of the 
right-hand side. 

The left-hand side of (\ref{2.33}) is only a symbol reminding of the 
identity (\ref{2.12}) demonstrated under the more restrictive setting 
of Subsection 2.2. But this will change under the specification of 
Section 3. In fact, having proved the subsequent Lemmas \ref{Lemma3.2} 
and \ref{Lemma3.3}, the right-hand side of (\ref{2.33}) will become 
the object that plays the role of a stochastic integral.  
\section{The Density Formula} 
\setcounter{equation}{0}

Besides the proof of Theorem \ref{Theorem1.11}, this section contains 
a number of technical preparations. There are two major technical ideas 
in order to prove Theorem \ref{Theorem1.11}. The first one is to re-map 
$X$ to $W$, to carry out the analysis on Wiener process related objects 
and then to map the result obtained there back to a result in terms of 
$X$. The second idea is to work on orthogonal projections of paths of 
$W$ rather than on the trajectories of $W$ itself. In Subsection 3.1 we 
will specify the flow $W^\rho:=W_{\cdot +\rho}+(A_\rho-Y_\rho)\1$, ${ 
\rho\in {\Bbb R}}$, and provide technical calculations relative to the 
finite dimensional projections of this flow. In Subsection 3.2 we will 
be interested in the stochastic integral relative to (\ref{2.33}). Then, 
in Subsection 3.3, we will collect a number of technical details of the 
approximation of the flow $W^\rho=W_{\cdot +\rho}+(A_\rho -Y_\rho)\1$, 
${\rho\in {\Bbb R}}$. The actual proof of Theorem \ref{Theorem1.11} 
will be presented in Subsection 3.4. 

Subsections 3.1 and 3.4 are entirely preparatory for the proof of 
Theorem \ref{Theorem1.11} whereas Subsections 3.2 and 3.3 provide some 
more insight on the related stochastic calculus. 

\subsection{Specification of $f$} 

One goal of this subsection is to specify $f$ in (\ref{2.5}), and with 
this also $h$, to those functionals we are going to use in the proof of 
Theorem \ref{Theorem1.11}. In particular, we will consider the flow 
\begin{eqnarray}\label{3.1}
X(W^\rho)=u(W^\rho)=u(W_{\cdot +\rho}+(A_\rho-Y_\rho)\1)=X(W_{\cdot +\rho} 
+(A_\rho-Y_\rho)\1)\, , \quad \rho\in {\Bbb R}.  
\end{eqnarray} 
Correspondingly, we also will suppose 
\begin{eqnarray}\label{3.2}
f(\rho,W):=W_{\cdot+\rho}+(A_\rho-Y_\rho)\1-W  
\end{eqnarray} 
which includes that $W^\rho=W+f(\rho,W)=W_{\cdot +\rho}+(A_\rho-Y_\rho) 
\1$, $\rho\in {\Bbb R}$, is a flow, cf. Remark (3) of Section 1. 
In the Lemmas \ref{Lemma3.1} and \ref{Lemma3.2} we will assume 
(\ref{3.1}) and (\ref{3.2}). Let us show that $f$ is compatible with 
the analysis of Subsection 2.3 and, in particular, satisfies the 
conditions (i) and (ii) of Proposition \ref{Proposition2.7}. Let us 
also recall that the space $D_{p',1}\equiv D_{p',1}(Q^{(m,r)}_{\sbnu} 
)$ is meaningfully defined for $p'\ge Q/(Q-1)$, cf. Definition 
\ref{Definition1.1} and Proposition \ref{Proposition6.3}. 
\begin{lemma}\label{Lemma3.1} 
Assume conditions (1'), (2), and (3) of Section 1 which include in 
particular (\ref{3.1}) as well as (\ref{3.2}). Furthermore, assume 
$A=A^1$ and $Y=Y^1$. Let $q$ be the number introduced in condition 
(1') of Section 1 and suppose $Y_\rho\in D_{q,1}\equiv D_{q,1}(Q^{ 
(m,r)}_{\sbnu})$, $\rho\in {\Bbb R}$.  
\medskip 

\nid
(a) For $s\in {\cal R}$, $i\in I(m,r)$, and $\rho\in {\Bbb R}$, we 
have $W_s,\langle H_i,d\dot{W}\rangle_{L^2}\in\bigcap_{p'\in [Q/(Q-1), 
\infty)}D_{p',1}$ and $W_{s+\rho}+A_\rho-Y_\rho\in D_{q,1}$ where $q$ 
is given by condition (1') of Section 1. Furthermore, it holds that 
\begin{itemize} 
\item[{}] $\displaystyle DW_s=\sum_{i\in I(m,r)}\int_0^s H_i(v)\, dv 
\cdot (H_i,0)+(0,{\sf e})$ for $Q^{(m,r)}_{\sbnu}$-a.e. $W\in\Omega$  
and that 
\item[{}] $D\langle H_i,d\dot{W}\rangle_{L^2}=\sum_{j\in I(m,r)}\langle 
H_i ,dH_j\rangle_{L^2}\cdot (H_j,0)$ for $Q^{(m,r)}_{\sbnu}$-a.e. $W 
\in\Omega$. 
\end{itemize} 
(b) $\left(W^\rho\right)_{\rho\in {\Bbb R}}$ satisfies (\ref{2.6}) and 
(\ref{2.7}). \\ 
(c) Let $1/p+1/q=1$. If $(\dot{A}_\rho -\dot{Y}_\rho)\in L^p (\Omega 
,Q^{(m,r)}_{\sbnu};F)$ then we have (i) and (ii) of Proposition 
\ref{Proposition2.7}, the latter demonstrating the compatibility of 
the weak mixed derivative in (\ref{2.7}) with the integral of 
Definition \ref{Definition2.6}. \\ 
(d) We have (iv) of Proposition \ref{Proposition2.7}. 
\end{lemma}
Proof. (a) This is straight forward by the following. By Lemma 
\ref{Lemma2.2} (a) we have
\begin{eqnarray*} 
DW_s=\sum_{i\in I(m,r)}\frac{\partial W_s}{\partial\kappa_i}\cdot ( 
H_i,0)+\left(0,\nabla_{W_0}W_s\right)=\sum_{i\in I(m,r)}\int_0^s H_i 
(v)\, dv\cdot (H_i,0)+(0,{\sf e})\, .   
\end{eqnarray*}
Recalling (\ref{2.2}) we get under the measure $Q^{(m,r)}_{\sbnu}$  
\begin{eqnarray*} 
D\langle H_i ,d\dot{W}\rangle_{L^2}=\sum_{j\in I(m,r)}\frac{\partial 
\langle H_i ,d\dot{W}\rangle_{L^2}}{\partial\kappa_j}\cdot (H_j ,0)= 
\sum_{j\in I(m,r)}\langle H_i ,dH_j\rangle_{L^2}\cdot (H_j,0)\, . 
\end{eqnarray*}
Thus $W_s,\langle H_i,d\dot{W}\rangle_{L^2}\in\bigcap_{p'\in [Q/(Q-1) 
,\infty)}D_{p',1}$ is now clear, in particular by the bounded support 
of $m$, cf. Subsection 1.2.

$W_{s+\rho}+A_\rho-Y_\rho\in D_{q,1}$ is a consequence of $A_\rho= 
X_\rho-W_\rho$ as well as $X_\rho\in D_{q,1}$ by condition (1') of 
Section 1 and $Y_\rho\in D_{q,1}$ by hypothesis. 
\medskip 

\nid  
(b) Let us verify (\ref{2.6}). By definition, $W^\sigma=W^\rho$ for 
some $\sigma\neq\rho$ would imply $W_{\cdot +\sigma}-W_{\cdot +\rho} 
=\left(A_\rho(W)-Y_\rho(W)-A_\sigma(W)+Y_\sigma(W)\right)\1$. Assuming 
this, the right-hand side was $Q^{(m,r)}_{\sbnu}$-a.e. constant in 
time but the left-hand side was not. 

Let us focus on (\ref{2.7}). We have $f(\rho,W)(s)=W_{s+\rho}+(A_\rho 
-Y_\rho)\1(s)-W_s$ which means in the setting of (\ref{2.7}) 
\begin{eqnarray*} 
\begin{array}{rl}
f^0(\rho,W)(s)&=W_{s+\rho}-W_s-W_\rho+W_0\, , \vphantom{\D\left( 
\dot{f}\right)} \\ 
B_\rho&=W_\rho-W_0+A_\rho -Y_\rho\, . \vphantom{\D\left(\dot{f} 
\right)}
\end{array}
\end{eqnarray*}
Keeping (\ref{2.8}) in mind we shall verify (\ref{2.7}) with 
\begin{eqnarray}\label{3.3}  
\dot{f}(\rho,W)(s)&&\hspace{-.5cm}:=\sum_{i\in I(m,r)}\langle H_i 
,dW\rangle_{L^2}\cdot H_i(s+\rho)+\dot{A}_\rho-\dot{Y}_\rho 
\nonumber \\ 
&&\hspace{-.5cm}=\sum_{i\in I(m,r)}\langle H_i,dW\rangle_{L^2} 
\cdot\left(H_i(s+\rho)-\frac12\left(H_i(\rho-)+H_i(\rho)\vphantom 
{l^1}\right)\right)+\dot{W}_\rho +\dot{A}_\rho-\dot{Y}_\rho\, , 
\nonumber \\ 
\end{eqnarray}
$s\in {\cal R}$, $\rho\in {\Bbb R}$. We have  
\begin{eqnarray*}
&&\hspace{-.5cm}\frac{1}{\lambda}\left(\left\langle g,df^0(\rho+ 
\lambda,W)\right\rangle_{L^2}-\left\langle g,df^0(\rho,W)\right 
\rangle_{L^2}\right)+\frac{1}{\lambda}\left(\left\langle x,B_{\rho 
+\lambda}\right\rangle_F-\left\langle x,B_\rho\right\rangle_F 
\right) \\ 
&&\hspace{.5cm}=\sum_{i\in I(m,r)}\langle H_i,dW\rangle_{L^2}\cdot 
\frac{1}{\lambda}\left(\vphantom{\dot{f}}\left\langle g,H_i(\cdot+ 
\rho +\lambda)-\frac12\left(H_i(\lambda+\rho-)+H_i(\lambda+\rho) 
\vphantom{l^1}\right)\cdot\1\right\rangle_{L^2}\right. 
 \\ 
&&\hspace{7cm}\left.-\left\langle g,H_i(\cdot+\rho)-\frac12\left 
(H_i(\rho-)+H_i(\rho)\vphantom{l^1}\right) 
\cdot\1\vphantom{f^0}\right\rangle_{L^2}\vphantom{\dot{f}}\right) 
 \\ 
&&\hspace{1cm}+\frac{1}{\lambda}\left(\vphantom{\dot{f}}\left 
\langle x,\left(\vphantom{f^0}W_{\rho+\lambda}+A_{\rho+\lambda}- 
Y_{\rho+\lambda}\right)-\left(W_{\rho}+A_{\rho}-Y_{\rho}\vphantom 
{f^0}\right)\right\rangle_F\right) 
\end{eqnarray*}
and therefore 
\begin{eqnarray}\label{3.4}
&&\hspace{-.5cm}\frac{d^\pm}{d\rho}\left\langle g,df^0(\rho,W) 
\right\rangle_{L^2}+\frac{d^\pm}{d\rho}\left\langle x,B_{\rho} 
\right\rangle_F\nonumber \\ 
&&\hspace{.5cm}=\sum_{i\in I(m,r)}\langle H_i,dW\rangle_{L^2} 
\cdot\left\langle g,dH_i(\cdot+\rho)\right\rangle_{L^2}+\left 
\langle x,\dot{W}_\rho+\dot{A}_\rho-\dot{Y}_\rho\right\rangle_F 
\vphantom{\int}\nonumber \\ 
&&\hspace{.5cm}=\left\langle g,d\dot{f}(\rho,W)\right\rangle_{ 
L^2}+\left\langle x,\dot{W}_\rho+\dot{A}_\rho-\dot{Y}_\rho 
\right\rangle_F\, , 
\end{eqnarray}
$k=(g,x)\in H$, $g:{\cal R}\to F$ cadlag of finite jump variation, 
$x\in F$, $\rho\in {\Bbb R}$ $Q^{(m,r)}_{\sbnu}$-a.e., as a direct 
consequence of (\ref{3.3}). In other words, we have (\ref{2.7}).  
\medskip

\nid  
(c) For condition (i) of Proposition \ref{Proposition2.7}, $\dot 
{f}(\rho,\cdot)\in I^p$, we take into consideration (\ref{3.3}). 
In particular, $(\dot{f}(\rho,\cdot))_a=(\dot{W}_\rho +\dot{A 
}_\rho-\dot{Y}_\rho)\1$. Recall also $(\dot{A}_\rho -\dot{Y}_\rho 
)\in L^p(\Omega ,Q^{(m,r)}_{\sbnu};F)$ by hypothesis. We get 
\begin{eqnarray*}
j^{-1}\left(\dot{f}(\rho,\cdot)\right)_a=(0,\dot{W}_\rho +\dot{A 
}_\rho-\dot{Y}_\rho)\in L^p (\Omega ,Q^{(m,r)}_{\sbnu};H)\, . 
\end{eqnarray*}
Furthermore, $\sum_{i\in I(m,r)}\langle H_i,dW\rangle_{L^2}\cdot 
\left(H_i(\cdot +\rho)-\frac12\left(H_i(\rho-)+H_i(\rho)\vphantom 
{l^1}\right)\right)$ is a cadlag pure jump process with $Q^{(m,r) 
}_{\sbnu}$-a.e. finitely many jumps whose jump times are 
non-random and whose magnitudes in the sense of Definition 
\ref{Definition1.4} (iii) belong to $L^p(\Omega,Q^{(m,r)}_{\sbnu}) 
$. Thus we have verified that $\dot{f}(\rho,\cdot)\in I^p$. 

For condition (ii) of Proposition \ref{Proposition2.7}, we recall 
now (\ref{3.4}) and the calculation before (\ref{3.4}) which now 
say 
\begin{eqnarray*}
\left\langle k,j^{-1}\dot{f}(\rho ,W)\right\rangle_H 
&&\hspace{-.5cm}=\frac{d^\pm}{d\rho}\left\langle g,df^0(\rho,W) 
\right\rangle_{L^2}+\left\langle x,\dot{B}_\rho\right\rangle_F\, ,  
\end{eqnarray*}
the left-hand side in the sense of Definition \ref{Definition2.6}. 
For compatibility with this definition, we restrict ourselves here 
to all test elements $k=(g,x)\in H$ with $g=g_c+g_j$ where $g_c\in 
C({\cal R} ;F)$ and $g_j:{\cal R}\to F$ being a cadlag pure jump 
function of finite jump variation and $x\in F$. 
\medskip

\nid  
(d) This follows directly from 
\begin{eqnarray*}
j^{-1}f(\rho,\cdot)(s)&&\hspace{-.5cm}=\left(\sum_{i\in I(m,r)} 
\langle H_i,dW\rangle_{L^2}\cdot\left(H_i(s+\rho)-H_i(s)\vphantom{ 
\dot{f}}\right)\, ,\, W_\rho +A_\rho -Y_\rho-W_0\right)\, , 
\end{eqnarray*}
hypotheses $A=A^1$ and $Y=Y^1$ together with conditions (2) (i) and 
(3) (iii) of Section 1, and $A_0=Y_0$. 
\qed 
\bigskip 

In order to prepare the crucial equalities in the proof of Theorem 
\ref{Theorem1.11}, namely (\ref{3.25}) and (\ref{3.26}) below, let us 
prove the subsequent lemma. 
\medskip 

For this, let $\vp$ be a cylindrical function on $C({\Bbb R};F)$ of the 
form $\vp (W)=f_0(W_0)\cdot f_1(W_{t_1}-W_0,\ldots ,W_{t_k}-W_0)$, $W\in 
\Omega$, where $f_0\in C_b^1(F)$, $f_1\in C_b^1(F^k)$, and $k\in {\Bbb 
N}$ such that $t_l\in {\cal R}\cap\{z\cdot t/2^m:z\in {\Bbb Z}\setminus 
\{0\}\}$, $l\in \{1,\ldots ,k\}$. We will also use the abbreviation 
$(D_F-D)(W_\sigma +A_\sigma -Y_\sigma):=D_F(W_\sigma +A_\sigma -Y_\sigma 
)-D(W_\sigma +A_\sigma -Y_\sigma)$ and the following integrals. Let 
$v\in [0,t]$. 
\begin{eqnarray*} 
L_1&&\hspace{-.5cm}:=\int\int_{\sigma=0}^v\left\langle D\vp\circ (W_{ 
\cdot +\sigma}+(A_\sigma -Y_\sigma)\1)\, ,\, j^{-1}\left(\dot{W}+ 
\left(\dot{A}_0-\dot{Y}_0\right)\1\right)\right\rangle_H\, d\sigma\, 
Q^{(m,r)}_{\sbnu}(dW) \\  
L_2&&\hspace{-.5cm}:=\int\int_{\sigma=0}^v\left(\vphantom{\int}f_1(W_{ 
t_1+\sigma}-W_\sigma,\ldots , W_{t_k+\sigma}-W_\sigma)\left\langle 
\vphantom{\dot{f}}\left\langle\nabla f_0(W_\sigma+ A_\sigma -Y_\sigma) 
\, ,\, \vphantom{l^1}\right.\right.\right. \\ 
&&\hspace{.5cm}\left.\left.\left.\vphantom{l^1} (D_F-D)(W_\sigma + 
A_\sigma -Y_\sigma)\right\rangle_{F\to F}\, ,\, j^{-1}\left(\dot{W}+ 
\left(\dot{A}_0-\dot{Y}_0\right)\1\right)\right\rangle_H\vphantom 
{\int}\right)\, d\sigma\, Q^{(m,r)}_{\sbnu}(dW)\vphantom{\int_{ 
\sigma=0}^\infty} 
\end{eqnarray*}
\begin{eqnarray*} 
R_1&&\hspace{-.5cm}:=\sum_{i\in I(r)}\int_{\sigma=0}^v\int\vp\left( 
[W_{\cdot+\sigma}+(A_\sigma -Y_\sigma)\1]\circ\pi_{m,r}\vphantom{l^1} 
\right)\cdot\left\langle H_i,d\left([W]\circ\pi_{m,r}\right)^\cdot 
\right\rangle_{L^2}\times \\ 
&&\hspace{.5cm}\vphantom{\int_{\sigma=0}^v}\times\langle H_i\, ,\, dW 
\rangle_{L^2}\, Q_{\sbnu}(dW)\, d\sigma \\ 
R_2&&\hspace{-.5cm}:=-\sum_{i\in I(r)}\int_{\sigma=0}^v\int\vp\left( 
[W_{\cdot +\sigma}+(A_\sigma-Y_\sigma)\1]\circ\pi_{m,r}\vphantom{l^1} 
\right)\times \\ 
&&\hspace{.5cm}\vphantom{\int_{\sigma=0}^v}\times\left\langle (H_i,0) 
\, ,\, D_F\left\langle H_i,d\left([W]\circ\pi_{m,r}\right)^\cdot\right 
\rangle_{L^2}\vphantom{\dot{f}}\right\rangle_H\, Q_{\sbnu}(dW)\, d 
\sigma 
\end{eqnarray*}
\begin{eqnarray*} 
J&&\hspace{-.5cm}:=\int_{\sigma=0}^v\int\left\langle D\vp\circ\left( 
W_{\cdot +\sigma}+(A_\sigma -Y_\sigma)\1\vphantom{l^1}\right)\, ,\, 
\left(0,\dot{W}_0+\dot{A}_0-\dot{Y}_0\right)\right\rangle_H\, Q^{(m,r) 
}_{\sbnu}(dW)\, d\sigma\, . 
\end{eqnarray*}
Well-definiteness of these integrals will be demonstrated in parts 
(d)-(f) of the following lemma. 
\begin{lemma}\label{Lemma3.2} 
Assume conditions (1'), (2), and (3) of Section 1 which in particular 
imply (\ref{3.1}) as well as (\ref{3.2}). Let $q$ be the number 
introduced in condition (1') of Section 1. Furthermore, suppose 
$Y_\rho\in D_{q,1}(Q^{(m,r)}_{\sbnu})$, $\rho\in {\Bbb R}$. Also assume 
(iii) of Proposition \ref{Proposition2.7} for $X$ as well as for $X$ 
replaced with $Y$. 
\medskip 

\nid 
(a) We have 
\begin{eqnarray*} 
\vp\circ\left(W_{\cdot +\sigma}+(A_\sigma-Y_\sigma)\1\vphantom{l^1}\right) 
\in D_{q,1}\, , \quad\sigma\in [0,t], 
\end{eqnarray*}
with respect to the measure $Q^{(m,r)}_{\sbnu}$. 
\medskip 

\nid 
(b) Assume $A=A^1$ and $Y=Y^1$ $Q^{(m,r)}_{\sbnu}$-a.e. For $v\in [0,t 
]$ we have $Q^{(m,r)}_{\sbnu}$-a.e. 
\begin{eqnarray*}
&&\hspace{-.5cm}\vp\left(W_{\cdot +v}+(A_v-Y_v)\1\vphantom{l^1}\right) 
-\vp(W) \\ 
&&\hspace{.5cm}=\int_{\sigma=0}^v\left\langle(D\vp)\left(\pi_{m,r} 
\left(W_{\cdot +\sigma}+(A_\sigma-Y_\sigma)\1\vphantom{l^1}\right) 
\right),j^{-1}\left(\pi_{m,r}\left(W_{\cdot +\sigma}+(A_\sigma-Y_\sigma) 
\1\vphantom{l^1}\right)\right)\, \dot{}\, \vphantom{\dot{f}}\right 
\rangle_H\, d\sigma\, . 
\end{eqnarray*} 
(c) Assume $A=A^1$ and $Y=Y^1$ $Q^{(m,r)}_{\sbnu}$-a.e. and let $\sigma 
\in [0,t]$. We have $Q^{(m,r)}_{\sbnu}$-a.e. 
\begin{eqnarray*} 
&&\hspace{-.5cm}\left\langle(D\vp)\left(\pi_{m,r}\left(W_{\cdot +\sigma} 
+(A_\sigma-Y_\sigma)\1\vphantom{l^1}\right)\right),j^{-1}\left(\pi_{m,r} 
\left(W_{\cdot +\sigma}+(A_\sigma-Y_\sigma)\1\vphantom{l^1}\right)\right 
)\, \dot{}\, \vphantom{\dot{f}}\right\rangle_H \\ 
&&\hspace{.5cm}=\left\langle D\vp\circ (W_{\cdot +\sigma}+(A_\sigma 
-Y_\sigma)\1)\, ,\, j^{-1}\left(\dot{W}+\left(\dot{A}_0-\dot{Y}_0\right) 
\1\right)\right\rangle_H \\ 
&&\hspace{1cm}+f_1(W_{t_1+\sigma}-W_\sigma,\ldots , W_{t_k+\sigma}- 
W_\sigma)\left\langle\vphantom{\dot{f}}\left\langle\nabla f_0(W_\sigma+ 
A_\sigma -Y_\sigma)\, ,\right.\right.\vphantom{\sum} \\ 
&&\hspace{1cm}\left.\left. (D_F-D)(W_\sigma +A_\sigma-Y_\sigma)\right 
\rangle_{F\to F}\, ,\, j^{-1}\left(\dot{W}+\left(\dot{A}_0-\dot{Y}_0\right 
)\1\right)\right\rangle_H\, . 
\end{eqnarray*}
For the well-definiteness of these integral terms recall Definition 
\ref{Definition2.6}. 
\medskip 

\nid 
%
(d) Let the assumptions of part (c) be in force. Suppose 
$(\dot{A}_\rho -\dot{Y}_\rho)\in L^p(\Omega,Q^{(m,r)}_{\sbnu};F)$, 
$\rho\in {\Bbb R}$, where $1/p+1/q=1$.

For $\sigma\in [0,t]$, let $\nabla_HX_\sigma\in L^q(\Omega ,Q^{(m,r) 
}_{\sbnu};H)$ and $\nabla_HY_\sigma\in L^q(\Omega ,Q^{(m,r)}_{\sbnu}; 
H)$. Suppose that the corresponding norms belong, with respect to 
$\sigma\in [0,t]$, to $L^1([0,t])$. The above defined integrals $L_1, 
L_2,R_1,R_2,J$ are well-defined and we have 
\begin{eqnarray*} 
L_1+L_2=R_1+R_2+J\, . 
\end{eqnarray*}
(e) For all $i\in I(m,r)$, we have 
\begin{eqnarray*} 
\vp\left(W_{\cdot +\sigma}+(A_\sigma-Y_\sigma)\1\vphantom{l^1}\right) 
\langle H_i,d\dot{W}\rangle_{L^2}\in D_{q',1}\, ,\quad\sigma\in [0,t] 
\vphantom{\sum^k_l}
\end{eqnarray*}
if $Q/(Q-1)<q'<q<\infty$. Furthermore, $(H_i,0)\in {\rm Dom}_{p'}(\delta 
)$ for all $i\in I(m,r)$ and all $p<p'<Q$ and $\delta (H_i,0)=\langle 
H_i,dW\rangle_{L^2}\in\bigcap_{p'\in [1,\infty)}L^{p'}(\Omega ,Q^{(m,r) 
}_{\sbnu})$ as well as 
\begin{eqnarray*} 
R_1&&\hspace{-.5cm}=\sum_{i\in I(m,r)}\int_{\sigma=0}^v\int\vp\left( 
W_{\cdot+\sigma}+(A_\sigma-Y_\sigma)\1\vphantom{l^1}\right)\cdot\left 
\langle H_i,d\dot{W}\right\rangle_{L^2}\langle H_i\, ,\, dW\rangle_{ 
L^2}\, Q^{(m,r)}_{\sbnu}(dW)\, d\sigma
\end{eqnarray*}
and
\begin{eqnarray*} 
R_2&&\hspace{-.5cm}=-\sum_{i\in I(m,r)}\int_{\sigma=0}^v\int\vp\left 
(W_{\cdot +\sigma}+(A_\sigma-Y_\sigma)\1\vphantom{l^1}\right)\cdot 
\left\langle D\left\langle H_i,d\dot{W}\right\rangle_{L^2},(H_i,0) 
\vphantom{\dot{f}}\right\rangle_H\, Q^{(m,r)}_{\sbnu}(dW)\, d\sigma 
\, . 
\end{eqnarray*}
(f) Assume $A=A^1$ and $Y=Y^1$ $Q^{(m,r)}_{\sbnu}$-a.e. and let $Q$ be 
the number that is used in the definition of the density $m$. Let $1/p 
+1/q=1$ with $p\ge q$ and assume $(\dot{A}_\rho -\dot{Y}_\rho)\in L^p 
(\Omega,Q^{(m,r)}_{\sbnu};F)$, $\rho\in {\Bbb R}$. Suppose 
the existence of the derivative  
\begin{eqnarray*} 
\nabla_{d,W_0}\left(\dot{A}_0-\dot{Y}_0\vphantom{l^1}\right)\quad\mbox{ 
\it in}\quad L^q(\Omega,Q^{(m,r)}_{\sbnu};F)\, . 
\end{eqnarray*}
We have 
$(0,\dot{W}_0+\dot{A}_0-\dot{Y}_0)\in {\rm Dom}_w(\delta)$ with $w= 
\min\left(q,(1/p+1/Q)^{-1}\right)$ in the sense of Subsection 1.2 and 
\begin{eqnarray*}
\delta\left(0,\dot{W}_0+\dot{A}_0-\dot{Y}_0\right)=-\left\langle\dot 
{W}_0+\dot{A}_0-\dot{Y}_0,\frac{\nabla m(W_0)}{m(W_0)}\right\rangle_F 
-\left\langle {\sf e},\nabla_{d,W_0}\left(\dot{A}_0-\dot{Y}_0\right) 
\right\rangle_F 
\end{eqnarray*} 
$\in L^w(\Omega,Q^{(m,r)}_{\sbnu})$. Equivalently, we have 
\begin{eqnarray}\label{3.5}
&&\hspace{-.5cm}\int\left\langle\dot{W}_0+\dot{A}_0-\dot{Y}_0,(D\psi 
)_2(W)\right\rangle_F\, Q^{(m,r)}_{\sbnu}(dW)\nonumber \\ 
&&\hspace{.5cm}=-\int\left(\left\langle\dot{W}_0+\dot{A}_0-\dot{Y}_0, 
\frac{\nabla m(W_0)}{m(W_0)}\right\rangle_F+\left\langle {\sf e}, 
\nabla_{d,W_0}\left(\dot{A}_0-\dot{Y}_0\right)\right\rangle_F\right) 
\cdot\psi(W)\, Q^{(m,r)}_{\sbnu}(dW)\nonumber \\ 
\end{eqnarray} 
for all $\psi\in D_{v,1}$ where $1/v+1/w=1$. In particular, 
\begin{eqnarray*}
J=\int\int_{\sigma=0}^v\vp\left(W_{\cdot +\sigma}+(A_\sigma-Y_\sigma) 
\1\vphantom{l^1}\right)\delta\left(0,\dot{W}_0+\dot{A}_0-\dot{Y}_0 
\right)\, d\sigma\, Q^{(m,r)}_{\sbnu}(dW)\, . 
\end{eqnarray*} 
\end{lemma}
Proof. (a) Actually, everything to prove can be taken from the 
differential calculus  of real functions in finite dimensional 
domains. However, the following more extensive presentation will 
also be a preparation of the proof of part (c). Keeping Definition 
\ref{Definition1.1}, (\ref{6.5}), (\ref{6.6}), Lemma \ref{Lemma2.2} 
(a), and the special choice of the $t_l$ in mind, we obtain  
\begin{eqnarray}\label{3.6}
\left\langle (D\vp)(W),j^{-1}\rho\right\rangle_H&&\hspace{-.5cm}= 
f_0(W_0)\cdot\sum_{l=1}^k\langle\nabla_l f_1(W_{t_1}-W_0,\ldots , 
W_{t_k}-W_0)\, , \, \rho_{t_l}-\rho_0\rangle_F\nonumber \\ 
&&\hspace{-.0cm}+f_1(W_{t_1}-W_0,\ldots ,W_{t_k}-W_0)\cdot\langle 
\nabla f_0(W_0),\rho_0\rangle_F\vphantom{\displaystyle\sum} 
\end{eqnarray} 
with $\rho\equiv (\rho -\rho_0\cdot\1,\rho_0)\in\{j(f,x):(f,x)\in 
H\}$ and $\nabla_l$ denoting the gradient with respect to the $l 
$-th entry. Accordingly,  
\begin{eqnarray}\label{3.7}
&&\hspace{-0.5cm}D\vp\circ (W_{\cdot +\sigma}+(A_\sigma-Y_\sigma) 
\1)=f_0(W_\sigma +A_\sigma-Y_\sigma)\times\nonumber \\ 
&&\hspace{2cm}\times\sum_{l=1}^k\langle\nabla_l f_1(W_{t_1+\sigma 
}-W_\sigma,\ldots ,W_{t_k+\sigma}-W_\sigma),D(W_{t_l+\sigma}- 
W_\sigma)\rangle_{F\to F} \nonumber \\ 
&&\hspace{1cm}+f_1(W_{t_1+\sigma}-W_\sigma,\ldots ,W_{t_k+\sigma} 
-W_\sigma)\times\vphantom{\sum_1}\nonumber \\ 
&&\hspace{2cm}\times\left\langle\nabla f_0(W_\sigma+A_\sigma- 
Y_\sigma),D(W_{\sigma}+A_\sigma-Y_\sigma)\right\rangle_{F\to F}\, .     
\end{eqnarray}
%
Now, by Lemma \ref{Lemma3.1} (a), we have $\vp\circ (W_{\cdot + 
\sigma}+(A_\sigma-Y_\sigma)\1)\in D_{q,1}$. 
\medskip 

\nid
(b) We use the fact that, because of its special form, we have 
$\vp\circ\pi_{m,r}=\vp$. Thus $Q^{(m,r)}_{\sbnu}$-a.e. it holds that 
\begin{eqnarray*} 
&&\hspace{-.5cm}\vp\left(W_{\cdot +v}+(A_v-Y_v)\1\vphantom{l^1}\right) 
-\vp(W)\vphantom{\int} \\ 
&&\hspace{.5cm}=\vp\circ\pi_{m,r}\left(W_{\cdot +v}+(A_v-Y_v)\1 
\vphantom{l^1}\right)-\vp(W) \vphantom{\int} \\ 
&&\hspace{.5cm}=\int_{\sigma=0}^v\frac{d^\pm}{d\sigma}\vp\circ\pi_{m, 
r}\left(W_{\cdot +\sigma}+(A_\sigma -Y_\sigma)\1\vphantom{l^1}\right) 
\, d\sigma \\ 
&&\hspace{.5cm}=\int_{\sigma=0}^v\left\langle(D\vp)\left(\pi_{m,r} 
\left(W_{\cdot +\sigma}+(A_\sigma-Y_\sigma)\1\vphantom{l^1}\right) 
\right),j^{-1}\left(\pi_{m,r}\left(W_{\cdot +\sigma}+(A_\sigma- 
Y_\sigma)\1\vphantom{l^1}\right)\right)\, \dot{}\, \vphantom{\dot{f}} 
\right\rangle_H\, d\sigma\, .  
\end{eqnarray*} 
Here we have also applied Lemma \ref{Lemma2.2} (a), the special form 
of $H_i$, $i\in I(m,r)$, and the differential calculus similar as in 
the proof of Proposition \ref{Proposition2.7}. 
\medskip 

\nid
(c) By Lemma \ref{Lemma3.1} (a) we have $Q^{(m,r)}_{\sbnu}\mbox 
{\rm -a.e.}$
\begin{eqnarray*} 
D\left(W_{t_l+\sigma}-W_\sigma\right)&&\hspace{-.5cm}=\sum_{i\in 
I(m,r)}\frac{\partial\left(W_{t_l+\sigma}-W_\sigma\right)}{\partial 
\kappa_i}\cdot (H_i,0) \\ 
&&\hspace{-.5cm}=\sum_{i\in I(m,r)}\int_\sigma^{t_l+\sigma}H_i(s)
\, ds\cdot (H_i,0)\, ,\quad l\in \{1,\ldots ,k\}, 
\end{eqnarray*}
which is by the definition of $I(m,r)$ also the orthogonal projection 
in $H$ of $({\sf e}\1_{[\sigma,t_l+\sigma)},0)$ to the subspace spanned 
by the functions $\left(\1_{[z\cdot t/2^m,(z+1)\cdot t/2^m)}\cdot e_j,0 
\right)$, $z\in {\Bbb Z}$, $j\in\{1,\ldots ,n\cdot d\}$. With 
\begin{eqnarray*} 
\sigma^-(m)\equiv\sigma^-(m,\sigma):=\max\{z\cdot t/2^m:z\in {\Bbb Z} 
,\ z\le\sigma\} 
\end{eqnarray*}
we obtain $Q^{(m,r)}_{\sbnu}$-a.e. for $l\in\{1,\ldots ,k\}$ 
\begin{eqnarray*} 
&&\hspace{-.5cm}D\left(W_{t_l+\sigma}-W_\sigma\right)=\left({\sf e} 
\1_{[\sigma^-(m),t_l+\sigma^-(m))},0\right)\vphantom{\frac12} \\ 
&&\hspace{0.5cm}+2^m\, \frac{\sigma-\sigma^-(m)}{t}\cdot\left(\left 
({\sf e}\1_{[t_l+\sigma^-(m),t_l+\sigma^-(m)+t/2^m)},0\right)-\left( 
{\sf e}\1_{[\sigma^-(m),\sigma^-(m)+t/2^m)},0\right)\vphantom{\dot 
{f}}\right)\, . 
\end{eqnarray*}
Recalling now the definition of the {\it weak mixed derivative} and 
taking Definition \ref{Definition2.6} into consideration it turns 
out that $Q^{(m,r)}_{\sbnu}$-a.e.  
\begin{eqnarray*} 
&&\hspace{-.5cm}\left\langle D\left(W_{t_l+\sigma}-W_\sigma\right), 
j^{-1}\left(\dot{W}+\left(\dot{A}_0-\dot{Y}_0\right)\1\right)\right 
\rangle_{H\to F}=\left\langle D\left(W_{t_l+\sigma}-W_\sigma\right) 
,j^{-1}\dot{W}\right\rangle_{H\to F} \vphantom{\frac12} \\ 
&&\hspace{.5cm}=\frac12\left(W'_{t_l+\sigma^-(m)}-W'_{t_l+\sigma^- 
(m)-t/2^m}\right)-\frac12\left(W'_{\sigma^-(m)}-W'_{\sigma^-(m)- 
t/2^m}\right) \\ 
&&\hspace{1cm}+2^{m-1}\, \frac{\sigma-\sigma^-(m)}{t}\cdot\left( 
W'_{t_l+\sigma^-(m)+t/2^m}-2W'_{t_l+\sigma^-(m)}+W'_{t_l+\sigma^- 
(m)-t/2^m}\right) \\ 
&&\hspace{1cm}-2^{m-1}\, \frac{\sigma-\sigma^-(m)}{t}\cdot\left( 
W'_{\sigma^-(m)+t/2^m}-2W'_{\sigma^-(m)}+W'_{\sigma^-(m)-t/2^m} 
\right) \\ 
&&\hspace{.5cm}=\frac12\left(\left(p_{m,r}W'_{\cdot +\sigma}\right 
)_{t_l}-\left(p_{m,r}W'_{\cdot +\sigma}\right)_{t_l-t/2^m}\right)- 
\frac12\left(\left(p_{m,r}W'_{\cdot +\sigma}\right)_0-\left(p_{m,r} 
W'_{\cdot +\sigma}\right)_{-t/2^m}\right) \\ 
&&\hspace{.5cm}=:\frac12\Delta\left(p_{m,r}\dot{W}_{\cdot+\sigma} 
\right)_{t_l}-\frac12\Delta\left(p_{m,r}\dot{W}_{\cdot +\sigma} 
\right)_0\, ,\quad l\in\{1,\ldots ,k\},  
\end{eqnarray*}
here, as in Subsection 1.2 introduced, $W'$ denoting the right 
derivative. It follows from Proposition \ref{Proposition2.7}, Lemma 
\ref{Lemma3.1} (c), (d), and condition (3) of Section 1 that  
\begin{eqnarray*}
\dot{W}_\sigma+\dot{A}_\sigma&&\hspace{-.5cm}\equiv\left.\frac{d} 
{d\rho}\right|^\pm_{\rho=0}X_{\sigma+\rho}(W)=\left.\frac{d}{d\rho} 
\right|^\pm_{\rho=0}X_\sigma (W^\rho) \\ 
&&\hspace{-.5cm}=\left\langle D_F(W_\sigma +A_\sigma),j^{-1}\left( 
\dot{W}+\left(\dot{A}_0-\dot{Y}_0\right)\1\right)\right\rangle_{H 
\to F}\vphantom{\sum}\quad Q^{(m,r)}_{\sbnu}\mbox{\rm -a.e.} 
\end{eqnarray*}
Furthermore, by condition (3) of Section 1 it holds that $Y_{\sigma 
+\rho}(W)=A_0(W^{\sigma +\rho})=A_0((W^\rho)^\sigma)$ $=Y_\sigma 
(W^\rho)$, $\rho\in {\Bbb R}$, $ ,\sigma\in [0,t]$. Similar to the 
above we get therefore 
\begin{eqnarray*}
\dot{Y}_\sigma&&\hspace{-.5cm}\equiv\left.\frac{d}{d\rho}\right 
|^\pm_{\rho=0}Y_{\sigma+\rho}(W)=\left.\frac{d}{d\rho}\right|^\pm_{ 
\rho=0}Y_\sigma (W^\rho) \\ 
&&\hspace{-.5cm}=\left\langle D_F(Y_\sigma),j^{-1}\left(\dot{W}+\left 
(\dot{A}_0-\dot{Y}_0\right)\1\right)\right\rangle_{H\to F}\quad Q^{ 
(m,r)}_{\sbnu}\mbox{\rm -a.e.}\vphantom{\sum}
\end{eqnarray*}
After these preparations we obtain  
\begin{eqnarray*}
&&\hspace{-.5cm}\sum_{l=1}^kf_0(W_\sigma +A_\sigma -Y_\sigma)\cdot 
\left\langle\vphantom{\dot{f}}\left\langle\nabla_lf_1(W_{t_1+\sigma}- 
W_\sigma,\ldots , W_{t_k+\sigma}-W_\sigma),\right.\right. \\ 
&&\hspace{4cm}\left.\left. D(W_{t_l+\sigma}-W_{\sigma})\right 
\rangle_{F\to F},j^{-1}\left(\dot{W}+\left(\dot{A}_0-\dot{Y}_0\right) 
\1\right)\vphantom{\dot{f}}\right\rangle_H \vphantom{\sum_{l=1^1}} \\ 
&&\hspace{1cm}+f_1(W_{t_1+\sigma}-W_\sigma,\ldots , W_{t_k+\sigma}- 
W_\sigma)\times\vphantom{\sum_{l=1^1}} \\ 
&&\hspace{1cm}\times\left\langle\vphantom{\dot{f}}\left\langle\nabla 
f_0(W_\sigma+A_\sigma -Y_\sigma)\, ,\, D_F(W_\sigma +A_\sigma-Y_\sigma 
)\right\rangle_{F\to F},j^{-1}\left(\dot{W}+\left(\dot{A}_0-\dot{Y}_0 
\right)\1\right)\right\rangle_H \\ 
&&\hspace{.5cm}=\sum_{l=1}^k\left\langle f_0(W_\sigma +A_\sigma- 
Y_\sigma)\cdot\nabla_lf_1(W_{t_1+\sigma}-W_\sigma,\ldots ,W_{t_k+ 
\sigma}-W_\sigma)\, ,\vphantom{\frac12}\right. \\ 
&&\hspace{2cm}\left.\frac12\Delta\left(p_{m,r}\dot{W}_{\cdot+\sigma} 
\right)_{t_l}-\frac12\Delta\left(p_{m,r}\dot{W}_{\cdot +\sigma}\right 
)_0\right\rangle_F \\ 
&&\hspace{1cm}+\left\langle \nabla f_0(W_\sigma +A_\sigma -Y_\sigma) 
\cdot f_1(W_{t_1+\sigma}-W_\sigma,\ldots , W_{t_k+\sigma}-W_\sigma), 
\dot{W}_\sigma+\dot{A}_\sigma-\dot{Y}_{\sigma}\right\rangle_F \\ 
&&\hspace{.5cm}=\left\langle(D\vp)\left(\pi_{m,r}\left(W_{\cdot + 
\sigma}+(A_\sigma-Y_\sigma)\1\vphantom{l^1}\right)\right),j^{-1}\left 
(\pi_{m,r}\left(W_{\cdot +\sigma}+(A_\sigma-Y_\sigma)\1\vphantom{l^1} 
\right)\right)\, \dot{}\, \vphantom{\dot{f}}\right\rangle_H\vphantom 
{\sum_{l=1}^k} 
\end{eqnarray*}
$Q^{(m,r)}_{\sbnu}$-a.e. The claim follows now from (\ref{3.6}), 
(\ref{3.7}), and Definition \ref{Definition2.6}. 
\medskip 

\nid 
(d) Using the notation $(D_{F,1})_s\xi(W)=\sum_{i\in I(r)}\frac 
{\partial\xi(W)}{\partial\kappa_i}\cdot H_i(s)$, $s\in {\cal R}$, the 
following calculation is straight forward. 
\begin{eqnarray*} 
&&\hspace{-.5cm}L_1+L_2=\int\int_{\sigma=0}^v\left\langle D\vp\circ 
(W_{\cdot +\sigma}+(A_\sigma -Y_\sigma)\1)\, ,\, j^{-1}\left(\dot{W}+ 
\left(\dot{A}_0-\dot{Y}_0\right)\1\right)\right\rangle_H\, d\sigma\, 
Q^{(m,r)}_{\sbnu}(dW) \\ 
&&\hspace{1cm}+\int\int_{\sigma=0}^v\left(\vphantom{\int}f_1(W_{t_1 
+\sigma}-W_\sigma,\ldots ,W_{t_k+\sigma}-W_\sigma)\left\langle\vphantom 
{\dot{f}}\left\langle\nabla f_0(W_\sigma+ A_\sigma -Y_\sigma)\, ,\, 
\vphantom{l^1}\right.\right.\right. \\ 
&&\hspace{1.5cm}\left.\left.\left.\vphantom{l^1}(D_F-D)(W_\sigma + 
A_\sigma -Y_\sigma)\right\rangle_{F\to F}\, ,\, j^{-1}\left(\dot{W}+ 
\left(\dot{A}_0-\dot{Y}_0\right)\1\right)\right\rangle_H\vphantom{\int} 
\right)\, d\sigma\, Q^{(m,r)}_{\sbnu}(dW)\vphantom{\int_{\sigma=0 
}^\infty} \\ 
&&\hspace{.5cm}=\int\int_{\sigma=0}^v\sum_{l=1}^kf_0(W_\sigma +A_\sigma 
-Y_\sigma)\cdot\left\langle\vphantom{\dot{f}}\left\langle\nabla_lf_1(W_{ 
t_1+\sigma}-W_\sigma,\ldots , W_{t_k+\sigma}-W_\sigma),\right.\right. \\ 
&&\hspace{7cm}\left.\left.D_1(W_{t_l+\sigma}-W_{\sigma})\right\rangle_{ 
F\to F},d\dot{W}\vphantom{\dot{f}}\right\rangle_{L^2}\, d\sigma\, Q^{(m 
,r)}_{\sbnu}(dW)\vphantom{\int_{\sigma=0}^v} \\ 
&&\hspace{1cm}+\int\int_{\sigma=0}^vf_1(W_{t_1+\sigma}-W_\sigma,\ldots 
,W_{t_k+\sigma}-W_\sigma)\times\vphantom{\int_{\sigma=0}^v} \\ 
&&\hspace{1.5cm}\times\left\langle\vphantom{\dot{f}}\left\langle\nabla 
f_0(W_\sigma+A_\sigma -Y_\sigma)\, ,\, D_{F,1}(W_\sigma +A_\sigma - 
Y_\sigma)\right\rangle_{F\to F},d\dot{W}\right\rangle_{L^2}\, d\sigma\, 
Q^{(m,r)}_{\sbnu}(dW)\vphantom{\int_{\sigma=0}^v} \\ 
&&\hspace{1cm}+\int\int_{\sigma=0}^v\left\langle D\vp\circ\left(W_{\cdot 
+\sigma}+(A_\sigma -Y_\sigma)\1\vphantom{l^1}\right)\, ,\, \left(0,\dot{ 
W}_0+\dot{A}_0-\dot{Y}_0\right)\right\rangle_H\, d\sigma\, Q^{(m,r)}_{ 
\sbnu}(dW)\, . 
\end{eqnarray*}
All integrals are well-defined by part (a), the hypotheses to the present 
part (d), and Proposition \ref{Proposition2.7}. 

We will now regard the first two integrals on the right-hand side as 
integrals with respect to the measure $Q_{\sbnu}$ rather than $Q^{(m,r) 
}_{\sbnu}$. As compensation the integrands will be considered under the 
projection $\pi_{m,r}$. We obtain 
\begin{eqnarray*} 
&&\hspace{-.5cm}L_1+L_2=\int\int_{\sigma=0}^v\sum_{l=1}^kf_0([W_\sigma + 
A_\sigma -Y_\sigma]\circ\pi_{m,r})\cdot\left\langle\vphantom{\dot{f}} 
\left\langle\nabla_lf_1([W_{t_1+\sigma}-W_\sigma]\circ\pi_{m,r},\ldots, 
\right.\right. \\ 
&&\hspace{1.5cm}\left.\left.[W_{t_k+\sigma}-W_\sigma]\circ\pi_{m,r})\, , 
\, D_1([W_{t_l+\sigma}-W_{\sigma}]\circ\pi_{m,r})\right\rangle_{F\to F} 
\, ,\, \vphantom{\dot{f}}\right.\vphantom{\int_{\sigma=0}^v} \\ 
&&\hspace{1.5cm}\left. d\left([W]\circ\pi_{m,r}\right)^\cdot\vphantom{ 
\dot{f}}\right\rangle_{L^2}\, d\sigma\, Q_{\sbnu}(dW)\vphantom{\int_{ 
\sigma=0}^v} \\ 
&&\hspace{1cm}+\int\int_{\sigma=0}^vf_1\left([W_{t_1+\sigma}-W_\sigma] 
\circ\pi_{m,r},\ldots ,[W_{t_k+\sigma}-W_\sigma]\circ\pi_{m,r}\right) 
\times\vphantom{\int_{\sigma=0}^v} \\ 
&&\hspace{1.5cm}\times\left\langle\vphantom{\dot{f}}\left\langle\nabla 
f_0([W_\sigma +A_\sigma -Y_\sigma ]\circ\pi_{m,r})\, ,\, D_{F,1}\left([ 
W_\sigma +A_\sigma -Y_\sigma ]\circ\pi_{m,r}\right)\right\rangle_{F\to 
F}\, ,\, \right.\vphantom{\int_{\sigma=0}^v} \\ 
&&\hspace{1.5cm}\left. d\left([W]\circ\pi_{m,r}\right)^\cdot\vphantom 
{\dot{f}}\right\rangle_{L^2}\, d\sigma\, Q_{\sbnu}(dW)+J\vphantom 
{\int_{\sigma=0}^v} \\ 
&&\hspace{0.5cm}=\int\int_{\sigma=0}^v\sum_{i\in I(r)}\left\langle D_{F, 
1}\vp\circ\left([W_{\cdot +\sigma}+(A_\sigma -Y_\sigma)\1]\circ\pi_{m,r} 
\vphantom{l^1}\right)\, ,\, H_i\right\rangle_{L^2}\times \\ 
&&\hspace{1.5cm}\times\left\langle H_i\, ,\, d\left([W]\circ\pi_{m,r} 
\right)^\cdot\right\rangle_{L^2}\, d\sigma\, Q_{\sbnu}(dW)\vphantom 
{\int_{\sigma=0}^v}+J  \\ 
&&\hspace{0.5cm}=:r_1+J\vphantom{\int_{\sigma=0}^v} 
%
\end{eqnarray*}
We note that $r_1$ is independent of the measure $\bnu$. Thus, $r_1$ 
is an integral with respect to the classical Wiener measure $Q(dV):= 
\int_{W_0\in D^n}\, Q_{\sbnu}(d(V+W_0\1))$ where $W=V+W_0\1\in\Omega$ 
and $V\in \Omega_0:=\{U\in\Omega :U_0=0\}$. Therefore, $D_{F,1}$ is 
the Malliavin derivative with respect to $Q$ on $(\Omega_0,{\cal B}( 
\Omega_0))$ and the corresponding Sobolev norm in $D_{q,1}$. We obtain 
\begin{eqnarray*}
&&\hspace{-.5cm}L_1+L_2=\sum_{i\in I(r)}\int_{\sigma=0}^v\int\vp\left( 
[W_{\cdot +\sigma}+(A_\sigma -Y_\sigma)\1]\circ\pi_{m,r}\vphantom{l^1} 
\right)\cdot\left\langle H_i,d\left([W]\circ\pi_{m,r}\right)^\cdot\right 
\rangle_{L^2}\times \\ 
&&\hspace{1.5cm}\vphantom{\int_{\sigma=0}^v}\times\langle H_i\, ,\, dW 
\rangle_{L^2}\, Q_{\sbnu}(dW)\, d\sigma \\ 
&&\hspace{1cm}-\sum_{i\in I(r)}\int_{\sigma=0}^v\int\vp\left([W_{\cdot 
+\sigma}+(A_\sigma-Y_\sigma)\1]\circ\pi_{m,r}\vphantom{l^1}\right)\times 
 \\ 
&&\hspace{1.5cm}\vphantom{\int_{\sigma=0}^v}\times\left\langle (H_i,0)\, 
,\, D_F\left\langle H_i,d\left([W]\circ\pi_{m,r}\right)^\cdot\right 
\rangle_{L^2}\vphantom{\dot{f}}\right\rangle_H\, Q_{\sbnu}(dW)\, d\sigma 
+J\, . 
\end{eqnarray*}
\medskip 

\nid 
(e) For $\langle H_i,d\dot{W}\rangle_{L^2}\in\bigcap_{r\in [Q/(Q-1) 
,\infty)}D_{r,1}$, recall Lemma \ref{Lemma3.1} (a), $i\in I(m,r)$. 
Together with part (a) this implies $\vp\circ (W_{\cdot +\sigma}+ 
(A_\sigma-Y_\sigma)\1)\langle H_i,d\dot{W}\rangle_{L^2}\in D_{q',1}$, 
$\sigma\in [0,t]$, for $Q/(Q-1)<q'<q$ and $i\in I(m,r)$. We have 
\begin{eqnarray*} 
R_1&&\hspace{-.5cm}=\sum_{i\in I(m,r)}\int_{\sigma=0}^v\int\vp\left( 
W_{\cdot+\sigma}+(A_\sigma-Y_\sigma)\1\vphantom{l^1}\right)\cdot\left 
\langle H_i,d\dot{W}\right\rangle_{L^2}\langle H_i\, ,\, dW\rangle_{ 
L^2}\, Q^{(m,r)}_{\sbnu}(dW)\, d\sigma \\ 
&&\hspace{.0cm}+\sum_{i\in I(r)\setminus I(m,r)}\int_{\sigma=0}^v\int 
\vp\left([W_{\cdot+\sigma}+(A_\sigma -Y_\sigma)\1]\circ\pi_{m,r} 
\vphantom{l^1}\right)\cdot\left\langle H_i,d\left([W]\circ\pi_{m,r} 
\right)^\cdot\right\rangle_{L^2}\times \\ 
&&\hspace{.5cm}\vphantom{\int_{\sigma=0}^v}\times\langle H_i\, ,\, dW 
\rangle_{L^2}\, Q_{\sbnu}(dW)\, d\sigma  \\ 
&&\hspace{-.5cm}=\sum_{i\in I(m,r)}\int_{\sigma=0}^v\int\vp\left( 
W_{\cdot+\sigma}+(A_\sigma-Y_\sigma)\1\vphantom{l^1}\right)\cdot\left 
\langle H_i,d\dot{W}\right\rangle_{L^2}\langle H_i\, ,\, dW\rangle_{ 
L^2}\, Q^{(m,r)}_{\sbnu}(dW)\, d\sigma
\end{eqnarray*}
since $\vp\left([W_{\cdot+\sigma}+(A_\sigma -Y_\sigma)\1]\circ\pi_{m, 
r}\vphantom{l^1}\right)\cdot\left\langle H_i,d\left([W]\circ\pi_{m,r} 
\right)^\cdot\right\rangle_{L^2}$ is for $i\in I(r)\setminus I(m,r)$ 
independent of $\langle H_i\, ,\, dW\rangle_{L^2}$ and $\int\langle 
H_i\, ,\, dW\rangle_{L^2}\, Q_{\sbnu}(dW)=0$. Furthermre, for $i\in 
I(r)\setminus I(m,r)$ we have $\left\langle H_i,d\left([W]\circ 
\pi_{m,r}\right)^\cdot\right\rangle_{L^2}=0$. Thus 
\begin{eqnarray*} 
R_2&&\hspace{-.5cm}=-\sum_{i\in I(m,r)}\int_{\sigma=0}^v\int\vp\left 
([W_{\cdot +\sigma}+(A_\sigma-Y_\sigma)\1]\circ\pi_{m,r}\vphantom 
{l^1}\right)\times \\ 
&&\hspace{4cm}\vphantom{\int_{\sigma=0}^v}\times\left\langle (H_i,0) 
\, ,\, D\left\langle H_i,d\left([W]\circ\pi_{m,r}\right)^\cdot\right 
\rangle_{L^2}\vphantom{\dot{f}}\right\rangle_H\, Q_{\sbnu}(dW)\, d 
\sigma \\ 
&&\hspace{-.5cm}=-\sum_{i\in I(m,r)}\int_{\sigma=0}^v\int\vp\left( 
W_{\cdot +\sigma}+(A_\sigma-Y_\sigma)\1\vphantom{l^1}\right)\times 
 \\ 
&&\hspace{4cm}\vphantom{\int_{\sigma=0}^v}\times\left\langle (H_i,0) 
\, ,\, D\left\langle H_i,d\dot{W}\right\rangle_{L^2}\vphantom{\dot{f}} 
\right\rangle_H\, Q^{(m,r)}_{\sbnu}(dW)\, d\sigma\, . 
\end{eqnarray*}
\medskip 

\nid 
(f) For $\vp$ cylindrical in the sense of this lemma, we have 
\begin{eqnarray*} 
&&\hspace{-.5cm}\int\left\langle\left(0,\dot{W}_0+\dot{A}_0- 
\dot{Y}_0\right),D\vp(W)\vphantom{\dot{f}}\right\rangle_H\, 
Q^{(m,r)}_{\sbnu}(dW) \\ 
&&\hspace{.5cm}=\int\int\left\langle\dot{W}_0+\dot{A}_0-\dot 
{Y}_0,\nabla_{W_0}\vp(W)\vphantom{\dot{f}}\right\rangle_F\, 
Q^{(m,r)}_{W_0}(dW)\, m(W_0)\, dW_0 \\ 
&&\hspace{.5cm}=\int\left\langle\int\left(\dot{W}_0+\dot{A}_0 
-\dot{Y}_0\right)\cdot f_1\left(W_{t_1}-W_0,\ldots ,W_{t_k}- 
W_0\vphantom{l^1}\right)\, Q^{(m,r)}_{W_0}(dW)\, ,\right.\\ 
&&\hspace{2cm}\left.\nabla_{W_0}f_0(W_0)\vphantom{\int}\right 
\rangle_F\, m(W_0)\, dW_0 \, . 
\end{eqnarray*} 
We observe that $f_1\left(W_{t_1}-W_0,\ldots ,W_{t_k}-W_0 
\vphantom{l^1}\right)\, Q^{(m,r)}_{W_0}(dW)$ is independent 
of $W_0$ and obtain 
\begin{eqnarray*} 
&&\hspace{-.5cm}\int\left\langle\left(0,\dot{W}_0+\dot{A}_0- 
\dot{Y}_0\right),D\vp(W)\vphantom{\dot{f}}\right\rangle_H\, 
Q^{(m,r)}_{\sbnu}(dW) \\ 
&&\hspace{.5cm}=-\int\left(\left\langle \dot{W}_0+\dot{A}_0- 
\dot{Y}_0,\frac{\nabla m(W_0)}{m(W_0)}\right\rangle_F+\left 
\langle {\sf e},\nabla_{d,W_0}\left(\dot{A}_0-\dot{Y}_0\right) 
\right\rangle_F\right)\cdot\vp(W)\, Q^{(m,r)}_{\sbnu}(dW) 
\end{eqnarray*} 
where we have taken into consideration $\nabla_{W_0}\dot{W}_0= 
0$. For the well-definiteness we refer to the definition of the 
density $m$, condition (1') of Subsection 1.2, as well as to 
$(\dot{A}_\rho -\dot{Y}_\rho)\in L^p(\Omega ,Q^{(m,r)}_{\sbnu}; 
F)$ where $p\ge q$ by hypothesis. Moreover, we recall the 
existence of $\nabla_{d,W_0}(\dot{A}_0-\dot{Y}_0)$ in $L^q 
(\Omega,Q^{(m,r)}_{\sbnu};F)$ and $q\ge w$ by hypothesis. The 
last chain of equation yields $(0,\dot{W}_0+\dot{A}_0-\dot{Y 
}_0)\in {\rm Dom}_w(\delta)$ and the stochastic integral 
$\delta (0,\dot{W}_0+\dot{A}_0-\dot{Y}_0)$ in the sense of 
Subsection 1.2. 
\qed 

\subsection{Related stochastic integral} 

The following lemma explains the stochastic integral in (\ref{2.33}) 
and provides detailed insight in the crucial step of the proof of 
Theorem \ref{Theorem1.11}, relation (\ref{3.26}). 
\begin{lemma}\label{Lemma3.3} 
Let $Q$ be the number appearing in the definition of $m$ in Subsection 
1.2, let $q$ be the number of condition (1') in Section 1, and let 
$1/p+1/q=1$ with $p\ge q$.  

Assume conditions (1'), (2), and in particular, 
\begin{eqnarray}\label{3.19}
W^\rho=f(\rho ,W)+W=W_{\cdot +\rho}+(A_\rho-Y_\rho)\1\, ,\quad\rho\in 
{\Bbb R}. 
\end{eqnarray} 
Assume furthermore  
\begin{itemize} 
\item[(i)] $\left(\dot{A}_0-\dot{Y}_0\right)\in L^p(\Omega,Q^{(m,r) 
}_{\sbnu};F)$ and $\nabla_{d,W_0}\left(\dot{A}_0-\dot{Y}_0\right)\in 
L^q(\Omega,Q^{(m,r)}_{\sbnu};F)$, 
\item[(ii)] $A=A^1$, $Y=Y^1$. 
\end{itemize} 
We have $\left(0,\dot{W}_0+\dot{A}_0-\dot{Y}_0\right)\in {\rm Dom}_w 
(\delta)$ with $w=\min\left(q,(1/p+1/Q)^{-1}\right)$ in the sense of 
Subsection 1.2 and for the term (\ref{2.33})
\begin{eqnarray*}
\hat{\delta}\left(j^{-1}\dot{f}(0,\cdot )\right)&&\hspace{-.5cm}= 
\delta\left(0,\dot{W}_0+\dot{A}_0-\dot{Y}_0\right)\nonumber \\ 
&&\hspace{-.5cm}=-\left\langle\dot{W}_0+\dot{A}_0-\dot{Y}_0,\frac{ 
\nabla m(W_0)}{m(W_0)}\right\rangle_F-\left\langle {\sf e},\nabla_{ 
d,W_0}\left(\dot{A}_0-\dot{Y}_0\right)\right\rangle_F 
\end{eqnarray*} 
for $Q^{(m,r)}_{\sbnu}$-a.e. $W\in\Omega$. 
\end{lemma}
{\bf Remarks. (1)} For the term (\ref{2.33}) we had assumed $\dot{f} 
(0,\cdot)(s)\in D_{p,1}$, $s\in {\cal R}$, following the context of 
Section 2. Now that we have specified $f$ and $\dot{f}$ (cf. 
(\ref{3.2}) and (\ref{3.3})) it is natural to suppose $\left(\dot{A}_0 
-\dot{Y}_0\right)\in L^p(\Omega,Q^{(m,r)}_{\sbnu};F)$. 
\medskip 

\nid 
{\bf (2)} By condition (1') in Subsection 1.2 and  $p\ge q$ we have 
$w>1$. By (i) of this lemma and the definition of $m$ in Subsection 
1.2 the first item of the right-hand side of the above representation 
of $\hat{\delta}\left(j^{-1}\dot{f}(0,\cdot )\right)$ belongs to $L^w 
(\Omega,Q^{(m,r)}_{\sbnu})$. By hypothesis we have $q\ge w$. Hence 
the second item of $\hat{\delta}\left(j^{-1}\dot{f}(0,\cdot )\right)$ 
belongs also to $L^w(\Omega,Q^{(m,r)}_{\sbnu})$. 
\medskip 

\nid
Proof. We show that the right-hand side of (\ref{2.33}) forms a 
stochastic integral in the sense of Subsection 1.2. First, we observe 
that the right-hand side of (\ref{2.33}), 
\begin{eqnarray*}
&&\hspace{-.5cm}I(W):=\sum_{i\in I(m,r)}\left\langle H_i,d\dot{f}(0, 
W)\right\rangle_{L^2}\cdot\langle H_i,dW\rangle_{L^2}-\sum_{i\in I(m 
,r)}\left\langle (H_i,0),D\left\langle H_i,d\dot{f}(0,W)\right 
\rangle_{L^2}\right\rangle_H\hspace{-.2cm}\nonumber \\ 
&&\hspace{1.5cm}-\left\langle\dot{f}(0,W)(0),\frac{\nabla m(W_0)}{m 
(W_0)}\right\rangle_F-\left\langle {\sf e},\nabla_{d,W_0}\dot{f}(0,W) 
(0)\right\rangle_F\, ,  
\end{eqnarray*} 
simplifies significantly because of 
\begin{eqnarray}\label{3.20}
\sum_{i\in I(m,r)}\left\langle H_i,d\dot{f}(0,W)\right\rangle_{L^2} 
\cdot\langle H_i,dW\rangle_{L^2}&&\hspace{-.5cm}=\sum_{i\in I(m,r)} 
\left\langle H_i,dW\right\rangle_{L^2}\cdot\left\langle H_i ,d\dot{W} 
\right\rangle_{L^2}=0\, ,\qquad
\end{eqnarray} 
and 
\begin{eqnarray}\label{3.21}
-\sum_{i\in I(m,r)}\left\langle (H_i,0),D\left\langle H_i ,d\dot{f} 
(0,W)\right\rangle_{L^2}\right\rangle_H&&\hspace{-.5cm}=-\sum_{i\in 
I(m,r)}\left\langle (H_i,0),D\left\langle H_i ,d\dot{W}\right 
\rangle_{L^2}\right\rangle_H \nonumber \\ 
&&\hspace{-.5cm}=-\sum_{i\in I(m,r)}\langle H_i,dH_i\rangle_{L^2} 
\nonumber \\ 
&&\hspace{-.5cm}=0\quad\mbox{\rm for }Q^{(m,r)}_{\sbnu} 
\mbox{\rm -a.e. }W\in\Omega.  
\end{eqnarray} 
In fact, relation (\ref{3.20}) follows from (\ref{3.19}) and 
(\ref{2.3}) and relation (\ref{3.21}) is a consequence of (\ref{3.19}), 
Lemma \ref{Lemma3.1} (a), and the definition of the integral $\langle 
H_i,dH_i\rangle_{L^2}$ in Subsection 2.1. We obtain 
\begin{eqnarray}\label{3.22}
I(W)&&\hspace{-.5cm}=-\left\langle\dot{f}(0,W)(0),\frac{\nabla m(W_0 
)}{m(W_0)}\right\rangle_F-\left\langle {\sf e},\nabla_{d,W_0}\dot{f} 
(0,W)(0)\right\rangle_F \nonumber \\ 
&&\hspace{-.5cm}=-\left\langle \dot{W}_0+\dot{A}_0-\dot{Y}_0,\frac{ 
\nabla m(W_0)}{m(W_0)}\right\rangle_F-\left\langle {\sf e},\nabla_{ 
d,W_0}\left(\dot{A}_0-\dot{Y}_0\right)\right\rangle_F=0 
\end{eqnarray} 
for $Q^{(m,r)}_{\sbnu}$-a.e. $W\in\Omega$, where we have taken into 
consideration $\nabla_{d,W_0}\dot{W}_0=0$. The lemma follows now 
from (\ref{3.22}) and Lemma \ref{Lemma3.2} (f), especially relation 
(\ref{3.5}). 
\qed 

\subsection{Technical details of the approximation} 

In this section we carry out all calculations for the approximation of 
trajectories with jumps, $A$, $X$, and $Y$, by continuous trajectories 
$A_n$, $X_n$, and $Y_n$. Furthermore, we present the technicalities 
for the extension of the analysis on $[-rt,(r-1)t]$ of the trajectories
to the analysis of the trajectories on ${\Bbb R}$. For the readers 
convenience we collect all claims in the subsequent Lemma 
\ref{Lemma3.4}. The proofs of all individual statements follow then. 

\begin{lemma}\label{Lemma3.4} 
Assume conditions (1)-(5) of Section 1. For all $W\in\Omega$, there 
exist sequences of processes $A_n(W)$ and $Y_n(W)$ such that with $ 
X_n(W)$ $:=W+A_n(W)$, $n\in {\Bbb N}$, the following holds. 
\begin{itemize} 
\item[(i)] $A_n=A_n^1$, $Y_n=Y_n^1$, $n\in {\Bbb N}$, and moreover 
$A_n=A_n^1$ and $Y_n=Y_n^1$ are continuously differentiable. 
\item[(ii)] Condition (iii) of Proposition \ref{Proposition2.7} holds 
for $X_n$ as well as $Y_n$. Moreover, for $s\in {\Bbb R}$, 
\begin{eqnarray*}
\nabla_HX_{n,s}\in \bigcap_{1\le p<\infty}L^p(\Omega,Q^{(m,r)}_{\sbnu} 
;H)\quad\mbox{\rm and}\quad\nabla_H Y_{n,s}\in\bigcap_{1\le p<\infty} 
L^p(\Omega ,Q^{(m,r)}_{\sbnu};H)\, . 
\end{eqnarray*} 
In addition, we have $A_{n,s},Y_{n,s}\in D_{q,1}(Q^{(m,r)}_{\sbnu})$ 
for all $q$ with $1/q+1/Q<1$ where $Q$ is the number introduced in the 
definition of the density $m$ in Subsection 1.2. 
\item[(iii)] We have 
\begin{eqnarray*}
\dot{A}_{n,s}-\dot{Y}_{n,s}\in \bigcap_{1\le p<\infty}L^p\left(\Omega, 
Q^{(m,r)}_{\sbnu};F\right)\cap L^p\left(\Omega,Q^{(m)}_{\sbnu};F\right) 
\, , \quad s\in {\Bbb R}. 
\end{eqnarray*} 
Furthermore, $\dot{A}_n(W)-\dot{Y}_n(W)$ is bounded on finite 
subintervals of ${\Bbb R}$ for all $W\in\bigcup_m\{\pi_m V:V\in\Omega 
\}$. 
\item[(iv)] For all $p\in [1,\infty)$, the derivative $\nabla_{W_0} 
\left(\dot{A}_{n,0}-\dot{Y}_{n,0}\right)$ exists in $L^p(\Omega,Q^{(m,
r)}_{\sbnu};F)$. In addition, for all $p\in [1,\infty)$, $\nabla_{W_0} 
\left(\dot{A}_{n,s}-\dot{Y}_{n,s}\right)$, seen as a limit in $L^p( 
\Omega,Q^{(m)}_{\sbnu};F)$, exists uniformly for all $s$ belonging to 
any finite subinterval of ${\Bbb R}$. In particular, $\nabla_{W_0} 
\left(\dot{A}_n(W)-\dot{Y}_n(W)\right)$ is bounded on finite 
subintervals ${\Bbb R}$ for all $W\in\bigcup_m\{\pi_m V:V\in\Omega\}$. 
Furthermore, 
\begin{eqnarray*} 
\int_S\left|\left\langle {\sf e},\nabla_{W_0}\left(\dot{A}_{n,s}- 
\dot{Y}_{n,s}\right)\right\rangle_F\right|\, ds\in \bigcap_{1\le p< 
\infty}L^p(\Omega,Q^{(m)}_{\sbnu})  
\end{eqnarray*} 
for every finite subinterval $S\subset {\Bbb R}$. 
\item[(v)] On $\bigcup_m\{\pi_m W:W\in\Omega\}$, $\nabla_{W_0}\left( 
A_n-Y_n\right)$ is continuously differentiable on ${\Bbb R}$ and we 
have $\nabla_{W_0}\left(\dot{A}_{n,s}-\dot{Y}_{n,s}\right)=\left( 
\nabla_{W_0}\left(A_{n,s}-Y_{n,s}\right)\right)^\cdot$. 
\item[(vi)] With $W^{n,u}:=W_{\cdot +u}+\left(A_{n,u}(W)-Y_{n,u}(W) 
\vphantom{l^1}\right)\1$, $u\in {\Bbb R}$, the process $X_n=W+A_n(W)$ 
is temporally homogeneous on $W\in\Omega$; i.e., we 
have 
\begin{eqnarray*} 
W^{n,0}=W\, , \quad A_{n,0}\left((W^{n,u})^v\vphantom{l^1}\right)= 
A_{n,0}\left(W^{n,u+v}\vphantom{l^1}\right)\, , 
\end{eqnarray*} 
and 
\begin{eqnarray*} 
X_{n,\cdot+v}(W)=X_n\left(W^{n,v}\vphantom{l^1}\right)\, ,\quad v\in 
{\Bbb R}. 
\end{eqnarray*} 
Furthermore, $Y_{n,v}(W)=A_{n,0}(W^{n,v})$, 
\begin{eqnarray*} 
Y_{n,\cdot+v}(W)=Y_n\left(W^{n,v}\vphantom{l^1}\right)\, ,\quad\mbox 
{\it and}\quad\dot{A}_{n,\cdot +v}(W)=\dot{A}_n(W^{n,v})\, ,\quad v 
\in {\Bbb R}. 
\end{eqnarray*} 
In addition, $A_{n,\cdot}$ and $Y_{n,\cdot}$ are constant on intervals  
where $W^{n,\cdot}_0\not\in\left\{x\in F:|x-z|<\frac{1}{n^3},\right.$ 
$\left.z\in D^n\vphantom{W^{n,\cdot}}\right\}$. 
\end{itemize} 
We have the following convergences.  
\medskip 
 
\nid
(a) $A_{n,s}(W)-Y_{n,s}(W)\stack{n\to\infty}{\lra}A_s(W)-Y_s(W)$ for 
all $W\in\bigcup_m\{\pi_m V:V\in\Omega\}$ with $W_0\not\in G(W-W_0\1)$, 
cf. Subsection 1.2, whenever neither $0$ nor $s$ is a jump time for $X$. 
\medskip 

\nid 
(b) Let $W\in\{\pi_m V:V\in\Omega\}$ with $W_0\not\in G(W-W_0\1)$. Then 
\begin{eqnarray*}
A_{n,s}(\pi_{m,r}W)-Y_{n,s}(\pi_{m,r}W)\stack{r\to\infty}{\lra}A_{n,s} 
(W)-Y_{n,s}(W)\, ,  
\end{eqnarray*} 
\begin{eqnarray*}
\dot{A}_{n,s}(\pi_{m,r}W)-\dot{Y}_{n,s}(\pi_{m,r}W)\stack{r\to\infty} 
{\lra}\dot{A}_{n,s}(W)-\dot{Y}_{n,s}(W)\, ,  
\end{eqnarray*} 
as well as 
\begin{eqnarray*}
\nabla_{W_0}\dot{A}_{n,s}(\pi_{m,r}W)-\nabla_{W_0}\dot{Y}_{n,s}(\pi_{m 
,r}W)\stack{r\to\infty}{\lra}\nabla_{W_0}\dot{A}_{n,s}(W)-\nabla_{W_0} 
\dot{Y}_{n,s}(W)\, , \quad s\in {\Bbb R}. 
\end{eqnarray*} 
Furthermore, 
\begin{eqnarray*}
\sup_{r\in {\Bbb N}}\left|\dot{A}_{n,s}(\pi_{m,r}W)-\dot{Y}_{n,s}(\pi_{ 
m,r}W)\right|\in\bigcap_{1\le p<\infty} L^p(\Omega,Q^{(m)}_{\sbnu};F) 
\end{eqnarray*} 
as well as 
\begin{eqnarray*}
\sup_{r\in {\Bbb N}}\left|\nabla_{W_0}\dot{A}_{n,s}(\pi_{m,r}W)-\nabla_{ 
W_0}\dot{Y}_{n,s}(\pi_{m,r}W)\right|\in L^1(\Omega,Q^{(m)}_{\sbnu};F)\, 
, \quad s\in {\Bbb R}, 
\end{eqnarray*} 
where the absolute values $|\, \cdot\, |$ and the $\sup_{r\in {\Bbb N}}$ 
are taken individually for all coordinates of $F$. 
\medskip 

\nid 
(c) For all $-\infty<l<r<\infty$ it holds that 
\begin{eqnarray*}
&&\hspace{-.5cm}\nabla_{W_0}\left((A_{n,l}-Y_{n,l})-(A_{n,r}-Y_{n,r}) 
\vphantom{l^1}\right)\vphantom{\int} \\ 
&&\stack{n\to\infty}{\lra}\left(\nabla_{W_0}A_l-\nabla_{W_0}Y_l 
\vphantom{l^1}\right)-\left(\nabla_{W_0}A_r-\nabla_{W_0}Y_r\vphantom 
{l^1}\right) 
\end{eqnarray*} 
uniformly bounded $Q^{(m)}_{\sbnu}$-a.e. on $\{\pi_mW:W\in\Omega\}$. 
\end{lemma} 
{\bf Remarks. (3)} For part (c), we recall that according to 
conditions (2) and (3) of Subsection 1.2, $A$ and $Y$ jump only on 
$Q^{(m)}_{\sbnu}$-zero sets. Following condition (1) of Subsection 
1.2 it turns out that $Q^{(m)}_{\sbnu}$-a.e. on $\{\pi_mW:W\in\Omega 
\}$ the set of all of jump times of $\nabla_{W_0}A$ and $\nabla_{W_0 
}Y$ is a subset of $\{\tau_k:k\in {\Bbb Z}\setminus\{0\}\}$. 
\medskip 

\nid
{\bf (4)} In order to ease the notation in the proof below we will use 
the symbol $|\, \cdot\, |$ to abbreviate $\langle\, \cdot\, ,\, \cdot\, 
\rangle_F^{1/2}$ in both cases, $F={\Bbb R}^{n\cdot d}$ as well as 
$F={\Bbb R}^{n\cdot d}\otimes {\Bbb R}^{n\cdot d}$. Furthermore, the 
symbol $|\, \cdot\, |$ will also be used for the coordinate wise 
absolute value. It will either be clear from the context or explicitly 
be mentioned which meaning the symbol has. 
\medskip 

\nid
{\bf (5)} For the sake of comprehensibility we first carry out the 
proofs of the parts (a) and (c) for the particular case that, for given 
$W\in\Omega$ with $W_0=0$ and $x\in F$, the jump times of $X=u(W+x\1)$ 
are independent of $x\in F$. In a separate paragraph at the end of 
this subsection, in we prove parts (a) and (c) of Lemma \ref{Lemma3.4} 
if parallel trajectories in $\Omega$ generate no longer identical jump 
times for $X$. There we consider jump times as introduced in Definition 
\ref{Definition1.7} of Subsection 1.2. 
\medskip 

{\bf Definition of $A_n$ and $Y_n$ and verification of (i) and (vi). } 
{\it Step 1 } We define $A_n$ and $Y_n$ and construct the crucial 
system of ODEs. 

Recall from Subsection 1.2 the definitions of $g_n$, $\gamma_n$, and 
\begin{eqnarray*} 
A_s(\, \cdot\, ,\gamma_n)(W)&&\hspace{-.5cm}:=\int_F\left\langle 
A_s(W+x\1)\, ,\, \gamma_n(x)\vphantom{l^1}\right\rangle_{F\to F}\, 
dx\, . 
\end{eqnarray*} 
Introduce  
\begin{eqnarray*} 
A_s(g,\gamma_n)(W):=\int_{\Bbb R}\left\langle A_{s-v}(\, \cdot 
\, ,\gamma_n)(W)\, ,\, g(v)\vphantom{l^1}\right\rangle_{F\to F} 
\, dv\, ,  
\end{eqnarray*} 
$s\in {\Bbb R}$, $g\in\left\{f{\sf e}:f\in C({\Bbb R})\right\}$, $n 
\in {\Bbb N}$, $W\in\Omega$. Likewise define $Y_s(g,\gamma_n)$ and 
$Y_s(\, \cdot\, ,\gamma_n)$, $A'_s(g,\gamma_n)$ and $A'_s(\, \cdot\, 
,\gamma_n),\ldots\, $. Furthermore, set 
\begin{eqnarray}\label{3.8} 
B_{n,s}:=A_s\left(g_n,\gamma_n\right) 
\end{eqnarray} 
and 
\begin{eqnarray}\label{3.9*} 
C_{n,s}:=A_s\left(g_n,\gamma_n\right)-Y_s\left(g_n,\gamma_n\right) 
\, ,\quad s\in {\Bbb R},\ n\in {\Bbb N}. 
\end{eqnarray} 

Let $W\in\Omega$. Let us consider the following system of first order 
ODEs 
\begin{eqnarray}\label{3.10*} 
\begin{array}{rl}
\dot{\vp}(s)\ =&\dot{C}_{n,0}\left(\vp(s)\, \1+W_{\cdot +s}\vphantom{ 
\displaystyle l^1}\right)=:F(s,\vp(s))\, ,\quad s\in {\Bbb R}, 
\vphantom{\left(\dot{f}\right)} \\ 
\vp(0)\ =&0\, . \vphantom{\left(\dot{f}\right)}
\end{array} 
\end{eqnarray} 
In particular, we have 
\begin{eqnarray}\label{3.11*}
F(s,x)&&\hspace{-.5cm}=\int_{y\in F}\int_{v\in\left(-\frac1n, 
\frac1n\right)}\left\langle\vphantom{\dot{f}}\left\langle A_{-v} 
\left(W_{\cdot +s}+y\, \1\vphantom{l^1}\right)-Y_{-v}\left( W_{ 
\cdot +s}+y\, \1\vphantom{l^1}\right)\, ,\, \right.\right. 
\nonumber \\ 
&&\hspace{5.0cm}\left.\left.g'_n(v)\vphantom{l^1}\right\rangle_{F 
\to F}\, ,\gamma_n(y-x)\vphantom{\dot{f}}\right\rangle_{F\to F}\, 
dv\, dy\vphantom{\int_{v\in\left(-\frac1n,\frac1n\right)}} 
\nonumber \\ 
&&\hspace{-.5cm}=\int_{y\in F}\int_{v\in\left(-\frac1n,\frac1n 
\right)}\left\langle\vphantom{\dot{f}}\left\langle A'_{-v}\left( 
W_{\cdot +s}+y\, \1\vphantom{l^1}\right)-Y'_{-v}\left(W_{\cdot +s} 
+y\, \1\vphantom{l^1}\right)\, ,\, \right.\right.\nonumber \\ 
&&\hspace{5.0cm}\left.\left.g_n(v)\vphantom{l^1}\right\rangle_{F 
\to F}\, ,\gamma_n(y-x)\vphantom{\dot{f}}\right\rangle_{F\to F}\, 
dv\, dy\vphantom{\int_{v\in\left(-\frac1n,\frac1n\right)}} 
\nonumber \\ 
&&\hspace{.0cm}+\int_F\sum_{k\in {\Bbb Z}\setminus\{0\}}\left 
\langle\vphantom{\dot{f}}\left\langle\Delta A_{\tau_k}\left(W_{ 
\cdot +s}+y\, \1\vphantom{l^1}\right)-\Delta Y_{\tau_k}\left(W_{ 
\cdot +s}+y\, \1\vphantom{l^1}\right)\, ,\, \right.\right. 
\nonumber \\ 
&&\hspace{3.0cm}\left.\left.g_n\left(\tau_k\circ u\left(W_{\cdot 
+s}+y\, \1\vphantom{l^1}\right)\right)\right\rangle_{F\to F}\, , 
\, \vphantom{\dot{f}}\gamma_n(y-x)\right\rangle_{F\to F}\, dy\, , 
\qquad 
\end{eqnarray} 
$s\in {\Bbb R}$, $x\in F$. It follows from Remark (6) of Section 1 
and the definition of $g_n$ as well as $\gamma_n$ that, for fixed 
$s_0>0$, $F(s,x)$ is bounded on $(s,x)\in [-s_0,s_0]\times F$. 
Continuity of $F(s,x)$ on $(s,x)\in [-s_0,s_0]\times F$ is a 
consequence of condition (4) (i) of Section 1. Furthermore, $F(s, 
\, \cdot\, )$ is continuously differentiable with uniformly bounded 
gradient for all $s\in [-s_0,s_0]$, cf. again Remark (6) of Section 
1. Thus, the Picard-Lindel\"of theorem can be applied to all 
intervals of the form $[-s_0,s_0]$. It follows that there exists a 
global unique solution $\vp\equiv\vp(\, \cdot\, ;n,W)$ to the 
equation (\ref{3.10*}). 
\medskip 

\nid 
{\it Step 2 } We verify (i) and (vi). Equation (\ref{3.10*}) is 
equivalent to 
\begin{eqnarray}\label{3.12*} 
\vp(s)&&\hspace{-.5cm}=\int_0^s\dot{C}_{n,0}\left(\vp(v)\, \1+W_{\cdot 
+v}\vphantom{l^1}\right)\, dv 
\end{eqnarray} 
and also to 
\begin{eqnarray}\label{3.13*} 
W_{\cdot +s}+\vp(s)\, \1&&\hspace{-.5cm}=W_{\cdot +s}+\int_0^s\dot{C 
}_{n,0}\left(\vp(v)\, \1+W_{\cdot +v}\vphantom{l^1}\right)\, dv\, \1 
\, , 
\end{eqnarray} 
$s\in {\Bbb R}$. We observe that by (\ref{3.13*}) 
\begin{eqnarray}\label{3.14*}
W^{n,s}:=W_{\cdot +s}+\vp (s)\, \1\, , \quad s\in {\Bbb R}, 
\end{eqnarray} 
is 
the unique solution to 
\begin{eqnarray}\label{3.15*}
W^{n,s}&&\hspace{-.5cm}=W_{\cdot +s}+\int_0^s\dot{C}_{n,0}\left( 
W^{n,v}\right)\, dv\, \1\, , \quad s\in {\Bbb R}. 
\end{eqnarray} 

As a direct consequence, $W^{n,s}$, $s\in {\Bbb R}$, is a flow 
which is continuous in the topology of uniform convergence on compact 
sets. Let us introduce 
\begin{eqnarray}\label{3.16*}
D_{n,s}(W):=\int_0^s\dot{C}_{n,0}\left(W^{n,v}\right)\, dv\, , \quad 
s\in {\Bbb R}. 
\end{eqnarray} 
From the flow property of $W^{n,s}$, $s\in {\Bbb R}$, it follows that 
\begin{eqnarray}\label{3.17*}
D_{n,s+v}(W)=D_{n,s}(W)+D_{n,v}\left(W^{n,s}\right)\vphantom{\left( 
\dot{f}\right)}\, , \quad v\in {\Bbb R}.    
\end{eqnarray} 
We define 
\begin{eqnarray}\label{3.18*}
A_{n,s}(W):=D_{n,s}(W)+B_{n,0}\left(W^{n,s}\right)\quad\mbox{\rm as 
well as}\quad Y_{n,s}(W):=B_{n,0}\left(W^{n,s}\right)\, , 
\end{eqnarray} 
and   
\begin{eqnarray*} 
X_{n,s}(W):=W_s+A_{n,s}(W)\, , \quad s\in {\Bbb R}.  
\end{eqnarray*} 
The following becomes now obvious. By (\ref{3.15*}) and (\ref{3.18*}) 
we have 
\begin{eqnarray}\label{3.19*}
W^{n,s}_v=W_{s+v}+A_{n,s}-Y_{n,s}\, , \quad s,v\in {\Bbb R}. 
\end{eqnarray} 
Furthermore, we obtain from (\ref{3.17*}) and (\ref{3.18*}) and 
the flow property of $W^{n,s}$, $s\in {\Bbb R}$, 
\begin{eqnarray}\label{3.20*}
Y_{n,s}(W^{n,v})=B_{n,0}\left(W^{n,s+v}\right)=Y_{n,v}(W^{n,s})  
\end{eqnarray} 
as well as 
\begin{eqnarray}\label{3.21*}
A_{n,s+v}(W)-A_{n,s}(W)&&\hspace{-.5cm}=D_{n,s+v}(W)+B_{n,0} 
\left(W^{n,s+v}\right)-D_{n,s}(W)-B_{n,0}\left(W^{n,s}\right) 
\vphantom{\left(\dot{f}\right)}\nonumber \\ 
&&\hspace{-.5cm}=D_{n,v}\left(W^{n,s}\right)+B_{n,0}\left(W^{n,v} 
\left(W^{n,s}\right)\vphantom{l^1}\right)-B_{n,0}\left(W^{n,s} 
\right)\nonumber \\ 
&&\hspace{-.5cm}=A_{n,v}\left(W^{n,s}\right)-A_{n,0}\left(W^{n,s} 
\right)\vphantom{\left(\dot{f}\right)}\, , \quad s,v\in {\Bbb R}.  
\end{eqnarray} 
As a consequence of (\ref{3.18*}) we obtain $A_{n,0}\left( 
\displaystyle W^{n,s}\right)=Y_{n,s}(W)$, $s\in {\Bbb R}$. 
Together with (\ref{3.19*}) and (\ref{3.21*}) this leads to 
\begin{eqnarray}\label{3.22*}
X_{n,v}\left(W^{n,s}\right)&&\hspace{-.5cm}=W^{n,s}_v+A_{n,v}\left( 
W^{n,s}\right)\vphantom{\left(\dot{f}\right)}\nonumber \\ 
&&\hspace{-.5cm}=W^{n,s}_v+A_{n,0}\left(W^{n,s}\right)+A_{n,s+v}(W)- 
A_{n,s}(W)\vphantom{\left(\dot{f}\right)}\nonumber \\ 
&&\hspace{-.5cm}=W_{s+v}+A_{n,s}(W)-Y_{n,s}(W)+A_{n,0}\left(W^{n,s} 
\right)+A_{n,s+v}(W)-A_{n,s}(W)\vphantom{\left(\dot{f}\right)}
\nonumber \\ 
&&\hspace{-.5cm}=W_{v+s}+A_{n,v+s}(W)\vphantom{\left(\dot{f}\right)} 
\nonumber \\ 
&&\hspace{-.5cm}=X_{n,v+s}(W)\, , \quad s,v\in {\Bbb R}. \vphantom 
{\left(\dot{f}\right)}  
\end{eqnarray} 

Recall the definitions (\ref{3.14*}), (\ref{3.19*}). Taking in 
particular into consideration that $\vp(s)=A_{n,s}-Y_{n,s}$ is the 
unique solution to (\ref{3.10*}), it turns out that $A_{n,s}-Y_{n,s}$ 
is continuously differentiable in $s\in {\Bbb R}$. Let us now 
demonstrate that $Y_{n,s}$ and hence also $A_{n,s}$ is continuously 
differentiable in $s\in {\Bbb R}$. Similar to (\ref{3.10*}) the 
system of first order ODEs 
\begin{eqnarray*} 
\begin{array}{rl}
\dot{\xi}(s)\ =&\dot{B}_{n,0}\left((\eta(s)-\xi(s))\, \1+W_{\cdot +s} 
\vphantom{\displaystyle l^1}\right)\vphantom{\left(\dot{f}\right)} \\ 
\dot{\eta}(s)\ =&\dot{C}_{n,0}\left((\eta(s)-\xi(s))\, \1+W_{\cdot +s} 
\right)+\dot{B}_{n,0}\left((\eta(s)-\xi(s))\, \1+W_{\cdot +s}\vphantom 
{\displaystyle l^1}\right)\, ,\quad s\in {\Bbb R},\vphantom{\left(\dot 
{f}\right)} \\ 
\xi(0)\ =&B_{n,0}(W) \vphantom{\left(\dot{f}\right)} \\ 
\eta(0)\ =&C_{n,0}(W)+B_{n,0}(W)\, . \vphantom{\left(\dot{f}\right)}
\end{array} 
\end{eqnarray*} 
has for $W\in\Omega$ a unique solution $(\eta,\xi)\in C^1({\Bbb R};F^2 
)$. We remark that the point wise derivative $\dot{B}_{n,0}$ exists by 
the definition (\ref{3.8}) and the convolutions in the terms above 
(\ref{3.8}). 

By (\ref{3.10*}) and (\ref{3.14*}) as well as (\ref{3.19*}), we have 
$\eta(s)-\xi(s)=\vp(s)=A_{n,s}-Y_{n,s}$, $s\in {\Bbb R}$. The first 
line of the above system and (\ref{3.18*}) as well as (\ref{3.20*}) 
imply 
\begin{eqnarray*} 
\dot{\xi}(s)=\dot{B}_{n,0}\left(W^{n,s}\vphantom{\displaystyle l^1} 
\right)=\dot{Y}_{n,0}\left(W^{n,s}\vphantom{\displaystyle l^1}\right) 
=\frac{d}{ds}Y_{n,s}(W)\, , \quad s\in {\Bbb R}. 
\end{eqnarray*} 
Thus $Y_{n,s}\in C^1({\Bbb R};F)$ and therefore also $A_{n,s}\in C^1 
({\Bbb R};F)$. As a consequence, (\ref{3.21*}) results in 
\begin{eqnarray*} 
\dot{A}_{n,s+v}(W)=\dot{A}_{n,s}(W^{n,v})\, , \quad s,v\in {\Bbb R}. 
\end{eqnarray*} 
Below, we will frequently use the calculation rules (\ref{3.20*}) and  
(\ref{3.21*}) together with its differential form displayed in the last 
relation. 

Furthermore, we recall the above system of first order ODEs for $\xi$ 
and $\eta$ and the hypotheses $A(W)\equiv 0$ and $Y(W)\equiv 0$ 
whenever $W_0\not\in D^n$, cf. the definition of the process $X$ in 
Subsection 1.2 and condition (3) (iii) of Section 1. It follows that 
$A_{n,\cdot}$ and $Y_{n,\cdot}$ are constant on intervals where 
$W^{n,\cdot}_0\not\in\left\{x\in F:|x-z|<\frac{1}{n^3},\ z\in D^n 
\vphantom{W^{n,\cdot}}\right\}$.
\medskip 

We have defined $A_n$ and $Y_n$ by (\ref{3.8}), (\ref{3.9*}), 
(\ref{3.16*}), and (\ref{3.18*}) and verified (i) and (vi) of the lemma. 
\qed 
\medskip 

{\bf Proof of part (a). } {\it Step 1 } We intend to apply 
Gronwall's inequality for this. The inequality is stated in the end 
of Step 3. In this step we prepare the inequality by analyzing the 
structure of $A_{n,s}-Y_{n,s}-(A_s-Y_s)$. 

In order to ease the notation we will frequently write $(A-Y)(\, 
\cdot\, )$ for $A(\, \cdot\, )-Y(\, \cdot\, )$, $(A-Y)'$ for $A'-Y'$, 
and $\Delta (A-Y)$ for $\Delta A-\Delta Y$. Fix $W\in\{\pi_m V:V\in 
\Omega\}$ with $W_0\not\in G(W-W_0\1)$. It follows from (\ref{3.9*}), 
(\ref{3.15*}), (\ref{3.16*}), (\ref{3.18*}), and (\ref{3.19*}) that  
\begin{eqnarray}\label{3.23*}
&&\hspace{-.5cm}A_{n,s}(W)-Y_{n,s}(W)\vphantom{\int}=\int_{r=0}^s 
\dot{C}_{n,0}\left(W^{n,r}\right)\, dr\nonumber \\ 
&&\hspace{.0cm}=\int_{r=0}^s\int_F\int\left\langle\vphantom{\dot 
{f}}\left\langle (A-Y)'_{-v}(W^{n,r}+y\1)\, ,\, g_n(v)\vphantom{l^1} 
\right\rangle_{F\to F}\, ,\, \vphantom{\dot{f}}\gamma_n(y)\right 
\rangle_{F\to F}\, dv\, dy\vphantom{\sum_{v\in (-\frac1n,\frac1n)}} 
\, dr\nonumber \\ 
&&\hspace{.5cm}+\int_{r=0}^s\int_F\sum_{k\in {\Bbb Z}\setminus\{0\}} 
\left\langle\vphantom{\dot{f}}\left\langle\Delta (A-Y)_{\tau_k}\left 
(W^{n,r}+y\1\vphantom{l^1}\right)\, ,\, \right.\right.\nonumber \\ 
&&\hspace{4.5cm}\left.\left.g_n\left(\tau_k\circ u\left(W^{n,r}+y\1 
\vphantom{l^1}\right)\right)\vphantom{l^1}\right\rangle_{F\to F}\, 
,\, \vphantom{\dot{f}}\gamma_n(y)\right\rangle_{F\to F}\, dy\, dr\, 
, \qquad 
\end{eqnarray} 
$s\in {\Bbb R}$. We take into consideration Remark (5) of the present 
subsection which, in particular, implies $\tau_k\circ u\left(W^{n,r}+ 
y\1\vphantom{l^1}\right)=\tau_k\circ u(W^r)$ and set 
\begin{eqnarray}\label{3.24*}
\psi(r)&&\hspace{-.5cm}\equiv\psi(r;n,W):=A_{n,r}(W)-Y_{n,r}(W)-(A 
-Y)_r(W)\, ,\vphantom{\int^r_v} \nonumber \\ 
K(s,x)&&\hspace{-.5cm}\equiv K(s,x;n,W)\vphantom{\int}\nonumber \\ 
&&\hspace{-.5cm}:=\int_F\int\left\langle\vphantom{\dot{f}}\left 
\langle (A-Y)'_{-v}\left(W^s+(x+y)\1\vphantom{l^1}\right)-(A-Y)'_{-v} 
(W^s)\, ,\, \right.\right.\nonumber \\ 
&&\hspace{4.5cm}\left.\left.g_n(v)\vphantom{l^1}\right\rangle_{F\to 
F}\, ,\, \vphantom{\dot{f}}\gamma_n(y)\right\rangle_{F\to F}\, dv\, 
dy\vphantom{\sum_{v\in(-\frac1n,\frac1n)}}\nonumber \\  
&&\hspace{-.0cm}+\int_F\sum_{k\in {\Bbb Z}\setminus\{0\}}\left\langle 
\vphantom{\dot{f}}\left\langle\Delta (A-Y)_{\tau_k}\left(W^s+(x+y)\1 
\vphantom{l^1}\right)-\Delta (A-Y)_{\tau_k}(W^s)\, ,\, \right.\right. 
\nonumber \\ 
&&\hspace{4.5cm}\left.\left.g_n\left(\tau_k\circ u(W^s)\right) 
\vphantom{l^1}\right\rangle_{F\to F}\, ,\, \vphantom{\dot{f}}\gamma_n 
(y)\right\rangle_{F\to F}\, dy\vphantom{\sum_{v}}\nonumber \\ 
%
&&\hspace{-.5cm}\equiv K_1(s,x)+K_2(s,x)\, ,\vphantom{\sum^s_v} \\  
\rho(s)&&\hspace{-.5cm}\equiv\rho(s;n,W)\vphantom{\int}\nonumber 
 \\ 
&&\hspace{-.5cm}:=\int_{r=0}^s\int\left\langle (A-Y)'_{-v}(W^r)\, 
,\, g_n(v)\vphantom{l^1}\right\rangle_{F\to F}\, dv\, dr\vphantom 
{\sum_{v\in (-\frac1n,\frac1n)}}\nonumber \\  
&&\hspace{.0cm}+\int_{r=0}^s\sum_{k\in {\Bbb Z}\setminus\{0\}}\left 
\langle\Delta (A-Y)_{\tau_k}(W^r)\, ,\, g_n\left(\tau_k\circ u(W^r) 
\right)\vphantom{l^1}\right\rangle_{F\to F}\, dy\, dr\nonumber \\  
&&\hspace{.0cm}-(A-Y)_s(W)\vphantom{\int^s}
\end{eqnarray} 
relation (\ref{3.23*}) results in 
\begin{eqnarray}\label{3.25*} 
\psi(s)=\int_0^s K(v,\psi(v))\, dv+\rho(s)\, , \quad s\in {\Bbb R}. 
\end{eqnarray} 
{\it Step 2 } We construct estimates on $\left|\int_{r=0}^s K_i(r, 
\psi(r))\, dr\right|$, $i=1,2$, and $\left|\int_{r=0}^s K(r,\psi(r 
))\, dr\right|$ $+|\rho(s)|$. Since $\nabla_{W_0}A'_s$ and 
$\nabla_{W_0}Y'_s$ exist in the sense of conditions (2) and (3) of 
Section 1,
\begin{eqnarray*}
\alpha'(S;W):=\sup\left|\nabla_x(A-Y)'_{-v}\left(W^s+x\, \1 
\vphantom{l^1}\right)\vphantom{\dot{f}}\right|<\infty 
\end{eqnarray*} 
where the supremum is taken over all $(v,s,x)\in (-1,1)\times S 
\times F$ and $S\subset {}\Bbb R$ is any finite subinterval. It 
follows that 
\begin{eqnarray}\label{3.27*}
\left|\int_{r=0}^sK_1(r,\psi(r))\, dr\right|\le\alpha'(S;W)\int_{ 
r=0}^s\left(|\psi(r)|+\textstyle{\frac{1}{n^3}}\cdot{\sf e}\right) 
\, dr\, , \quad s\in S,  
\end{eqnarray} 
where the absolute value has been taken coordinate wise. 
\medskip 

Since $\Delta A$ and $\Delta Y$ are piecewise continuously 
differentiable on $F$, cf. conditions (2) and (3) of Subsection 
1.2, 
\begin{eqnarray*}
\beta'(S;W):=\sup\left|\nabla_x\Delta (A-Y)_{\tau_k}\left(W^s+x\, 
\1\vphantom{l^1}\right)\vphantom{\dot{f}}\right|<\infty  
\end{eqnarray*} 
where the supremum is taken over all $(k,s,x)\in\{l\in{\Bbb Z} 
\setminus\{0\}:\tau_l\in (-1,1)\}\times S\times F$ and $S\subset 
{\Bbb R}$ is any finite subinterval. Now, recall again that $\Delta 
A$ and $\Delta Y$ are piecewise continuously differentiable on $F$. 
In particular for any $k\in {\Bbb Z}\setminus\{0\}$, there exist 
finitely many mutually exclusive open sets $F_1\equiv F_1(k),F_2 
\equiv F_2(k),\ldots\, $ with piecewise $C^1$-boundary and 
$\overline{F}=\bigcup_i\overline{F_i}$ such that $\Delta A\left. 
(W+\, \cdot\, \1)\right|_{F_i}$, $\Delta Y\left.(W+\, \cdot\, \1) 
\right|_{F_i}\in C_b(F_i)$. 

The last relation together with $\beta'(S;W)<\infty$ and  condition 
(3) (ii) of Section 1 implies the existence of $\gamma'(S;W)< 
\infty$ such that 
\begin{eqnarray*}
\left|\Delta (A-Y)_{\tau_k}\left(W^s+x\1\vphantom{l^1}\right)-\Delta 
(A-Y)_{\tau_k}(W^s)\vphantom{\dot{f}}\right|\le\gamma'(S;W)\cdot |x| 
\vphantom{\sum}
\end{eqnarray*} 
for $(k,s,x)\in\{l\in{\Bbb Z}\setminus\{0\}:\tau_l\in (-1,1)\}\times 
S\times F$ and any finite subinterval $S\subset {\Bbb R}$. Also here, 
the absolute value has been taken coordinate wise. It follows that 
\begin{eqnarray}\label{3.28*} 
&&\hspace{-0.5cm}\left|\int_{r=0}^sK_2(r,\psi(r))\, dr\right|\nonumber 
 \\ 
&&\hspace{0.5cm}\le\gamma'(S;W)\sum_{k\in {\Bbb Z\setminus\{0\}}}\int_{ 
r=0}^s\left\langle\vphantom{\dot{f}}|\psi(r)|+\textstyle{\frac{1}{n^3}} 
\cdot{\sf e}\, ,\, g_n\left(\tau_k\circ u(W^r)\right)\vphantom{\dot{f}} 
\right\rangle_{F\to F}\, dr\, ,\quad s\in S. \qquad\qquad
\end{eqnarray} 

By $X(W^r)=X_{\cdot +r}$, cf. condition (3) of Section 1, we have 
\begin{eqnarray*}
\tau_k\circ u(W^r)=\tau_{k'}\circ u(W)-r 
\end{eqnarray*} 
for $k\in {\Bbb Z}\setminus\{0\}$ and some $k'\equiv k'(k;W,r)\in 
{\Bbb Z}\setminus\{0\}$ such that $k-k'$ is a constant $c\equiv 
c(W,r)$ if $k$ and $k'$ are both positive or both are negative. 
Moreover with the same constant $c$, $k-k'=c+1$ if $k>0>k'$ and 
$k-k'=c-1$ if $k<0<k'$, cf. the definition of the jumps in 
Subsection 1.2. 
\medskip 

Relations (\ref{3.24*})-(\ref{3.25*}) together with the estimates 
on $\int_{r=0}^sK_i(r,\psi(r))\, dr$, $i=1,2$, in (\ref{3.27*}) 
and (\ref{3.28*}), result in 
\begin{eqnarray}\label{3.29*} 
|\psi(s)|&&\hspace{-.5cm}\le\left|\int_0^s K(v,\psi(v))\, dv 
\right|+|\rho(s)|\vphantom{\sum_{k'\in {\Bbb Z}}}\nonumber \\ 
&&\hspace{-.5cm}\le\alpha'(S;W)\int_{r=0}^s\left(|\psi(r)|+ 
\textstyle{\frac{1}{n^3}}\cdot{\sf e}\right)\, dr\vphantom{ 
\sum_{k'\in {\Bbb Z}}}\nonumber \\ 
&&\hspace{.0cm}+\gamma'(S;W)\sum_{k\in {\Bbb Z}\setminus\{0\}} 
\int_{r=0}^s\left\langle|\psi(r)|+\textstyle{\frac{1}{n^3}} 
\cdot{\sf e}\, ,\, g_n\left(\tau_k\circ u(W^r)\right)\vphantom 
{\dot{f}}\right\rangle_{F\to F}\, dr+|\rho(s)|\vphantom{\sum_{ 
k'\in {\Bbb Z}\setminus\{0\}}}\nonumber \\ 
&&\hspace{-.5cm}=\int_{r=0}^s\left\langle|\psi(r)|\, ,\, \beta 
(r;n,W)\vphantom{\dot{f}}\right\rangle_{F\to F}\, dr+R(s;n,W)\, 
,\quad s\in S,\ n\in {\Bbb N}, 
\end{eqnarray} 
for any finite subinterval $S\subset {\Bbb R}$ where all absolute 
values are coordinate wise, 
\begin{eqnarray*}
R(s;n,W):=&&\hspace{-.5cm}\frac{\alpha'(S;W)\cdot s}{n^3}\cdot 
{\sf e}+\frac{\gamma'(S;W)}{n^3}\sum_{k'\in {\Bbb Z\setminus\{0 
\}}}\int_{r=0}^sg_n\left(\tau_{k'}-r\right)\, dr+|\rho(s)|\, ,  
\end{eqnarray*} 
and 
\begin{eqnarray*}
\beta(r;n,W):=\alpha'(S;W)\cdot {\sf e}+\gamma'(S;W)\sum_{k' 
\in {\Bbb Z\setminus\{0\}}}g_n\left(\tau_{k'}-r\right)\, . 
\end{eqnarray*} 

\nid 
{\it Step 3 } We apply Gronwall's inequality  to carry out the 
proof of part (a) of Lemma \ref{Lemma3.4}. From condition (3) 
and Remark (3) of Section 1 we get 
\begin{eqnarray}\label{3.30*}
&&\hspace{-.5cm}\rho(s;n,W)=\int_{r=0}^s\int_{v\in\left(-\frac1n, 
\frac1n\right)}\left\langle (A-Y)'_{-v}(W^r)\, ,\, g_n(v)\vphantom 
{\dot{f}}\right\rangle_{F\to F}\, dv\, dr\vphantom{\int_{u=0}^s} 
\nonumber \\ 
&&\hspace{1.0cm}+\int_{r=0}^s\sum_{v\in\left(-\frac1n,\frac1n\right 
)}\left\langle\Delta (A-Y)_{v}(W^r)\, ,\, g_n(v)\vphantom{\dot{f}} 
\right\rangle_{F\to F}\, dr-(A-Y)_s(W)\vphantom{\int_{u=0}^s} 
\nonumber \\ 
&&\hspace{.5cm}=\int_{r=0}^s\int_{v\in\left(-\frac1n,\frac1n\right)} 
\left\langle (A-Y)'_{r-v}(W)\, ,\, g_n(v)\vphantom{\dot{f}}\right 
\rangle_{F\to F}\, dv\, dr\vphantom{\int_{u=0}^s}\nonumber \\ 
&&\hspace{1.0cm}+\int_{r=0}^s\sum_{v\in\left(-\frac1n,\frac1n\right 
)}\left\langle\Delta (A-Y)_{r-v}(W)\, ,\, g_n(v)\vphantom{\dot{f}} 
\right\rangle_{F\to F}\, dr-(A-Y)_s(W)\vphantom{\int_{u=0}^s} 
\nonumber \\ 
&&\hspace{.5cm}=\int_{v\in\left(-\frac1n,\frac1n\right)}\left\langle 
(A-Y)_{s-v}(W)-(A-Y)_{-v}(W)\, ,\, g_n(v)\vphantom{\dot{f}}\right 
\rangle_{F\to F}\, dv\nonumber \\ 
&&\hspace{1.0cm}-(A-Y)_s(W)\vphantom{\int_{u=0}^s} 
\end{eqnarray} 
which implies 
\begin{eqnarray*}
\rho(s;n,W)\stack{n\to\infty}{\lra}0 
\end{eqnarray*} 
if neither $0$ nor $s$ is a jump time for $X$. With the last 
equality sign of (\ref{3.30*}) we get 
\begin{eqnarray*}
&&\hspace{-0.5cm}\left\langle\rho(s;n,W)\, ,\, \beta(s;n,W) 
\vphantom{l^1}\right\rangle_{F\to F}\vphantom{\int_v} \\ 
&&\hspace{0.5cm}=\alpha'(S;W)\int_{v\in {\Bbb R}}\left\langle 
\vphantom{\dot{f}}(A-Y)_{s-v}(W)-(A-Y)_s(W)\, ,\, g_n(v) 
\vphantom{\dot{f}}\right\rangle_{F\to F}\, dv \\ 
&&\hspace{1.0cm}+\gamma'(S;W)\sum_{k'\in {\Bbb Z}\setminus\{0\} 
}\int_{v\in {\Bbb R}}\left\langle\vphantom{\dot{f}}(A-Y)_v(W) 
\left\langle g_n(s-v)\, ,\, g_n\left(\tau_{k'}-s\right)\vphantom{ 
l^1}\right\rangle_{F\to F}\vphantom{\dot{f}}\right\rangle_{F\to 
F}\, dv \\ 
&&\hspace{1.0cm}-\gamma'(S;W)\sum_{k'\in {\Bbb Z\setminus\{0\}}} 
\left\langle\vphantom{\dot{f}}(A-Y)_s(W)\, , \, g_n\left(\tau_{ 
k'}-s\right)\vphantom{\dot{f}}\right\rangle_{F\to F} \\ 
&&\hspace{1.0cm}-\left\langle\int_{v\in {\Bbb R}}\left\langle 
(A-Y)_{-v}(W)\, ,\, g_n(v)\vphantom{\dot{f}}\right\rangle_{F\to 
F}\, dv\, ,\, \beta(s;n,W)\right\rangle_{F\to F}\, .  
\end{eqnarray*} 
Keeping in mind $A_0=Y_0$ we obtain 
\begin{eqnarray*} 
\left\langle\rho(\, \cdot\, ;n,W)\, ,\, \beta(\, \cdot\, ;n,W) 
\vphantom{l^1}\right\rangle_{F\to F}\stack{n\to\infty}{\lra}0 
\quad\mbox{\rm in }L^1([0,s];F) 
\end{eqnarray*} 
if, without loss of generality, $0\le s$ and neither $0$ nor $s$ is 
a jump time for $X$. It is a consequence of the definition of $R(\, 
\cdot\, ;n,W)$ and the fact that 
\begin{eqnarray*}
\sum_{k'\in {\Bbb Z\setminus\{0\}}}\int_{r=0}^sg_n\left(\tau_{k'} 
-r\right)\, dr 
\end{eqnarray*} 
is on every component uniformly bounded in $n$ by the number of 
jumps of $X=u(W)$ on $(-1,s+1]$ that 
\begin{eqnarray}\label{3.31*} 
\begin{array}{rcl}
R(s;n,W)&\stack{n\to\infty}{\lra}&0\,  \\ & & \\ 
\left\langle R(\, \cdot\, ;n,W)\, ,\, \beta(\, \cdot\, ;n,W)\vphantom 
{l^1}\right\rangle_{F\to F}&\stack{n\to\infty}{\lra}&0\quad\mbox{\rm 
in }L^1([0,s];F) 
\end{array}
\end{eqnarray} 
if, without loss of generality, $0\le s$ and neither $0$ nor $s$ are 
jump times for $X$. It follows now from (\ref{3.29*}) and Gronwall's 
inequality that 
\begin{eqnarray}\label{3.32*}
&&\hspace{-.5cm}\left|A_{n,s}(W)-Y_{n,s}(W)-\left(A_s(W)-Y_s(W) 
\vphantom{l^1}\right)\vphantom{\dot{f}}\right|=|\psi(s)|\nonumber 
 \\ 
&&\hspace{.5cm}\le R(s;n,W)+\int_{v=0}^s\left\langle R(v;n,W)\, ,\, 
\beta(v;n,W)\vphantom{l^1}\right\rangle_{F\to F} \cdot\exp{\left\{ 
\int_{r=v}^s|\beta(r;n,W)|\, dr\right\}}\, dv\nonumber \\ 
&&\hspace{.5cm}\le R(s;n,W)+\exp{\left\{\int_{r=0}^s|\beta(r;n,W)| 
\, dr\right\}}\cdot\int_{v=0}^s\left\langle R(v;n,W)\, ,\, \beta(v; 
n,W)\vphantom{l^1}\right\rangle_{F\to F}\, dv\, , \nonumber \\ 
\end{eqnarray} 
$n\in {\Bbb N}$. Here the absolute values in the first line are 
coordinate wise, the absolute value on $\beta$ is not, cf. Remark 
(4) of the present section. Part (a) follows now on the one hand 
from (\ref{3.31*}). On the other hand we use again the fact that 
$\sum_{k'\in {\Bbb Z}\setminus\{0\}}\int_{r=0}^sg_n\left(\tau_{k'} 
-r\right)\, dr$ is on every coordinate uniformly bounded in $n$. 
This implies that $\int_{r=0}^s|\beta(r;n,W)|\, dr$  is uniformly 
bounded in $n\in {\Bbb N}$ and we have accomplished the proof of (a). 
\medskip 

As an addendum to this step but also as a preparation to other  
proofs of the present subsection we obtain the following. By 
(\ref{3.9*}) and (\ref{3.18*}) we have 
\begin{eqnarray*}
&&\hspace{-.5cm}\dot{A}_{n,s}(W)-\dot{Y}_{n,s}(W)=\int_{\Bbb R} 
\left\langle A_{-v}(\, \cdot\, , \gamma_n)(W^{n,s})-Y_{-v}(\,  
\cdot\, ,\gamma_n)(W^{n,s})\, ,\, g'_n(v)\vphantom{l^1}\right 
\rangle_{F\to F}\, dv 
\end{eqnarray*} 
and therefore 
\begin{eqnarray}\label{3.33*}
\left|\dot{A}_{n,s}(W)-\dot{Y}_{n,s}(W)\right|&&\hspace{-.5cm}\le 
\frac{2\|g'_n\|}{n}\sup_{v\in\left(-\frac1n,\frac1n\right)}\left| 
A_{-v}(\, \cdot\, , \gamma_n)(W^{n,s})-Y_{-v}(\, \cdot\, , 
\gamma_n)(W^{n,s})\vphantom{\dot{f}}\right|\qquad 
\end{eqnarray} 
$s\in {\Bbb R}$, $n\in {\Bbb N}$, where the absolute value $|\, 
\cdot\, |$ has been taken individually for each coordinate of $F$. 

We have $\sup_{s,W}|Y_s(W)|<\infty$ by condition (3) of Section 
1 and $D^n$ is bounded by hypothesis. If $W^{n,s}_0+x\not\in D^n$ 
then, by hypothesis of Subsection 1.2, $A(W^{n,s}+x\1)\equiv 0$. 
On the other hand if $W^{n,s}_0+x\in D^n$ then, again by 
hypothesis of Subsection 1.2, we may assume that $X_\cdot (W^{n,s} 
+x\1)\in D^n$ and take into consideration that $A_\cdot (W^{n,s}+ 
x\1)=X_\cdot(W^{n,s}+x\1)-(W^{n,s}_0+x)-\left(W_{\cdot +s}-W_s 
\right)$. 

In any case, by (\ref{3.33*}) there exists a positive constant 
$c\equiv c(n)$ which may depend on $n\in {\Bbb N}$ but is 
independent of $s\in {\Bbb R}$ and $W\in \{\pi_m V:V\in\Omega\}$ 
such that 
\begin{eqnarray}\label{3.34*} 
\left|\dot{A}_{n,s}(W)-\dot{Y}_{n,s}(W)\right|&&\hspace{-.5cm} 
\le c(n)\, {\sf e}+\frac{2\|g'_n\|}{n}\sup_{v\in\left(-\frac1n, 
\frac1n\right)}\left|W_{s-v}-W_s\right|\, .  
\end{eqnarray} 
\medskip 
\qed

{\bf Proof of parts (c) and (v). } Let $s\in {\Bbb R}$, $n\in {\Bbb 
N}$, $W\in\{\pi_m V:V\in\Omega\}$ with $W_0\not\in G(W-W_0\1)$, and  
\begin{eqnarray}\label{3.35*}  
\kappa(s)\equiv\kappa(s;n,W)\vphantom{\left(\dot{f}\right)}:= 
\nabla_{W_0}\rho(s;n,W)\, ,   
\end{eqnarray} 
for the well-definiteness of this expression for all $s\in {\Bbb 
R}$ recall (\ref{3.30*}) and Remark (7) of Section 1. Let us also 
recall the existence of a local spatial gradient for $A^1$ and 
$Y^1$ and the piecewise continuous differentiability of the jump 
magnitudes, cf. conditions (2) as well as (3) and Remark (2) of 
Section 1. As a consequence of the definition of $K$ in (\ref{3.24*}) 
the integral $\int_0^s$ in (\ref{3.25*}) commutes with the derivative 
$\nabla_{W_0}$. Keeping these arguments in mind and noting that 
\begin{eqnarray*}
\dot{C}_{n,0}=\int_{v\in\left(-\frac2n,0\right)}\left\langle A_{
-v}(\, \cdot\, , \gamma_n)-Y_{-v}(\, \cdot\, , \gamma_n)\, , \, 
g_n'\left(v+{\textstyle\frac1n}\right)\vphantom{A_{\frac1n}}\right 
\rangle_{F\to F}\, dv
\end{eqnarray*} 
we verify also that the integral $\int_0^s$ in (\ref{3.12*}) commutes 
with the derivative $\nabla_{W_0}$. It follows now from (\ref{3.12*}) 
that $\nabla_{W_0}\vp\equiv\nabla_{W_0}\vp(\, \cdot\, ;n,W)$ satisfies, 
whenever it exists,  
\begin{eqnarray}\label{3.39*}
\begin{array}{rl}
\left(\nabla_{W_0}\vp(s)\right)\dot{}\ =&\left(\nabla_{W_0} 
\dot{C}_{n,0}\right)\left(W^{n,s}\vphantom{\displaystyle l^1}\right) 
\cdot\left(\nabla_{W_0}\vp(s)+{\sf e}\right)\, , \quad s\in {\Bbb R}, 
\vphantom{\left(\dot{f}\right)} \\ 
\nabla_{W_0}\vp(0)\ =&0 \vphantom{\left(\dot{f}\right)}
\end{array}\, .  
\end{eqnarray} 
But existence of $\nabla_{W_0}\vp$ follows by the Picard-Lindel\"of 
theorem. Together with (\ref{3.12*}) this yields, en passant, (v). 
More precisely, we note that $\nabla_{W_0}\vp$ is the unique solution 
to (\ref{3.39*}) by the representation 
\begin{eqnarray}\label{3.40*}
\nabla_{W_0}\dot{C}_{n,0}&&\hspace{-.5cm}=\int_{v\in\left(-\frac 
1n,\frac1n\right)}\left\langle\nabla_{W_0}\left(A_{-v}(\, \cdot\, , 
\gamma_n)\vphantom{l^1}\right)-\nabla_{W_0}\left(Y_{-v}(\, \cdot\, 
, \gamma_n)\vphantom{l^1}\right)\, ,\, g_n'(v)\vphantom{A_{\frac1n} 
}\right\rangle_{F\to F}\, dv\qquad
\end{eqnarray} 
and 
\begin{eqnarray}\label{3.41*} 
\nabla_{W_0}\left(A_{-v}(\, \cdot\, ,\gamma_n)(W^{n,s})\vphantom 
{l^1}\right)=-\int_F\left\langle A_{-v}(W^{n,s}+x\1)\, {\sf e}\, , 
\, \nabla_x\gamma_n(x)\vphantom{l^1}\right\rangle_{F\to F}\, dx\, ,
\end{eqnarray} 
the multiplication $\langle\ , \, \rangle_{F\to F}$ taken here 
adequately coordinate wise; a similar representation holds for 
$\nabla_{W_0}\left(Y_{-v}(\, \cdot\, ,\gamma_n)(W^{n,s})\D 
\vphantom{l^1}\right)$. In fact, the Picard-Lindel\"of theorem 
yields uniqueness by the continuity of $s\to W^{n,s}$ demonstrated 
in the proof of (i) of Lemma \ref{Lemma3.4} which together with 
condition (4) (i) of Subsection 1.2 implies continuity of 
\begin{eqnarray*}
s\to \left(\nabla_{W_0}\dot{C}_{n,0}\right)\left(W^{n,s}\vphantom 
{\displaystyle l^1}\right)\, . 
\end{eqnarray*} 

It follows from (\ref{3.25*}) that 
\begin{eqnarray}\label{3.36*} 
\nabla_{W_0}\psi(s)=\int_0^s\left(\nabla_xK\vphantom{l^1}\right) 
(v,\psi(v))\cdot \nabla_{W_0}\psi(v)\, dv+\kappa(s)\, , \quad s 
\in {\Bbb R},  
\end{eqnarray} 
where, for the well-definiteness of $\nabla_{W_0}\psi(s)$ for all 
$s\in {\Bbb R}$, we first recall the well-definiteness of 
$\nabla_{W_0}\vp$. Then, for $\nabla_{W_0}\psi$ we keep Remark (7) 
of Section 1 in mind. According to (\ref{3.24*}) and Remark (5) of 
the present section we have also 
\begin{eqnarray*}
K(s,x)&&\hspace{-.5cm}=\int_F\int\left\langle\vphantom{\dot{f}} 
\left\langle (A-Y)'_{-v}(W^s+(x+y)\1)-(A-Y)'_{-v}(W^s)\, ,\, 
\right.\right.\nonumber \\ 
&&\hspace{5cm}\left.\left.g_n(v)\vphantom{l^1}\right\rangle_{ 
F\to F}\, ,\, \vphantom{\dot{f}}\gamma_n(y)\right\rangle_{F 
\to F}\, dv\, dy\vphantom{\sum_{v\in (-\frac1n,\frac1n)}} \\ 
&&\hspace{-.0cm}+\int_F\sum_{k\in {\Bbb Z}\setminus\{0\}}\left 
\langle\vphantom{\dot{f}}\left\langle\Delta (A-Y)_{\tau_k}\left 
(W^s+(x+y)\1\vphantom{l^1}\right)-\Delta (A-Y)_{\tau_k}(W^s)\, 
,\, \right.\right.\nonumber \\ 
&&\hspace{5cm}\left.\left.g_n\left(\tau_{k'}-s\right)\vphantom 
{l^1}\right\rangle_{F\to F}\, ,\, \vphantom{\dot{f}}\gamma_n(y) 
\right\rangle_{F\to F}\, dy\vphantom{\int_F}\, ,\quad x\in F.   
\end{eqnarray*} 
By similar arguments as in the previous step, using conditions 
(1)-(3) of Section 1, (\ref{3.36*}) and this representation of 
$K$ yield 
\begin{eqnarray}\label{3.37*}
&&\hspace{-.5cm}\left|\nabla_{W_0}\psi(s)\right|\le\alpha'(S;W) 
\int_{r=0}^s |\nabla_{W_0}\psi(r)|\, dr\vphantom{\sum_{k'\in {\Bbb 
Z}}}\nonumber \\ 
&&\hspace{.3cm}+\beta'(S;W)\sum_{k'\in {\Bbb Z\setminus\{0\}}} 
\int_{r=0}^s|\nabla_{W_0}\psi(r)|\cdot\left|g_n\left(\tau_{k'} 
-r\right)\right|\, dr+|\kappa(s)|\qquad 
\end{eqnarray} 
for $s\in S$, $n\in {\Bbb N}$, and any finite subinterval $S\subset 
{\Bbb R}$. With 
\begin{eqnarray*}
\delta(r;n,W):=\alpha'(S;W)+\beta'(S;W)\sum_{k'\in {\Bbb Z\setminus 
\{0\}}}\left|g_n\left(\tau_{k'}-r\right)\right|
\end{eqnarray*} 
we obtain as in (\ref{3.32*}) from (\ref{3.37*}) and Gronwall's 
inequality 
\begin{eqnarray}\label{3.38*} 
&&\hspace{-.5cm}\left|\vphantom{\dot{f}}\nabla_{W_0}\left(A_{n,s}- 
Y_{n,s}\vphantom{l^1}\right)-\left(\nabla_{W_0}A_s-\nabla_{W_0}Y_s 
\vphantom{l^1}\right)\vphantom{\dot{f}}\right|=\left|\nabla_{W_0} 
\psi(s)\right|\nonumber \\ 
&&\hspace{.5cm}\le |\kappa(s;n,W)|+\exp{\left\{\int_{r=0}^s\delta 
(r;n,W)\, dr\right\}}\cdot\int_{r=0}^s|\kappa(r;n,W)|\cdot\delta 
(r;n,W)\, dr\, , \qquad 
\end{eqnarray} 
$n\in {\Bbb N}$, $s\in S$ for any finite, without loss of generality 
non-negative, subinterval $S\subset {\Bbb R}$. Adjusting the arguments 
of (\ref{3.30*}) and (\ref{3.31*}) to $\nabla_{W_0}A$ and $\nabla_{W_0} 
Y$ and in particular noticing the right continuity of $\nabla_{W_0}A$ 
and $\nabla_{W_0}Y$, cf. Remark (7) of Section 1, we get from 
(\ref{3.35*}) 
\begin{eqnarray*} 
\begin{array}{rcl}
\kappa(s;n,W)&\stack{n\to\infty}{\lra}&0 \\ 
& & \\ 
|\kappa(\, \cdot\, ;n,W)|\cdot\delta(\, \cdot\, ;n,W)&\stack{n\to 
\infty}{\lra}&0\quad\mbox{\rm in }L^1([0,s]) 
\end{array}
\end{eqnarray*} 
if, without loss of generality, $0\le s$ and neither $0$ nor $s$ is a
jump time for $X$. Together with (\ref{3.38*}) we get from here that 
\begin{eqnarray*}
&&\hspace{-.5cm}\nabla_{W_0}\left((A_{n,l}-Y_{n,l})-(A_{n,r}-Y_{n,r}) 
\vphantom{l^1}\right)\vphantom{\int} \\ 
&&\stack{n\to\infty}{\lra}\left(\nabla_{W_0}A_l-\nabla_{W_0}Y_l 
\vphantom{l^1}\right)-\left(\nabla_{W_0}A_r-\nabla_{W_0}Y_r\vphantom 
{l^1}\right)\, ,\quad l<r,  
\end{eqnarray*} 
if neither $l$ nor $r$ is a jump time for $X$. 

Now we take special attention to condition (1) (iv), the uniform 
boundedness of $\nabla_{W_0}A'$ in condition (2), the pendant for 
$\nabla_{W_0}Y'$ according to condition (3) in Section 1, and 
Remark (7) of Section 1. Reviewing the definitions of $\delta$ and 
$\kappa$ via the definitions of $\alpha'$, $\beta'$, and (\ref{3.30*}) 
we verify that the latter limit even holds uniformly bounded $Q^{(m)
}_{\sbnu}$-a.e. on $\{\pi_mW:W\in\Omega\}$.  We have shown (c) where, 
for its formulation, we recall that according to Remark (3) of this 
section $\nabla_{W_0}A$ and $\nabla_{W_0}Y$ jump only on $Q^{(m)}_{ 
\sbnu}$-zero sets. 
\medskip 
\qed

{\bf Proof of part (b). }  {\it Step 1 } In this step we verify the 
claims relative to $A-Y$ and $\dot{A}-\dot{Y}$. We apply Gronwall's 
inequality for the proofs of the limits. 

Let $W\in\{\pi_mV:V\in\Omega\}$ with $W_0\not\in G(W-W_0\1)$ and 
$W^{(r)}:=\pi_{m,r}W$. As above, let $\vp$ denote the solution to 
(\ref{3.10*}) but let $\vp^{(r)}$ denote the solution to (\ref{3.10*}) 
for $W$ replaced with $W^{(r)}$. Let $s\in {\Bbb R}$, $\ve >0$, and 
choose $r\in {\Bbb N}$ such that $|s|+1<rt$ and 
\begin{eqnarray}\label{3.42*}
&&\hspace{-.5cm}\int_{v\in\left(-\frac1n,\frac1n\right)}\left|A_{-v} 
(\, \cdot\, ,\gamma_n)\left(W_{\cdot +s}+\vp^{(r)}(s)\1\vphantom{l^1} 
\right)\vphantom{\left(W^{(r)}_{\cdot +s}\right)}-Y_{-v}(\, \cdot\, , 
\gamma_n)\left(W_{\cdot +s}+\vp^{(r)}(s)\1\vphantom{l^1}\right)\right. 
\nonumber \\ 
&&\hspace{1cm}-\left.\left(A_{-v}(\, \cdot\, ,\gamma_n)\left(W^{(r) 
}_{\cdot +s}+\vp^{(r)}(s)\1\vphantom{l^1}\right)-Y_{-v}(\, \cdot\, ,  
\gamma_n)\left(W^{(r)}_{\cdot +s}+\vp^{(r)}(s)\1\vphantom{l^1}\right) 
\right)\right|\, dv\vphantom{\int_{v\in\left(-\frac1n,\frac1n\right)}} 
\nonumber \\ 
&&\hspace{.5cm}<\frac{n\ve}{2\|g'_n\|}\, . 
\end{eqnarray} 
Such a choice of $r\in {\Bbb N}$ is possible according to condition 
(4) (ii) of Section 1. It follows from Remark (6) of Section 1 that 
there is some $\beta >0$ not depending on $r$ such that 
\begin{eqnarray}\label{3.43*}
&&\hspace{-.5cm}\left|\vphantom{\dot{f}}A_{-v}(\, \cdot \, ,\gamma_n) 
\left(W_{\cdot +s}+\vp(s)\, \1\vphantom{l^1}\right)-A_{-v}(\, \cdot\, 
,\gamma_n)\left(W_{\cdot +s}+\vp^{(r)}(s)\, \1\vphantom{l^1}\right) 
\vphantom{\left|\dot{f}\right|}\right.\vphantom{\int}\nonumber \\ 
&&\hspace{1.0cm}-\left(\vphantom{\dot{f}}Y_{-v}(\, \cdot\, ,\gamma_n 
)\left(W_{\cdot +s}+\vp(s)\1\vphantom{l^1}\right)\vphantom{\left|\dot 
{f}\right|}\left.-Y_{-v}(\, \cdot\, ,\gamma_n)\left(W_{\cdot +s}+\vp^{ 
(r)}(s)\, \1\vphantom{l^1}\right)\vphantom{\dot{f}}\right)\vphantom{ 
\dot{f}}\right|\vphantom{\int}\nonumber \\ 
&&\hspace{.5cm}<\frac{n\beta}{2\|g'_n\|}\cdot\left|\vp(s)-\vp^{(r)}(s) 
\vphantom{\dot{f}}\right| 
\end{eqnarray} 
for all $v\in\left(-\frac1n,\frac1n\right)$, $|s|+1<rt$. It follows 
from (\ref{3.10*}) and (\ref{3.11*}) that 
\begin{eqnarray*}
\dot{\vp}(s)-\dot{\vp^{(r)}}(s)\vphantom{\int}&&\hspace{-.5cm}=\int_{ 
\left(-\frac1n,\frac1n\right)}\left\langle A_{-v}(\, \cdot\, ,\gamma_n) 
\left(W_{\cdot+ s}+\vp(s)\, \1\vphantom{l^1}\right)\vphantom{\left(W^{ 
(r)}_{\cdot +s}\right)}\right.\nonumber \\ 
&&\hspace{.5cm}\left.-A_{-v}(\, \cdot\, ,\gamma_n)\left(W^{(r)}_{\cdot 
+s}+\vp^{(r)}(s)\, \1\right)\, ,\, g'_n(v)\right\rangle_{F\to F}\, dv 
\vphantom{\int_{\Bbb R}}\nonumber \\ 
&&\hspace{.0cm}-\int_{\left(-\frac1n,\frac1n\right)}\left\langle Y_{-v} 
(\, \cdot\, ,\gamma_n)\left(W_{\cdot+s}+\vp(s)\, \1\vphantom{l^1}\right) 
\vphantom{\left(W^{(r)}_{\cdot +s}\right)}\right.\nonumber \\ 
&&\hspace{.5cm}\left.-Y_{-v}(\, \cdot\, ,\gamma_n)\left(W^{(r)}_{\cdot 
+s}+\vp^{(r)}(s)\, \1\right)\, ,\, g'_n(v)\right\rangle_{F\to F}\, dv 
\vphantom{\int_{\Bbb R}}\, . 
\end{eqnarray*} 
Together with relation (\ref{3.42*}) and (\ref{3.43*}), this implies 
\begin{eqnarray}\label{3.44*}
\left|\dot{\vp}(s)-\dot{\vp^{(r)}}(s)\vphantom{\dot{f}}\right| 
&&\hspace{-.5cm}<\ve+\beta\cdot\left|\vp(s)-\vp^{(r)}(s)\vphantom 
{\dot{f}}\right| \vphantom{\int_{\Bbb R}}\, . 
\end{eqnarray} 

By (\ref{3.10*}) we have $\vp(0)=0=\vp^{(r)}(0)$. Moreover, because of 
(\ref{3.14*}) and (\ref{3.19*}) it holds that $\vp(s)=A_{n,s}(W)-Y_{n,s} 
(W)$ and $\vp^{(r)}(s)=A_{n,s}\left(W^{(r)}\right)-Y_{n,s}\left(W^{(r)} 
\right)$. By means of Gronwall's inequality, (\ref{3.44*}) implies now 
\begin{eqnarray}\label{3.45*}
&&\hspace{-.5cm}\left|\dot{A}_{n,s}(W)-\dot{Y}_{n,s}(W)-\left(\dot{A 
}_{n,s}\left(W^{(r)}\right)-\dot{Y}_{n,s}\left(W^{(r)}\right)\right) 
\vphantom{\dot{f}}\right|=\left|\dot{\vp}(s)-\dot{\vp^{(r)}}(s) 
\vphantom{\dot{f}}\right|\vphantom{\int} \nonumber \\ 
&&\hspace{.5cm}\le\ve\cdot e^{\beta |s|}\, , \quad |s|\le rt. 
\end{eqnarray} 
We have thus verified $\dot{A}_{n,\cdot}(\pi_{m,r}W)-\dot{Y}_{n,\cdot} 
(\pi_{m,r}W)\stack{r\to\infty}{\lra}\dot{A}_{n,\cdot}(W)-\dot{Y}_{n, 
\cdot}(W)$, uniformly on every finite subinterval $S\subset{\Bbb R}$. 
Relation $A_{n,\cdot}(\pi_{m,r}W)-Y_{n,\cdot}(\pi_{m,r}W)\stack{r\to 
\infty}{\lra}A_{n,\cdot}(W)-Y_{n,\cdot}(W)$, uniformly on every finite 
subinterval $S\subset{\Bbb R}$ follows now from (\ref{3.45*}) together 
with $A_{n,0}-Y_{n,0}=0$, cf. (\ref{3.18*}). The latter also implies 
\begin{eqnarray}\label{3.46*}
\left(W^{(r)}_v\right)^{n,s}\stack{r\to\infty}{\lra}W^{n,s}_v  
\end{eqnarray} 
uniformly on $(s,v)\in S^2$ for every finite subinterval $S\subset{\Bbb 
R}$. Furthermore, taking the absolute value $|\, \cdot\, |$ and the 
$\sup_{r\in {\Bbb N}}$ individually for each coordinate of $F$, 
\begin{eqnarray}\label{3.47*}
\sup_{r\in {\Bbb N}}\left|\dot{A}_{n,s}(\pi_{m,r}W)-\dot{Y}_{n,s} 
(\pi_{m,r}W)\right|\in\bigcap_{1\le p<\infty} L^p(\Omega,Q^{(m)}_{ 
\sbnu};F)\, ,\quad s\in {\Bbb R},  
\end{eqnarray} 
is a consequence of (\ref{3.34*}) and the finite L\'evy-Ciesielsky 
construction of $\pi_{m,r}W$ in (\ref{2.1}). Note that the $\langle H_i, 
dW\rangle_{L^2}$, $i\in I(m,r)$, in (\ref{2.1}) are independent $N(0,1) 
$-random variables and that by the $\sup_{v\in\left(-\frac1n,\frac1n 
\right)}$ in (\ref{3.34*}) we use just a selection of them which is 
independent of $r$. 
\medskip 

\nid
{\it Step 2 } We prove the claims relative to $\nabla_{W_0}\dot{A}- 
\nabla_{W_0}\dot{Y}$. Again we use Gronwall's inequality. 

Let $W\in\{\pi_m V:V\in\Omega\}$ with $W_0\not\in G(W-W_0\1)$. Recalling 
(\ref{3.18*}), (\ref{3.19*}), as well as (\ref{3.39*}) and (\ref{3.40*}) 
we obtain  
\begin{eqnarray}\label{3.48*}
&&\hspace{-.5cm}\left(\nabla_{W_0}A_{n,s}\vphantom{\dot{f}}\right)^\cdot 
\left(W^{(r)}\right)-\left(\nabla_{W_0}Y_{n,s}\vphantom{\dot{f}}\right 
)^\cdot\left(W^{(r)}\right)=\nabla_{W_0}\dot{A}_{n,s}\left(W^{(r)}\right) 
-\nabla_{W_0}\dot{Y}_{n,s}\left(W^{(r)}\right) \nonumber \\ 
&&\hspace{.5cm}=\nabla_{W_0}\left(\dot{C}_{n,0}\circ\left(W^{(r)}\right 
)^{n,s}\right) \nonumber \\ 
&&\hspace{.5cm}=\left(\nabla_{W_0}\dot{C}_{n,0}\right)\left(\left(W^{(r)} 
\right)^{n,s}\right)\cdot\left(\nabla_{W_0}A_{n,s}\left(W^{(r)}\right)- 
\nabla_{W_0}Y_{n,s}\left(W^{(r)}\right)+{\sf e}\vphantom{\dot{f}}\right) 
\, , \qquad 
\end{eqnarray} 
$s\in {\Bbb R}$. Let us now consider (\ref{3.48*}) a first order system 
\begin{eqnarray}\label{3.49*}
\begin{array}{rl}
\dot{\xi(s)}\ =&\left(\nabla_{W_0}\dot{C}_{n,0}\right)\left(\left(W^{(r)} 
\right)^{n,s}\right)\cdot\left(\xi(s)+{\sf e}\right)\, ,\quad s\in {\Bbb R}, 
\vphantom{\left(\dot{f}\right)} \\ 
\xi(0)\ =&0\vphantom{\left(\dot{f}\right)}
\end{array} 
\end{eqnarray} 
with solution 
\begin{eqnarray*} 
\xi(s)=\nabla_{W_0}A_{n,s}\left(W^{(r)}\right)-\nabla_{W_0}Y_{n,s} 
\left(W^{(r)}\right)\, , \quad |s|\le rt,  
\end{eqnarray*} 
which is unique by the arguments for (\ref{3.39*})-(\ref{3.41*}). 
Together with representations (\ref{3.40*}), (\ref{3.41*}), relation 
(\ref{3.46*}), and condition (4) (i) of Section 1 we get 
\begin{eqnarray*}
\left(\nabla_{W_0}\dot{C}_{n,0}\right)\left(\left(W^{(r)}\right)^{n,s} 
\right)\stack{r\to\infty}{\lra}\left(\nabla_{W_0}\dot{C}_{n,0}\right) 
\left(W^{n,s}\right) 
\end{eqnarray*} 
uniformly on every finite subinterval $s\in S\subset{\Bbb R}$. Now 
it follows from (\ref{3.39*}) and (\ref{3.49*}) that 
\begin{eqnarray}\label{3.50*}
\nabla_{W_0}\dot{A}_{n,s}\left(W^{(r)}\right)-\nabla_{W_0}\dot{Y}_{ 
n,s}\left(W^{(r)}\right)\stack{r\to\infty}{\lra}\nabla_{W_0}\dot{A 
}_{n,s}(W)-\nabla_{W_0}\dot{Y}_{n,s}(W)\, , \quad s\in {\Bbb R}. 
\end{eqnarray} 
By (\ref{3.39*})-(\ref{3.41*}) we have 
\begin{eqnarray}\label{3.51*}
&&\hspace{-.5cm}\nabla_{W_0}\dot{A}_{n,s}(W)-\nabla_{W_0}\dot{Y 
}_{n,s}(W)\vphantom{\int}\nonumber \\ 
&&\hspace{.5cm}=\left(\int_{\left(-\frac1n,\frac1n\right)}\left 
\langle\int_F\left\langle\left(A_{-v}\left(\, \cdot\, +x\1 
\vphantom{l^1}\right)-Y_{-v}\left(\, \cdot\, +x\1\vphantom{l^1} 
\right)\vphantom{\dot{f}}\right){\sf e}\, ,\, \right.\right. 
\right.\nonumber \\ 
&&\hspace{1.5cm}\left.\left.\left.\nabla_{W_0}\gamma_n(x) 
\vphantom{Y_{-v}}\vphantom{\dot{f}}\right\rangle_{F\to F}\, dx\, 
,\, g'_n(v)\vphantom{\int_F}\right\rangle_{F\to F}\, dv\vphantom 
{\int_{\left(-\frac1n,\frac1n\right)}}\right)\left(W^{n,s}\right) 
\times\nonumber \\ 
&&\hspace{1cm}\times\left(\nabla_{W_0}A_{n,s}(W)-\nabla_{W_0}Y_{n 
,s}(W)\vphantom{\dot{f}}+{\sf e}\right)\vphantom{\left(\int_{\Bbb 
R}\right)}\, ,\quad s\in {\Bbb R}. 
\end{eqnarray} 
Let us set 
\begin{eqnarray*}
b(n):=\frac{2\|g'_n\|\|\nabla \gamma_n\|\lambda_F\left(B_{\frac{1} 
{n^3}}\right)}{n} 
\end{eqnarray*} 
where $B_{\frac{1}{n^3}}\subset F$ denotes any ball with radius 
$\frac{1}{n^3}$. As for the verification of (\ref{3.34*}) we obtain 
from (\ref{3.51*}) 
\begin{eqnarray}\label{3.52*} 
&&\hspace{-0.5cm}\left|\nabla_{W_0}\dot{A}_{n,s}(W)-\nabla_{W_0} 
\dot{Y}_{n,s}(W)\right|\vphantom{\left(\int_R\right)}\nonumber \\ 
&&\hspace{0.5cm}\le\left(a(n)+b(n)\sup_{v\in\left(-\frac1n,\frac1n 
\right)}\left|W_{-v+s}-W_s\vphantom{l^1}\right|\right)\times 
\nonumber \\ 
&&\hspace{5.0cm}\times\left(1+\left|\nabla_{W_0}A_{n,s}(W)-\nabla_{ 
W_0}Y_{n,s}(W)\vphantom{\dot{f}}\right|\vphantom{\int}\right)  
\end{eqnarray} 
where $a(n)$ is some  positive constant which may depend on $n\in 
{\Bbb N}$ but is independent of $s\in {\Bbb R}$ and $W\in\{\pi_m 
V:V\in\Omega\}$. The absolute values on the left-hand side and in 
the second factor of the right-hand side have been taken coordinate 
wise. Taking into consideration $\nabla_{W_0}A_{n,0}-\nabla_{W_0} 
Y_{n,0}=0$, Gronwall's inequality shows now that 
\begin{eqnarray}\label{3.53*} 
&&\hspace{-0.5cm}\left|\nabla_{W_0}A_{n,s}(W)-\nabla_{W_0}Y_{n,s} 
(W)\vphantom{\dot{f}}\right|\vphantom{\left(\int_R\right)}\le s 
\cdot\left(a(n)+2\, b(n)\sup_{u\in\left(-\frac1n,s+\frac1n\right)} 
\left|W_u-W_0\vphantom{l^1}\right|\right)\times\nonumber \\ 
&&\hspace{0.5cm}\times\int_{u=0}^s\exp\left\{a(n)+b(n)\sup_{v\in 
\left(-\frac1n,\frac1n\right)}\left|W_{-v+u}-W_u\vphantom{l^1} 
\right|\right\}du 
\end{eqnarray} 
for, without loss of generality, $s\ge 0$. In order to verify that 
all exponential moments of $\sup_{v\in\left(-\frac1n,\frac1n\right 
)}c\cdot\left|W_{-v+u}-W_u\vphantom{l^1}\right|$, $c>0$, are 
finite with respect to the measure $Q^{(m)}_{\sbnu}$ and bounded 
with respect to $u\in [0,s]$ one may read, for example \cite{P05}, 
V.8, the proof of the lemma following Theorem 43 therein, and 
adjust this to the measure $Q^{(m)}_{\sbnu}$. Applying the 
arguments of (\ref{3.47*}) to the composition of (\ref{3.52*}) 
with (\ref{3.53*}) we obtain  
\begin{eqnarray}\label{3.54*}
\sup_{r\in {\Bbb N}}\left|\nabla_{W_0}\dot{A}_{n,s}(\pi_{m,r}W) 
-\nabla_{W_0}\dot{Y}_{n,s}(\pi_{m,r}W)\right|\in L^1(\Omega,Q^{ 
(m)}_{\sbnu};F)\, ,\quad s\in {\Bbb R}, 
\end{eqnarray} 
the absolute value $|\, \cdot\, |$ and the $\sup_{r\in {\Bbb N}}$ 
has been taken coordinate wise. 

Everything stated in (b) is now in (\ref{3.45*}) as well as the 
sentence underneath, (\ref{3.47*}), (\ref{3.50*}), and 
(\ref{3.54*}). 
\qed 
\medskip 

{\bf Proof of the remaining parts (ii)-(iv). } Item (iv) of 
the lemma is a direct consequence of (\ref{3.52*}) together with 
(\ref{3.53*}). As an immediate consequence of (\ref{3.34*}) we get 
\begin{eqnarray*}
\dot{A}_{n,s}-\dot{Y}_{n,s}\in \bigcap_{1\le p<\infty}L^p\left( 
\Omega,Q^{(m,r)}_{\sbnu};F\right)\cap L^p\left(\Omega,Q^{(m)}_{ 
\sbnu};F\right)\, , \quad s\in {\Bbb R},  
\end{eqnarray*} 
and verify that $\dot{A}_n(W)-\dot{Y}_n(W)$ is bounded on finite 
subintervals of ${\Bbb R}$ for all $W\in\{\pi_m V:V\in\Omega\}$. 
In other words, we have (iii) of the present lemma. 

The remainder of the proof focuses on (ii). In the first two steps 
we prove the statements about $\nabla_HX$ and $\nabla_HY$. 
\medskip 

\nid 
{\it Step 1 } In this step we construct path wise estimates on the 
individual directional derivatives of $A_{n,s}$ and $Y_{n,s}$ with 
respect to coordinates $\kappa_i$ and $\lambda_i$. Furthermore, we 
prove belongingness of such derivatives to the stated $L^p$-spaces. 

According to (\ref{3.8}), (\ref{3.9*}), and 
(\ref{3.16*})-(\ref{3.18*}) we have 
\begin{eqnarray}\label{3.541*}
A_{n,s}(W)-Y_{n,s}(W)=\int_0^s\left(A_0\left(g'_n,\gamma_n\right) 
\left(W^{n,u}\right)-Y_0\left(g'_n,\gamma_n\right)\left(W^{n,u} 
\right)\vphantom{\dot{f}}\right)\, du 
\end{eqnarray} 
and 
\begin{eqnarray}\label{3.542*}
Y_{n,s}(W)=A_0\left(g_n,\gamma_n\right)(W^{n,s})\, , \quad s\in 
{\Bbb R}. 
\end{eqnarray} 

In the following calculations we will use the notation $\lambda_j 
=e_j\1(\cdot)$, for all coordinates $j\in F$, here $F={\Bbb R}^{n 
\cdot d}$. In order to simplify the notation, let us write $(A-Y 
)_0\left(g'_n,\gamma_n\right)$ instead of $A_0\left(g'_n,\gamma_n 
\right)-Y_0\left(g'_n,\gamma_n\right)$. According to condition (5) 
(ii) of Section 1 the equation 
\begin{eqnarray}\label{3.55*} 
\begin{array}{rl}
Z_{i,s}\ =&\D\int_0^s\left\langle\left(\nabla_H(A-Y)_0 
\left(g'_n,\gamma_n\right)\vphantom{\dot{f}}\right)\left(W^{n,u} 
\right)\vphantom{\frac{\partial}{\partial\kappa_i}}\, ,\, \left( 
H_i(\cdot +u),Z_{i,u}\right)\right\rangle_{H\to F}\, du \\ 
Z_{i,0}\ =&0\vphantom{\sum_{I}}
\end{array} 
\end{eqnarray} 
is well-defined on $jH$. By the continuity of $u\to W^{n,u}$ 
in the sense of (\ref{3.15*}) and the Picard-Lindel\"of theorem 
there exists a unique solution in some neighborhood $U$ of $s=0$. 
It follows also from the Picard-Lindel\"of theorem that $U$ can 
be chosen so that it is independent of $i\in I(r)$. Recalling 
again condition (5) (ii) of Section 1 we notice that 
\begin{eqnarray*}  
&&\hspace{-.5cm}\left\langle\left(\nabla_H(A-Y)_0\left(g'_n, 
\gamma_n\right)\vphantom{\dot{f}}\right)\left(W^{n,u}\right) 
\vphantom{\frac{\partial}{\partial\kappa_i}}\, ,\, \left (H_i 
(\cdot +u),Z_{i,u}\right)\right\rangle_{H\to F}\vphantom{\sum_{ 
i'\in I(m,r)}} \\ 
&&\hspace{.5cm}=\sum_{i'\in I(r)}\left(\frac{\partial} 
{\partial\kappa_{i'}}\dot{C}_{n,0}\right)\left(W^{n,u}\right) 
\cdot\frac{\partial}{\partial\kappa_i}\left\langle H_{i'},dW^{ 
n,u}\right\rangle_{L^2} \\ 
&&\hspace{1cm}+\sum_{j'}\left(\frac{\partial}{\partial\lambda_{ 
j'}}\dot{C}_{n,0}\right)\left(W^{n,u}\right)\cdot\left\langle 
e_{j'},Z_{i,u}\right\rangle_F\, . 
\end{eqnarray*} 
Now we obtain from the facts that $W^{n,\cdot}$ is the unique 
solution to (\ref{3.15*}) and that condition (5) (ii) provides 
the chain rule the existence of the derivative $\frac{\partial} 
{\partial\kappa_i}\left(A_{n,s}-Y_{n,s}\vphantom{l^1}\right)$ and 
\begin{eqnarray*}  
Z_{i,s}=\frac{\partial}{\partial\kappa_i}\left(A_{n,s}-Y_{n,s} 
\vphantom{l^1}\right) 
\end{eqnarray*} 
for all $s$ in some neighborhood $U$ of $0$. This implies the 
existence of the derivative $\frac{\partial}{\partial\kappa_i} 
Y_{n,s}$ for all $s$ in the same neighborhood $U$ of $0$ and 
\begin{eqnarray}\label{3.56*} 
\frac{\partial}{\partial\kappa_i}Y_{n,s}
&&\hspace{-.5cm}=\left\langle\left(\nabla_HA_0\left(g_n,\gamma_n 
\right)\vphantom{\dot{f}}\right)\left(W^{n,s}\right)\, ,\, \left( 
H_i(\cdot +s),\frac{\partial}{\partial\kappa_i}\left(A_{n,s}-Y_{n 
,s}\vphantom{l^1}\right)\right)\right\rangle_{H\to F}\, .  
\qquad
\end{eqnarray} 
Let $s\in {\cal R}\cap U\equiv {\cal R}(r)\cap U$ and, without 
loss of generality, $s\ge 0$. 
With condition (5) (i) of Section 1 it follows from (\ref{3.55*}) 
that 
\begin{eqnarray}\label{3.57*} 
&&\hspace{-.5cm}\left|\frac{\partial}{\partial\kappa_i}A_{ 
n,s}-\frac{\partial}{\partial\kappa_i}Y_{n,s}\right|=\left| 
\int_0^s\lim_{\ve\to 0}\frac{1}{\ve}\left(\left((A-Y)_0\left( 
g'_n,\gamma_n\right)\vphantom{\dot{f}}\right)\left(W^{n,u}+\ve 
\cdot j\left(H_i(\cdot +u),Z_{i,u}\vphantom{l^1}\right)\right) 
\vphantom{\int}\right.\right.\nonumber \\ 
&&\hspace{1.2cm}\left.\vphantom{\int_0^s}\left.\vphantom{\int} 
-\left((A-Y)_0\left(g'_n,\gamma_n\right)\vphantom{\dot{f}} 
\right)\left(W^{n,u}\right)\right)\, du\right|\nonumber \\ 
&&\hspace{.5cm}\le\int_0^s\left(\left|\left(\nabla_{W_0}A_0 
\left(|g'_n|,\gamma_n\right)\vphantom{\dot{f}}\right)\left( 
W^{n,u}\right)\right|\vphantom{\left(\frac{\partial}{\partial 
\kappa_i}\right)}+\left|\left(\nabla_{W_0}Y_0\left(|g'_n|, 
\gamma_n\right)\vphantom{\dot{f}}\right)\left(W^{n,u}\right) 
\right|\right)\times\nonumber \\ 
&&\hspace{1.5cm}\times\left(\left|\left\langle H_i(\cdot +u), 
h_u\right\rangle_{L^2}\right|+\left|\frac{\partial}{\partial 
\kappa_i}A_{n,u}-\frac{\partial}{\partial\kappa_i}Y_{n,u} 
\right|\right)\, du 
\end{eqnarray} 
and from (\ref{3.56*}) that 
\begin{eqnarray}\label{3.58*}
&&\hspace{-.5cm}\left|\frac{\partial}{\partial\kappa_i}Y_{n 
,s}\right|\le\left|\left(\nabla_{W_0}A_0\left(g_n,\gamma_n 
\right)\vphantom{\dot{f}}\right)\left(W^{n,s}\right)\right| 
\times\nonumber \\ 
&&\hspace{1.5cm}\times\left(\left|\left\langle H_i(\cdot +s) 
,h_s\right\rangle_{L^2}\right|+\left|\frac{\partial}{\partial 
\kappa_i}A_{n,s}-\frac{\partial}{\partial\kappa_i}Y_{n,u} 
\right|\right) 
\end{eqnarray} 
where $h_\cdot\in H$ is defined in condition (5) (i) of Section 1. 


For $g=|g'_n|$ or $g=g_n$ and $W\in jH$ 
we obtain as in (\ref{3.34*}) or (\ref{3.51*}), (\ref{3.52*}) 
\begin{eqnarray}\label{3.59*} 
\left|\left(\nabla_{W_0}A_0\left(g,\gamma_n\right)\vphantom{\dot{f 
}}\right)\left(W^{n,u}\right)\right|&&\hspace{-.5cm}=\left|\left( 
\int_{\left(-\frac1n,\frac1n\right)}\left\langle\int_F\left\langle 
A_{-v}\left(\, \cdot\, +\left(x-W_0\right)\1\vphantom{l^1}\right) 
{\sf e}\, ,\, \right.\right.\right.\right.\nonumber \\ 
&&\hspace{.5cm}\left.\left.\left.\left.\nabla_{W_0}\gamma_n\left( 
x-W_0\vphantom{l^1}\right)\vphantom{Y_{-v}}\vphantom{\dot{f}}\right 
\rangle_{F\to F}\, dx\, ,\, g(v)\vphantom{\int_F}\right\rangle_{F 
\to F}\, dv\vphantom{\int_{\left(-\frac2n,0\right)}}\right)\left( 
W^{n,u}\right)\right|\nonumber \\ 
&&\hspace{-.5cm}\le d(n)\left(1+\sup_{v\in\left(-\frac1n,\frac1n 
\right)}\left|W_{-v+u}-W_u\right|\right) 
\end{eqnarray} 
as well as 
\begin{eqnarray}\label{3.60*} 
\left|\left(\nabla_{W_0}Y_0\left(g,\gamma_n\right)\vphantom{\dot 
{f}}\right)\left(W^{n,u}\right)\right|\le d(n)\, ,  
\end{eqnarray} 
$u\in [0,s]$, where $d(n)$ is some  positive constant which may 
depend on $n\in {\Bbb N}$ but is independent of $u\in [0,s]$ 
and $W\in jH$. 
Relations (\ref{3.57*}), (\ref{3.59*}), and (\ref{3.60*}) yield 
\begin{eqnarray}\label{3.61*} 
\left|\frac{\partial}{\partial\kappa_i}A_{n,s}-\frac{\partial}{ 
\partial\kappa_i}Y_{n,s}\right|&&\hspace{-.5cm}\le\int_0^s d(n)
\left(2+\sup_{v\in\left(-\frac1n,\frac1n\right)}\left|W_{-v+u}- 
W_u\vphantom{l^1}\right|\vphantom{\dot{f}}\right)\times\nonumber 
 \\ 
&&\hspace{.5cm}\times\left(\left|\left\langle H_i(\cdot +u) 
,h_u\right\rangle_{L^2}\right|+\left|\frac{\partial}{\partial 
\kappa_i}A_{n,u}-\frac{\partial}{\partial\kappa_i}Y_{n,u}\right| 
\right)\, du\, . \qquad
\end{eqnarray} 
Recalling (\ref{3.55*}) we verify that $u\to\frac{\partial} 
{\partial\kappa_i}A_{n,u}-\frac{\partial}{\partial\kappa_i}Y_{ 
n,u}$ is a continuous function on $u\in [0,s]$. Applying now 
Gronwall's inequality, relation (\ref{3.61*}) shows that with 
\begin{eqnarray}\label{3.6109*}
c(W):=2d(n)\sup_{u\in\left(-\frac1n,s+\frac1n\right)}\left(1+ 
\left|W_u-W_0\vphantom{l^1}\right|\vphantom{\dot{f}}\right) 
\end{eqnarray} 
we have 
\begin{eqnarray}\label{3.611*} 
\left|\frac{\partial}{\partial\kappa_i}A_{n,s}-\frac{\partial} 
{\partial\kappa_i}Y_{n,s}\right|&&\hspace{-.5cm}\le c(W)\int_0^s 
\left|\left\langle H_i(\cdot +u),h_u\right\rangle_{L^2}\right| 
\, du\nonumber \\ 
&&\hspace{-.0cm}+c(W)\cdot e^{s\cdot c(W)}\int_0^s\int_0^u\left| 
\left\langle H_i(\cdot +v),h_v\right\rangle_{L^2}\right|\, dv\, 
du\, .\qquad
\end{eqnarray} 
Together with (\ref{3.58*}) and (\ref{3.59*}) we obtain from 
(\ref{3.611*})
\begin{eqnarray}\label{3.612*} 
\left|\frac{\partial}{\partial\kappa_i}Y_{n,s}(W)\right|&&\hspace 
{-.5cm}\le c(W)\cdot\left|\left\langle H_i(\cdot +s),h_s\right 
\rangle_{L^2}\right|+(c(W))^2\int_0^s\left|\left\langle H_i(\cdot 
+u),h_u\right\rangle_{L^2}\right|\, du\nonumber \\ 
&&\hspace{-.0cm}+(c(W))^2\cdot e^{s\cdot c(W)}\int_0^s\int_0^u 
\left|\left\langle H_i(\cdot +v),h_v\right\rangle_{L^2}\right|\, 
dv\, du\qquad
\end{eqnarray} 
and 
\begin{eqnarray}\label{3.613*} 
\left|\frac{\partial}{\partial\kappa_i}A_{n,s}(W)\right|&& 
\hspace{-.5cm}\le c(W)\cdot\left|\left\langle H_i(\cdot +s),h_s 
\right\rangle_{L^2}\right|+\left(c(W)+(c(W))^2\right)\int_0^s 
\left|\left\langle H_i(\cdot +u),h_u\right\rangle_{L^2}\right|\, 
du\nonumber \\ 
&&\hspace{-.0cm}+\left(c(W)+(c(W))^2\right)\cdot e^{s\cdot c(W) 
}\int_0^s\int_0^u\left|\left\langle H_i(\cdot +v),h_v\right 
\rangle_{L^2}\right|\, dv\, du\qquad
\end{eqnarray} 
for all $W\in jH$ and all $i\in I(r)$. Recalling the argument after 
(\ref{3.53*}) this gives 
\begin{eqnarray}\label{3.62*} 
\frac{\partial}{\partial\kappa_i}A_{n,s},\frac{\partial}{\partial 
\kappa_i}Y_{n,s}\in\bigcap_{1\le p<\infty} L^p(\Omega,Q^{(m,r)}_{ 
\sbnu};F)\, , \quad i\in I(r). 
\end{eqnarray} 
Similarly, we get 
\begin{eqnarray}\label{3.63*} 
\frac{\partial}{\partial\lambda_j}A_{n,s},\frac{\partial}{\partial 
\lambda_j}Y_{n,s}\in\bigcap_{1\le p<\infty} L^p(\Omega,Q^{(m,r)}_{ 
\sbnu};F) 
\end{eqnarray} 
by examining 
\begin{eqnarray}\label{3.64*} 
\begin{array}{rl}
z_{j,s}\ =&\D\int_0^s\left\langle\left(\nabla_H\left((A-Y)_0\left( 
g'_n,\gamma_n\right)\right)\vphantom{\dot{f}}\right)\left(W^{n,u} 
\right)\vphantom{\frac{\partial}{\partial\kappa_i}}\, ,\, \left(0, 
e_j+z_{j,u}\right)\right\rangle_{H\to F}\, du \\ 
z_{j,0}\ =&0\vphantom{\sum_{I}} 
\end{array} 
\end{eqnarray} 
and 
\begin{eqnarray}\label{3.65*} 
\frac{\partial}{\partial\lambda_j}Y_{n,s}&&\hspace{-.5cm}=\left 
\langle\left(\nabla_HA_0\left(g_n,\gamma_n\right)\vphantom{\dot 
{f}}\right)\left(W^{n,s}\right)\, , \, \left(0,e_j+\frac{\partial} 
{\partial\lambda_j}\left(A_{n,s}-Y_{n,s}\vphantom{l^1}\right) 
\right)\right\rangle_{H\to F}\, .  
\end{eqnarray} 
Indeed, as above, it follows that $z_{j,s}=\frac{\partial}{\partial 
\lambda_j}\left(A_{n,s}-Y_{n,s}\vphantom{l^1}\right)$ is the unique 
solution to (\ref{3.64*}) for all $s$ in some neighborhood of $0$. 
\medskip 

\nid 
{\it Step 2 } In this step we represent $\nabla_HX_{n,s}$ and 
$\nabla_HY_{n,s}$ by means of the directional derivatives and prove 
belongingness to the stated $L^p$-spaces. As a byproduct of Steps 1 
and 2 we verify that $A_{n,s}$ and $Y_{n,s}$ belong to the stated 
$D_{q,1}$-spaces. 

Let 
\begin{eqnarray*} 
k&&\hspace{-.5cm}=\sum_{i\in I(r)}\langle k,(H_i,0)\rangle_H\cdot 
(H_i,0)+\sum_j\langle k,(0,e_j)\rangle_H\cdot (0,e_j)\nonumber \\ 
&&\hspace{-.5cm}=:k^{(1)}+k^{(2)}\equiv \left(k^{(1)},k^{(2)}\right) 
\in H
\end{eqnarray*} 
and 
\begin{eqnarray*} 
\kappa^{(1)}:=\int_0^\cdot k^{(1)}(v)\, dv\, , \quad\kappa^{(2)} 
:=k^{(2)}\1(\cdot)\, , \quad\kappa:=\kappa^{(1)}+\kappa^{(2)}\, . 
\end{eqnarray*} 
On the one hand, keeping condition (5) (ii) in mind, as in 
(\ref{3.55*}) and (\ref{3.56*}) it turns out that the derivative 
$\frac{\partial}{\partial\kappa}\left(A_{n,s}-Y_{n,s}\right)$ 
exists on $jH$. 
In fact,  
\begin{eqnarray*}
Z_s:=\frac{\partial}{\partial\kappa}\left(A_{n,s}-Y_{n,s}\right) 
\end{eqnarray*} 
is the unique solution to 
\begin{eqnarray}\label{3.66*} 
\begin{array}{rl}
Z_s\ =&\D\int_0^s\left\langle\left(\nabla_H\left((A-Y)_0\left(g'_n, 
\gamma_n\right)\right)\vphantom{\dot{f}}\right)\left(W^{n,u}\right) 
\left(k^{(1)}(\cdot+u),k^{(2)}+Z_u\right)\right\rangle_{H\to F}\, du 
 \\ 
Z_0\ =&0 
\end{array} 
\end{eqnarray} 
for all $s$ in some neighborhood of $0$ and  
\begin{eqnarray}\label{3.67*} 
\frac{\partial}{\partial\kappa}Y_{n,s}&&\hspace{-.5cm}=\left\langle 
\left(\nabla_HA_0\left(g_n,\gamma_n\right)\vphantom{\dot{f}}\right) 
\left(W^{n,s}\right)\, ,\, \vphantom{\frac{\partial}{\partial\kappa 
}}\right.\nonumber \\ 
&&\hspace{2cm}\left.\left(k^{(1)}(\cdot+s),k^{(2)}+\frac{\partial} 
{\partial\kappa}\left(A_{n,s}-Y_{n,s}\vphantom{l^1}\right)\right) 
\right\rangle_{H\to F}\, . 
\end{eqnarray} 
On the other hand, by (\ref{3.612*}) we get 
\begin{eqnarray*} 
\left(\frac{\partial}{\partial\kappa_i}Y_{n,s}(W)\right)^2&& 
\hspace{-.5cm}\le 3(c(W))^2\cdot\left\langle H_i(\cdot +s),h_s 
\right\rangle_{L^2}^2+3(c(W))^4rt\int_0^s\left\langle H_i(\cdot 
+u),h_u\right\rangle_{L^2}^2\, du\nonumber \\ 
&&\hspace{-.0cm}+3(c(W))^4\cdot e^{2s\cdot c(W)}(rt)^2\int_0^s 
\int_0^u\left\langle H_i(\cdot +v),h_v\right\rangle_{L^2}^2\, 
dv\, du
\end{eqnarray*} 
where we recall the definition of $c(W)$ in (\ref{3.6109*}). 
With $h_s$ extending to a function of type ${\Bbb R}\to F$ by 
setting $h_s=0$ on ${\Bbb R}\setminus {\cal R}$ and using the 
abbreviation $c\equiv c(W)$ the latter implies 
\begin{eqnarray}\label{3.671*} 
\sum_{i\in I(r)}\left(\frac{\partial}{\partial\kappa_i}Y_{n,s} 
(W)\right)^2&&\hspace{-.5cm}\le 3(c(W))^2\cdot\|h_s(\cdot -s) 
\|_{L^2}^2+3(c(W))^4rt\int_0^s\|h_u(\cdot -u)\|_{L^2}^2\, du 
\nonumber \\ 
&&\hspace{-.0cm}+3(c(W))^4\cdot e^{2s\cdot c(W)}(rt)^2\int_0^s 
\int_0^u\|h_v(\cdot -v)\|_{L^2}^2\, dv\, du\nonumber \\ 
&&\hspace{-.5cm}\le\left(3c^2+3c^4\left((rt)^2+(rt)^4e^{2rt 
\cdot c}\right)\vphantom{\int}\right)\cdot\sup_{s\in {\cal R}} 
\|h_s\|^2\, . 
\end{eqnarray} 
Similarly, by (\ref{3.613*}) 
\begin{eqnarray}\label{3.672*} 
\sum_{i\in I(r)}\left(\frac{\partial}{\partial\kappa_i}A_{n,s} 
(W)\right)^2\le\left(3c^2+3(c+c^2)^2\left((rt)^2+(rt)^4e^{2rt 
\cdot c}\right)\vphantom{\int}\right)\cdot\sup_{s\in {\cal R}} 
\|h_s\|^2\, . 
\end{eqnarray} 
In other words, we obtain by linear combination of (\ref{3.55*}) 
and (\ref{3.56*}) as well as (\ref{3.64*}) and (\ref{3.65*}) and 
from the uniqueness of (\ref{3.66*}) and (\ref{3.67*})
\begin{eqnarray}\label{3.68*} 
\frac{\partial}{\partial\kappa}A_{n,s}&&\hspace{-.5cm}=\sum_{i 
\in I(r)}\langle k,(H_i,0)\rangle_H\cdot\frac{\partial}{\partial 
\kappa_i}A_{n,s}+\sum_j\langle k,(0,e_j)\rangle_H\cdot\frac{ 
\partial}{\partial\lambda_j}A_{n,s} 
\end{eqnarray} 
and 
\begin{eqnarray}\label{3.69*} 
\frac{\partial}{\partial\kappa}Y_{n,s}&&\hspace{-.5cm}=\sum_{i 
\in I(r)}\langle k,(H_i,0)\rangle_H\cdot\frac{\partial}{\partial 
\kappa_i}Y_{n,s}+\sum_j\langle k,(0,e_j)\rangle_H\cdot\frac{ 
\partial}{\partial\lambda_j}Y_{n,s} 
\end{eqnarray} 
in some neighborhood of $s=0$ on $jH$. Now we replace in 
(\ref{3.55*}), (\ref{3.64*}), and (\ref{3.66*}) the argument $W^{ 
n,u}$ by $W^{n,u+s_0}$ as well as in (\ref{3.56*}), (\ref{3.65*}), 
and (\ref{3.67*}) $W^{n,s}$ by $W^{n,s+s_0}$. By adjusting the 
initial conditions in (\ref{3.55*}), (\ref{3.64*}), and 
(\ref{3.66*}) we verify (\ref{3.68*}) and (\ref{3.69*}) for all 
$s=s_0\in {\Bbb R}$. Together with (\ref{3.671*}), (\ref{3.672*}) 
and $c\equiv c(W)$ in (\ref{3.6109*}) this says 
\begin{eqnarray}\label{3.70*} 
\nabla_HA_{n,s}&&\hspace{-.5cm}=\sum_{i\in I(r)}\frac{\partial} 
{\partial\kappa_i}A_{n,s}\cdot (H_i,0)+\sum_j\frac{\partial} 
{\partial\lambda_j}A_{n,s}\cdot (0,e_j)\, , 
\end{eqnarray} 
\begin{eqnarray}\label{3.71*} 
\nabla_HY_{n,s}&&\hspace{-.5cm}=\sum_{i\in I(r)}\frac{\partial} 
{\partial\kappa_i}Y_{n,s}\cdot (H_i,0)+\sum_j\frac{\partial} 
{\partial\lambda_j}Y_{n,s}\cdot (0,e_j)\, ,\quad s\in {\Bbb R},  
\end{eqnarray} 
and 
\begin{eqnarray*}
\nabla_HX_{n,s}\in\bigcap_{1\le p<\infty} L^p(\Omega,Q^{(m,r)}_{ 
\sbnu};H)\quad\mbox{\rm and}\quad\nabla_H Y_{n,s}\in\bigcap_{1\le 
p<\infty} L^p(\Omega ,Q^{(m,r)}_{\sbnu};H)\, . 
\end{eqnarray*} 

We recall (\ref{3.541*}) and (\ref{3.542*}), as well as the 
paragraph below (\ref{3.33*}) in order to verify $A_{n,s},Y_{n,s} 
\in\bigcap_{1\le p<\infty} L^p(\Omega,Q^{(m,r)}_{\sbnu};F)$. It 
follows now from (\ref{3.62*}), (\ref{3.63*}) and 
(\ref{3.69*}) and Theorem 1.11 in \cite{KO84} that $A_{n,s},Y_{n, 
s}\in D_{q,1}(Q^{(m,r)}_{\sbnu})$ for all $q$ with $1/q+1/Q<1$. 
\medskip 

\nid 
{\it Step 3 } It remains to verify property (iii) of Proposition 
\ref{Proposition2.7} for $X_n$ as well as $Y_n$. Let us take a 
look at (\ref{3.70*}) and (\ref{3.71*}). It turns out that 
\begin{eqnarray*} 
D_GA_{n,s}(W)&&\hspace{-.5cm}=\sum_{i\in I(r)}\frac{\partial} 
{\partial\kappa_i}A_{n,s}(W)\cdot (H_i,0)+\sum_j\frac{\partial} 
{\partial\lambda_j}A_{n,s}(W)\cdot (0,e_j)
\end{eqnarray*} 
is the G\^ateaux derivative in the space $jH\equiv jH^{(r)}$ of 
the function $A_{n,s}$, $n\in {\Bbb N}$, $s\in {\Bbb R}$, at the 
point $W\in jH$. Similarly we get the G\^ateaux derivative of 
$Y_{n,s}$. 

We observe that all $\frac{\partial}{\partial\kappa_i}A_{n,s}$, 
$\frac{\partial}{\partial\kappa_i}Y_{n,s}$, $i\in I(r)$, and all 
$\frac{\partial}{\partial\lambda_j}A_{n,s}$, $\frac{\partial} 
{\partial\lambda_j}Y_{n,s}$, $s\in {\cal R}$, are continuous on 
$jH$, i.e., with respect to the coordinates $\kappa_i$, $i\in 
I(r)$, and $\lambda_j$. In fact, we just use (\ref{3.55*}), 
(\ref{3.56*}) as well as (\ref{3.64*}), (\ref{3.65*}) together 
with conditions (4) (i) and (5) (i) of Section 1 where we 
recall from Subsection 1.2 the definitions of the jump times, 
the relation $\bnu(G(W))=0$, $W\in\Omega$, and the definitions 
of $A(\, \cdot\, ,\gamma_n)$ as well as $Y(\, \cdot\, ,\gamma_n 
)$. Furthermore, we take into consideration the continuity 
$jH^{(r)}\ni W\to W^{n,u}$, $u\in[0,s]$. The latter we derive 
from (\ref{3.9*}), (\ref{3.10*}), and (\ref{3.14*}) and again 
condition (4) (i) of Section 1. 

With (\ref{3.671*}), (\ref{3.672*}), and (\ref{3.6109*}) it 
follows now that $D_GA_{n,s}$ and $D_GY_{n,s}$ are continuous 
in the space $jH$. We conclude that $D_GA_{n,s}(W)$ and $D_G 
Y_{n,s}(W)$ are at the same time the Fr\'echet derivatives 
$D_FA_{n,s}(W)$ and $D_FY_{n,s}(W)$ of $A_{n,s}$ and $Y_{n,s 
}$ in the space $jH$. 
\medskip

Next we show differentiability of $(D_FA_{n,s}(W))_1$ and $(D_F 
Y_{n,s}(W))_1$ for $W\in jH$ and $s\in {\cal R}$. According to 
condition (5) (iii), $\nabla_G(A_0-Y_0)(g_n',\gamma_n)$ is 
continuously differentiable which we symbolize by $(\dot{\nabla 
}_G)_\cdot(A_0-Y_0)(g_n',\gamma_n)$. For fixed $r\in {\cal R}$ 
and variable $s\in {\cal R}$, we consider the first order ODE 
\begin{eqnarray*} 
\begin{array}{rl}
\D\zeta_s(r)\ =&\D\int_0^s\left(\vphantom{\int}\left((\dot{ 
\nabla}_G)_{r-u}(A_0-Y_0)(g_n',\gamma_n)\vphantom{\dot{f}}\right) 
(W^{n,u})\right. \\ 
 & \hspace{1cm}\D\left.+\left\langle\left(\nabla_{W_0}(A_0-Y_0) 
(g_n',\gamma_n)\vphantom{\dot{f}}\right)(W^{n,u})\, ,\, \zeta_u 
(r)\vphantom{\dot{f}}\right\rangle_F\vphantom{\int}\right)\, du 
 \\ 
\D\zeta_0(r)\ =&0\vphantom{\sum_{I}}\, . 
\end{array} 
\end{eqnarray*} 
Using the arguments for (\ref{3.55*}), this equation has a unique 
solution. Furthermore, 
\begin{eqnarray}\label{3.72*} 
Z_s(r):=\int_0^r\zeta_s(p)\, dp
\end{eqnarray} 
is the unique solution to 
\begin{eqnarray}\label{3.73*} 
\begin{array}{rl}
\D Z_s(r)\ =&\D\int_0^s\left(\vphantom{\int}\left(\vphantom{\dot{f}} 
(\nabla_G)_{r-u}(A_0-Y_0)(g_n',\gamma_n)\right)(W^{n,u})\right.\\ 
 & \hspace{1cm}\D\left.+\left\langle\left(\vphantom{\dot{f}}\nabla_{ 
W_0}(A_0-Y_0)(g_n',\gamma_n)\right)(W^{n,u})\, ,\, Z_u(r)\vphantom 
{\dot{f}}\right\rangle_F\vphantom{\int}\right)\, du \\ 
\D Z_0\ =&0\vphantom{\sum_{I}}\, . 
\end{array} 
\end{eqnarray} 
For $W\in jH$ we have 
\begin{eqnarray*} 
&&\hspace{-.5cm}\vphantom{\int}\left(\vphantom{\dot{f}}\left( 
\nabla_G\right)_{\cdot -u}(A_0-Y_0)(g_n',\gamma_n)\right)(W^{n,u}) 
 \\ 
&&\hspace{1cm}+\left\langle\left(\vphantom{\dot{f}}\nabla_{W_0}(A_0 
-Y_0)(g_n',\gamma_n)\right)(W^{n,u})\, ,\, D_{F,1}(A_{n,u}-Y_{n,u}) 
(W)\vphantom{\dot{f}}\right\rangle_F\vphantom{\int} \\ 
&&\hspace{.5cm}=\sum_{i\in I(r)}H_i\left\langle\left(\vphantom{\dot 
{f}}\nabla_G (A_0-Y_0)(g_n',\gamma_n)\right)(W^{n,u})\, ,\, H_i 
(\cdot +u)\vphantom{\dot{f}}\right\rangle_{L^2}\vphantom{\int} \\ 
&&\hspace{1cm}+\left\langle\left(\vphantom{\dot{f}}\nabla_{W_0}(A_0 
-Y_0)(g_n',\gamma_n)\right)(W^{n,u})\, ,\, D_{F,1}(A_{n,u}-Y_{n,u}) 
(W)\vphantom{\dot{f}}\right\rangle_F\vphantom{\int} \\ 
&&\hspace{.5cm}=\frac{d}{du}D_{F,1}(A_{n,u}-Y_{n,u})(W) 
\end{eqnarray*} 
where, for the last equality sign we have used (\ref{3.55*}) and 
condition (5) (ii) of Section 1. This and $D_{F,1}(A_{n,0}-Y_{n,}) 
(W)=0$ on ${\cal R}$, (\ref{3.73*}) imply $Z_s=D_{F,1}(A_{n,s}-Y_{ 
n,s})(W)$. Relation (\ref{3.72*}) says now that $D_{F,1}A_{n,s}(W) 
-D_{F,1}Y_{n,s}(W)$ is differentiable. 

As a consequence of (\ref{3.56*}) and condition (5) (ii) of Section 
1, for $W\in jH$ it holds that 
\begin{eqnarray*} 
D_{F,1}Y_{n,s}&&\hspace{-.5cm}=\sum_{i\in I(r)}H_i\left\langle 
\left(\vphantom{\dot{f}}\nabla_G A_0(g_n,\gamma_n)\right)(W^{n,s}) 
\, ,\, H_i(\cdot +s)\vphantom{\dot{f}}\right\rangle_{L^2}\vphantom 
{\int} \\ 
&&\hspace{0cm}+\left\langle\left(\vphantom{\dot{f}}\nabla_{W_0}A_0 
(g_n,\gamma_n)\right)(W^{n,s})\, ,\, D_{F,1}(A_{n,s}-Y_{n,s})(W) 
\vphantom{\dot{f}}\right\rangle_F\vphantom{\int} \\ 
&&\hspace{-.5cm}=\left(\vphantom{\dot{f}}\left(\nabla_G\right)_{ 
\cdot -s}A_0(g_n,\gamma_n)\right)(W^{n,s})\vphantom{\int} \\ 
&&\hspace{0cm}+\left\langle\left(\vphantom{\dot{f}}\nabla_{W_0}A_0 
(g_n,\gamma_n)\right)(W^{n,s})\, ,\, D_{F,1}(A_{n,s}-Y_{n,s})(W) 
\vphantom{\dot{f}}\right\rangle_F\vphantom{\int}\, .  
\end{eqnarray*} 
Differentiability of $D_{F,1}Y_{n,s}(W)$, $W\in jH$, follows 
now from condition (5) (iii) of Section 1 and differentiability 
of $D_{F,1}A_{n,s}(W)-D_{F,1}Y_{n,s}(W)$. Thus, $D_{F,1}A_{n,s} 
(W)$ and $D_{F,1}Y_{n,s}(W)$ are differentiable. As an immediate 
consequence we get (\ref{2.30}) for $X_n$ and $Y_n$ instead of 
$X$. 
\medskip 

For the second order Fr\'echet derivatives we repeat the whole 
part from (\ref{3.55*}) till the existence of first order 
Fr\'echet derivatives $D_FA_{n,s}(W)$ and $D_FY_{n,s}(W)$ of 
$A_{n,s}$ and $Y_{n,s}$ in the space $jH$. This gives direction 
how to verify the existence of the corresponding second order 
Fr\'echet derivatives which we sketch below. 

The difference is now that we deal with the mixed second order 
directional derivatives instead of the first order ones. For 
example, for the derivatives $\frac{\partial^2}{\partial\kappa_i 
\partial\kappa_{i'}} A_{n,s}$ and $\frac{\partial^2}{\partial 
\kappa_i\partial\kappa_{i'}} Y_{n,s}$ we modify relations 
(\ref{3.55*}) and (\ref{3.56*}) to 
\begin{eqnarray*} 
\begin{array}{rl}
Z_{i,i',s}=&\D\int_0^s\left\langle\left(\nabla^2_H(A-Y)_0 
\left(g'_n,\gamma_n\right)\vphantom{\dot{f}}\right)\left(W^{ 
n,u}\right)\vphantom{\frac{\partial}{\partial\kappa_i}}\, ,\, 
\right.\vphantom{\left(\left(\int\right)^1_1\right)} \\ 
 & \hspace{.0cm}\D\left.\vphantom{\D\frac{\partial}{\partial 
\kappa_i}}\left(\left(H_i(\cdot +u),\frac{\partial}{\partial 
\kappa_i}(A_{n,u} -Y_{n,u})\right),\left(H_{i'}(\cdot +u), 
\frac{\partial}{\partial\kappa_{i'}}(A_{n,u} -Y_{n,u})\right) 
\vphantom{l^1}\right)\right\rangle_{H\otimes H\to F}\, du 
\vphantom{\left(\left(\int\right)^1_1\right)} \\ 
&+\D\int_0^s\left\langle\left(\nabla_H(A-Y)_0\left(g'_n, 
\gamma_n\right)\vphantom{\dot{f}}\right)\left(W^{n,u}\right) 
\vphantom{\frac{\partial}{\partial\kappa_i}}\, ,\, \left(0, 
Z_{i,i',u}\vphantom{l^1}\right)\right\rangle_{H\to F}\, du 
\vphantom{\left(\left(\int\right)^1_1\right)} \\ 
Z_{i,i',0}\ =&0\vphantom{\left(\left(\int\right)^1_1\right)}
\end{array} 
\end{eqnarray*} 
with unique solution $Z_{i,i',s}=\frac{\partial^2}{\partial 
\kappa_i\partial\kappa_{i'}}(A_{n,s} -Y_{n,s})$ and 
\begin{eqnarray*} 
\frac{\partial^2}{\partial\kappa_i\partial\kappa_{i'}}Y_{n,s}
&&\hspace{-.5cm}=\left\langle\left(\nabla^2_HA_0\left(g_n, 
\gamma_n\right)\vphantom{\dot{f}}\right)\left(W^{n,s}\right) 
\, ,\, \vphantom{\frac{\partial^2}{\partial\kappa_i\partial 
\kappa_{i'}}}\right.\\ 
&&\hspace{.0cm}\left.\left(\left(H_i(\cdot +s),\frac{\partial} 
{\partial\kappa_i}(A_{n,s}-Y_{n,s})\right),\left(H_{i'}(\cdot 
+s),\frac{\partial}{\partial\kappa_{i'}}(A_{n,s} -Y_{n,s}) 
\right)\right)\right\rangle_{H\otimes H\to F} \\ 
&&\hspace{.0cm}+\int_0^s\left\langle\left(\nabla_HA_0\left(g_n 
,\gamma_n\right)\vphantom{\dot{f}}\right)\left(W^{n,u}\right) 
\vphantom{\frac{\partial}{\partial\kappa_i}}\, ,\, \left(0, 
\frac{\partial^2}{\partial\kappa_i\partial\kappa_{i'}}(A_{n,s} 
-Y_{n,s})\right)\right\rangle_{H\to F}\, .  
\end{eqnarray*} 
For the well-finiteness of these relations we use condition (5) 
(ii) of Section 1. Moreover, for the convergence in $H\otimes 
H$ of sums that correspond to (\ref{3.671*}) and (\ref{3.672*}), 
we pay particular attention to the crucial step for this, namely 
the counterpart to (\ref{3.57*}). In order to handle the second 
order partial derivatives which appear now we apply 
\begin{eqnarray*}
\left\langle\nabla^2_H\xi,(x,y)\right\rangle_{H\otimes H}=\lim_{ 
\ve,\delta\to 0}\frac{\xi(\cdot +\ve x+\delta y)+\xi(\cdot -\ve x 
-\delta y)-\xi(\cdot +\ve x-\delta y)-\xi(\cdot -\ve x+\delta y)} 
{4\ve\delta}
\end{eqnarray*} 
if $x\perp y$, $\left\langle\nabla^2_H\xi,(x,y)\right\rangle_{H 
\otimes H}=\lim_{\ve\to 0}\ve^{-2}\left(\xi(\cdot +\ve x)+\xi( 
\cdot -\ve x)-2\xi\right)$ if $x=y$, as well as linear 
combinations of these two formulas for general $x,y$, and we use 
condition (5) (i) of Section 1. Recalling Remark (2) of Section 
2, we get the remaining part of condition (iii) of Proposition 
\ref{Proposition2.7} for $X_n$ and $Y_n$ instead of $X$. 
\qed
\bigskip 

{\bf Parts (a) and (c) of Lemma \ref{Lemma3.4} if parallel 
trajectories in $\Omega$ do not necessarily generate identical 
jump times for $X$. } Let $W\in\{\pi_m V:V\in\Omega\}$ with 
$W_0\not\in G(W-W_0\1)$, as in the above proofs of the parts (a) 
and (c) of Lemma \ref{Lemma3.4} under the restriction of Remark 
(5) of this section. We concentrate on part (a). Part (c) can be 
treated in a similar way. Let $\sigma(s,x)\equiv\sigma(s,x;W,x_0 
)$, $s\in {\Bbb R}$, $x\in F$, be the time re-parametrization 
with respect to $W-W_0\1$ and $W_0$ introduced in Remark (1) of 
Section 1. 
\medskip 

\nid 
{\it Step 1 } We define a certain time substitution based on 
$A_{n,\cdot}$, $Y_{n,\cdot}$, and $\sigma$. According to 
condition (2) (i) of Section 1 we may assume that $s=0$ is not 
a jump time. We set $\eta(0):=0$, $\bar{A}_{n,s}(W):=A_0\left(( 
W^{n,s})^{-s}\right)$, and $\bar{Y}_{n,s}(W):=Y_0\left((W^{n,s} 
)^{-s}\right)$, $s\in {\Bbb R}$. We differentiate (\ref{3.23*}) 
with respect to $s$ and multiply this coordinate wise with 
\begin{eqnarray}\label{1}
d\eta(s)\equiv d\eta(s;n):=2ds-d\sigma(s,W_0+\bar{A}_{n,s}-\bar 
{Y}_{n,s}) 
\end{eqnarray} 
where $\eta\equiv\eta(\cdot;n):{\Bbb R}\to {\Bbb R}$ with $\eta 
(0)=0$ is supposed to be right continuous. Integrating the result 
we get 
\begin{eqnarray*}
&&\hspace{-.5cm}\int_{r=0}^s\left(\dot{A}_{n,r}(W)-\dot{Y}_{n,r} 
(W)\vphantom{l^1}\right)\, d\eta(r) \\ 
&&\hspace{.5cm}=\int_{r=0}^s\dot{C}_{n,0}\left(W^{n,r}\right)\,  
d\eta(r) \\ 
&&\hspace{.5cm}=\int_{r=0}^s\int_F\int\left\langle\vphantom{\dot 
{f}}\left\langle (A-Y)'_{-v}\left(W^{n,r}+y\1\vphantom{l^1} 
\right)\, ,\, g_n(v)\vphantom{l^1}\right\rangle_{F\to F}\, ,\, 
\vphantom{\dot{f}}\gamma_n(y)\right\rangle_{F\to F}\, dv\, dy\, 
d\eta(r) \\  
&&\hspace{1.0cm}+\sum_{k\in{\Bbb Z}\setminus\{0\}}\int_{r=0}^s 
\int_F\left\langle\vphantom{\dot{f}}\left\langle\Delta(A-Y)_{ 
\tau_k}\left(W^{n,r}+y\1\vphantom{l^1}\right)\, ,\right.\right. 
 \\ 
&&\hspace{5.5cm}\vphantom{\int_F}\left.\vphantom{\dot{f}}\left. 
g_n\left(\tau_k\circ u\left(W^{n,r}+y\1\vphantom{l^1}\right) 
\right)\right\rangle_{F\to F}\gamma_n(y)\right\rangle_{F\to F}\, 
dy\, d\eta(r)\, , 
\end{eqnarray*} 
$s\in {\Bbb R}$, where we note that the integrands with respect 
to $r$ of the above integrals are continuous by mollifying with 
$\gamma_n$. 

Next, we compare the definitions of both mollifier functions 
$g_n$ and $\gamma_n$ in Subsection 1.2 preparing condition (5) 
therein. We give now a motivation for the order $n^3$ in the 
definition of $\gamma_n$. We recall that for $y<1/n^3$ we have 
coordinate wise $|g_n(v+y)-g_n(v)|\le n^2\max |g_1'|\cdot |y| 
\stack{n\to\infty}{\lra}0$. This implies the existence of a 
sequence of real measurable functions ${\Bbb R}\ni s\to c_n(s) 
\equiv c_n(s;W)$, $n\in {\Bbb N}$, with $c_n\stack{n\to\infty} 
{\lra}0$ uniformly on every finite subinterval of ${\Bbb R}$  
such that 
\begin{eqnarray}\label{2}
&&\hspace{-.5cm}\int_{r=0}^s\left(\dot{A}_{n,r}(W)-\dot{Y}_{n,r} 
(W)\vphantom{l^1}\right)\, d\eta(r)\nonumber \\ 
&&\hspace{.5cm}=\int_{r=0}^s\int\left\langle (A-Y)'_{-v}(W^{n,r} 
)\, ,\, g_n(v)\vphantom{l^1}\right\rangle_{F\to F}\, dv \, d\eta 
(r)\nonumber \\ 
&&\hspace{1.0cm}+\sum_{k\in{\Bbb Z}\setminus\{0\}}\int_{r=0}^s 
\left\langle\Delta(A-Y)_{\tau_k}(W^{n,r})\, ,\, g_n\left(\tau_k 
\circ u(W^{n,r})\right)\vphantom{l^1}\right\rangle_{F\to F}\, d 
\eta(r)\nonumber \\ 
&&\hspace{1.0cm}+\int_{r=0}^s c_n(r)\, d\eta(r)\cdot {\sf e}\, , 
\quad s\in {\Bbb R}. 
\end{eqnarray} 
Here, the integrand with respect to $r$ of the second 
integral of the right-hand side is no longer necessarily 
continuous. However, by Definition \ref{Definition1.7} and the 
continuity of $r\to W^{n,r}$ in the sense of (\ref{3.15*}) the 
functions $r\to\Delta (A-Y)_{\tau_k}(W^{n,r})$ and $r\to g_n 
\left(\tau_k\circ u(W^{n,r})\right)$ are piecewise continuous. 
By the asymptotic behavior of the mollifier function $\gamma_n 
$ as $n\to\infty$ the increase of the second integral of the 
right-hand side at jump times is now appropriately defined as 
follows. Form coordinate wise the product of the increase of 
$\eta$ times the average of the left and the right limit of 
the integrand at jump points. Note that these limits exist by 
Definition \ref{Definition1.7} (j). 
\medskip 

For the moment, fix $n\in {\Bbb N}$ and, without loss of 
generality, $s>0$. By (\ref{1}) and property (i) of Remark 
(1) in Section 1, the function $\eta(r)$ is piecewise 
continuously differentiable in $r\in [0,s)$. Thus the function 
$\eta$ has on $[0,s)$ countably many disjoint real intervals 
$(a_l^+,b_l^+)\subseteq [0,s)$ of strict continuous increase 
resp. strict continuous decrease $(a_l^-,b_l^-)\subseteq [0,s) 
$, $l\in {\Bbb N}$. On these intervals $\eta(r)$ is  
continuously differentiable. 

In addition $\eta$ has on $[0,s)$ countably many positive jumps, 
say at $J_i^+\in [0,s)$ and countably many negative jumps, say 
at $J_i^-\in [0,s)$. By property (i) of Remark (1) in Section 
1, $\eta$ has left and right limits at these jump times. 

For each of the intervals of strict continuous increase or 
decrease, the function $\eta$ has a well-defined {\it local} 
inverse $\eta^{-1}_{l,+}$ resp. $\eta^{-1}_{l,-}$. In other 
words, there exist countably many not necessarily disjoint 
intervals $(\alpha_l^+,\beta_l^+)$ resp. $(\alpha_l^-, 
\beta_l^-)$, $l\in {\Bbb N}$, such that $\eta^{-1}_{l,+}$ maps 
$(\alpha_l^+,\beta_l^+)$ onto $(a_l^+,b_l^+)$ resp. $\eta^{-1 
}_{l,-}$ maps $(\alpha_l^-,\beta_l^-)$ onto $(a_l^-,b_l^-)$. 
Furthermore, there are real bounded open intervals $I_i^+$ 
resp. $I_i^-$ such that $\eta(J_i^+)-=\inf I_i^+$ and $\eta( 
J_i^+)+=\sup I_i^+$ resp. $\eta(J_i^-)-=\sup I_i^-$ and $\eta 
(J_i^-)+=\inf I_i^-$. We introduce corresponding {\it local} 
inverse maps $\eta_{J_i^+}^{-1}$ and $\eta_{J_i^-}^{-1}$ by 
$J_i^+=\eta_{J_i^+}^{-1}(r)$, $r\in I_i^+$, and $J_i^-=\eta_{ 
J_i^-}^{-1}(r)$, $r\in I_i^+$. To simplify the notation, we 
summarize the intervals and points 
\begin{eqnarray*}
{\cal I}:=\left\{\mbox{all}\ (a_l^+,b_l^+),\, (a_l^-,b_l^-)\ 
\mbox{and all}\ J_i^+, J_i^-\right\} 
\end{eqnarray*} 
and denote for each $I\in {\cal I}$ the above local inverse map 
by $\bar{\eta}^{-1}\equiv\bar{\eta}^{-1}_I$. If $I\in {\cal I}$ 
is an interval, the inverse $\bar{\eta}$ to $\bar{\eta}^{-1}_I$ 
is the restriction of $\eta$ to $I$. If $I\in {\cal I}$ is a 
point, let the {\it formal} inverse $\bar{\eta}$ to $\bar{\eta 
}^{-1}_I$ be defined by $\bar{\eta}(J_i^+):=\sup I_i^+=\eta( 
J_i^+)+$ and $\bar{\eta}(J_i^-):=\inf I_i^-=\eta(J_i^-)+$. 
\medskip 

\nid
{\it Step 2 } We discover some of properties of the time 
substitution $\eta$. Let $r\in [0,s)$ and $k\in {\Bbb Z} 
\setminus\{0\}$. We take into consideration that by definition 
$\eta(r)=2r-\sigma(r,W_0+\bar{A}_{n,r}-\bar{Y}_{n,r})$ and 
\begin{eqnarray}\label{25}
&&\hspace{-.5cm}\tau_k\circ u(W^{n,r})+r=\tau_{k'}\circ u\left( 
\vphantom{l^1}W+(\bar{A}_{n,r}-\bar{Y}_{n,r})\1\right)\vphantom 
{\dot{f}}\nonumber \\ 
&&\hspace{.5cm}=\sigma\left(\tau_{k'}\circ u(W),W_0+\bar{A}_{n, 
r}-\bar{Y}_{n,r}\vphantom{l^1}\right)\, ,\quad r\in[0,s), 
\vphantom{\dot{f}}
\end{eqnarray} 
for some $k'\equiv k'(k;W,r)\in {\Bbb Z}\setminus\{0\}$ which 
we specify below, cf. also Remark (1) (v) of Section 1. 
According to Definition \ref{Definition1.7} (jjj), for all $k 
\in {\Bbb Z}\setminus\{0\}$ such that $\tau_{k'}\circ u(W)\in 
[0,s)$ it holds that 
\begin{eqnarray*}
&&\hspace{-.5cm}\tau_k\circ u(W^{n,r})+r=\sigma\left(\tau_{k'} 
\circ u(W),W_0+\bar{A}_{n,r}-\bar{Y}_{n,r}\vphantom{l^1}\right 
)\vphantom{\frac1n} \\ 
&&\hspace{.5cm}=\sigma\left(r,W_0+\bar{A}_{n,r}-\bar{Y}_{n,r} 
\vphantom{l^1}\right)\vphantom{\frac1n} \\ 
&&\hspace{1.0cm}+\sigma\left(\tau_{k'}\circ u(W),W_0+\bar{A}_{ 
n,r}-\bar{Y}_{n,r}\vphantom{l^1}\right)-\sigma\left(r,W_0+\bar 
{A}_{n,r}-\bar{Y}_{n,r}\vphantom{l^1}\right)\vphantom{\frac1n} 
 \\ 
&&\hspace{.5cm}=\sigma\left(r,W_0+\bar{A}_{n,r}-\bar{Y}_{n,r} 
\vphantom {l^1}\right)+\tau_{k'}\circ u(W)-r\vphantom{\frac1n} 
 \\ 
&&\hspace{.5cm}=2r-\eta(r)+\tau_{k'}\circ u(W)-r=r-\eta(r)+ 
\tau_{k'}\circ u(W)\quad\mbox{\rm if}\quad\left|\tau_k\circ 
u(W^{n,r})\right|\le\frac1n 
\end{eqnarray*} 
for sufficiently large $n\in {\Bbb N}$. Summarizing the last 
paragraph, we get  
\begin{eqnarray}\label{3}
\tau_k\circ u(W^{n,r})=-\eta(r)+\tau_{k'}\circ u(W)\quad\mbox{ 
\rm if}\quad\left|\tau_k\circ u(W^{n,r})\right|\le\frac1n 
\end{eqnarray} 
for sufficiently large $n\in {\Bbb N}$. We mention that in this 
case $\tau_k\circ u(W^r)=-r+\tau_{k'}\circ u(W)$ which means 
that $k$ and $k'$ are related as in Step 2 of the proof of part 
(a) of Lemma \ref{Lemma3.4} under the restriction of Remark (5) 
of the present section. 

It follows from (v) of Remark (1) of Section 1 and (jjj) of 
Definition \ref{Definition1.7} that for $k'$ as above 
\begin{eqnarray}\label{35}
&&\hspace{-.5cm}2\left(\tau_{k'}\circ u(W+(x-W_0)\1)+\delta 
\vphantom{l^1}\right)-\sigma\left(\tau_{k'}\circ u(W+(x-W_0)\1) 
+\delta,x\vphantom{l^1}\right)\vphantom{\dot{f}}\nonumber \\ 
&&\hspace{.5cm}=\tau_{k'}\circ u(W)+\delta\, ,\quad x\in D^n, 
\vphantom{\dot{f}}
\end{eqnarray} 
where we keep in mind that we take here $W-W_0\1$ for $W$ in 
Definition \ref{Definition1.7} and $W_0$ here for $x_0$ there. 
We pick up the situation in (\ref{3}), i. e., we assume that 
$r$ is chosen such that $\left|\tau_k\circ u(W^{n,r})\right|\le 
\frac1n$ for sufficiently large $n\in {\Bbb N}$. We set $\delta: 
=-\tau_k\circ u(W^{n,r})$. It follows from (\ref{25}) that $r= 
\tau_{k'}\circ u\left(W+\left(\bar{A}_{n,r}-\bar{Y}_{n,r}\right) 
\1\right)+\delta$. With $x:=W_0+\bar{A}_{n,r}-\bar{Y}_{n,r}$ we 
obtain from (\ref{35}) 
\begin{eqnarray}\label{36}
&&\hspace{-.5cm}\eta\left(\tau_{k'}\circ u\left(W+\left(\bar{A 
}_{n,r}-\bar{Y}_{n,r}\right)\1\right)+\delta\right)=\eta(r) 
\vphantom{\dot{f}}\nonumber \\ 
&&\hspace{.5cm}=2r-\sigma\left(r,W_0+\bar{A}_{n,r}-\bar{Y}_{n, 
r}\right)\vphantom{\dot{f}}\nonumber \\ 
&&\hspace{.5cm}=\tau_{k'}\circ u(W)+\delta\vphantom{\dot{f}} 
\nonumber \\ 
&&\hspace{.5cm}=\eta\left(\tau_{k'}\circ u\left(W+\left(\bar{A 
}_{n,\tau_{k'}}-\bar{Y}_{n,\tau_{k'}}\right)\1\right)\right)+ 
\delta\vphantom{\dot{f}}
\end{eqnarray} 
where $\tau_{k'}$ without any further argument stands until the 
end of the present Step 2 for $\tau_{k'}\circ u\left(W+\left( 
\bar{A}_{n,\tau_{k'}}-\bar{Y}_{n,\tau_{k'}}\right)\1\right)$, 
i. e., the last equality sign is (\ref{3}) for $r\equiv r_0=: 
\tau_{k'}\circ u\left(W+\left(\bar{A}_{n,\tau_{k'}}-\bar{Y}_{n, 
\tau_{k'}}\right)\1\right)$. We note that for this particular 
$r_0$ we have $\tau_k\circ u(W^{n,r_0})=0$. For $I\in {\cal I}$ 
being an interval and $\tau_{k'}\circ u\left(W+\left(\bar{A}_{n, 
\tau_{k'}}-\bar{Y}_{n,\tau_{k'}}\right)\1\right)\in I$ we obtain 
\begin{eqnarray*}
\bar{\eta}_I^{-1}\left(\tau_{k'}\circ u(W)\right)=\tau_{k'}\circ 
u\left(W+\left(\bar{A}_{n,\tau_{k'}}-\bar{Y}_{n,\tau_{k'}}\right 
)\1\right)\vphantom{\dot{f}}
\end{eqnarray*} 
and 
\begin{eqnarray*}
\bar{\eta}_I^{-1}\left(\tau_{k'}\circ u(W)+\delta\right)=r= 
\tau_{k'}\circ u\left(W+\left(\bar{A}_{n,r}-\bar{Y}_{n,r}\right) 
\1\right)+\delta\, . 
\end{eqnarray*} 
Recalling $\delta\equiv\delta_k(r)=-\tau_k\circ u(W^{n,r})$ this 
leads to 
\begin{eqnarray}\label{37}
&&\hspace{-.5cm}\frac{d}{d\delta}\bar{\eta}_I^{-1}\left(\tau_{k'} 
\circ u(W)+\delta\right)=\frac{dr}{d\delta}
=1+\frac{d\tau_{k'}\circ u\left(W+\left(\bar{A}_{n,r}-\bar{Y}_{n, 
r}\right)\1\right)}{dr}\left/\ \frac{d\delta_k(r)}{dr}\right. 
\qquad 
\end{eqnarray} 
for sufficiently large $n\in {\Bbb N}$ which means that we assume 
$\left|\tau_k\circ u(W^{n,r})\right|\le \frac1n$. 
\medskip 

\nid 
{\it Step 3 } Setting $p=\tau_{k'}\circ u(W)+\delta$ and looking 
at (\ref{37}), the goal of this step is an estimate on $|\bar{d\eta 
}_I^{-1}(p)/dp -1|$. For the subsequent calculations we recall 
condition (3) (i) and Remark (3) of Section 1. Furthermore, we 
keep Lemma \ref{Lemma3.4} (vi) in mind. In particular it holds that 
\begin{eqnarray}\label{375}
&&\hspace{-.5cm}\bar{A}_{n,r}(W)-\bar{Y}_{n,r}(W)=A_{n,r}(W)-Y_{n 
,r}(W)+(A-Y)_{-r}\left(W^{n,r}\right)\vphantom{\int}\nonumber \\ 
&&\hspace{.5cm}=\left(A_{n,r}(W)-Y_{n,r}(W)-(A-Y)_r(W)\vphantom{ 
\dot{f}}\right)+\left((A-Y)_{-r}\left(W^{n,r}\right)-(A-Y)_{-r}( 
W^r)\vphantom{\dot{f}}\right)\, ,\vphantom{\int}\nonumber \\ 
\end{eqnarray} 
As in (\ref{3.24*}) let $\psi(r)=A_{n,r}(W)-Y_{n,r}(W)-(A-Y)_r(W 
)$. From the convention that $0$ is not a jump time for $X$ it 
follows that  $-r$ is not a jump time for $u(W^r)$. For 
sufficiently large $n\in {\Bbb N}$ we may because of $\left| 
\tau_k\circ u(W^{n,r})\right|\le \frac1n$ also assume that $-r$ 
is not a jump time for $u\left(W^{n,r}\right)$. Using (\ref{375}) 
we obtain for sufficiently large $n\in {\Bbb N}$ and $r\in [0,s)$ 
not being a jump time for $u(W)=X$ such that $\left|\tau_k\circ 
u(W^{n,r})\right|\le\frac1n$ 
\begin{eqnarray*}
&&\hspace{-.5cm}\frac{d}{dr}\left(\bar{A}_{n,r}-\bar{Y}_{n,r} 
\right)(W)=\frac{d}{dr}\left(\left(A_{n,r}-Y_{n,r}\right)(W)-(A-Y 
)_{-r}\left(W^{n,r}\right)\vphantom{\dot{f}}\right)\vphantom{\int} 
 \\ 
&&\hspace{.5cm}=\left(\left(\dot{A}_{n,r}-\dot{Y}_{n,r}\right)(W)
-(A-Y)'_r(W)\right)+\left((A-Y)'_{-r}\left(W^{n,r}\right)-(A-Y)'_{ 
-r}\left(W^r\right)\vphantom{l^1}\right)\vphantom{\int} \\ 
&&\hspace{1.0cm}-\left(\lim_{h\to 0}\frac1h\left((A-Y)_{-r}\left( 
W^{n,r+h}\right)-(A-Y)_{-r}\left(W^{n,r}\right)\right)\right. 
\vphantom{\int} \\ 
&&\hspace{5.0cm}\left.-\lim_{h\to 0}\frac1h\left((A-Y)_{-r}\left( 
W^{r+h}\right)-(A-Y)_{-r}\left(W^r\right)\right)\right)\, . 
\vphantom{\int} 
\end{eqnarray*} 
We treat the three differences on the right-hand side, separated by 
parentheses, separately. For the first one we form the difference 
with the left-hand side and recall how we have applied (\ref{3.23*}) 
in Step 1 of this proof. Furthermore, we take into consideration that 
$(A-Y)'_r(W)=(A-Y)'_0(W^r)$. For the second one, we use the idea of 
(\ref{3.27*}). In order to treat the third difference we recall the 
way we have used Proposition \ref{Proposition2.7} in Lemma 
\ref{Lemma3.2} (b) and (c). Again with the idea of estimating used 
in (\ref{3.27*}) we get 
\begin{eqnarray}\label{38}
&&\hspace{-.5cm}\left|\frac{d}{dr}\left(\bar{A}_{n,r}-\bar{Y}_{n,r} 
\right)(W)-\sum_{k\in{\Bbb Z}\setminus\{0\}}\left\langle\Delta(A-Y 
)_{\tau_k}(W^{n,r})\, ,\, g_n\left(\tau_k\circ u\left(W^{n,r}\right) 
\right)\vphantom{l^1}\right\rangle_{F\to F}\right|\nonumber\vphantom 
{\int} \\ 
&&\hspace{.5cm}\le\left|(A-Y)'_0\left(W^{n,r}\right)-(A-Y)'_0(W^r) 
\vphantom{\sum_{k\in{\Bbb Z}\setminus\{0\}}}\right|+|c_n(r)|+\alpha' 
(S;W)\cdot|\psi(r)|+\delta'(S;W)\cdot|\psi(r)|\vphantom{\int} 
\nonumber \\ 
\end{eqnarray} 
for some positive constant $\delta'(S;W)<\infty$ where $S=[0,s)$. 
In (\ref{38}) and in the next relation the absolute values are taken 
coordinate wise. Furthermore by (\ref{375}),   
\begin{eqnarray}\label{39}
&&\hspace{-.5cm}\left|\left(\nabla_{W_0}\tau_{k'}\circ u\vphantom 
{l^1}\right)\left(W+\left(\bar{A}_{n,r}(W)-\bar{Y}_{n,r}(W)\vphantom 
{l^1}\right)\1\right)-\left(\nabla_{W_0}\tau_{k'}\circ u\vphantom 
{l^1}\right)(W)\right|\vphantom{\int}\nonumber \\  
&&\hspace{.5cm}\le\kappa'(S;W)\cdot|\psi(r)|\vphantom{\int}
\end{eqnarray} 
for some positive constant $\kappa'(S;W)<\infty$ by (j) of 
Definition \ref{Definition1.7} and condition (3) (ii) of Section 1. 
For sufficiently large $n\in {\Bbb N}$ and $r\in [0,s)$ not being a 
jump time for $u(W)=X$ such that $\left|\tau_k\circ u(W^{n,r})\right 
|\le\frac1n$ the subsequent estimate uses (\ref{38}) and (\ref{39}) 
and 
\begin{eqnarray*}
\left\langle\tau_{k'}\circ u\left(W+\left(\bar{A}_{n,r}(W)-\bar{Y}_{ 
n,r}(W)\vphantom{l^1}\right)\1\right)\, ,\, \Delta(A-Y)_{\tau_k}\left 
(W^{n,r}\right)\right\rangle_F=0\, ,\vphantom{\int}  
\end{eqnarray*} 
recall $\left(W+\left(\bar{A}_{n,r}(W)-\bar{Y}_{n,r}(W)\vphantom{l^1 
}\right)\1\right)^r=W^{n,r}$, $\tau_{k'}\circ u\left(\left(W+\left( 
\bar{A}_{n,r}(W)-\bar{Y}_{n,r}(W)\vphantom{l^1}\right)\1\right)^r 
\right)=\tau_k\circ u\left(W^{n,r}\right)$, and (jv) of Definition 
\ref{Definition1.7}. We get the estimate 
\begin{eqnarray*}
&&\hspace{-.5cm}\left|\frac{d}{dr}\tau_{k'}\circ u\left(W+\left( 
\bar{A}_{n,r}(W)-\bar{Y}_{n,r}(W)\vphantom{l^1}\right)\1\right) 
\right|\vphantom{\int}\nonumber \\ 
&&\hspace{.5cm}=\left|\left\langle\left(\nabla_{W_0}\tau_{k'}\circ 
u\vphantom{l^1}\right)\left(W+\left(\bar{A}_{n,r}(W)-\bar{Y}_{n,r} 
(W)\vphantom{l^1}\right)\1\right)\, ,\vphantom{\sum_{k\in{\Bbb Z} 
}}\right.\right.\nonumber \\ 
&&\hspace{1.5cm}\left.\frac{d}{dr}\left(\bar{A}_{n,r}-\bar{Y}_{n,r} 
\right)(W)-\sum_{k\in{\Bbb Z}\setminus\{0\}}\left\langle\Delta(A-Y 
)_{\tau_k}(W^{n,r})\, ,\, g_n\left(\tau_k\circ u(W^{n,r})\right) 
\vphantom{l^1}\right\rangle_{F\to F}\vphantom{\dot{f}}\right 
\rangle_F\nonumber \\ 
&&\hspace{1.0cm}+\sum_{k\in{\Bbb Z}\setminus\{0\}}\left\langle 
\tau_{k'}\circ u\left(W+\left(\bar{A}_{n,r}(W)-\bar{Y}_{n,r}(W) 
\vphantom{l^1}\right)\1\right)\, ,\vphantom{\sum}\right. 
\nonumber \\ 
&&\hspace{1.5cm}\left.\left.\left\langle\Delta(A-Y)_{\tau_k}(W^{ 
n,r})\, ,\, g_n\left(\tau_k\circ u(W^{n,r})\right)\vphantom{l^1} 
\right\rangle_{F\to F}\vphantom{\sum}\right\rangle_F\vphantom{ 
\sum_{k\in{\Bbb Z}}}\right|\nonumber \\ 
&&\hspace{.5cm}\le\left\langle\left|\left(\nabla_{W_0}\tau_{k'} 
\circ u\vphantom{l^1}\right)(W)\right|\, ,\, |c_n(r)|+2\alpha'(S; 
W)\cdot|\psi(r)|+\delta'(S;W)\cdot|\psi(r)|\right\rangle_F 
\vphantom{\int}\nonumber \\ 
&&\hspace{1.0cm}+\left\langle\kappa'(S;W)\cdot\left|\psi(r)(W) 
\right|\, ,\, |c_n(r)|+2\alpha'(S;W)\cdot|\psi(r)|+\delta'(S;W) 
\cdot|\psi(r)|\vphantom{l^1}\right\rangle_F\, .\vphantom{\int} 
\end{eqnarray*} 
For this we recall also that for the whole proof we hold $W\in\{\pi_m 
V:V\in\Omega\}$ with $W_0\not\in G(W-W_0\1)$ fixed. Summing up, we 
have 
\begin{eqnarray}\label{40}
&&\hspace{-.5cm}\left|\frac{d}{dr}\tau_{k'}\circ u\left(W+\left(\bar 
{A}_{n,r}(W)-\bar{Y}_{n,r}(W)\vphantom{l^1}\right)\1\right)\right| 
\le C_1\cdot|\psi(r)|+C_2\cdot|\psi(r)|^2+d_n(r)\qquad 
\end{eqnarray} 
for some positive constants $C_1\equiv C_1(S;W)<\infty$, $C_2\equiv 
C_2(S;W)<\infty$ and $d_n(r)\stack{n\to\infty}{\lra}0$ boundedly on 
$r\in [0,s)\equiv S$. As above we mention that (\ref{40}) holds for 
sufficiently large $n\in {\Bbb N}$ and for $r\in [0,s)$ such that 
$\left|\tau_k\circ u(W^{n,r})\right|\le\frac1n$. The restriction not 
being a jump time for $u (W)=X$ can now be neglected by looking at 
$r=\tau_{k'}\circ u\left(W+\left(\bar{A}_{n,r}-\bar{Y}_{n,r}\right) 
\1\right)+\delta$ again and recalling that $\delta\equiv\delta_k(r) 
=-\tau_k\circ u(W^{n,r})$.

In order to prepare the conclusions of this step we recall (\ref{3}) 
which in particular says that $p=\delta+\tau_{k'}\circ u(W)=\eta(r)$. 
We note that $\psi(r)\equiv\psi(r;n,W)$ is bounded in $n\in {\Bbb N}$ 
on $r\in [0,s)=S$. This is a consequence of the last sentence in 
Lemma \ref{Lemma3.4} (vi). 

We observe that, without loss of generality, we may assume $|\liminf_r 
d\delta_k(r)/dr|>0$ for $r$ such that $\left|\tau_k\circ u(W^{n,r}) 
\right|\le\frac1n$ and $n\in {\Bbb N}$ sufficiently large. This 
assumption is justified by the fact that supposing $d\delta_k(r)/dr=-d 
\tau_k\circ u\left(W^{n,r}\right)/dr\to 0$ on some subsequence $n_k\to 
\infty$ and some $r\equiv r_{n_k}$ with the above hypotheses on $r$, 
we would by (\ref{3.23*}) obtain $|\psi(r)|\to\infty$ on this 
subsequence. However we have already given an argument that this is 
impossible. 
\medskip

For some $C_3\equiv C_3(W;S)$ independent of $n\in {\Bbb N}$ we obtain 
ow from (\ref{37}) and (\ref{40})
\begin{eqnarray}\label{41}
\left|\frac{d\bar{\eta}_I^{-1}(p)}{dp}-1\right|\le C_3\cdot\left| 
\psi\left(\bar{\eta}_I^{-1}(p)\right)\right|+d_n\left(\bar{\eta 
}_I^{-1}(p)\right) 
\end{eqnarray} 
for all $p\in\{\eta(r'):r'\in I\}\cap\left[\tau_{k'}\circ u(W)-\frac 
1n\, ,\, \tau_{k'}\circ u(W)+\frac1n\right]$ and all $k'\in {\Bbb Z} 
\setminus\{0\}$ such that $\tau_{k'}\circ u\left(W+\left(A_{n,\tau_{ 
k'}}-Y_{n,\tau_{k'}}\right)\1\right)\in I$ where $I\in {\cal I}$ is 
an interval. 

Moreover, there exist $b_n(I)>0$ for all intervals $I\in {\cal I}$ 
such that the sum of the $b_n(I)$ over all such $I$ is finite and 
and $b_n(I)\stack{n\to\infty}{\lra}0$ boundedly and 
\begin{eqnarray}\label{42}
\int\left|\frac{d\bar{\eta}_I^{-1}(p)}{dp}-1\right|\, dp\le b_n(I) 
\, .
\end{eqnarray} 
for sufficiently large $n\in N$. Here the integral ranges over all 
$p\in\{\eta(r'):r'\in I\}\setminus\left[\tau_{k'}\circ u(W)-\frac 
1n\, ,\, \tau_{k'}\circ u(W)+\frac1n\right]$. This is a consequence 
of the following observation. For time intervals between two jumps 
of the form $\tau_k\circ u(W^{n,r})$, (\ref{3.23*}) can be treated 
in a similar manner, no matter whether or not we require the 
restriction of Remark (5) of this section. 
\medskip 

\nid 
{\it Step 4 } We introduce a decomposition similar to 
(\ref{3.23*})-(\ref{3.25*}). Let $I\in {\cal I}$, $q\in I$, and 
if $I$ is just a point, set $\inf I:=I-$. For the sake of 
formality, we write also $\bar{\eta}(\inf{I})$ for $\eta(\inf{I 
})$. Using (\ref{2}) and (\ref{3}) we obtain for $k'\equiv k'(k 
;W,\bar{\eta}^{-1}(r))$ and sufficiently large $n\in {\Bbb N}$ 
\begin{eqnarray}\label{4}
&&\hspace{-.5cm}\int_{\inf{I}}^q\left(\dot{A}_{n,r}(W)-\dot{Y 
}_{n,r}(W)\vphantom{l^1}\right)\, d\eta(r)=\int_{\bar{\eta}(\inf 
{I})}^{\bar{\eta}(q)}\left(\dot{A}_{n,\bar{\eta}^{-1}(r)}(W)- 
\dot{Y}_{n,\bar{\eta}^{-1}(r)}(W)\vphantom{l^1}\right)\, dr 
\nonumber \\ 
&&\hspace{.5cm}=\int_{\bar{\eta}(\inf{I})}^{\bar{\eta}(q)}\int 
\left\langle (A-Y)'_{-v}(W^{n,\bar{\eta}^{-1}(r)})\, ,\, g_n(v) 
\vphantom{l^1}\right\rangle_{F\to F}\, dv\, dr\nonumber \\ 
&&\hspace{1.0cm}+\sum_{k\in{\Bbb Z}\setminus\{0\}}\int_{\bar{ 
\eta}(\inf{I})}^{\bar{\eta}(q)}\left\langle\Delta(A-Y)_{\tau_k} 
(W^{n,\bar{\eta}^{-1}(r)})\, ,\, g_n\left(\tau_{k'}\circ u(W)-r 
\right)\vphantom{l^1}\right\rangle_{F\to F}\, dr\nonumber \\ 
&&\hspace{1.0cm}+\int_{\bar{\eta}(\inf{I})}^{\bar{\eta}(q)}c_n 
\left(\bar{\eta}^{-1}(r)\right)\, dr\cdot {\sf e} 
\end{eqnarray} 
where we note that the variable substitution is compatible with 
the above mode of integration with respect to $d\eta$ at jump 
times. We set 
\begin{eqnarray*}
\hat{r}(r)&&\hspace{-.5cm}\equiv\hat{r}(r;n,W):=\int\left\langle 
(A-Y)'_{-v}\left(W^{\bar{\eta}^{-1}(r)}\right)\, ,\, g_n(v) 
\vphantom{l^1}\right\rangle_{F\to F}\, dv\vphantom{\sum_{v\in (- 
\frac1n,\frac1n)}}\nonumber \\  
&&\hspace{.0cm}+\sum_{k\in{\Bbb Z}\setminus\{0\}}\left\langle 
\Delta (A-Y)_{\tau_k}\left(W^{\bar{\eta}^{-1}(r)}\right)\, ,\, 
g_n\left(\tau_{k'}\circ u(W)-r\right)\vphantom{l^1}\right 
\rangle_{F\to F}+c_n\left(\bar{\eta}^{-1}(r)\right)\cdot {\sf e} 
\vphantom{\sum^s_v}\, .  
\end{eqnarray*} 
Differentiating (\ref{4}) and re-integrating then with respect to 
$\bar{\eta}^{-1}$ we get for $I\in {\cal I}$ and $q\in I$ 
\begin{eqnarray}\label{5}
&&\hspace{-.5cm}(A-Y)_{n,q}(W)-(A-Y)_{n,\inf I}(W)\vphantom{\sum_{ 
k\in {\Bbb Z}}}\nonumber \\  
&&\hspace{5.0cm}-\left((A-Y)_q(W)-(A-Y)_{\inf I}(W)\vphantom{l^1} 
\right)\vphantom{\sum_{k\in {\Bbb Z}}}\nonumber \\
&&\hspace{.5cm}=\int_{\bar{\eta}(\inf I)}^{\bar{\eta}(q)}\left(\dot 
{A}_{n,\bar{\eta}^{-1}(r)}(W)-\dot{Y}_{n,\bar{\eta}^{-1}(r)}(W) 
\right)\, d\bar{\eta}^{-1}(r)\vphantom{\sum_{k\in {\Bbb Z}}}\nonumber 
 \\  
&&\hspace{5.0cm}-\left((A-Y)_q(W)-(A-Y)_{\inf I}(W)\right)\vphantom 
{\sum_{k\in {\Bbb Z}}}\nonumber \\ 
&&\hspace{.5cm}=\int_{\bar{\eta}(\inf I)}^{\bar{\eta}(q)}\int\left 
\langle (A-Y)'_{-v}\left(W^{n,\bar{\eta}^{-1}(r)}\vphantom{l^1} 
\right)-(A-Y)'_{-v}\left(W^{\bar{\eta}^{-1}(r)}\right)\, ,\, 
\right.\vphantom{\sum_{k\in {\Bbb Z}}}\nonumber \\  
&&\hspace{5.0cm}\left.g_n(v)\vphantom{\left(W^{\bar{\eta}^{-1}(r)} 
\right)}\right\rangle_{F\to F}\, dv\, d\bar{\eta}^{-1}(r)\vphantom{ 
\sum_{k\in {\Bbb Z}}}\nonumber \\  
&&\hspace{1.0cm}+\sum_{k\in {\Bbb Z}\setminus\{0\}}\int_{\bar{\eta} 
(\inf I)}^{\bar{\eta}(q)}\left\langle\Delta (A-Y)_{\tau_k}\left(W^{ 
n,\bar{\eta}^{-1}(r)}\vphantom{l^1}\right)-\Delta (A-Y)_{\tau_k} 
\left(W^{\bar{\eta}^{-1}(r)}\right)\, ,\, \right.\nonumber \\ 
&&\hspace{5.0cm}\left.g_n\left(\tau_{k'}-r\right)\vphantom{W^{\bar{ 
\eta}^{-1}(s)}}\right\rangle_{F\to F}\, d\bar{\eta}^{-1}(r) 
\vphantom{\sum_{k\in {\Bbb Z}}}\nonumber \\  
&&\hspace{1.0cm}+\int_{\bar{\eta}(\inf I)}^{\bar{\eta}(q)}\hat{r}(r 
)\, d\bar{\eta}^{-1}(r)-\left((A-Y)_q(W)-(A-Y)_{\inf I}(W) 
\vphantom{l^1}\right)\vphantom{\sum_{k\in {\Bbb Z}}}  
\end{eqnarray} 
where, in case that $I$ is just a point, we read the integral 
$\int_{\bar{\eta}(\inf I)}^{\bar{\eta}(q)}(\, \cdot\, )(\bar{ 
\eta}^{-1}(r))\, d\bar{\eta}^{-1}(r)$ as $\int_{\inf I}^q(\, 
\cdot\, )(r)\, dr=0$. Now we recall that the time 
re-parametrization $\sigma(s,x)$, $s\in {\Bbb R}$, $x\in D^n$, 
with respect to $W-W_0\1$ and $W_0$ of Remark (1) in Section 1 
is on every finite subinterval of ${\Bbb R}$ just in the time 
points $\tau_{k'}\circ u(W)$ (and all $x\in D^n$) uniquely 
defined, cf. (v) of Remark (1) in Section 1. 
\medskip 

\nid 
{\it Step 5 } We construct estimates on the items of (\ref{5}). 
Using (\ref{1}) and (\ref{36}) we observe the following. Except 
for neighborhoods of the finitely many jumps times $\tau_{k'} 
\circ u\left(W+\left(\bar{A}_{n,\tau_{k'}}-\bar{Y}_{n,\tau_{k'}} 
\right)\1\right)\in [0,s)$ such that they correspond via $\eta$ 
to $\ve$-neighborhoods of $\tau_{k'}\circ u(W)$ in the sense of 
(\ref{36}), the functions $\eta\equiv\eta(\cdot;n)$ can be 
slightly altered on $[0,s)$. For $\ve>0$ let {\boldmath${\ve}$ 
}$(S)$ denote the union of these neighborhoods of $\tau_{k'} 
\circ u\left(W+\left(\bar{A}_{n,\tau_{k'}}-\bar{Y}_{n,\tau_{k'}} 
\right)\1\right)$ intersected with $S=[0,s)$. 

%

In particular, we may assume that for no $n\in {\Bbb N}$ there 
is some interval of constancy on $[0,s)\setminus${\boldmath${\ve 
}$}$([0,s))$. Without loss of generality, we may also assume that 
for all $I\in{\cal I}$ being intervals the restriction of $\bar 
{\eta}_I^{-1}$ to $\{\bar{\eta}_I(r):r\in I\setminus${\boldmath$ 
{\ve}$}$([0,s))\}$ is continuously differentiable. 
\medskip

We set 
\begin{eqnarray*}
\hat{\psi}(r)\equiv\hat{\psi}(r\, ;n,W):=A_{n,\bar{\eta}^{-1}(r)} 
(W)-Y_{n,\bar{\eta}^{-1}(r)}(W)-(A-Y)_{\bar{\eta}^{-1}(r)}(W)\, . 
\end{eqnarray*} 
Recalling how we proceeded in Steps 2 and 3 of the proof of part 
(a) of Lemma \ref{Lemma3.4} under the restriction of Remark (5) of 
this section, we arrive with (\ref{5}) at 
\begin{eqnarray}\label{8} 
&&\hspace{-.5cm}(A-Y)_{n,q}(W)-(A-Y)_{n,\inf I}(W)-\left((A-Y)_q 
(W)-(A-Y)_{\inf I}(W)\vphantom{l^1}\right)\vphantom{\sum}\nonumber 
 \\
&&\hspace{.5cm}\le\int_{\bar{\eta}(\inf I)}^{\bar{\eta}(q)}\left( 
\alpha'(S;W)\cdot\left|\hat{\psi}(r)\right|+\gamma'(S;W)\sum_{k\in 
{\Bbb Z}\setminus\{0\}}\left\langle\left|\hat{\psi}(r)\right|\, , 
\, g_n\left(\tau_{k'}-r\right)\vphantom{\dot{f}}\right\rangle_{F 
\to F}\right)\, d\bar{\eta}^{-1}(r)\vphantom{\sum_{k\in {\Bbb Z}}} 
\nonumber \\  
&&\hspace{1.0cm}+\left|\int_{\bar{\eta}(\inf I)}^{\bar{\eta}(q)} 
\int\left\langle(A-Y)'_{-v}\left(W^{\bar{\eta}^{-1}(r)}\right)\, , 
\, g_n(v)\vphantom{l^1}\right\rangle_{F\to F}\, dv\, d\bar{\eta}^{ 
-1}(r)\right.\vphantom{\sum_{k\in {\Bbb Z}}}\nonumber \\  
&&\hspace{1.0cm}+\sum_{k\in{\Bbb Z}\setminus\{0\}}\int_{\bar{\eta} 
(\inf I)}^{\bar{\eta}(q)}\left\langle\Delta (A-Y)_{\tau_k}\left( 
W^{\bar{\eta}^{-1}(r)}\right)\, ,\, g_n\left(\tau_{k'}-r\right) 
\vphantom{l^1}\right\rangle_{F\to F}\, d\bar{\eta}^{-1}(r)\nonumber 
 \\ 
&&\hspace{1.0cm}\left.-\left((A-Y)_q(W)-(A-Y)_{\inf I}(W)\vphantom 
{l^1}\right)\vphantom{\int_{\inf{I}-}^s}\right|+\int_{\inf I}^q 
\left|c_n(r)\right|\, dr\cdot {\sf e}\vphantom{\sum_{k\in {\Bbb Z}} 
}\nonumber \\ 
&&\hspace{0.5cm}=:R_1+|T_1|+R_2\vphantom{\sum_{k\in {\Bbb Z}}}
\end{eqnarray} 
where $R_1\equiv R_1(n;I)$ is the term in the first line after 
the `` $\le$ " sign, $T_1\equiv T_1(n;I)$ is the term inside the 
coordinate wise absolute value signs that follows $R_1$. $R_2 
\equiv R_2(n;I)$ is the last integral on the right-hand side. 
Because of (\ref{3}) we have for all $I\in {\cal I}$ and 
sufficiently large $n\in {\Bbb N}$ 
\begin{eqnarray*}
\tau_k\circ u(W^{n,\bar{\eta}^{-1}(r)})=-r+\tau_{k'}\circ u(W)= 
\tau_k\circ u(W^r)
\end{eqnarray*} 
or, equivalently, 
\begin{eqnarray*}
\tau_k\circ u(W^{n,\bar{\eta}^{-1}(r)})-\tau_k\circ u(W^{\bar{ 
\eta}^{-1}(r)})=\tau_k\circ u(W^r)-\tau_k\circ u(W^{\bar{\eta}^{ 
-1}(r)})
\end{eqnarray*} 
if $\left|\tau_k\circ u(W^{n,\bar{\eta}^{-1}(r)})\right|\le\frac 
1n$ which implies $\left|-r+\tau_{k'}\circ u(W)\right|\le\frac1n$. 
This yields 
\begin{eqnarray*}
\left|\Delta (A-Y)_{\tau_k}\left(W^{\bar{\eta}^{-1}(r)}\right)- 
\Delta (A-Y)_{\tau_k}\left(W^r\right)\right|\le\gamma'(S;W)\cdot 
\left|\hat{\psi}(r)\right|
\end{eqnarray*} 
and therefore 
\begin{eqnarray}\label{9} 
T_1&&\hspace{-.5cm}\le\left|\int_{\inf I}^q\int\left\langle(A-Y 
)'_{-v}\left(W^r\right)\, ,\, g_n(v)\vphantom{l^1}\right\rangle_{F 
\to F}\, dv\, dr\right.\vphantom{\sum_{k\in {\Bbb Z}}}\nonumber \\ 
&&\hspace{.0cm}+\sum_{k\in{\Bbb Z}\setminus\{0\}}\int_{\bar{\eta} 
(\inf I)}^{\bar{\eta}(q)}\left\langle\Delta (A-Y)_{\tau_k}\left(W^r 
\right)\, ,\, g_n\left(\tau_{k'}-r\right)\vphantom{l^1}\right 
\rangle_{F\to F}\, d\bar{\eta}^{-1}(r)\nonumber \\ 
&&\hspace{.0cm}-\left.\left((A-Y)_q(W)-(A-Y)_{\inf I}(W)\vphantom 
{l^1}\right)\vphantom{\int_{\bar{\eta}(\inf I)}^{\bar{\eta}(q)}} 
\right|\nonumber \\ 
&&\hspace{.0cm}+\gamma'(S;W)\sum_{k\in{\Bbb Z}\setminus\{0\}} 
\int_{\bar{\eta}(\inf I)}^{\bar{\eta}(q)}\left\langle\left|\hat{ 
\psi}(r)\right|\, ,\, g_n\left(\tau_{k'}-r\right)\vphantom{l^1} 
\right\rangle_{F\to F}\, d\bar{\eta}^{-1}(r)\nonumber \\ 
&&\hspace{-.5cm}=:|T_2|+R_3\vphantom{\sum_{k\in {\Bbb Z}}}
\end{eqnarray} 
where $T_2\equiv T_2(n;I)$ is the term inside the coordinate wise 
absolute value signs after the `` $\le$ " sign. $R_3\equiv R_3 
(n;I)$ is the last line on the right-hand side. Furthermore, if 
$I$ is an interval, we have 
\begin{eqnarray*}
&&\hspace{-0.5cm}\int_{\bar{\eta}(\inf I)}^{\bar{\eta}(q)}\int_{v 
\in\left(-\frac1n,\frac1n\right)}\left\langle (A-Y)'_{r-v}(W)\, , 
\, g_n(v)\vphantom{\dot{f}}\right\rangle_{F\to F}\, dv\, dr 
\vphantom{\left|\sum_{k\in {\Bbb Z}}\right|} \\ 
&&\hspace{1.0cm}+\int_{\bar{\eta}(\inf I)}^{\bar{\eta}(q)}\sum_{v 
\in\left(-\frac1n,\frac1n\right)}\left\langle\Delta (A-Y)_{r-v}(W) 
\, ,\, g_n(v)\vphantom{\dot{f}}\right\rangle_{F\to F}\, dr 
\vphantom{\left|\sum_{k\in {\Bbb Z}}\right|} \\ 
&&\hspace{.5cm}=\int_{v\in\left(-\frac1n,\frac1n\right)}\left 
\langle(A-Y)_{\bar{\eta}(q)-v}(W)-(A-Y)_{\bar{\eta}(\inf I)-v}(W) 
\, ,\, g_n(v)\vphantom{\dot{f}}\right\rangle_{F\to F}\, dv 
\vphantom{\left|\sum_{k\in {\Bbb Z}}\right|}
\end{eqnarray*} 
and consequently  
\begin{eqnarray}\label{10}
T_2&&\hspace{-0.5cm}\le\left|\left(\int_{v\in\left(-\frac1n,\frac1n 
\right)}\left\langle(A-Y)_{\bar{\eta}(q)-v}(W)-(A-Y)_{\bar{\eta}(\inf 
I)-v}(W)\, ,\, g_n(v)\vphantom{\dot{f}}\right\rangle_{F\to F}\, dv 
\right.\right.\nonumber \\ 
&&\hspace{.0cm}\left.-\left((A-Y)_{\bar{\eta}(q)}(W)-(A-Y)_{\bar{ 
\eta}(\inf I)}\vphantom{l^1}\right)\vphantom{\int_{\bar{\eta}(\inf 
I)}^{\bar{\eta}(q)}}\right)\nonumber \\ 
&&\hspace{.0cm}+\left(\int_{v\in\left(-\frac1n,\frac1n\right)}\left 
\langle\int_{\inf I}^q(A-Y)'_{r-v}(W)\, dr-\int_{\bar{\eta}(\inf I) 
}^{\bar{\eta}(q)}(A-Y)'_{r-v}(W)\, dr\, ,\, g_n(v)\vphantom{\dot{f}} 
\right\rangle_{F\to F} dv\right.\nonumber \\ 
&&\hspace{.0cm}+\int_{\inf I}^q\sum_{v\in\left(-\frac1n,\frac1n\right 
)}\left\langle\Delta (A-Y)_{r-v}\, ,\, g_n(v)\vphantom{\dot{f}}\right 
\rangle_{F\to F}\, dr\nonumber \\ 
&&\hspace{4cm}-\int_{\bar{\eta}(\inf I)}^{\bar{\eta}(q)}\sum_{v\in 
\left(-\frac1n,\frac1n\right)} \left\langle\Delta (A-Y)_{r-v}\, dr\, , 
\, g_n(v)\vphantom{\dot{f}}\right\rangle_{F\to F}\nonumber \\ 
&&\hspace{.0cm}\left.-\left((A-Y)_q(W)-(A-Y)_{\inf I}(W)\vphantom{l^1} 
\right)+\left((A-Y)_{\bar{\eta}(q)}(W)-(A-Y)_{\bar{\eta}(\inf I)}(W) 
\vphantom{l^1}\right)\vphantom{\int_{v\in\left(-\frac1n,\frac1n\right)} 
}\right)\nonumber \\ 
&&\hspace{.0cm}-\left(\vphantom{\int_{v\in\left(-\frac1n,\frac1n\right) 
}}\right.\int_{\inf I}^q\sum_{v\in\left(-\frac1n,\frac1n\right)}\left 
\langle\Delta (A-Y)_{r-v}\, ,\, g_n(v)\vphantom{\dot{f}}\right\rangle_{ 
F\to F}\, dr\nonumber \\ 
&&\hspace{4cm}-\int_{\bar{\eta}(\inf I)}^{\bar{\eta}(q)}\sum_{v\in 
\left(-\frac1n,\frac1n\right)} \left\langle\Delta (A-Y)_{r-v}\, dr\, , 
\, g_n(v)\vphantom{\dot{f}}\right\rangle_{F\to F}\left.\left.\vphantom 
{\int_{v\in\left(-\frac1n,\frac1n\right)}}\right)\right|+R_4 
\nonumber \\ 
&&\hspace{-.5cm}=|T_{31}+T_{32}-Q|+R_4\vphantom{\int_v}
\end{eqnarray} 
where $T_{31}\equiv T_{31}(n;I)$, $T_{32}\equiv T_{32}(n;I)$, and 
$Q\equiv Q(n;I)$ are the terms inside the coordinate wise absolute 
value signs separated by parentheses and 
\begin{eqnarray}\label{11}
R_4&&\hspace{-.5cm}\equiv R_4(n;I):= 
\left|\sum_{k\in{\Bbb Z}\setminus\{0\}}\int_{ 
\bar{\eta}(\inf I)}^{\bar{\eta}(q)}\left\langle\Delta (A-Y)_{ 
\tau_k}(W^r)\, ,\, g_n(\tau_{k'}-r)\vphantom{\dot{f}}\right 
\rangle_{F\to F}\, d\left(\bar{\eta}^{-1}(r)-r\right)\right| 
\vphantom{\left|\sum_{k\in {\Bbb Z}}\right|}\, .\vphantom{\left| 
\sum_{k\in {\Bbb Z}}\right|}\nonumber \\ 
\end{eqnarray} 

\nid 
{\it Step 6 } In this step we look at (\ref{8})-(\ref{11}) 
simultaneously for all $I\in {\cal I}$. By definition, the most 
right $I\in {\cal I}$ is an interval. For that $I$ we choose $q: 
=s$. For all other intervals $I\in {\cal I}$, let $q:=\sup I-$. 
If $I$ is just a point then let $q:=I$. In this case we remind of 
the convention $\int_{\bar{\eta}(\inf I)}^{\bar{\eta}(q)}(\,  
\cdot\, )(\bar{\eta}^{-1}(r))\, d\bar{\eta}^{-1}(r)=$ $\int_{\inf 
I}^q(\, \cdot\, )(r)\, dr=0$. Moreover, if $I=q$ is just a point, 
we have $R_1=R_2=R_3=0$ and $T_1=T_2=-\Delta(A-Y)_q(W)$. If $I$ 
is just a point, we set $T_{31}:=-\left((A-Y)_{\bar{\eta}(q)}(W)- 
(A-Y)_{\bar{\eta}(\inf I)}\vphantom{l^1}\right)$, $T_{32}:=-\left 
((A-Y)_q(W)-(A-Y)_{\inf I}\vphantom{l^1}\right)+\left((A-Y)_{\bar 
{\eta}(q)}(W)-(A-Y)_{\bar{\eta}(\inf I)}\vphantom{l^1}\right)$, 
and consequently $R_4:=0$. 

Recalling (\ref{41}) we get with some positive constant $C_4\equiv 
C_4(S;W)<\infty$ not depending on $n\in {\Bbb N}$ the estimate
\begin{eqnarray}\label{12} 
&&\hspace{-.5cm}\sum_{I\in {\cal I}\ \mbox{\rm\scriptsize is 
interval}}R_4(n;I)\nonumber \\ 
&&\hspace{.5cm}\le C_4\cdot\hspace{-.5cm}\sum_{I\in {\cal I}\mbox{ 
\rm\scriptsize is interval}}\ \sum_{k\in {\Bbb Z}\setminus\{0\}} 
\int_{\bar{\eta}(\inf I)}^{\bar{\eta}(\sup I)}\left\langle g_n( 
\tau_{k'}-p)\, ,\, \left(\left|\hat{\psi}(p)\right|+d_n\left(\bar 
{\eta}^{-1}(p)\right)\cdot {\sf e}\right)\right\rangle_{F\to F}\, 
dp\nonumber \\ 
\end{eqnarray} 
where for the usage of (\ref{41}) we have taken into 
consideration that $\hat{\psi}(p)=\psi\left(\bar{\eta}^{-1}(p) 
\right)$. 
\medskip 

We complete the proof by referring to the similarities in 
(\ref{3.29*})-(\ref{3.31*}). First we take a look at $\sum_{I\in 
{\cal I}}T_{31}(n;I)$ and $\sum_{I\in {\cal I}}T_{32}(n;I)$, 
where $T_{31}$ and $T_{32}$ are defined in (\ref{10}). Here we 
point out similarities to (\ref{3.30*}) if $I$ is an interval. 
For the choice of $q=I$ in case that $I$ is just a point recall 
the first paragraph of the present step. As in (\ref{3.30*}) it 
turns out that  
\begin{eqnarray}\label{13} 
\hat{\rho}_1(s;n,W):=\sum_{I\in{\cal I}}T_{31}(n;I)\stack{n\to 
\infty}{\lra}0\, . 
\end{eqnarray} 
In addition, assuming without loss of generality $\eta(s)\le s$,  
\begin{eqnarray}\label{14} 
&&\hspace{-.5cm}\hat{\rho}_2(s;n,W):=\sum_{I\in{\cal I}}T_{32} 
(n;I)=\int_{v\in\left(-\frac1n,\frac1n\right)}\left\langle\int_{ 
\eta(s)}^s(A-Y)'_{r-v}(W)\, dr\, ,\, g_n(v)\vphantom{\dot{f}} 
\right\rangle_{F\to F}\, dv\nonumber \\ 
&&\hspace{.5cm}+\int_{\eta(s)}^s\sum_{v\in\left(-\frac1n,\frac1n 
\right)}\left\langle\Delta (A-Y)_{r-v}(W)\, ,\, g_n(v)\vphantom 
{\dot{f}}\right\rangle_{F\to F}\, dr\nonumber \\ 
&&\hspace{.5cm}-\left((A-Y)_s(W)-(A-Y)_{\eta(s)}(W)\vphantom{l^1} 
\right)\stack{n\to\infty}{\lra}0\vphantom{\left|\sum_{k\in{\Bbb 
Z}}\right|}\, . 
\end{eqnarray} 
Moreover, with $\tau_{k'_s}:=\max\{k':\tau_{k'}<s\}$ and $d_s:=s- 
\tau_{k'_s}$, we have for sufficiently large $n\in {\Bbb N}$
\begin{eqnarray*} 
&&\hspace{-.5cm}\sum_{I\in{\cal I}}Q(n;I)=\int_{\eta(s)}^s\sum_{ 
v\in\left(-\frac1n,\frac1n\right)}\left\langle\Delta (A-Y)_{r-v}(W) 
\, ,\, g_n(v)\vphantom{\dot{f}}\right\rangle_{F\to F}\, dr\nonumber 
 \\ 
&&\hspace{.5cm}\le\frac{C_5}{d_s}\left|\textstyle{\left(s-\frac{d_s 
}{2}\right)-\eta\left(s-\frac{d_s}{2}\right)}\right|=\frac{C_5}{d_s} 
\left|\textstyle{\sigma\left(s-\frac{d_s}{2},W_0+(\bar{A}-\bar{Y})_{ 
n,s-\frac{d_s}{2}}\right)-\left(s-\frac{d_s}{2}\right)}\right| 
\vphantom{\sum_{v\in\left(-\frac1n,\frac1n\right)}}
\end{eqnarray*} 
since we assume that $s$ is not a jump time for $X$. Here $C_5 
\equiv C_5(S;W)<\infty$ is some positive constant not depending 
on $n\in {\Bbb N}$. The absolute values are take coordinate wise. 
From Remark (1) (i), (iii) of Section 1 we get the existence of 
some positive constant $C_6\equiv C_6(S;W)<\infty$ not depending 
on $n\in {\Bbb N}$ such that $\sum_{I\in{\cal I}}Q(n;I)\le C_6/ 
d_s|(\bar{A}-\bar{Y})_{n,s-\frac{d_s}{2}}|$. Recalling how we 
have estimated in (\ref{3.27*}), by (\ref{375}) there are constants 
$C_7$ and $C_8$ with the properties of $C_6$ such that 
\begin{eqnarray}\label{15} 
&&\hspace{-.5cm}\sum_{I\in{\cal I}}Q(n;I)\le\frac{C_6(1+\alpha'( 
(-s,s);W)))}{d_s}\left|\textstyle{\psi\left(s-\frac{d_s}{2}\right)} 
\right|=\frac{C_7}{d_s}\left|\textstyle{\psi\left(s-\frac{d_s}{2} 
\right)}\right|\vphantom{\sum_{v\in\left(-\frac1n,\frac1n\right)}} 
\nonumber \\ 
&&\hspace{.5cm}=\frac{C_7}{d_s}\left|\textstyle{\hat{\psi}\left( 
\eta(s-\frac{d_s}{2})\right)}\right|\le\frac{C_8}{d_s}\int_{r=0 
}^{\eta(s)}\left\langle|\hat{\psi}(r)|\, ,\, g_n\left(\textstyle{ 
\eta(s-\frac{d_s}{2})-r}\right)\right\rangle_{F\to F}\, dr\vphantom 
{\sum_{v\in\left(-\frac1n,\frac1n\right)}}\, .\quad 
\end{eqnarray} 
For the last line we have also taken into consideration that 
because of (\ref{42}) it holds that $\eta(s-\frac{d_s}{2})<\eta(s)$ 
for sufficiently large $n\in {\Bbb N}$. 

Now we screen the terms $R_1$, $R_2$, $R_3$, and $R_4$. We take into 
consideration relations (\ref{41}) as well as (\ref{42}) and the 
fact that, in case that $I\in {\cal I}$ is just a point, we have 
$R_1=R_2=R_3=R_4=0$. We introduce 
\begin{eqnarray*} 
\hat{\beta}_1(r;n,W):=\alpha'(S;W)\cdot{\sf e}+2\gamma'(S;W)\sum_{ 
k\in {\Bbb Z}\setminus\{0\}}g_n\left(\tau_{k'}-r\right)+\frac{C_8} 
{d_s}g_n\left(\textstyle{\eta(s-\frac{d_s}{2})-r}\right)\, , 
\end{eqnarray*} 
which is the sum of the coefficients for $|\hat{\psi}(r)|$ in $Q$ 
and $R_1$ as well as $R_3$ for integration with respect to $dr$. 
In addition, using (\ref{41}), we get some coefficient $\hat{\beta 
}_2(r;n,W)$ for $|\hat{\psi}(r)|$ coming from $R_1$, $R_3$, and 
$R_4$ for integration with respect to $d\bar{\eta}^{-1}(r)-dr$. 
Relations (\ref{11}) and (\ref{12}) give the idea how to construct 
$\hat{\beta}_2(r;n,W)$. We set $\hat{\beta}(r;n,W):=\hat{\beta}_1 
(r;n,W)+\hat{\beta}_2(r;n,W)$. 

Furthermore, let $\hat{R}(s;n,W)$ be the sum of $\hat{\rho}_1(s;n,W) 
$, $\hat{\rho}_2(s;n,W)$ (cf. (\ref{13}) and (\ref{14})), $\sum_{I 
\in {\cal I}}R_2(n;I)$, and all other terms in $\sum_{I\in {\cal I}} 
R_1(n;I)$, $\sum_{I\in {\cal I}}R_3(n;I)$, $\sum_{I\in {\cal I}}R_4( 
n;I)$ tending by (\ref{41}) and (\ref{42}) to zero as $n\to\infty$. 
Also here, relations (\ref{11}) and (\ref{12}) give the idea. We 
arrive at 
\begin{eqnarray*} 
\left|\hat{\psi}(\eta(s))\right|=|\psi(s)|&&\hspace{-.5cm}\le\int_{ 
r=0}^{\eta(s)}\left\langle\left|\hat{\psi}(r)\right|\, ,\, \hat{ 
\beta}(r;n,W)\vphantom{\dot{f}}\right\rangle_{F\to F}\, dr+\hat{R}( 
s;n,W)\, ,\quad n\in {\Bbb N}, 
\end{eqnarray*} 
where all absolute values are coordinate wise. This is the 
counterpart to (\ref{3.29*})-(\ref{3.31*}). The claim follows now as 
in (\ref{3.32*}). 
\qed

\subsection{Proof of Theorem \ref{Theorem1.11}} 

{\bf Proof of Theorem \ref{Theorem1.11} } Let $W^{n,s}$, $s\in {\Bbb 
R}$, denote the flow defined in Lemma \ref{Lemma3.4} (vi). In Step 1, 
we will establish an equation for the measures $Q^{(m)}_{\sbnu}\circ 
W^{n,-v}$, $v\in (0,t]$, assuming $t>0$. 

Using this equation we will demonstrate in Step 2 that the 
Radon-Nikodym derivatives 
\begin{eqnarray*}
\omega^{(m)}_{n,-v}:=\frac{dQ^{(m)}_{\sbnu}\circ W^{n,-v}}{dQ^{(m)}_{ 
\sbnu}} 
\end{eqnarray*}
exist in a certain sense. Furthermore, we will provide a representation 
of $\omega^{(m)}_{n,-v}$. In other words, we will derive a pre-version 
(\ref{1.1}) on projections to $H_i$, $i\in J(m)$, and $e_j$, $j\in\{1, 
\ldots ,n\cdot d\}$, as far as jumps are not involved. 

In Step 3, we are going to incorporate the jumps, i. e., we are going to 
determine the weak limit $\lim_{n\to\infty}\omega^{(m)}_{n,-v}Q^{(m)}_{ 
\sbnu}$. In this way we will actually prove the existence of densities of 
the form 
\begin{eqnarray*}
\omega^{(m)}_{-v}:=\frac{dQ^{(m)}_{\sbnu}\circ W^{-v}}{dQ^{(m)}_{\sbnu}} 
\, ,\quad v\in {\Bbb R}. 
\end{eqnarray*}

In Step 4 we will prepare the verification of the conditions of Theorem 
\ref{Theorem2.1}. Moreover, in Step 5, we shall verify (i) and (iii) of 
Theorem \ref{Theorem2.1} and carry out the approximation $\lim_{m\to 
\infty}\omega^{(m)}_{-t}$ which will lead to the representation 
(\ref{1.1}). In Step 6, we are going to verify the remaining condition 
(ii) of Theorem \ref{Theorem2.1} which will prove the theorem.
\medskip 
 
\nid 
{\it Step 1 } In this step, we shall derive an equation for the measures 
$Q^{(m)}_{\sbnu}\circ W^{n,-v}$, $v\in (0,t]$. 

Let $\vp$ be a cylindrical function on $C({\Bbb R};F)$ of the form $\vp 
(W)=f_0(W_0)\cdot f_1(W_{t_1}-W_0,\ldots ,W_{t_l}-W_0)$ where $f_0\in 
C_b^1(F)$, $f_1\in C_b^1(F^l)$, and $t_i\in \{z\cdot t/2^m:z\in {\Bbb Z} 
\setminus\{0\}\}$, $i\in \{1,\ldots ,l\}$, $l\in {\Bbb N}$. We notice 
that if we consider these test $\vp$ functions under the measure $Q^{(m, 
r)}_{\sbnu}$ then we deal exactly with the test functions used in Lemma 
\ref{Lemma3.2}. 
\medskip 

Let $X_n=W+A_n$ and $Y_n$, $n\in {\Bbb N}$, be the sequence of 
processes introduced in Lemma \ref{Lemma3.4} and let $W^{n,s}:= 
W_{\cdot +s}+(A_{n,s}-Y_{n,s})\1$, $s\in {\Bbb R}$, be the flow 
which has been defined in Lemma \ref{Lemma3.4} (vi) on $\left\{ 
(\pi_m W)_{\cdot +w}:W\in\Omega\, , \ w\in\left[0,\frac{1}{2^m} 
\cdot t\right)\right\}$. For the sake of definiteness we extend 
this definition to $\Omega$ by setting $A_{n,\cdot}(W):=A(W)$ 
and $Y_{n,\cdot}(W):=Y(W)$ if $W\not\in\left\{(\pi_m W)_{\cdot 
+w}:W\in\Omega\, ,\ w\in\left[0,\frac{1}{2^m}\cdot t\right) 
\right\}$. Note that now $W^{n,s}$, $s\in {\Bbb R}$, is a flow 
on $\Omega$ with $W^{n,0}=W$. We get 
\begin{eqnarray}\label{3.24}
&&\hspace{-.5cm}\int\vp\, ~ dQ^{(m)}_{\sbnu}\circ W^{n,-v}-\int\vp 
\,~ dQ^{(m)}_{\sbnu} \nonumber \\ 
&&\hspace{.5cm}=\int\vp(W^{n,v})\, Q^{(m)}_{\sbnu}(dW)-\int\vp(W) 
\, Q^{(m)}_{\sbnu}(dW) \nonumber \\ 
&&\hspace{.5cm}=\int\left(\vp\left(W_{\cdot +v}+(A_{n,v}-Y_{n,v}) 
\1\vphantom{l^1}\right)-\vp(W)\vphantom{\dot{f}}\right)\, Q^{(m) 
}_{\sbnu}(dW)\nonumber \\ 
&&\hspace{.5cm}=\lim_{r\to\infty}\int\left(\vp\left(W_{\cdot +v}+ 
(A_{n,v}-Y_{n,v})\1\vphantom{l^1}\right)-\vp(W)\vphantom{\dot{f}} 
\right)\, Q^{(m,r)}_{\sbnu}(dW)
\end{eqnarray} 
where the limit $\lim_{r\to\infty}$ follows from Lemma \ref{Lemma3.4} 
(b) and Lemma \ref{Lemma2.2} (c) with $\psi\equiv\psi(A_n-Y_n)$. For 
this limit we take also into consideration condition (4) (i) of 
Section 1 and that $A_n=A_n^1$ and $Y_n=Y_n^1$ by Lemma \ref{Lemma3.4} 
(i).

Recalling Definition \ref{Definition2.6} (b) and using Lemma 
\ref{Lemma3.2} (b), we obtain 
\begin{eqnarray*} 
&&\hspace{-.5cm}\int\left(\vp\left(W_{\cdot +v}+(A_{n,v}-Y_{n,v}) 
\1\vphantom{l^1}\right)-\vp(W)\vphantom{\dot{f}}\right)\, Q^{(m,r) 
}_{\sbnu}(dW) \\ 
&&\hspace{.5cm}=\int\int_{\sigma=0}^v\left\langle (D\vp)\left(\pi_{ 
m,r}\left(W_{\cdot +\sigma}+(A_{n,\sigma}-Y_{n,\sigma})\1\vphantom{ 
l^1}\right)\right)\, ,\, \vphantom{j^{-1}(\dot{W}_{\cdot +\sigma}+ 
\dot{A}_{n,\sigma}}\right. \\ 
&&\hspace{2.5cm}\left. j^{-1}\left(\pi_{m,r}\left(W_{\cdot +\sigma} 
+(A_{n,\sigma}-Y_{n,\sigma})\1\vphantom{l^1}\right)\right)\, \dot{} 
\, \vphantom{\dot{f}}\right\rangle_H\, d\sigma\, Q^{(m,r)}_{\sbnu} 
(dW) 
\end{eqnarray*}
where the ``$\ \dot{}\ $" indicates the weak mixed derivative with 
respect to $\sigma\in [0,t]$. Because of Lemma \ref{Lemma3.2} (c)  
this is thus equal to 
\begin{eqnarray*}
&&\hspace{-.5cm}\int\left(\vp\left(W_{\cdot +v}+(A_{n,v}-Y_{n,v})\1 
\vphantom{l^1}\right)-\vp(W)\vphantom{\dot{f}}\right)\, Q^{(m,r)}_{ 
\sbnu}(dW) \\
&&\hspace{.5cm}=\int\int_{\sigma=0}^v\left\langle D\vp\circ\left( 
W_{\cdot +\sigma}+(A_{n,\sigma}-Y_{n,\sigma})\1\vphantom{l^1}\right) 
\, ,\, \vphantom{\left(\dot{f}\right)}\right. \\ 
&&\hspace{2.5cm}\left.\vphantom{\left(\dot{f}\right)}j^{-1}\left( 
\dot{W}+\left(\dot{A}_{n,0}-\dot{Y}_{n,0}\right)\1\right)\right 
\rangle_H\, d\sigma\, Q^{(m,r)}_{\sbnu}(dW)\, . 
\end{eqnarray*} 
Summarizing everything from (\ref{3.24}) on, we get 
\begin{eqnarray}\label{3.25}
&&\hspace{-.8cm}\int\vp\, ~ dQ^{(m)}_{\sbnu}\circ W^{n,-v}-\int\vp 
\,~ dQ^{(m)}_{\sbnu}  \nonumber \\
&&\hspace{.2cm}=\lim_{r\to\infty}\int\int_{\sigma =0}^v\left\langle 
D\vp\circ \left(W_{\cdot +\sigma}+(A_{n,\sigma}-Y_{n,\sigma})\1 
\vphantom{l^1}\right)\, ,\, \vphantom{\left(\dot{f}\right)}\right. 
\nonumber \\ 
&&\hspace{3.2cm}\left. \vphantom{\left(\dot{f}\right)}j^{-1}\left( 
\dot{W}+\left(\dot{A}_{n,0}-\dot{Y}_{n,0}\right)\1\right)\right 
\rangle_H\, d\sigma\, Q^{(m,r)}_{\sbnu}(dW) \nonumber \\ 
&&\hspace{.2cm}=\lim_{r\to\infty}\int\int_{\sigma=0}^v\left\langle 
\left(D\vp\circ\left(W_{\cdot +\sigma}+(A_{n,\sigma}-Y_{n,\sigma}) 
\1\vphantom{l^1}\right)\right)_1\, , \, \vphantom{\left(\dot{f} 
\right)} d\dot{W}\right\rangle_{L^2}\, d\sigma\, Q^{(m,r)}_{\sbnu} 
(dW)\nonumber \\ 
&&\hspace{.6cm}+\lim_{r\to\infty}\int\int_{\sigma=0}^v\left\langle 
\nabla_{W_0}\vp\left(W_{\cdot +\sigma}+(A_{n,\sigma}-Y_{n,\sigma}) 
\1\vphantom{l^1}\right)\, , \, \vphantom{\left(\dot{f}\right)}\dot 
{W}+\dot{A}_{n,0}-\dot{Y}_{n,0}\right\rangle_F\, d\sigma\, Q^{(m,r 
)}_{\sbnu}(dW)\, .  \nonumber \\ 
\end{eqnarray}

In order to apply Lemma \ref{Lemma3.2} (d) we note again that the 
weak mixed derivative defined in Subsection 2.3 is compatible with 
the integral of Definition \ref{Definition2.6} (a), cf. also Lemma 
\ref{Lemma3.1} (b) and its proof. For the following this means in 
particular that 
\begin{eqnarray*}
\left\langle (H_i,0)\, ,\, j^{-1}\left(\dot{W}+(\dot{A}_{n,0}- 
\dot{Y}_{n,0})\1\vphantom{l^1}\right)\right\rangle_H=\left\langle 
H_i,d\dot{W}\right\rangle_{L^2}\, . 
\end{eqnarray*}
By Lemma \ref{Lemma3.2} (d)-(f) and by (\ref{2.3}) as well as 
Lemma \ref{Lemma3.3} we verify now  
\begin{eqnarray}\label{3.26}
&&\hspace{-.5cm}\int\left\langle D\vp\circ\left(W_{\cdot +\sigma} 
+(A_{n,\sigma}-Y_{n,\sigma})\vphantom{l^1}\1\right)\, ,\, j^{-1} 
\left(\dot{W}+(\dot{A}_{n,0}-\dot{Y}_{n,0})\1\vphantom{l^1}\right) 
\right\rangle_H\, Q^{(m,r)}_{\sbnu}(dW) \nonumber \\ 
&&\hspace{.5cm}=\sum_{i\in I(m,r)}\int\left\langle D\vp\circ 
\left(W_{\cdot +\sigma}+(A_{n,\sigma}-Y_{n,\sigma})\vphantom{l^1} 
\1\right)\, ,\, (H_i,0)\right\rangle_H\cdot\left\langle H_i,d\dot 
{W}\right\rangle_{L^2}\, Q^{(m,r)}_{\sbnu}(dW) \nonumber \\ 
&&\hspace{1cm}+\int\left\langle D\vp\circ\left(W_{\cdot +\sigma} 
+(A_{n,\sigma}-Y_{n,\sigma})\vphantom{l^1}\1\right),\left(0,\dot 
{W}_0+\dot{A}_{n,0}-\dot{Y}_{n,0}\right)\right\rangle_H\, Q^{(m, 
r)}_{\sbnu}(dW)\vphantom{\int}\nonumber \\ 
&&\hspace{.5cm}=\int\vp\left(W_{\cdot +\sigma}+(A_{n,\sigma}-Y_{ 
n,\sigma})\1\vphantom{l^1}\right)\times\nonumber \\ 
&&\hspace{1cm}\times\left(\sum_{i\in I(m,r)}\left(\vphantom{\frac 
12}\left\langle H_i,d\dot{W}\right\rangle_{L^2}\cdot\delta (H_i 
,0)-\left\langle D\left\langle H_i,d\dot{W}\right\rangle_{L^2}, 
(H_i,0)\right\rangle_H\vphantom{\frac12}\right)\vphantom{\left( 
\sum_{i\in I(m,r)}\right)}\right.\nonumber \\ 
&&\hspace{1cm}\left.-\left\langle\frac{\nabla m(W_0)}{m(W_0)}\, 
, \, \dot{W}_0+\dot{A}_{n,0}-\dot{Y}_{n,0}\right\rangle_F-\left 
\langle {\sf e},\nabla_{d,W_0}\left(\dot{A}_{n,0}-\dot{Y}_{n,0} 
\right)\right\rangle_F\vphantom{\sum_{i\in I(m,r)}}\right)\, Q^{ 
(m,r)}_{\sbnu}(dW)\nonumber \\ 
&&\hspace{.5cm}=-\int\vp\left(W_{\cdot +\sigma}+(A_{n,\sigma}-Y_{ 
n,\sigma})\1\vphantom{l^1}\right)\left(\left\langle\frac{\nabla m 
(W_0)}{m(W_0)}\, ,\, \dot{W}_0+\dot{A}_{n,0}-\dot{Y}_{n,0}\right 
\rangle_F\right.\nonumber \\ 
&&\hspace{6.5cm}+\left.\left\langle {\sf e},\nabla_{d,W_0}\left( 
\dot{A}_{n,0}-\dot{Y}_{n,0}\right)\right\rangle_F\vphantom{\int} 
\right)\, Q^{(m,r)}_{\sbnu}(dW) \nonumber \\ 
&&\hspace{.1cm}\stack{r\to\infty}{\lra}-\int\vp\left(W_{\cdot + 
\sigma}+(A_{n,\sigma}-Y_{n,\sigma})\1\vphantom{l^1}\right)\left( 
\left\langle\frac{\nabla m(W_0)}{m(W_0)}\, , \, \dot{W}_0+\dot 
{A}_{n,0}-\dot{Y}_{n,0}\right\rangle_F\right.\nonumber \\ 
&&\hspace{6.5cm}+\left.\left\langle {\sf e},\nabla_{d,W_0}\left( 
\dot{A}_{n,0}-\dot{Y}_{n,0}\right)\right\rangle_F\vphantom{\int} 
\right)\, Q^{(m)}_{\sbnu}(dW)\, . 
\end{eqnarray}
Here the $\lim_{r\to\infty}$ follows from Lemma \ref{Lemma2.2} 
(c) with $\psi\equiv\psi(A_n-Y_n,\dot{A}_n-\dot{Y}_n,\nabla_{d,W_0} 
\dot{A}_n-\nabla_{d,W_0}\dot{Y}_n)$ and Lemma \ref{Lemma3.4} (b). In 
the following, let $(\cdot )^{-1}_c$ denote the coordinate wise and 
inverse and $\times_c$ denote the coordinate wise product. The 
$\log$ of a vector will denote the coordinate wise logarithm and the 
absolute value of a vector is going to be the vector of the absolute 
values. Recalling (\ref{3.24})-(\ref{3.26}), it holds that 
\begin{eqnarray}\label{3.27}
&&\hspace{-.5cm}\int\vp\, ~ dQ^{(m)}_{\sbnu}\circ W^{n,-v}-\int 
\vp\,~ dQ^{(m)}_{\sbnu}  \nonumber \\ 
&&\hspace{.5cm}=-\int_{\sigma=0}^v\int\vp(W^{n,\sigma})\left 
(\left\langle\frac{\nabla m(W_0)}{m(W_0)}\, , \, \dot{W}_0+\dot 
{A}_{n,0}-\dot{Y}_{n,0}\right\rangle_F\right.\nonumber \\ 
&&\hspace{4.5cm}+\left.\left\langle {\sf e},\nabla_{d,W_0}\left 
(\dot{A}_{n,0}-\dot{Y}_{n,0}\right)\right\rangle_F\vphantom 
{\int}\right)\, Q^{(m)}_{\sbnu}(dW)\, d\sigma \nonumber \\ 
&&\hspace{.5cm}=-\int_{\sigma=0}^v\int\vp\cdot\left(\left\langle 
\frac{\nabla m(W^{n,-\sigma}_0)}{m(W^{n,-\sigma}_0)}\, , \, \dot 
{X}_{n,-\sigma}-\dot{Y}_{n,-\sigma}\right\rangle_F\right. 
\nonumber \\ 
&&\hspace{4.5cm}+\left.\left\langle {\sf e},\nabla_{d,W_0}\left 
(\dot{A}_{n,0}-\dot{Y}_{n,0}\right)\left(W^{n,-\sigma}\right) 
\right\rangle_F\vphantom{\int}\right)\, dQ^{(m)}_{\sbnu}\circ 
W^{n,-\sigma}\, d\sigma \nonumber \\ 
&&\hspace{.5cm}=-\int_{\sigma=0}^v\int\vp\cdot\left\langle\frac{ 
\nabla m(W^{n,-\sigma}_0)}{m(W^{n,-\sigma}_0)}\, ,\, \dot{X}_{n, 
-\sigma}-\dot{Y}_{n,-\sigma}\right\rangle_F\, dQ^{(m)}_{\sbnu} 
\circ W^{n,-\sigma}\, ds\nonumber \\ 
&&\hspace{1.0cm}-\int_{\sigma=0}^v\int\vp\cdot\left\langle {\sf 
e},\left(\nabla_{d,W_0}\left(A_{n,-\sigma}-Y_{n,-\sigma}\right) 
\vphantom{\dot{f}}\right)^{\cdot}\times_c\right.\nonumber \\ 
&&\hspace{4.5cm}\left.\times_c\left({\sf e}+\nabla_{d,W_0}\left( 
A_{n,-\sigma}-Y_{n,-\sigma}\right)\vphantom{\dot{f}}\right)^{-1 
}_c\right\rangle_F\, dQ^{(m)}_{\sbnu}\circ W^{n,-\sigma}\, d\sigma 
\nonumber \\ 
&&\hspace{.5cm}=-\int_{\sigma=0}^v\int\vp\cdot\left\langle\frac{ 
\nabla m(W^{n,-\sigma}_0)}{m(W^{n,-\sigma}_0)}\, ,\, \dot{X}_{n, 
-\sigma}-\dot{Y}_{n,-\sigma}\right\rangle_F\, dQ^{(m)}_{\sbnu} 
\circ W^{n,-\sigma}\, ds\nonumber \\ 
&&\hspace{1.0cm}-\int_{\sigma=0}^v\int\vp\cdot\left\langle {\sf 
e},\left(\log\left|{\sf e}+\nabla_{d,W_0}\left(A_{n,-\sigma}-Y_{ 
n,-\sigma}\right)\vphantom{\dot{f}}\right|\right)^{\cdot}\right 
\rangle_F\, dQ^{(m)}_{\sbnu}\circ W^{n,-\sigma}\, d\sigma\qquad 
\qquad
\end{eqnarray}
where for the second equality sign we have applied the flow property 
of $W^{n,\, \cdot}$ and in particular $W^{n,0}=W$. Furthermore, we 
have made use of $X_{n,\cdot-\sigma}(W)=X_n(W^{n,-\sigma})$, $Y_{n, 
\cdot-\sigma}(W)=Y_n(W^{n,-\sigma})$, and $\dot{A}_{n,\cdot-\sigma} 
(W)=\dot{A}_n(W^{n,-\sigma})$, cf. Lemma \ref{Lemma3.4} (vi). 
\medskip

\nid 
{\it Step 2 } Let ${\cal F}^{(m)}$ denote the $\sigma$-algebra on 
$\Omega$ which is generated by the projection operator $\pi_m$. In 
this step, we shall specify a subset of $[0,t]$ such that, for $v$ 
belonging to this subset, the Radon-Nikodym derivative $dQ^{(m)}_{ 
\sbnu}\circ W^{n,-v}/dQ^{(m)}_{\sbnu}$ exists on $(\Omega,{\cal F 
}^{(m)})$. We will also derive a representation of the density. 
Let us keep up with the following form of (\ref{3.27}), 
\begin{eqnarray*}
&&\hspace{-.5cm}\int\vp\, ~ dQ^{(m)}_{\sbnu}\circ W^{n,-v}-\int\vp 
\,~ dQ^{(m)}_{\sbnu} \\ 
&&\hspace{.5cm}=-\int_{\sigma=0}^v\int\vp\cdot\left(\left(\left 
\langle\frac{\nabla m(W_0)}{m(W_0)}\, , \, \dot{X}_{n,0}-\dot{Y}_{ 
n,0}\right\rangle_F\right.\right.\nonumber \\ 
&&\hspace{4.5cm}+\left.\left.\left\langle {\sf e},\nabla_{d,W_0} 
\left(\dot{A}_{n,0}-\dot{Y}_{n,0}\right)\right\rangle_F\vphantom 
{\int}\right)\, dQ^{(m)}_{\sbnu}\right)\circ W^{n,-\sigma}\, d 
\sigma\, .
\end{eqnarray*}
We recall that by $\dot{X}_{n,0}-\dot{Y}_{n,0}=\dot{W}+\dot{A}_{n,0} 
-\dot{Y}_{n,0}$, the definitions of the density $m$ and the exponent 
$q$ in Subsection 1.2, as well as by Lemma \ref{Lemma3.4} (iii) and 
(iv), 
\begin{eqnarray*}
\left(\left\langle\frac{\nabla m(W_0)}{m(W_0)}\, , \, \dot{X}_{n,0}- 
\dot{Y}_{n,0}\right\rangle_F+\left\langle {\sf e},\nabla_{d,W_0}\left 
(\dot{A}_{n,0}-\dot{Y}_{n,0}\right)\right\rangle_F\vphantom{\int} 
\right)\, dQ^{(m)}_{\sbnu} 
\end{eqnarray*}
is a finite signed measure. Using monotone convergence it turns out 
that (\ref{3.27}) holds for all indicator functions $\vp(W)=\1_{\{ 
W_{t_1}\in A_1,\ldots ,W_{t_l}\in A_l\}}$ where $A_i\in {\cal B}(F 
)$ are open sets, $t_i\in\{z\cdot t/2^m:z\in {\Bbb Z}\}$, $i\in\{1, 
\ldots ,l\}$, $l\in {\Bbb N}$. Using Sierpinski's monotone class 
theorem it follows now that (\ref{3.27}) even holds for all $\vp\in 
B_b(\Omega,{\cal F}^{(m)})$, the set of all bounded and with respect 
to ${\cal F}^{(m)}$ measurable functions on $\Omega$. 

Let us now establish a version of (\ref{3.27}) for time dependent 
test functions $\psi$ which are defined everywhere on $[0,t]\times 
\Omega$ and satisfy the following. 
\begin{eqnarray}\label{3.288} 
\psi(v,\, \cdot \, )\in B_b(\Omega,{\cal F}^{(m)})\, , \quad v\in 
[0,t]\, . 
\end{eqnarray}
There exists $(\partial\psi/\partial v)(v,\, \cdot\, )$ such that 
\begin{eqnarray}\label{3.289} 
&&\hspace{-.5cm}\lim_{h\to 0}\frac{\psi(v+h,\, \cdot\, )\circ W^{n, 
v+h}-\psi(v,\, \cdot\, )\circ W^{n,v+h}}{h}\nonumber \\ 
&&\hspace{.5cm}=\lim_{h\to 0}\frac{\psi(v+h,\, \cdot\, )\circ W^{n, 
v}-\psi(v,\, \cdot\, )\circ W^{n,v}}{h}\nonumber \\ 
&&\hspace{.5cm}=\frac{\partial\psi}{\partial v}(v,\, \cdot\, )\circ 
W^{n,v}\vphantom{\frac{W^h}{W_h}}\, , \quad v\in [0,t],\quad\mbox{ 
\rm in }\ L^1(\Omega,Q^{(m)}_{\sbnu}). 
\end{eqnarray}
It follows that, for all $v\in [0,t]$, there exists the limit 
\begin{eqnarray*} 
&&\hspace{-.5cm}\frac{d}{dv}\int\psi(v,\, \cdot\, )\, dQ^{(m)}_{\sbnu} 
\circ W^{n,-v} \\ 
&&\hspace{.5cm}=\lim_{h\to 0}\frac{\D\int\psi(v+h,\, \cdot\, )\, dQ^{ 
(m)}_{\sbnu}\circ W^{n,-v-h}-\int\psi(v,\, \cdot\, )\, dQ^{(m)}_{\sbnu} 
\circ W^{n,-v}}{h\vphantom{\dot{f}}}   
\end{eqnarray*}
whenever $\psi$ satisfies (\ref{3.288}) and (\ref{3.289}) and that 
\begin{eqnarray*}
&&\hspace{-.5cm}\frac{d}{dv}\int\psi(v,\, \cdot\, )\, dQ^{(m)}_{\sbnu} 
\circ W^{n,-v}-\int\frac{\partial\psi}{\partial v}(v,\, \cdot\, )\, d 
Q^{(m)}_{\sbnu}\circ W^{n,-v} \\ 
&&\hspace{.5cm}=\frac{d}{dv}\int\psi(v,\, \cdot\, )\, dQ^{(m)}_{\sbnu} 
\circ W^{n,-v}-\lim_{h\to 0}\int\frac{\psi(v+h,\, \cdot\, )-\psi(v,\, 
\cdot\, )}{h}\, dQ^{(m)}_{\sbnu}\circ W^{n,-v} \\ 
&&\hspace{.5cm}=\lim_{h\to 0}\frac1h\left(\int\psi(v+h,\, \cdot\, )\, d 
Q^{(m)}_{\sbnu}\circ W^{n,-v-h}-\int\psi(v+h,\, \cdot\, )\, dQ^{(m)}_{ 
\sbnu}\circ W^{n,-v}\right) \\ 
&&\hspace{.5cm}=\lim_{h\to 0}\frac1h\left(\int\psi(v,\, \cdot\, )\, d 
Q^{(m)}_{\sbnu}\circ W^{n,-v-h}-\int\psi(v,\, \cdot\, )\, dQ^{(m)}_{ 
\sbnu}\circ W^{n,-v}\right) \\ 
&&\hspace{1cm}+\lim_{h\to 0}\int\frac{\psi(v+h,\, \cdot\, )-\psi 
(v,\, \cdot\, )}{h}\, dQ^{(m)}_{\sbnu}\circ W^{n,-v-h} \\ 
&&\hspace{3cm}-\lim_{h\to 0}\int\frac{\psi(v+h,\, \cdot\, )-\psi(v,\, 
\cdot\, )}{h}\, dQ^{(m)}_{\sbnu}\circ W^{n,-v} \\ 
&&\hspace{.5cm}=\lim_{h\to 0}\frac1h\left(\int\psi(v,\, \cdot\, )\, d 
Q^{(m)}_{\sbnu}\circ W^{n,-v-h}-\int\psi(v,\, \cdot\, )\, dQ^{(m)}_{ 
\sbnu}\circ W^{n,-v}\right) \\ 
&&\hspace{.5cm}=-\int\psi(v,\, \cdot\, )\cdot\left(\left\langle\frac{ 
\nabla m(W^{n,-v}_0)}{m(W^{n,-v}_0)}\, , \, \dot{X}_{n,-v}-\dot{Y}_{n 
,-v}\right\rangle_F\right. \\ 
&&\hspace{2.5cm}+\left.\left\langle {\sf e},\nabla_{d,W_0}\left(\dot{A 
}_{n,0}-\dot{Y}_{n,0}\right)\left(W^{n,-v}\right)\right\rangle_F 
\vphantom{\int}\right)\, dQ^{(m)}_{\sbnu}\circ W^{n,-v}\, , \quad v\in 
[0,t], 
\end{eqnarray*}
the last line by (\ref{3.288}) of the present step and (\ref{3.27}). 
In addition, let us require that the subsequent condition on $\psi$ 
is satisfied. 
\begin{eqnarray}\label{3.290} 
\int_0^t\left|\frac{\partial\psi}{\partial v}(v,\, \cdot\, )\circ 
W^{n,v}\right|\, dv\in L^1(\Omega,Q^{(m)}_{\sbnu})\, . 
\end{eqnarray} 
Now the last equation has an integral version 
\begin{eqnarray}\label{3.291} 
&&\hspace{-.5cm}\int\psi(v,\, \cdot\, )\, dQ^{(m)}_{\sbnu}\circ 
W^{n,-v}-\int\psi(0,\, \cdot\, )\, dQ^{(m)}_{\sbnu}-\int_{\sigma= 
0}^v\int\frac{\partial\psi}{\partial\sigma}(\sigma,\, \cdot\, )\, 
dQ^{(m)}_{\sbnu}\circ W^{n,-\sigma}\, d\sigma\nonumber \\ 
&&\hspace{.5cm}=-\int_{\sigma=0}^v\int\psi(\sigma,\, \cdot\, )\cdot 
\left(\left\langle\frac{\nabla m(W^{n,-\sigma}_0)}{m(W^{n,-\sigma 
}_0)}\, ,\, \dot{X}_{n,-\sigma}-\dot{Y}_{n,-\sigma}\right\rangle_F 
\right.\nonumber \\ 
&&\hspace{4.5cm}+\left.\left\langle {\sf e},\nabla_{d,W_0}\left( 
\dot{A}_{n,0}-\dot{Y}_{n,0}\right)\left(W^{n,-\sigma}\right)\right 
\rangle_F\vphantom{\int}\right)\, dQ^{(m)}_{\sbnu}\circ W^{n,- 
\sigma}\, d\sigma\, ,\nonumber \\ 
\end{eqnarray}
$v\in [0,t]$. Let $\Psi(v,\, \cdot\, )\in B_b(\Omega,{\cal F}^{(m)} 
)$, $v\in [0,t]$. Until further specification, assume that the 
function 
\begin{eqnarray}\label{3.2910} 
\hat{\psi}(v,\, \cdot\, ):=\exp\left\{\int_{\sigma=0}^v\Psi(\sigma, 
\, \cdot\, )\, d\sigma\right\}\cdot\1_G\, , \quad v\in [0,t], 
\end{eqnarray}
is for any $G\in {\cal F}^{(m)}$ and $c>0$ a function satisfying 
(\ref{3.288})-(\ref{3.290}). Thus $\hat{\psi}$ satisfies also 
(\ref{3.291}). We obtain from (\ref{3.291}) 
\begin{eqnarray}\label{3.2911}  
&&\hspace{-.5cm}\int\hat{\psi}(v,\, \cdot\, )\, dQ^{(m)}_{\sbnu} 
\circ W^{n,-v}-Q^{(m)}_{\sbnu}(G)\nonumber \\ 
&&\hspace{.5cm}=\int\hat{\psi}(v,\, \cdot\, )\, dQ^{(m)}_{\sbnu} 
\circ W^{n,-v}-\int\hat{\psi}(0,\, \cdot\, )\, dQ^{(m)}_{\sbnu} 
\nonumber \\ 
&&\hspace{.5cm}=\int_{\sigma=0}^v\int\hat{\psi}(\sigma,\, \cdot\, 
)\cdot\left(\Psi(\sigma,\, \cdot\, )-\left\langle\frac{\nabla m 
(W^{n,-\sigma}_0)}{m(W^{n,-\sigma}_0)}\, ,\, \dot{X}_{n,-\sigma}- 
\dot{Y}_{n,-\sigma}\right\rangle_F\right.\nonumber \\ 
&&\hspace{2cm}-\left.\left\langle {\sf e},\nabla_{d,W_0}\left(\dot 
{A}_{n,0}-\dot{Y}_{n,0}\right)\left(W^{n,-\sigma}\right)\right 
\rangle_F\vphantom{\int}\right)\, dQ^{(m)}_{\sbnu}\circ W^{n,- 
\sigma}\, d\sigma\, ,\quad v\in [0,t]\, . \nonumber \\ 
\end{eqnarray}

The following observation is a consequence of the law of iterated 
logarithm and elementary geometric properties of piecewise linear 
functions with equidistant linearity intervals. For every $v\in 
[0,t]\setminus\{z\cdot t/2^m+t/2^{m+1}:z\in {\Bbb Z}\}$ and every 
$W\in\Omega$ with 
\begin{eqnarray*}
\lim_{z\to\pm\infty\atop z\in {\Bbb Z}}\frac1z\log\left|W_{z\cdot 
t/2^m}\right|=0 
\end{eqnarray*}
there exists exactly one $V\in\{\left(\pi_mU\right)_{\cdot +v}:U 
\in\Omega\}$ with $V_s=W_s$ for all $s\in\{z\cdot t/2^m:z\in {\Bbb 
Z}\}$, which also means that $\pi_mW=\pi_mV$, such that 
\begin{eqnarray*}
\lim_{z\to\pm\infty\atop z\in {\Bbb Z}}\frac1z\log\left|V_{z\cdot 
t/2^m+v}\right|=0\, .  
\end{eqnarray*}
In this notation we shall write $V=\pi_{m;v}W$. Recalling $W^{n,-v} 
=W_{\cdot -v}+\left(A_{-v}-Y_{-v}\right)\1$ it is obvious that 
\begin{eqnarray*}
\ \ \ Q_{\sbnu}^{(m)}\circ W^{n,-v}\left(V\in\{\pi_{m;v}W:W\in 
\Omega\}:\lim_{z\to\pm\infty\atop z\in {\Bbb Z}}\frac1z\log\left| 
V_{z\cdot t/2^m+v}\right|=0\right)=1\, , \quad v\in [0,t].  
\end{eqnarray*}
In other words, with $E_{\sbnu}^{(m)}\circ W^{n,-v}$ denoting the 
expectation with respect to the measure $Q_{\sbnu}^{(m)}\circ W^{n 
,-v}$, $v\in [0,t]$, and ${\cal F}^{(m)}\subset{\cal F}$ being the 
$\sigma$-algebra generated by $\pi_m$, it holds that  
\begin{eqnarray}\label{3.2930} 
E_{\sbnu}^{(m)}\circ W^{n,-v}\left.\left(\vphantom{l^1}\xi\right| 
{\cal F}^{(m)}\right)=\xi\circ\pi_{m;v}\, , \quad \xi\in L^1(\Omega 
,Q_{\sbnu}^{(m)}\circ W^{n,-v})\, , \vphantom{\int}
\end{eqnarray}
$v\in [0,t]\setminus\{z\cdot t/2^m+t/2^{m+1}:z\in {\Bbb Z}\}$. 
Furthermore, for the same $v$ we have 
\begin{eqnarray}\label{3.2931} 
\xi\circ\pi_{m;v}=\xi\, , \quad Q_{\sbnu}^{(m)}\circ W^{n,-v}\mbox{ 
\rm -a.e.}\vphantom{\int}
\end{eqnarray}

Next, let us specify $\hat{\psi}$ in (\ref{3.2910}). For $c>0$ we 
choose  
\begin{eqnarray}\label{3.2939} 
\Phi_c(v,\, \cdot\, )=\left(\left\langle\frac{\nabla m(W^{n,-v}_0)} 
{m(W^{n,-v}_0)}\, ,\, \dot{X}_{n,-v}-\dot{Y}_{n,-v}\right\rangle_F+ 
\left\langle{\sf e},\nabla_{d,W_0}\left(\dot{A}_{n,0}-\dot{Y}_{n,0} 
\right)\left(W^{n,-v}\right)\right\rangle_F\right)\wedge c\nonumber 
\\ 
\end{eqnarray}
and define 
\begin{eqnarray}\label{3.2933} 
\Psi_c(v,\, \cdot\, ):=\Phi_c(v,\, \cdot\, )\circ \pi_{m;v}  
\end{eqnarray}
as well as  
\begin{eqnarray*} 
\hat{\psi}_c(v,\, \cdot\, ):=\exp\left\{\int_{\sigma=0}^v\Psi_c( 
\sigma,\, \cdot\, )\, d\sigma\right\}\cdot\1_G\, ,  
\end{eqnarray*}
$v\in [0,t]\setminus\{z\cdot t/2^m+t/2^{m+1}:z\in {\Bbb Z}\}$. By 
(\ref{3.2930}) we see that $\hat{\psi}_c$ has in particular the form 
of $\hat{\psi}$ in (\ref{3.2910}). We observe that 
\begin{eqnarray}\label{3.294} 
\hat{\psi}_c(v,\, \cdot\, )&&\hspace{-.5cm}=\exp\left\{\int_{\sigma 
=0}^v\Phi_c(\sigma,\, \cdot\, )\, d\sigma\right\}\circ\pi_{m;v}\cdot 
\1_G\, , 
\end{eqnarray}
$v\in [0,t]\setminus\{z\cdot t/2^m+t/2^{m+1}:z\in {\Bbb Z}\}$. In 
order to apply (\ref{3.2911}) to $\hat{\psi}_c$ we have to verify 
(\ref{3.288})-(\ref{3.290}). Relation (\ref{3.288}) is clear by 
definition. Relation (\ref{3.289}) is on the one hand a consequence 
of the fact that for $k(h)$ either $h$ or $0$ we have 
\begin{eqnarray*} 
&&\hspace{-.5cm}\left|\frac{\hat{\psi}_c(v+h,\, \cdot\, )\circ W^{n, 
v+k(h)}-\hat{\psi}_c(v,\, \cdot\, )\circ W^{n,v+k(h)}}{h}\vphantom{ 
\left(\int_b^b\right)}\right| \\ 
&&\hspace{.5cm}=\left|\frac1h\int_{a=0}^h\Psi_c(v+a,W^{n,v+k(a)}) 
\cdot\exp\left\{\int_{\sigma=0}^{v+a}\Psi_c(\sigma,W^{n,v+k(a)})\, d 
\sigma\right\}\, da\cdot\1_G\right| \\ 
&&\hspace{.2cm}\stack{h\to 0}{\lra}\Psi_c(v,W^{n,v})\cdot\exp\left\{ 
\int_{\sigma=0}^v\Psi_c(\sigma,W^{n,v})\, d\sigma\right\}\cdot\1_G 
\quad Q_{\sbnu}^{(m)}\mbox{\rm -a.e.}
\end{eqnarray*}
by (\ref{3.2939}) and (\ref{3.2933}) as well as Lemma \ref{Lemma3.4} 
(i), (v), and (vi). On the other hand, 
\begin{eqnarray*} 
&&\hspace{-.5cm}\left\|\Psi_c(v+a,W^{n,v+k(a)})\cdot\exp\left\{\int_{ 
\sigma=0}^{v+a}\Psi_c(\sigma,W^{n,v+k(a)})\, d\sigma\right\}\right 
\|_{L^1(\Omega,Q_{\sbnu}^{(m)})}\le c\cdot e^{c(v+\ve)}
\end{eqnarray*}
independent of $a\in (-\ve,\ve)$ for some $\ve>0$, $v\in [0,t] 
\setminus\{z\cdot t/2^m+t/2^{m+1}:z\in {\Bbb Z}\}$, which completes 
the verification of (\ref{3.289}). The last relation yields also 
(\ref{3.290}) for $\hat{\psi}_c$. 
\medskip 

It follows now from (\ref{3.2911}) as well as (\ref{3.2931}) and 
(\ref{3.294}) that 
\begin{eqnarray}\label{3.295}
&&\hspace{-.5cm}\int_G\exp\left\{\int_{\sigma=0}^v\Phi_c(\sigma,\, 
\cdot\, )\, d\sigma\right\}\, dQ^{(m)}_{\sbnu}\circ W^{n,-v}-Q^{(m) 
}_{\sbnu}(G)\nonumber \\ 
&&\hspace{.5cm}=\int\hat{\psi}_c(v,\, \cdot\, )\, dQ^{(m)}_{\sbnu} 
\circ W^{n,-v}-Q^{(m)}_{\sbnu}(G)\nonumber \\ 
&&\hspace{.5cm}=\int_{\sigma=0}^v\int\hat{\psi}_c(\sigma,\, \cdot\, 
)\cdot\left(\Psi_c(\sigma,\, \cdot\, )-\left\langle\frac{\nabla m( 
W^{n,-\sigma}_0)}{m(W^{n,-\sigma}_0)}\, ,\, \dot{X}_{n,-\sigma}-\dot 
{Y}_{n,-\sigma}\right\rangle_F\right.\nonumber \\ 
&&\hspace{2cm}-\left.\left\langle {\sf e},\nabla_{d,W_0}\left(\dot{A 
}_{n,0}-\dot{Y}_{n,0}\right)\left(W^{n,-\sigma}\right)\right 
\rangle_F\vphantom{\int}\right)\, dQ^{(m)}_{\sbnu}\circ W^{n,-\sigma 
}\, d\sigma\nonumber \\ 
&&\hspace{.5cm}=\int_{\sigma=0}^v\int\hat{\psi}_c(\sigma,\, \cdot\, 
)\cdot\left(\Phi_c(\sigma,\, \cdot\, )-\left\langle\frac{\nabla m( 
W^{n,-\sigma}_0)}{m(W^{n,-\sigma}_0)}\, ,\, \dot{X}_{n,-\sigma}-\dot 
{Y}_{n,-\sigma}\right\rangle_F\right.\nonumber \\ 
&&\hspace{2cm}-\left.\left\langle {\sf e},\nabla_{d,W_0}\left(\dot 
{A}_{n,0}-\dot{Y}_{n,0}\right)\left(W^{n,-\sigma}\right)\right 
\rangle_F\vphantom{\int}\right)\, dQ^{(m)}_{\sbnu}\circ W^{n,-\sigma 
}\, d\sigma\, , 
\end{eqnarray}
$v\in [0,t]\setminus\{z\cdot t/2^m+t/2^{m+1}:z\in {\Bbb Z}\}$. In 
particular, (\ref{3.295}) together with (\ref{3.2939}) means that 
\begin{eqnarray*} 
\int_G\exp\left\{\int_{\sigma=0}^v\Phi_c(\sigma,\, \cdot\, )\, d\sigma 
\right\}\, dQ^{(m)}_{\sbnu}\circ W^{n,-v}\le Q^{(m)}_{\sbnu}(G)\, , 
\quad G\in {\cal F}^{(m)},  
\end{eqnarray*}
$v\in [0,t]\setminus\{z\cdot t/2^m+t/2^{m+1}:z\in {\Bbb Z}\}$. 
Extending definition (\ref{3.2939}) to $c=\infty$, monotone convergence 
yields for the same $v$ 
\begin{eqnarray}\label{3.296} 
\int_G\exp\left\{\int_{\sigma=0}^v\Phi_\infty(\sigma,\, \cdot\, )\, d 
\sigma\right\}\, dQ^{(m)}_{\sbnu}\circ W^{n,-v}\le Q^{(m)}_{\sbnu}(G) 
\, , \quad G\in {\cal F}^{(m)}. 
\end{eqnarray}
Plugging (\ref{3.294}) into (\ref{3.295}), taking into consideration 
(\ref{3.2931}), and plugging then (\ref{3.296}) into (\ref{3.295}) we 
obtain 
\begin{eqnarray}\label{3.297}
0&&\hspace{-.5cm}\le Q^{(m)}_{\sbnu}(G)-\int_G\exp\left\{\int_{\sigma 
=0}^v\Phi_c(\sigma,\, \cdot\, )\, d\sigma\right\}\, dQ^{(m)}_{\sbnu} 
\circ W^{n,-v}\nonumber \\ 
&&\hspace{-.5cm}=\int_{\sigma=0}^v\int_G\exp\left\{\int_{u=0}^\sigma 
\Phi_c(u,\, \cdot\, )\, du\right\}\cdot\left(\Phi_\infty(\sigma,\,  
\cdot\, )\vphantom{l^1}-\Phi_c(\sigma,\, \cdot\, )\vphantom{\int}\right) 
\, dQ^{(m)}_{\sbnu}\circ W^{n,-\sigma}\, d\sigma\nonumber \\ 
&&\hspace{-.5cm}\le\int_{\sigma=0}^v\int_G\exp\left\{\int_{u=0}^\sigma 
\left(\left(\Phi_\infty (u,\, \cdot\, )-\Phi_c(u,\, \cdot\, )\vphantom 
{l^1}\right)\vee c\vphantom{\dot{f}}\right)\, du\right\}\times\nonumber 
 \\ 
&&\hspace{5.5cm}\times\left(\left(\Phi_\infty (\sigma,\, \cdot\, )- 
\Phi_c(\sigma,\, \cdot\, )\vphantom{l^1}\right)\vee c\vphantom{\dot{f}} 
\right)\, dQ^{(m)}_{\sbnu}\, d\sigma\nonumber \\ 
&&\hspace{-.5cm}=\int_G\left(\exp\left\{\int_{\sigma=0}^v\left(\left( 
\Phi_\infty(\sigma,\, \cdot\, )-\Phi_c(\sigma,\, \cdot\, )\vphantom{l^1 
}\right)\vee c\vphantom{\dot{f}}\right)\, d\sigma\right\}-1\right)\, d 
Q^{(m)}_{\sbnu}\, ,\quad G\in {\cal F}^{(m)}, \qquad
\end{eqnarray}
$v\in [0,t]\setminus\{z\cdot t/2^m+t/2^{m+1}:z\in {\Bbb Z}\}$. Letting 
$c\to\infty$ in (\ref{3.297}) and recalling (\ref{3.296}), it follows 
again by monotone convergence that for such $v$ 
\begin{eqnarray*}
\int_G\exp\left\{\int_{\sigma=0}^v\Phi_\infty(\sigma,\, \cdot\, )\, d 
\sigma\right\}\, dQ^{(m)}_{\sbnu}\circ W^{n,-v}=Q^{(m)}_{\sbnu}(G)\, , 
\quad G\in {\cal F}^{(m)}. 
\end{eqnarray*}

Noting that 
\begin{eqnarray*} 
&&\hspace{-.5cm}-\int_{\sigma=0}^v\left\langle\frac{\nabla m\left(X_{n, 
-\sigma}-Y_{n,-\sigma}\right)}{m\left(X_{n,-\sigma}-Y_{n,-\sigma}\right) 
}\, ,\, d(X_{n,-\sigma}-Y_{n,-\sigma})\right\rangle_F\nonumber \\ 
&&\hspace{.5cm}=\log m\left(X_{n,v}-Y_{n,v}\right)-\log m\left(X_{n,0} 
-Y_{n,0}\right)\vphantom{\int}
\end{eqnarray*} 
is path wise just a Lebesgue-Stieltjes integral with no jumps in the 
integrator and recalling $X_{n,0}-Y_{n,0}=W_0+A_{n,0}-Y_{n,0}=W_0$ and 
Lemma \ref{Lemma3.4} (v), we have shown that on $(\Omega ,{\cal F}^{(m) 
})$ the Radon-Nikodym derivative $dQ^{(m)}_{\sbnu}\circ W^{n,-v}/dQ^{( 
m)}_{\sbnu}$ exists and that it coincides with 
\begin{eqnarray*}
\omega^{(m)}_{n,-v}&&\hspace{-.5cm}:=\frac{dQ^{(m)}_{\sbnu}\circ W^{ 
n,-v}}{dQ^{(m)}_{\sbnu}}=\exp\left\{-\int_{\sigma=0}^v\Phi_\infty( 
\sigma,\, \cdot\, )\, d\sigma\right\} \\ 
&&\hspace{-.5cm}=\exp\left\{-\int_{\sigma=0}^v\left\langle\frac{\nabla 
m(X_{n,-\sigma}-Y_{n,-\sigma})}{m(X_{n,-\sigma}-Y_{n,-\sigma})}\, ,\, 
\dot{X}_{n,-\sigma}-\dot{Y}_{n,-\sigma}\right\rangle_F d\sigma\right 
\}\times \\ 
&&\hspace{3cm}\times\exp\left\{-\int_{\sigma=0}^v\left\langle {\sf e}, 
\nabla_{d,W_0}\left(\dot{A}_{n,0}-\dot{Y}_{n,0}\right)(W^{n,-\sigma}) 
\right\rangle_F d\sigma\right\} \\ 
&&\hspace{-.5cm}=\frac{m(X_{n,-v}-Y_{n,-v})}{m(W_0)}\cdot\exp\left\{ 
-\int_{\sigma=0}^v\left\langle{\sf e}\, , \, \left(\log\left|{\sf e}+ 
\nabla_{d,W_0}\left(A_{n,-\sigma}-Y_{n,-\sigma}\right)\vphantom{\dot{ 
f}}\right|\right)^{\cdot}\right\rangle_F\, d\sigma\right\} \\ 
&&\hspace{-.5cm}=\frac{m(X_{n,-v}-Y_{n,-v})}{m(W_0)}\cdot\exp\left\{ 
\left\langle {\sf e},\log\left|{\sf e}+\nabla_{d,W_0}\left( A_{n,-v}- 
Y_{n,-v}\right)\vphantom{\dot{f}}\right|\right\rangle_F\right\} \\ 
&&\hspace{-.5cm}=\frac{m(X_{n,-v}-Y_{n,-v})}{m(W_0)}\cdot\prod_{i=1}^{ 
n\cdot d}\left|{\sf e}+\nabla_{d,W_0}\left( A_{n,-v}-Y_{n,-v}\vphantom 
{\dot{f}}\right)\right|_i
\end{eqnarray*}
$Q^{(m)}_{\sbnu}$-a.e., $v\in [0,t]\setminus\{z\cdot t/2^m+t/2^{m+1}: 
z\in {\Bbb Z}\}$. Recall also (\ref{3.27}). 
\medskip

\nid 
{\it Step 3 } Now we are interested in the weak limit $\lim_{n\to\infty} 
\omega^{(m)}_{n,-v}Q^{(m)}_{\sbnu}$ on $(\Omega,{\cal F}^{(m)})$, $v\in 
[0,t]\setminus\{z\cdot t/2^m+t/2^{m+1}:z\in {\Bbb Z}\}$. By Lemma 
\ref{Lemma3.4} (iv), (c), and (a) together with the fact that $Y,X$ jump 
for fixed $v$ only on a $Q_{\sbnu}^{(m)}$-zero set, cf. conditions (2) and 
(3) of Section 1, the boundedness and continuity of the density $m$, and 
again $A_{n,0}=Y_{n,0}$, we have 
\begin{eqnarray}\label{3.30}
\lim_{n\to\infty}\omega^{(m)}_{n,-v}&&\hspace{-.5cm}=\lim_{n\to\infty} 
\frac{m(X_{n,-v}-Y_{n,-v})}{m(W_0)}\cdot\prod_{i=1}^{n\cdot d}\left|{\sf 
e}+\nabla_{d,W_0}\left( A_{n,-v}-Y_{n,-v}\vphantom{\dot{f}}\right)\right 
|_i\nonumber \\ 
&&\hspace{-.5cm}=\frac{m(X_{-v}-Y_{-v})}{m(W_0)}\cdot\prod_{i=1}^{n\cdot 
d}\left|{\sf e}+\nabla_{d,W_0}A_{-v}-\nabla_{d,W_0}Y_{-v}\vphantom{\dot{ 
f}}\right|_i
\end{eqnarray} 
$Q^{(m)}_{\sbnu}$-a.e. on $(\Omega,{\cal F}^{(m)})$; recall also Remark 
(3) of this section. Let us turn to test functions $\vp$ of the form 
$\vp(W)=f(W_{t_1},\ldots ,W_{t_l})$ where $f\in C_b(F^l)$, and $t_i\in 
\{z\cdot t/2^m:z\in {\Bbb Z}\}$, $i\in \{1,\ldots ,l\}$, $l\in {\Bbb 
N}$. For $v\in [0,t]\setminus\{z\cdot t/2^m+t/2^{m+1}:z\in {\Bbb Z}\}$ 
we have 
\begin{eqnarray*} 
&&\hspace{-.5cm}\int\vp\,~ dQ^{(m)}_{\sbnu}\circ W^{-v} \\ 
&&\hspace{.5cm}=\int\vp(W^v)\, Q^{(m)}_{\sbnu}(dW) \\ 
&&\hspace{.5cm}=\int\vp\left(W_{\cdot +v}+(A_v-Y_v)\1\vphantom{l^1} 
\right)\, Q^{(m)}_{\sbnu}(dW) \\
&&\hspace{.5cm}=\lim_{n\to\infty}\int\vp\left(W_{\cdot +v}+(A_{n,v} 
-Y_{n,v})\1\vphantom{l^1}\right)\, Q^{(m)}_{\sbnu}(dW) \\
&&\hspace{.5cm}=\lim_{n\to\infty}\int\vp\,~ dQ^{(m)}_{\sbnu}\circ 
W^{n,-v} \\ 
&&\hspace{.5cm}=\lim_{n\to\infty}\int\vp\, \omega^{(m)}_{n,-v}\,~d 
Q^{(m)}_{\sbnu} \\ 
&&\hspace{.5cm}=\lim_{n\to\infty}\int\vp\, m(X_{n,-v}-Y_{n,-v})\cdot 
\prod_{i=1}^{n\cdot d}\left|{\sf e}+\nabla_{d,W_0}A_{n,-v}-\nabla_{ 
d,W_0}Y_{n,-v}\vphantom{\dot{f}}\right|_i\,~dQ^{(m)}_{\lambda_F} \\ 
&&\hspace{.5cm}=\int\vp\frac{m(X_{-v}-Y_{-v})}{m(W_0)}\cdot\prod_{ 
i=1}^{n\cdot d}\left|{\sf e}+\nabla_{d,W_0}A_{-v}-\nabla_{d,W_0} 
Y_{-v}\vphantom{\dot{f}}\right|_i\,~ dQ^{(m)}_{\sbnu}  
\end{eqnarray*} 
where the limit $\lim_{n\to\infty}$ in the fourth line is a 
consequence of Lemma \ref{Lemma3.4} (a) together with the fact that 
$A,Y$ jump for fixed $v$ only on a $Q_{\sbnu}^{(m)}$-zero set, cf. 
conditions (2) and (3) of Section 1. The last line follows from 
(\ref{3.30}) and Lemma \ref{Lemma3.4} (c). For $v\in [0,t]\setminus 
\{z\cdot t/2^m+t/2^{m+1}:z\in {\Bbb Z}\}$ we have shown 
\begin{eqnarray}\label{3.32}
&&\hspace{-.5cm}\omega^{(m)}_{-v}(W):=\frac{dQ^{(m)}_{\sbnu}\circ 
W^{-v}}{dQ^{(m)}_{\sbnu}}=\frac{m(X_{-v}-Y_{-v})}{m(W_0)}\cdot 
\prod_{i=1}^{n\cdot d}\left|{\sf e}+\nabla_{d,W_0}A_{-v}-\nabla_{d, 
W_0}Y_{-v}\vphantom{\dot{f}}\right|_i\qquad
\end{eqnarray} 
on $(\Omega,{\cal F}^{(m)})$. In the remainder of this step we 
show that (\ref{3.32}) holds also for negative $v$. To do so we 
still let $v\ge 0$ and take a look on $\omega_{-v}(W^v)$. We 
also recall Lemma \ref{Lemma3.4} (vi), Lemma \ref{Lemma3.4} (a) 
and (c), as well as conditions 2 (i) and (4) (i) of Section 1 and 
obtain $Q^{(m)}_{\sbnu}$-a.e. 
\begin{eqnarray*} 
&&\hspace{-.5cm}\omega^{(m)}_{-v}(W^v)=\frac{m\left(X_{-v}(W^v)- 
Y_{-v}(W^v)\vphantom{l^1}\right)}{m(W_0^v)}\cdot\prod_{i=1}^{n 
\cdot d}\left|{\sf e}+\nabla_{d,W_0}A_{-v}(W^v)-\nabla_{d,W_0}Y_{ 
-v}(W^v)\vphantom{\dot{f}}\right|_i \\ 
&&\hspace{.5cm}=\frac{m\left(W_0\vphantom{l^1}\right)}{m(W_0^v)} 
\cdot\exp\left\{\left\langle {\sf e},\log\left|{\sf e}+\nabla_{d 
,W_0}A_{-v}(W^v)-\nabla_{d,W_0}Y_{-v}(W^v)\vphantom{\dot{f}} 
\right|\right\rangle_F\right\} \\ 
&&\hspace{.5cm}=\lim_{n\to\infty}\left(\frac{m\left(W_0\vphantom 
{l^1}\right)}{m(W_0^{n,v})}\cdot\exp\left\{\left\langle {\sf e}, 
\log\left|{\sf e}+\nabla_{d,W_0}A_{n,-v}(W^{n,v})-\nabla_{d,W_0} 
Y_{n,-v}(W^{n,v})\vphantom{\dot{f}}\right|\right\rangle_F\right\} 
\right) \\ 
&&\hspace{.5cm}=\lim_{n\to\infty}\left(\frac{m\left(W_0\vphantom 
{l^1}\right)}{m(W_0^{n,v})}\cdot\exp\left\{-\int_{s=0}^v\left 
\langle{\sf e}\, ,\, \left(\log\left|{\sf e}+\nabla_{d,W_0}\left( 
A_{n,-s}-Y_{n,-s}\right)\vphantom{\dot{f}}\right|\right)^{\cdot} 
\times_c\right.\right.\right.  \\ 
&&\hspace{10.8cm}\left.\left.\left.\vphantom{\int_{s=0}^v}\times_c 
(W^{n,v})\right\rangle_F\, ds\vphantom{\dot{f}}\right\}\right) \\ 
&&\hspace{.5cm}=\lim_{n\to\infty}\left(\frac{m\left(W_0\vphantom 
{l^1}\right)}{m(W_0^{n,v})}\cdot\exp\left\{-\int_{s=0}^v\left 
\langle {\sf e},\nabla_{d,W_0}\left(\left(\dot{A}_{n,0}-\dot{Y}_{n 
,0}\right)(W^{n,-s})\right)\times_c\right.\right.\right.  \\ 
&&\hspace{5.7cm}\left.\left.\left.\times_c\left({\sf e}+\nabla_{d, 
W_0}\left(A_{n,-s}-Y_{n,-s}\right)\vphantom{\dot{f}}\right)^{-1 
}_c\right\rangle_F\, ds\circ W^{n,v}\right\}\right) \\ 
&&\hspace{.5cm}=\lim_{n\to\infty}\left(\frac{m\left(W_0\vphantom 
{l^1}\right)}{m(W_0^{n,v})}\cdot\exp\left\{-\int_{s=0}^v\left 
\langle {\sf e},\nabla_{d,W_0}\left(\dot{A}_{n,0}-\dot{Y}_{n,0} 
\right)(W^{n,-s})\right\rangle_F\, ds\circ W^{n,v}\right\}\right) 
 \\ 
&&\hspace{.5cm}=\lim_{n\to\infty}\left(\frac{m\left(W_0\vphantom 
{l^1}\right)}{m(W_0^{n,v})}\cdot\exp\left\{-\int_{s=-v}^0\left 
\langle {\sf e},\nabla_{d,W_0}\left(\dot{A}_{n,0}-\dot{Y}_{n,0} 
\right)(W^{n,-s})\right\rangle_F\, ds\right\}\right)\, . 
\end{eqnarray*}
Performing similar calculations, now the other way around, we 
obtain 
\begin{eqnarray*} 
&&\hspace{-.5cm}\frac{dQ^{(m)}_{\sbnu}\circ W^v}{dQ^{(m)}_{\sbnu}  
}=\frac{1}{\omega^{(m)}_{-v}(W^v)} \\ 
&&\hspace{.5cm}=\lim_{n\to\infty}\left(\frac{m(W_0^{n,v})}{m\left 
(W_0\vphantom{l^1}\right)}\cdot\exp\left\{-\int_{s=0}^{-v}\left 
\langle {\sf e},\nabla_{d,W_0}\left(\dot{A}_{n,0}-\dot{Y}_{n,0} 
\right)(W^{n,-s})\right\rangle_F\, ds\right\}\right) \\ 
&&\hspace{.5cm}=\lim_{n\to\infty}\left(\frac{m(W_0^{n,v})}{m\left 
(W_0\vphantom{l^1}\right)}\cdot\exp\left\{-\int_{s=0}^{-v}\left 
\langle{\sf e}\, ,\, \left(\log\left|{\sf e}+\nabla_{d,W_0}\left( 
A_{n,-s}-Y_{n,-s}\right)\vphantom{\dot{f}}\right|\right)^{\cdot} 
\right\rangle_F\, ds\right\}\right) \\ 
&&\hspace{.5cm}=\lim_{n\to\infty}\left(\frac{m(W_0^{n,v})}{m\left 
(W_0\vphantom{l^1}\right)}\cdot\exp\left\{\left\langle {\sf e}, 
\log\left|{\sf e}+\nabla_{d,W_0}A_{n,v}-\nabla_{d,W_0}Y_{n,v} 
\vphantom{\dot{f}}\right|\right\rangle_F\right\}\right) \\ 
&&\hspace{.5cm}=\frac{m(W_0^v)}{m\left(W_0\vphantom{l^1}\right)} 
\cdot\exp\left\{\left\langle {\sf e},\log\left|{\sf e}+\nabla_{d, 
W_0}A_v-\nabla_{d,W_0}Y_v\vphantom{\dot{f}}\right|\right\rangle_F 
\right\} \\ 
&&\hspace{.5cm}=\frac{m\left(X_v-Y_v\vphantom{l^1}\right)}{m(W_0) 
}\cdot\prod_{i=1}^{n\cdot d}\left|{\sf e}+\nabla_{d,W_0}A_v- 
\nabla_{d,W_0}Y_v\vphantom{\dot{f}}\right|_i\, . 
\end{eqnarray*} 
This implies that relation (\ref{3.32}) holds also for $v<0$. 
\medskip 

\nid
{\it Step 4 } In this step we will prepare the verification of the 
conditions of Theorem \ref{Theorem2.1}. 
\medskip 

Let us consider the L\'evy-Ciesielsky construction of the process 
$W_s=\eta+\sum_{i=1}^\infty\xi_i\cdot\int_0^s H_i(u)\, du$, $s\in 
{\Bbb R}$, with independent $N(0,t)$-distributed random variables 
$\xi_1 ,\xi_2 ,\ldots\, $ and $\eta$ being an $F$-valued random 
variable with distribution $\bnu$ independent of $\xi_1, \xi_2 , 
\ldots\, $. We recall that the right-hand side converges $Q_{\sbnu} 
$-a.e. uniformly in $s$ on compact subsets of ${\Bbb R}$ and that 
$\xi_i=\langle H_i,dW\rangle_{L^2}$, $i\in {\Bbb N}$, and $\eta = 
W_0$. In Subsection 2.1 we had introduced the projection of $W$ to 
the linear span of 
\begin{eqnarray*}
\left\{\int_0^{\, \cdot} H_i(u)\, du,\, e_j\cdot\1:i\in J(m),\, 
j\in\{1,\ldots ,n\cdot d\}\right\} 
\end{eqnarray*}
by considering the process $W$ under the measure $Q^{(m)}_{\sbnu}$. 
\medskip 

For a path $W\in\Omega$ with $W_s=y+\sum_{i=1}^\infty x_i\cdot 
\int_0^s H_i(u)\, du$, $s\in {\Bbb R}$, we will use the 
identifications 
\begin{eqnarray*}
Q_{\sbnu}(dW)\equiv Q_{\sbnu}\left(\eta\in dy,\left\{\xi_i\in dx_i: 
i\in {\Bbb N}\right\}\right) 
\end{eqnarray*}
and for 
\begin{eqnarray*}
W^{(m)}_s:=\eta+\sum_{i\in J(m)}\xi_i\cdot\int_0^s H_i(u)\, du 
\, , \quad s\in {\Bbb R},   
\end{eqnarray*}
we shall write 
\begin{eqnarray*}
Q_{\sbnu}^{(m)}(dW^{(m)})\equiv Q_{\sbnu}\left(\eta\in dy,\left\{ 
\xi_i\in dx_i:i\in J(m)\right\}\right)\, . 
\end{eqnarray*}

Define $A^{(m)}:=A(W^{(m)})$ and $Y^{(m)}:=Y(W^{(m)})$. Let $i\in 
{\Bbb N}$ and $i'\equiv i'(i)$ are related by $H_{i'}=H_i(\cdot +t)$. 
Note that $i\to i'$ is a bijection of type ${\Bbb N}\to {\Bbb N}$ or 
$J(m)\to J(m)$. It follows from the L\'evy-Ciesielsky construction of 
the path $W\in \Omega$ that   
\begin{eqnarray}\label{3.33}
Z_t(W)&&\hspace{-.5cm}:=W_{\cdot -t}+\left(A^{(m)}_{-t}-Y^{(m)}_{-t} 
\right)\1\vphantom{\int_0^\cdot}\nonumber \\ 
&&\hspace{-.5cm}=\left(\eta+A^{(m)}_{-t}-Y^{(m)}_{-t}\right)\1+\sum_{ 
i\in {\Bbb N}}x_i\cdot\int_0^{\, \cdot -t}H_i(u)\, du\nonumber 
\\ 
&&\hspace{-.5cm}=\left(\eta+A^{(m)}_{-t}-Y^{(m)}_{-t}\right)\1+\sum_{ 
i\in {\Bbb N}}x_{i'}\cdot\int_t^{\, \cdot}H_i(u)\, du\, .  
\end{eqnarray}
Similarly, we obtain 
\begin{eqnarray}\label{3.333}
W^{(m)}_{\cdot -t}+\left(A^{(m)}_{-t}-Y^{(m)}_{-t}\right)\1=\left( 
\eta+A^{(m)}_{-t}-Y^{(m)}_{-t}\right)\1+\sum_{i\in J(m)}x_{i'}\cdot 
\int_t^{\, \cdot}H_i(u)\, du\, .  
\end{eqnarray}

Let us demonstrate that (\ref{3.33}) implies that $W\to W_{\cdot -t} 
+\left(A^{(m)}_{-t}-Y^{(m)}_{-t}\right)\1$ is an injection with 
\begin{eqnarray}\label{3.331}
Q_{\sbnu}\left(\left\{ W_{\cdot -t}+\left(A^{(m)}_{-t}-Y^{(m)}_{-t} 
\right)\1:W\in\Omega\right\}\right)=1\, . 
\end{eqnarray}

To begin with let us verify injectivity. Given $Z$, the values of 
$x_i$, $i\in {\Bbb N}$, are uniquely determined by the linear 
combination $Z-Z_0\1=\sum_{i\in {\Bbb N}}x_{i'}\cdot\int_0^{\, \cdot} 
H_i(u)\, du$ and the bijectivity of $i\to i'$. It remains to 
demonstrate that for given $x_i$, $i\in {\Bbb N}$, and $\eta+A^{ 
(m)}_{-t}-Y^{(m)}_{-t}$, the value of $\eta$ is unique. Because of 
(\ref{3.333}) the term $W^{(m)}_{\cdot -t}+\left(A^{(m)}_{-t}-Y^{(m) 
}_{-t}\right)\1$ is now given and therefore 
\begin{eqnarray*}
X_{\cdot -t}(W^{(m)})=u\left(W^{(m)}_{\cdot -t}+\left(A^{(m)}_{-t}- 
Y^{(m)}_{-t}\right)\1\right)=\left(W^{(m)}\right)^{-t}+A\left(\left( 
W^{(m)}\right)^{-t}\right)  
\end{eqnarray*}
is given as well, cf. condition (3) of Section 1. By the injectivity 
of the map $u$, cf. Subsection 1.2, we get $Y^{(m)}_{-t}=Y_{-t}\left( 
W^{(m)}\right)=A_0\left(\left(W^{(m)}\right)^{-t}\right)$ uniquely. 
Furthermore, we get $X\left(W^{(m)}\right)=W^{(m)}+A^{(m)}$ and again 
by the injectivity of the map $u$ this yields $A_{-t}^{(m)}$ uniquely. 
This shows finally that $\eta$ is unique. Therefore, $W\to W_{\cdot 
-t}+\left(A^{(m)}_{-t}-Y^{(m)}_{-t}\right)\1$ is an injection. From 
$W^{(m)}_{\cdot -t}=(W_{\cdot -t})^{(m)}$ and $\vp\equiv 1$ in 
(\ref{3.25}) we obtain 
\begin{eqnarray*} 
Q^{(m)}_{\sbnu}\left(\left\{W^{(m)}_{\cdot -t}+\left(A^{(m)}_{-t} 
-Y^{(m)}_{-t}\right)\1:W\in\Omega\right\}\right)=1\, .    
\end{eqnarray*}  
Together with (\ref{3.33}), (\ref{3.333}), and the independence of 
$A^{(m)}$ of $\xi_i$, $i\not\in J(m)$, this means 
\begin{eqnarray*} 
1&&\hspace{-.5cm}=Q^{(m)}_{\sbnu}\left(\left\{W^{(m)}_{\cdot -t}+ 
\left(A^{(m)}_{-t}-Y^{(m)}_{-t}\right)\1:W\in\Omega\right\}\right)\cdot 
Q_{\sbnu}\left(\left\{\xi_i\in {\Bbb R}:i\not\in J(m)\right\}\right) 
 \\ 
&&\hspace{-.5cm}=Q_{\sbnu}\left(\left\{W_{\cdot -t}+\left(A^{(m)}_{ 
-t}-Y^{(m)}_{-t}\right)\1:W\in\Omega\right\}\right) 
\end{eqnarray*}  
which is (\ref{3.331}). The flow property of $W^v$, $v\in {\Bbb R}$, 
verified in Remark (3) of Section 1 implies that $W\to W_{\cdot -t} 
+\left(A_{-t}-Y_{-t}\right)\1=W^{-t}$ is an injection. Let us also 
demonstrate that 
\begin{eqnarray}\label{3.332} 
Q_{\sbnu}\left(\left\{W_{\cdot -t}+\left(A_{-t}-Y_{-t}\right)\1: W 
\in\Omega\right\}\right)=1\, .  
\end{eqnarray} 

By (\ref{3.331}), for $Q_{\sbnu}$-a.e. $V\in\Omega$ and every $m\in 
{\Bbb N}$, there is a $\rho(m,V)\in F$ such that with $W(m,V):=V_{ 
\cdot +t}+(\rho(m,V)-V_t)\1$, $A^{(m)}\equiv A^{(m)}(W(m,V))$, and 
$Y^{(m)}\equiv Y^{(m)}(W(m,V))$ we have $W_0(m,V)\in D^n$ and 
\begin{eqnarray*}
V=W_{\cdot -t}(m,V)+\left(A^{(m)}_{-t}-Y^{(m)}_{-t}\right)\1\, . 
\end{eqnarray*}
In particular, we have $\rho(m,V)=W_0(m,V)\in D^n$. We also recall 
that by the conventions of Subsection 2.1 $D$ is bounded. Therefore 
we may choose an accumulation point $\rho$ of the sequence $\rho(m, 
V)$, $m\in {\Bbb N}$, and a subsequence $\rho(m_k,V)$, $k\in {\Bbb 
N}$, with $\rho(m_k,V)\stack{k\to\infty}{\lra}\rho$. In other words, 
$W(m_k;V)\stack{k\to\infty}{\lra}W:=V_{\cdot +t}+(\rho -V_t)\1$. By 
the hypothesis that $A_{-t}-Y_{-t}$ jumps only on a set of $Q_{\sbnu 
}$-measure zero in condition (2) of Section 1 and by condition (4) 
(i) this implies $A_{-t}^{(m_k)}-Y_{-t}^{(m_k)}\stack {k\to\infty} 
{\lra}A_{-t}-Y_{-t}$. We have shown that for $Q_{\sbnu}$-a.e. 
$V\in\Omega$ there exist a $W\in\Omega$, namely $W=V_{\cdot +t}+ 
(\rho -V_t)\1=V_{\cdot +t}+\lim_{k\to\infty}(\rho(m_k,V)-V_t)\1$, 
such that 
\begin{eqnarray*}
V=W_{\cdot -t}+\left(A_{-t}-Y_{-t}\right)\1 
\end{eqnarray*}
which is an equivalent formulation of (\ref{3.332}). 

As a further preparation for the next step we mention that because of 
the injectivity of $W\to W_{\cdot -t}+\left(A_{-t}-Y_{-t}\right)\1$ 
for $Q_{\sbnu}$-a.e. $V\in\Omega$ there is just one accumulation point 
of the sequence $\rho(m,V)$, $m\in {\Bbb N}$, i. e., the limit $\rho= 
\lim_{m\to\infty}\rho(m,V)$ is well-defined. 
\medskip 

\nid 
{\it Step 5 } We shall verify conditions (i) and (iii) of Theorem 
\ref{Theorem2.1} and establish the limit 
\begin{eqnarray*}
\lim_{n\to\infty}\frac{dQ_{\sbnu}\circ f_n^{-1}}{dQ_{\sbnu}}\equiv 
\lim_{m\to\infty}\frac{dQ_{\sbnu}\circ (f^{(m)})^{-1}}{dQ_{\sbnu}} 
\end{eqnarray*}
in $L^1(\Omega,Q_{\sbnu})$. In the context of Theorem \ref{Theorem2.1} 
let $M:=\Omega$ and $\mu:=Q_{\sbnu}$. By the result of Step 4 we can 
introduce $f_n\equiv f^{(m)}$ and $f$ for $Q_{\sbnu}$-a.e. $W\in\Omega$ 
implicitly by their inverse functions  
\begin{eqnarray}\label{3.34} 
\begin{array}{c}
\D f^{-1}(W):=W_{\cdot -t}+\left(A_{-t}-Y_{-t}\right)\1\, ,  \\ \\ 
\D f_n^{-1}(W)\equiv \left(f^{(m)}\right)^{-1}(W):=W_{\cdot -t}+ 
\left(A_{-t}^{(m)}-Y_{-t}^{(m)}\right)\1\, ,  
\end{array}
\end{eqnarray}
where we stress that from here on $n\in {\Bbb N}$ has a different 
meaning than in Steps 1-3. Using the notation of Step 4 we have 
for $Q_{\sbnu}$-a.e. $V\in\Omega$ 
\begin{eqnarray*} 
\begin{array}{c}
\D f(V)=V_{\cdot +t}+(\rho-V_t)\1\, ,  \\ \\ 
\D f_n(V)\equiv f^{(m)}(V)=V_{\cdot +t}+(\rho(m,V)-V_t)\1\, .   
\end{array}
\end{eqnarray*}
and 
\begin{eqnarray*} 
f_n(V)\stack{n\to\infty}{\lra}f(V)\, , 
\end{eqnarray*}
i. e., we have condition (iii) of Theorem \ref{Theorem2.1}. 

Let $W\in\Omega$ with $W_s=y+\sum_{i=1}^\infty x_i\cdot\int_0^s H_i 
(u)\, du$, $s\in {\Bbb R}$. From (\ref{3.34}) we obtain 
\begin{eqnarray*}
f^{-1}_n(W)&&\hspace{-.5cm}=W_{\cdot -t}+\left(A^{(m)}_{-t}-Y_{-t 
}^{(m)}\right)\1 \\ 
&&\hspace{-.5cm}=(W_{\cdot -t})^{(m)}+\left(A^{(m)}_{-t}-Y_{-t}^{ 
(m)}\right)\1+\left(W_{\cdot -t}-(W_{\cdot -t})^{(m)}\right)\, .  
\end{eqnarray*}
From this and the fact that $A^{(m)}_{-t}=A_{-t}(W^{(m)})$ as well 
as $Y^{(m)}_{-t}=Y_{-t}(W^{(m)})$ are  by definition independent of 
$\xi_i$, $i\not\in J(m)$, we conclude that $Q_{\sbnu}\circ f_n^{-1} 
$ factorizes into 
\begin{eqnarray}\label{3.35} 
&&\hspace{-.5cm}dQ_{\sbnu}\circ f_n^{-1}(W)\equiv Q_{\sbnu}\left( 
df^{-1}_n(W)\right)\vphantom{\left(A^{(m)}_{-t}\right)}\nonumber 
 \\ 
&&\hspace{.5cm}=Q_{\sbnu}^{(m)}\left(d\left((W_{\cdot -t})^{(m)}+ 
\left(A_{-t}^{(m)}-Y_{-t}^{(m)}\right)\1\right)\right)\times Q_{ 
\sbnu}\left(\xi_i\in dx_{i'}:i\not\in J(m)\vphantom{l^1}\right) 
\, . \qquad 
\end{eqnarray}
Let us recall that $(W_{\cdot -t})^{(m)}=W^{(m)}_{\cdot -t}$ and 
therefore $dQ^{(m)}_{\sbnu}\circ f_n^{-1}(W)=Q^{(m)}_{\sbnu}\left( 
d f_n^{-1}(W)\right)=Q_{\sbnu}^{(m)}\left(d\left((W_{\cdot -t})^{ 
(m)}+\left(A^{(m)}_{-t}-Y_{-t}^{(m)}\right)\1\right)\right)$. For 
$W\in\Omega$ we obtain by using (\ref{3.35}) 
\begin{eqnarray}\label{3.37} 
&&\hspace{-.5cm}\omega^{(m)}_{-t}(W^{(m)})=\frac{Q_{\sbnu}^{(m)} 
\left(d\left((W_{\cdot -t})^{(m)}+\left(A^{(m)}_{-t}-Y_{-t}^{(m)} 
\right)\1\right)\right)}{Q_{\sbnu}^{(m)}\left(dW^{(m)}\right) 
\vphantom{\left(\dot{f}\right)}} \nonumber \\ 
&&\hspace{.5cm}=\frac{Q^{(m)}_{\sbnu}\left(d\left((W_{\cdot -t})^{ 
(m)}+\left(A_{-t}^{(m)}-Y_{-t}^{(m)}\right)\1\right)\right)\cdot 
Q_{\sbnu}\left(\xi_i\in dx_{i'}:i\not\in J(m)\vphantom{l^1}\right) 
}{Q^{(m)}_{\sbnu}\left(dW^{(m)}\right)\cdot Q_{\sbnu}\left(\xi_i 
\in dx_i:i\not\in J(m)\vphantom{l^1}\right)\vphantom{\left(\dot{f} 
\right)}} \nonumber \\ 
&&\hspace{.5cm}=\frac{Q_{\sbnu}\left(d\left(W_{\cdot -t}+\left(A_{ 
-t}^{(m)}-Y_{-t}^{(m)}\right)\1\right)\right)}{Q_{\sbnu}\left(dW 
\right)\vphantom{\left(\dot{f}\right)}} \nonumber \\ 
&&\hspace{.5cm}=\frac{dQ_{\sbnu}\circ f_n^{-1}}{dQ_{\sbnu}}\, (W) 
\, . 
\end{eqnarray}

We get immediately condition (i) of Theorem \ref{Theorem2.1}. Next 
recall conditions (1) (iv) and (4) (i) of Section 1. Furthermore, 
take into consideration that $X_{-t}$, $A_{-t}$, $Y_{-t}$ jump only 
on a set of $Q_{\sbnu}$-measure zero, cf. conditions (2) and (3) of 
Section 1, that the density $m$ is bounded and continuous, and that 
$\lambda_F(D^n)<\infty$ by the boundedness of $D$. Using the result 
of Step 3, namely (\ref{3.32}), relation (\ref{3.37}) shows now that 
\begin{eqnarray}\label{3.38}
\lim_{n\to\infty}\frac{dQ_{\sbnu}\circ f_n^{-1}}{dQ_{\sbnu}}&&\hspace 
{-.5cm}=\lim_{m\to\infty}\left(\frac{m(X_{-t}-Y_{-t})}{m(W_0)}\cdot 
\prod_{i=1}^{n\cdot d}\left|{\sf e}+\nabla_{d,W_0}A_{-t}-\nabla_{d, 
W_0}Y_{-t}\vphantom{\dot{f}}\right|_i\right)\circ\pi_m\nonumber \\ 
&&\hspace{-.5cm}=\frac{m(X_{-t}-Y_{-t})}{m(W_0)}\cdot\prod_{i=1}^{n 
\cdot d}\left|{\sf e}+\nabla_{d,W_0}A_{-t}-\nabla_{d,W_0}Y_{-t} 
\vphantom{\dot{f}}\right|_i
\end{eqnarray} 
in $L^1(\Omega,Q_{\sbnu})$ where $m\to\infty$ in the second line 
refers to the projection $\pi_m$, not to the density $m$. 

After having verified condition (ii) of Theorem \ref{Theorem2.1}, 
relation (\ref{3.38}) will yield (\ref{1.1}). 
\medskip 
 
\nid 
{\it Step 6 } We verify condition (ii) of Theorem \ref{Theorem2.1}. 
For this, we have to demonstrate that the sequence of densities $d 
\mu\circ f_n^{-1}/d\mu$, $n\in {\Bbb N}$, is uniformly integrable. 
By (\ref{3.32}) and (\ref{3.37}) this means that we have to verify 
uniform integrability with respect to $Q_{\sbnu}$ of the terms 
\begin{eqnarray*} 
\frac{dQ_{\sbnu}\circ f_n^{-1}}{dQ_{\sbnu}\vphantom{\left(\dot{f} 
\right)}}&&\hspace{-.5cm}=\left(\frac{m(X_{-t}-Y_{-t})}{m(W_0)}
\cdot\prod_{i=1}^{n\cdot d}\left|{\sf e}+\nabla_{d,W_0}A_{-t}- 
\nabla_{d,W_0}Y_{-t}\vphantom{\dot{f}}\right|_i\right)\circ\pi_m 
\vphantom{\int_\int} \\ 
&&\hspace{-.5cm}\le\frac{\|m\|}{m(W_0)}\cdot\sup_m\prod_{i=1}^{n 
\cdot d}\left|{\sf e}+\nabla_{d,W_0}A_{-t}-\nabla_{d,W_0}Y_{-t} 
\vphantom{\dot{f}}\right|_i\circ\pi_m\vphantom{\int^\int}
\end{eqnarray*}
where the right-hand side is in $L^1(\Omega,Q_{\sbnu})$ by 
condition (1) (iv) of Section 1. We have thus verified the 
uniform integrability with respect to $Q_{\sbnu}$ of the sequence 
$dQ_{\sbnu}\circ f_n^{-1}/dQ_{\sbnu}$, $n\in {\Bbb N}$.  
\qed

\section{Partial Integration} 
\setcounter{equation}{0} 

Let us assume (1)-(5) of Section 1. Then for $t\in {\Bbb R}$ we have 
according to Theorem \ref{Theorem1.11} and Corollary \ref{Corollary1.13} 
\begin{eqnarray*} 
\rho_{-t}(X)&&\hspace{-.5cm}=\frac{P_{\sbmu}(dX_{\cdot -t})}{P_{\sbmu} 
(dX)} \\ 
&&\hspace{-.5cm}=\frac{\D m\left(X_{-t}-Y_{-t}\circ u^{-1}\vphantom 
{l^1}\right)}{\D m\left(u^{-1}_0\vphantom{l^1}\right)}\cdot\prod_{i=1 
}^{n\cdot d}\left|{\sf e}+\nabla_{d,W_0}A_{-t}\circ u^{-1}-\nabla_{d, 
W_0}Y_{-t}\circ u^{-1}\vphantom{\dot{f}}\right|_i\, . 
\end{eqnarray*} 

In contrast to the previous sections, where $\bmu$ was used as a part 
of the symbol $P_{\sbmu}=Q_{\sbnu}\circ u^{-1}$, we give $\bmu$ now a 
meaning by $\bmu(A):=P_{\sbmu}(X_0\in A)$, $A\in {\cal B}(F)$. For $x 
\in D^n$, let $E_x$ denote the expectation relative to $P_x:=P_{\sbmu} 
(\, \cdot\, |X_0=x)$ and let $E_{\sbmu}$ denote the expectation relative 
to $P_{\sbmu}$. In addition to condition (4) (i) of Section 1 we may 
formulate the following condition. 
\begin{itemize}
\item[(6)] For all $W\in\Omega$, the jump times for $\nabla_{d,W_0}A(W)$ 
and $\nabla_{d,W_0}Y(W)$ are a subset of $\{\tau_k(u(W))$ $\equiv\tau_k 
(X):k\in {\Bbb Z}\setminus\{0\}\}$. 
\end{itemize}
We may also assume that 
\begin{itemize}
\item[(7)] 
\begin{eqnarray*} 
\frac{d\bmu}{d\bnu}\in C_b(D^n)\quad\mbox{\rm and}\quad\frac{d\bmu} 
{d\bnu}>0\, .
\end{eqnarray*} 
\end{itemize}
We recall that, by condition (2) (i) of Section 1, for any fixed $s 
\in {\Bbb R}$ and $\lambda_F$-a.e. $x\in D^n$, $A^2_s\circ u^{-1}$ 
jumps only on a $P_x$-zero set. As a consequence of this and conditions 
(6) and (7) as well as conditions (1) (iv) and (4) (i) of Section 1 we 
obtain 
\begin{eqnarray}\label{4.1} 
\lim_{t\downarrow 0}\frac{P_{\sbmu}(X_t\in dx)}{P_{\sbmu}(X_0\in dx)}  
=\lim_{t\downarrow 0}\int\rho_{-t}(X)\, P_x(dX)=\1 
\end{eqnarray} 
boundedly for $\lambda_F$-a.e. $x\in D^n$. Let us set 
\begin{eqnarray}\label{4.2}
r_{-t}\left(X,x\right)&&\hspace{-.5cm}:=E_x\left.\left(\frac{\D m\left( 
W_0^{-t}\circ u^{-1}\right)}{\D m\left(u^{-1}_0\vphantom{l^1}\right)}\, 
\prod_{i=1}^{n\cdot d}\left|{\sf e}+\nabla_{d,W_0}A_{-t}\circ u^{-1}- 
\nabla_{d,W_0}Y_{-t}\circ u^{-1}\vphantom{\dot{f}}\right|_i\vphantom{ 
\frac{\D m\left(\left(W^{-t}\circ u^{-1}\vphantom{l^1}\right)_0\right)} 
{\D m\left(u^{-1}_0\vphantom{l^1}\right)}}\right|X_{-t}\right)\nonumber 
 \\ 
&&\hspace{-.4cm}=E_{\sbmu}\left.\left(\frac{P_{\sbmu}(dX_{\cdot -t})} 
{P_{\sbmu}(dX)}\right|X_{-t},X_0=x\right)\, , \quad x\in D^n,\ t\in 
{\Bbb R}. 
\end{eqnarray} 
Moreover, for $x\in D^n$, let $E_{t,x}$ be the expectation relative to 
$P_{\sbmu}(\, \cdot \, |X_t=x)$. Let $D(A)$ and $D(A^\ast)$ denote the set 
of $f\in L^2(F,\bmu)$ for which, respectively, the limit 
\begin{eqnarray*} 
Af:=\lim_{t\downarrow 0}\frac1t\left(E_\cdot f(X_t)-f\right)\quad\mbox{ 
\rm and}\quad A^\ast f:=-\lim_{t\downarrow 0}\frac1t\left(E_{t,\cdot} f(X_0) 
-f\right)
\end{eqnarray*} 
exists in $L^2(F,\bmu)$. 

Let us write $\langle\, \cdot\, ,\, \cdot\, \rangle_{\sbmu}$ and $\|\cdot\|_{ 
\sbmu}$ for the inner product and the norm in ${L^2(F,\bmu)}$. For $F={\Bbb R 
}^{n\cdot d}$ and $f\equiv(f_1,\ldots ,f_{n\cdot d})\in C^1(\overline{D^n};F)$ 
let $\nabla_{d,x}f:=\left(\frac{\partial f_1}{\partial x_1},\ldots ,\frac{ 
\partial f_{n\cdot d}}{\partial x_{n\cdot d}}\right)$.

Furthermore, let $C^{1,1}\left({\Bbb R}\times\overline{D^n};F\right)$ denote 
the space of all functions $a\equiv a(s,x)$, $s\in{\Bbb R}$, $x\in\overline 
{D^n}$, for which $a,a'\equiv(\partial a)/(\partial s),\nabla_{d,x}a,(\nabla_{ 
d,x} a)'\in C\left({\Bbb R}\times\overline{D^n};F\right)$.  It follows that 
$\nabla_{d,x}(a')\in C\left({\Bbb R}\times\overline{D^n};F\right)$ and that 
$\nabla_{d,x} (a')=(\nabla_{d,x} a)'$. 

In the remainder of this section we will use the letter $Z$ to abbreviate 
$A-Z$. This also means that, for example, we write $\nabla_{d,W_0}Z_{-t} 
\circ u^{-1}$ for $\nabla_{d,W_0}A_{-t}\circ u^{-1}-\nabla_{d,W_0}Y_{-t} 
\circ u^{-1}$ and $\nabla_{d,W_0}Z'_0\circ u^{-1}$ for $\nabla_{d,W_0}A'_0 
\circ u^{-1}-\nabla_{d,W_0}Y'_0\circ u^{-1}$. 

Now we are able to state and prove the $L^2$-{\it version} of the integration 
by parts theorem. 
\begin{theorem}\label{Theorem4.1} 
Suppose that the hypotheses of Theorem \ref{Theorem1.11} and Corollary 
\ref{Corollary1.13} as well as conditions (6) and (7) of this section are 
satisfied. In addition, assume the following. 
\begin{itemize} 
\item[(i)] The derivative 
\begin{eqnarray*}
r^\ast:=\lim_{t\downarrow 0}\frac1t\left(E_\cdot r_{-t}(X,\cdot)-\1\right) 
\end{eqnarray*} 
exists in $L^2(F,\bmu)$. 
\end{itemize} 
(a) Then for $f\in D(A)$ and $g\in D(A^\ast)\cap B_b(D^n)$, we have 
\begin{eqnarray*} 
\langle Af\, , \, g\rangle_{\sbmu}+\langle f\, , \, A^\ast g\rangle_{\sbmu}= 
\left\langle f\, , \, r^\ast\cdot g\right\rangle_{\sbmu}\, .  
\end{eqnarray*} 
(b) Let, in addition to (i), 
\begin{itemize} 
\item[(ii)] the limit 
\begin{eqnarray*} 
B^\ast m:=-\lim_{t\downarrow 0}\frac1t\left(E_\cdot m\left(W_0^{-t}\circ 
u^{-1}\right)-m\right) 
\end{eqnarray*} 
exists in $L^2(D^n,\lambda_F)$ and 
\item[(iii)] $A\equiv A_s(W-W_0,W_0)\in C^{1,1}\left({\Bbb R}\times 
\overline{D^n};F\right)$ $Q_x$-a.e. and, for the restriction of $\nabla_{ 
d,W_0}A$ to $[-1,0]$, $\nabla_{d,W_0}A'\in L^\infty\left(\Omega,Q_x;C_b([ 
-1,0];F)\right)$ for a.e. $x\in D^n$, and the same for $Y$
\end{itemize} 
then 
\begin{eqnarray*} 
r^\ast=-\frac{B^\ast m}{m}-E_\cdot\left\langle {\sf e},\nabla_{d,W_0}Z'_0 
\circ u^{-1}\right\rangle_F\, .  
\end{eqnarray*} 
\end{theorem}
Proof. (a) By (\ref{4.1}) and condition (i) we have for $g\in D(A^\ast) 
\cap B_b(D^n)$ 
\begin{eqnarray}\label{4.3}
&&\hspace{-.5cm}\lim_{t\downarrow 0}\frac1t\left(\int g(X_{-t})\cdot 
r_{-t}\left(X,\cdot\right)\, dP_\cdot -g\right) \nonumber \\ 
&&\hspace{.5cm}=\lim_{t\downarrow 0}\frac1t\int\left(g(X_{-t})-g\right) 
\cdot r_{-t}\left(X,\cdot\right)\, dP_\cdot +\lim_{t\downarrow 0}\frac1t 
\int g\cdot\left(r_{-t}\left(X,\cdot\right)-\1\right)\, dP_\cdot 
\nonumber \\ 
&&\hspace{.5cm}=\lim_{t\downarrow 0}\frac1t\int\left(g(X_{-t})-g\right) 
\cdot E_{\sbmu}\left.\left(\frac{P_{\sbmu}(dX_{\cdot -t})}{P_{\sbmu}(dX)} 
\right|X_{-t},X_0=\, \cdot\, \right)\, dP_\cdot +r^\ast\cdot g\nonumber 
 \\ 
&&\hspace{.5cm}=\lim_{t\downarrow 0}\frac1t\int E_{\sbmu}\left.\left( 
\left(g(X_{-t})-g\right)\cdot\frac{P_{\sbmu}(dX_{\cdot -t})}{P_{\sbmu} 
(dX)}\right|X_{-t},X_0=\, \cdot\, \right)\, dP_\cdot +r^\ast\cdot g 
\nonumber \\ 
&&\hspace{.5cm}=\lim_{t\downarrow 0}\left(\frac1t\int\left(g(X_0)-g 
\right)\, P_{\sbmu}(dX|X_t=\cdot)\ \frac{P_{\sbmu}(X_t\in d\, \cdot)} 
{P_{\sbmu}(X_0\in d\, \cdot)}\right)+r^\ast\cdot g \nonumber \\ 
&&\hspace{.5cm}=-A^\ast g+r^\ast\cdot g\, \quad\mbox{\rm in } L^2(F, 
\bmu). \vphantom{\int}
\end{eqnarray}
Keeping (\ref{4.2}) and (\ref{4.3}) in mind, the following chain of 
equations is now self-explaining. It holds for $f\in D(A)$ and $g\in 
D(A^\ast)\cap B_b(D^n)$ that 
\begin{eqnarray}\label{4.4}
\langle Af\, , \, g\rangle_{\sbmu}&&\hspace{-.5cm}=\left\langle\left. 
\frac{d^+}{dt}\right|_{t=0}E_\cdot f(X_t)\, , \, g\right\rangle_{\sbmu} 
\vphantom{\left.\frac{d^+}{dt}\right|_{t=0}} \nonumber \\ 
&&\hspace{-1.5cm}=\left.\frac{d^+}{dt}\right|_{t=0}\left\langle E_\cdot 
f(X_t)\, , \, g\right\rangle_{\sbmu} \nonumber \\ 
&&\hspace{-1.5cm}=\lim_{t\downarrow 0}\frac1t\left(\int f(X_t)g(X_0)\, d 
P_{\sbmu}-\langle f\, , \, g\rangle_{\sbmu}\right) \nonumber \\ 
&&\hspace{-1.5cm}=\lim_{t\downarrow 0}\frac1t\left(\int f(X_0) g(X_{-t}) 
\, P_{\sbmu}(dX_{\cdot -t})-\langle f\, , \, g\rangle_{\sbmu}\right) 
\nonumber \\ 
&&\hspace{-1.5cm}=\lim_{t\downarrow 0}\frac1t\left(\int f(X_0)g(X_{-t}) 
\cdot E_{\sbmu}\left.\left(\frac{P_{\sbmu}(dX_{\cdot -t})}{P_{\sbmu}(d 
X_\cdot )}\right|X_{-t},X_0\right)\, dP_{\sbmu}-\langle f\, , \, g 
\rangle_{\sbmu}\right) \nonumber \\ 
&&\hspace{-1.5cm}=\lim_{t\downarrow 0}\frac1t\int f(X_0)\left(g(X_{-t}) 
\cdot E_{\sbmu}\left.\left(\frac{P_{\sbmu}(dX_{\cdot -t})}{P_{\sbmu}(d 
X_\cdot )}\right|X_{-t},X_0\right)-g(X_0)\right)\, dP_{\sbmu} \nonumber 
\\ 
&&\hspace{-1.5cm}=\lim_{t\downarrow 0}\frac1t\int f(x)\left(\int g(X_{- 
t})\cdot E_{\sbmu}\left.\left(\frac{P_{\sbmu}(dX_{\cdot -t})}{P_{\sbmu} 
(dX_\cdot)}\right|X_{-t},X_0=x\right)\, dP_x-g(x)\right)\, \bmu(dx) 
\nonumber \\ 
&&\hspace{-1.5cm}=\lim_{t\downarrow 0}\frac1t\int f(x)\left(\int g(X_{- 
t})\cdot r_{-t}\left(X,x\right)\, dP_x-g(x)\right)\, \bmu(dx) \nonumber 
\\ 
&&\hspace{-1.5cm}=-\langle f\, , \, A^\ast g\rangle_{\sbmu}+\left\langle 
f\, , \, r^\ast\cdot g\right\rangle_{\sbmu}\, .  
\vphantom{\int}
\end{eqnarray}
We obtain part (a). 
\medskip 

\nid 
(b)  By (\ref{4.2}) and (ii) and (iii) of the present theorem there is 
a sequence $t_n>0$, $n\in {\Bbb N}$, with $t_n\stack{n\to\infty}{\lra}0$ 
such that $\bmu$-a.e. 
\begin{eqnarray*}
&&\hspace{-.5cm}\left|E_\cdot\left(\frac{r_{-t_n}\left(X,\cdot\right)- 
\1}{t_n}\right)+\frac{B^\ast m}{m}+E_\cdot\left\langle {\sf e},\nabla_{ 
d,W_0}Z'_0\circ u^{-1}\right\rangle_F\right|\vphantom{\sum_1^1} \\  
&&\hspace{.5cm}=\left|\frac1{t_n}\left(E_\cdot\left(\frac{m\left(W_0^{- 
t_n}\circ u^{-1}\right)}{m(u_0^{-1})}\cdot\prod_{i=1}^{n\cdot d}\left|{ 
\sf e}+\nabla_{d,W_0}Z_{-t_n}\circ u^{-1}\vphantom{\dot{f}}\right|_i 
\right)-\1\right)\right.\vphantom{\sum_1^1} \\  
&&\hspace{1.0cm}\left.\vphantom{\prod_{i=1}^n}+\frac{B^\ast m}{m}+E_\cdot 
\left\langle {\sf e},\nabla_{d,W_0}Z'_0\circ u^{-1}\right\rangle_F\right| 
\vphantom{\sum_1^1} \\ 
&&\hspace{.5cm}\le\left|\frac1{t_n}E_\cdot\left(\frac{m\left(W_0^{-t_n} 
\circ u^{-1}\right)}{m(u_0^{-1})}\cdot\left(\prod_{i=1}^{n\cdot d}\left| 
{\sf e}+\nabla_{d,W_0}Z_{-t_n}\circ u^{-1}\vphantom{\dot{f}}\right|_i-\1 
\right)\right)\right. \\ 
&&\hspace{1.0cm}\left.+E_\cdot\left\langle {\sf e},\nabla_{d,W_0}Z'_0\circ 
u^{-1}\right\rangle_F\vphantom{\left(\prod_{i=1}^{n\cdot d}\right)}\right| 
+\left|\frac1{t_n}E_\cdot\left(\frac{m\left(W_0^{-t_n}\circ u^{-1}\right)} 
{m(u_0^{-1})}-\1\right) +\frac{B^\ast m}{m}\right| \\ 
&&\hspace{.0cm}\stack{n\to\infty}{\lra}0\, . 
\end{eqnarray*} 
Part (b) is now a consequence of condition (i) of this theorem. 
\qed 
\bigskip 

Let us turn to the {\it weak version} of the integration by parts theorem. 
For this, let $D(A_w)$ denote the set of all $f\in C_b(D^n)$ for which the 
limit 
\begin{eqnarray*} 
A_w f(g):=\lim_{t\downarrow 0}\frac1t\int g\left(E_\cdot f(X_t)-f\right)\, 
d\bmu 
\end{eqnarray*} 
exists for all $g\in C_b(D^n)$. Furthermore, let $D(A^\ast_w)$ be the set 
of all $g\in C_b(D^n)$ for which 
the limit 
\begin{eqnarray}\label{4.5} 
A^\ast_w g(f):=-\lim_{t\downarrow 0}\frac1t\int f\left(E_{t,\cdot} g(X_0) 
-g\right)\, d\bmu 
\end{eqnarray} 
exists for all $f\in C_b(D^n)$ such that $A^\ast_w g$ is a bounded linear 
functional on $C(\overline{D^n})$. Let us recall conditions (1) (iv) and 
(3) (iii) of Section 1 and introduce the following condition. 
\begin{itemize}
\item[(8)] We have $\sup\left|\nabla_{d,W_0}Z_{-t}(W)\right|<1$ where $| 
\cdot|$ refers to the maximal coordinate wise absolute value and the 
supremum is taken over all $W\in\bigcup_m\{\pi_{m}V:V\in\Omega\}$ and $t 
\in [0,1]$. 
\end{itemize}
By conditions 2 (i) and (4) (i) of Section 1 we obtain in a 
coordinate wise way 
\begin{eqnarray*}
\log\left({\sf e}+\nabla_{d,W_0}Z_{-t}\right)=\log\left|{\sf e}+\nabla_{d, 
W_0}Z_{-t}\right|\quad Q_{\sbnu}\mbox{\rm -a.e.} 
\end{eqnarray*} 
and thus also $Q_{\lambda_F}$-a.e., $t\in [0,1]$.  
\begin{theorem}\label{Theorem4.2} 
Suppose that the hypotheses of Theorem \ref{Theorem1.11} and Corollary 
\ref{Corollary1.13} as well as conditions (6)-(8) of this section are 
satisfied. \\ 
(a) In addition, assume the following. 
\begin{itemize}  
\item[(iv)] The limit  
\begin{eqnarray*} 
r^\ast_w(h):=\lim_{t\downarrow 0}\frac1t\int h\left(E_\cdot r_{-t}(X,\cdot) 
-\1\right)\, d\bmu 
\end{eqnarray*} 
exists for all $h\in C_b(D^n)$. 
\end{itemize} 
Then for $f\in D(A_w)$ and $g\in D(A^\ast_w)$ we have 
\begin{eqnarray*} 
A_w f(g)+A^\ast_w g(f)=r^\ast_w(fg)\, .  
\end{eqnarray*} 
(b) Assume the following.  
\begin{itemize} 
\item[(v)] The limit  
\begin{eqnarray*} 
B^\ast_wm(h):=-\lim_{t\downarrow 0}\frac1t\int h\left(E_\cdot m\left( 
W_0^{-t}\circ u^{-1}\right)-m\right)\, dx 
\end{eqnarray*} 
exists for all $h\in C_b(D^n)$. 
\item[(vi)] 
\begin{eqnarray*} 
\limsup_{t\downarrow 0}\frac1t\int\left(m\left(W_0^{-t}\circ u^{-1}\right) 
-m(u_0^{-1})\vphantom{l^1}\right)^2\, dP_{\lambda_F}<\infty\, . 
\end{eqnarray*} 
\item[(vii)] The limit  
\begin{eqnarray*} 
e^\ast_w(h):=\lim_{t\downarrow 0}\frac1t\int hE_\cdot\left\langle {\sf 
e},\nabla_{d,W_0}Z_{-t}\circ u^{-1}\right\rangle_F\, d\bmu 
\end{eqnarray*}  
exists for all $h\in C_b(D^n)$. 
\item[(viii)] For all ${\sf g}\in {\Bbb R}^{n\cdot d}$ with $\max_{j\in 
\{1,\ldots ,n\cdot d\}}|{\sf g}_j|\le 1$ and all $h\in C_b(D^n)$ there 
exists $d_h>0$ such that for all $k\ge 1$ and $t\in [0,1]$, 
\begin{eqnarray*} 
\left|\int hE_\cdot\left(\left\langle{\sf g},\nabla_{d,W_0}Z_{-t}\circ 
u^{-1}\right\rangle_F\vphantom{\dot{f}}\right)^k\, dx\right|\le d_ht^k\, . 
\end{eqnarray*} 
\end{itemize} 
Then we have (iv) and 
\begin{eqnarray*} 
r^\ast_w(h)=-B^\ast_wm(h)-e^\ast_w(h)\, , \quad h\in C_b(D^n).  
\end{eqnarray*} 
\end{theorem} 
Proof. (a) We modify the proof of Theorem \ref{Theorem4.1} (a). By 
(\ref{4.1}) and (iv) we verify the weak version of (\ref{4.3}), 
\begin{eqnarray*}
&&\hspace{-.5cm}\lim_{t\downarrow 0}\frac1t\int f\left(\int g(X_{-t}) 
\cdot r_{-t}\left(X,\cdot\right)\, dP_\cdot -g\right)\, d\bmu \\ 
&&\hspace{.5cm}=\lim_{t\downarrow 0}\frac1t\int f\left(\int\left(g( 
X_{-t})-g\right)\cdot r_{-t}\left(X,\cdot\right)\, dP_\cdot\right)\, 
d\bmu \\ 
&&\hspace{1.0cm}+\lim_{t\downarrow 0}\frac1t\int f\int g\cdot\left( 
r_{-t}\left(X,\cdot\right)-\1\right)\, dP_\cdot\, d\bmu \\ 
&&\hspace{.5cm}=\lim_{t\downarrow 0}\frac1t\int f(x)\left(E_{t,x}g( 
X_0)-g(x)\vphantom{l^1}\right)\frac{P_{\sbmu}(X_t\in dx)}{P_{\sbmu} 
(X_0\in dx)}\ \bmu(dx)+r^\ast_w(fg)\nonumber \\ 
&&\hspace{.5cm}=r^\ast_w(fg)-A^\ast_wg(f)\, ,\vphantom{\int}
\end{eqnarray*}
for $f\in D(A_w)$ and $g\in D(A^\ast_w)$. Anything else is now as in 
(\ref{4.4}). 
\medskip 

\nid 
(b) Let $h\in C_b(D^n)$. By the above condition (8) and assumptions (vii) 
and (viii) it turns out that 
\begin{eqnarray}\label{4.6}
&&\hspace{-.5cm}\frac1t\int hE_\cdot\left(\prod_{i=1}^{n\cdot d}\left|{\sf 
e}+\nabla_{d,W_0}Z_{-t}\circ u^{-1}\vphantom{\dot{f}}\right|-\1\right)\, d 
\bmu\stack{t\downarrow 0}{\lra}-e^\ast_w(h)\, . 
\end{eqnarray} 
Recalling condition (8) and that, by the above condition (7) it holds that 
$d\bmu/d\bnu\in C_b(D^n)$, we obtain by (vi) and (viii) 
\begin{eqnarray}\label{4.7}
&&\hspace{-.5cm}\left(\frac1t\int h(X_0)\left(\frac{m\left(W_0^{-t}\circ 
u^{-1}\right)}{m(u_0^{-1})}-\1\right)\cdot\left(\prod_{i=1}^{n\cdot d}\left 
|{\sf e}+\nabla_{d,W_0}Z_{-t}\circ u^{-1}\vphantom{\dot{f}}\right|-\1\right) 
\, dP_{\sbmu}\right)^2\nonumber \\
&&\hspace{.5cm}\le\frac1t\int\left({m\left(W_0^{-t}\circ u^{-1}\right)-m(u^{ 
-1}_0)}\vphantom{l^1}\right)^2\, dP_{\lambda_F}\times\nonumber \\
&&\hspace{1.0cm}\times\frac1t\int h(X_0)^2\ \frac{d\bmu}{d\bnu}(X_0)^2\left( 
\prod_{i=1}^{n\cdot d}\left|{\sf e}+\nabla_{d,W_0}Z_{-t}\circ u^{-1}\vphantom 
{\dot{f}}\right|-\1\right)^2\, dP_{\lambda_F}\nonumber \\ 
&&\hspace{-.0cm}\stack{t\downarrow 0}{\lra}0\, . \vphantom{\int} 
\end{eqnarray} 
It follows now from (\ref{4.6}) and (\ref{4.7}) that 
\begin{eqnarray}\label{4.8}
&&\hspace{-.5cm}\lim_{t\downarrow 0}\frac1t\int hE_\cdot\left(\frac{m\left 
(W_0^{-t}\circ u^{-1}\right)}{m}\cdot\left(\prod_{i=1}^{n\cdot d}\left|{\sf 
e}+\nabla_{d,W_0}Z_{-t}\circ u^{-1}\vphantom{\dot{f}}\right|-\1\right)\right) 
\, d\bmu\nonumber \\ 
&&\hspace{.0cm}=\lim_{t\downarrow 0}\frac1t\int hE_\cdot\left(\prod_{i=1}^{n 
\cdot d}\left|{\sf e}+\nabla_{d,W_0}Z_{-t}\circ u^{-1}\vphantom{\dot{f}}\right| 
-\1\right)\, d\bmu\nonumber \\ 
&&\hspace{.5cm}+\lim_{t\downarrow 0}\frac1t\int h(X_0)\left(\frac{m\left( 
W_0^{-t}\circ u^{-1}\right)}{m(u^{-1}_0)}-\1\right)\cdot\left(\prod_{i=1}^{ 
n\cdot d}\left|{\sf e}+\nabla_{d,W_0}Z_{-t}\circ u^{-1}\vphantom{\dot{f}} 
\right|-\1\right)\, dP_{\sbmu}\nonumber \\ 
&&\hspace{.0cm}=-e^\ast_w(h)\, . \vphantom{\int}
\end{eqnarray} 
Recalling (\ref{4.2}) and condition (v) we obtain from (\ref{4.8}) 
\begin{eqnarray*} 
&&\hspace{-.5cm}-B^\ast_wm(h)-e^\ast_w(h)=\lim_{t\downarrow 0}\frac1t 
\int hE_\cdot\left(\frac{m\left(W_0^{-t}\circ u^{-1}\right)}{m}-\1\right 
)\, d\bmu \\ 
&&\hspace{1.0cm}+\lim_{t\downarrow 0}\frac1t\int hE_\cdot\left(\frac{m 
\left(W_0^{-t}\circ u^{-1}\right)}{m}\cdot \left(\prod_{i=1}^{n\cdot d} 
\left|{\sf e}+\nabla_{d,W_0}Z_{-t}\circ u^{-1}\vphantom{\dot{f}}\right|- 
\1\right)\right)\, d\bmu \\ 
&&\hspace{.5cm}=\lim_{t\downarrow 0}\frac1t\int h\left(E_\cdot r_{-t}(X, 
\cdot)-\1\right)\, d\bmu \\ 
&&\hspace{.5cm}=r^\ast_w(h)\, , \quad h\in C_b(D^n). \vphantom{\int}
\end{eqnarray*} 
\qed 
\bigskip

Let $\cal C$ be the set of all $f\in C_b(D^n)$ such that for all $g\in 
D(A_w)$ the limit 
\begin{eqnarray*} 
R(f,g):=\lim_{t\downarrow 0}\frac1t\int\left(f(X_t)-f(X_0)\right)\left(g 
(X_t)-g(X_0)\right)\, dP_{\sbmu}
\end{eqnarray*} 
exists. Denote by $D(A^\ast_w,{\cal C})$ the set of all $g\in C_b(D^n)$ 
for which the limit (\ref{4.5}) exists for all $f\in {\cal C}$. 
\begin{proposition}\label{Proposition4.3} 
Suppose that the hypotheses of Theorem \ref{Theorem1.11} and Corollary 
\ref{Corollary1.13} as well as conditions (6)-(8) of this section are 
satisfied. \\ 
(a) Furthermore assume (v) and (viii) of Theorem \ref{Theorem4.2} and 
the following. 
\begin{itemize} 
\item[(ix)] $h\to B_w^\ast m(h)$ is a bounded linear functional on $C( 
\overline{D^n})$. 
\end{itemize}
It holds that $D(A_w)\subseteq D(A^\ast_w,{\cal C})$. For $f\in {\cal C}$ 
and $g\in D(A_w)$ we have 
\begin{eqnarray*} 
A_w^\ast g(f)=A_w g(f)+R(f,g)\, . \vphantom{\int}
\end{eqnarray*} 
(b) Assume in addition that 
\begin{itemize} 
\item[(x)] the process $Y$ is constant zero or, equivalently, $X_0=W_0$ 

and $m\in D(A_w)$. 
\end{itemize} 
Then for all $h\in C_b(D^n)$ such that $h/m\in {\cal C}$ we have $m\in 
D(A^\ast_w,{\cal C})$ and 
\begin{eqnarray*} 
B^\ast_wm(h)=A^\ast_wm(h/m)=A_wm(h/m)+R(h/m,m)\, . \vphantom{\int}
\end{eqnarray*} 
\end{proposition} 
Proof. (a) By (6) and (7) as well as its consequences outlined in the 
beginning of this section, we have 
\begin{eqnarray}\label{4.12} 
&&\hspace{-.5cm}\frac{P_{\sbmu}(X_t\in dx)}{P_{\sbmu}(X_0\in dx)}\, E_{ 
t,x}(g(X_0)-g(x))=\frac{P_{\sbmu}(X_t\in dx)}{P_{\sbmu}(X_0\in dx)}\, 
E_{t,x}\left(g(X_0)-g(X_t)\right)\vphantom{\int}\nonumber \\ 
&&\hspace{.5cm}=E_x(g(X_{-t})-g(X_0))\vphantom{\int}\nonumber \\ 
&&\hspace{1.0cm}+\int(g(X_{-t})-g(X_0))(\rho_{-t}(X)-\1)\, P_x(dX) 
\nonumber \\ 
&&\hspace{.5cm}=E_x(g(X_{-t})-g(X_0))+\frac{1}{m(x)}\int(g(X_{-t})-g( 
X_0))\left(m\left(W_0^{-t}\circ u^{-1}\right)\times\right.\vphantom 
{\prod_{i=1}^{n\cdot d}}\nonumber \\ 
&&\hspace{1.0cm}\left.\times\prod_{i=1}^{n\cdot d}\left|{\sf e}+\nabla_{ 
d,W_0}Z_{-t}\circ u^{-1}\vphantom{\dot{f}}\right|_i-m(x)\right)\, P_x(d 
X)\nonumber \\ 
&&\hspace{.1cm}\stack{t\downarrow 0}{\lra}0\vphantom{\int}
\end{eqnarray} 
for $g\in C_b(D^n)$ and $\lambda_F$-a.e. $x\in D^n$. Let $f,g\in 
C_b(D^n)$. We have also 
\begin{eqnarray}\label{4.13}
&&\hspace{-.5cm}\frac1t\int f\left(E_{t,\cdot}g(X_0)-g\right)(E_\cdot 
\rho_{-t}-\1)\, d\bmu\nonumber \\ 
&&\hspace{.5cm}=\frac1t\int f\left(E_{t,\cdot}g(X_0)-g\right)\cdot 
E_\cdot\left(\frac{m\left(W_0^{-t}\circ u^{-1}\right)}{m}\cdot\prod_{ 
i=1}^{n\cdot d}\left|{\sf e}+\nabla_{d,W_0}Z_{-t}\circ u^{-1}\vphantom 
{\dot{f}}\right|_i-\1\right)\, d\bmu\nonumber \\ 
&&\hspace{.5cm}=\int f\left(E_{t,\cdot}g(X_0)-g\right)\cdot E_\cdot 
\left(\frac{m\left(W_0^{-t}\circ u^{-1}\right)-m}{t}\right)\, dx 
\nonumber \\
&&\hspace{1.0cm}+\int f\left(E_{t,\cdot}g(X_0)-g\right)\times\nonumber 
\vphantom{\prod_{i=1}^{n\cdot d}} \\ 
&&\hspace{1.5cm}\times E_\cdot\left(\frac{m\left(W_0^{-t}\circ u^{-1} 
\right)-m}{t}\cdot\left(\prod_{i=1}^{n\cdot d}\left|{\sf e}+\nabla_{d, 
W_0}Z_{-t}\circ u^{-1}\vphantom{\dot{f}}\right|_i-\1\right)\right)\, 
dx\nonumber \\
&&\hspace{1.0cm}+\int f\left(E_{t,\cdot}g(X_0)-g\right)\cdot\frac1t 
\left(E_\cdot\prod_{i=1}^{n\cdot d}\left|{\sf e}+\nabla_{d,W_0}Z_{-t} 
\circ u^{-1}\vphantom{\dot{f}}\right|_i-\1\right)\, d\bmu\nonumber 
\\ 
&&\hspace{.5cm}=T_1(t)+T_2(t)+T_3(t)\, . 
\end{eqnarray} 
Now we take into consideration that because of (\ref{4.12}) we have 
$\lim_{t\downarrow 0}(E_{t,\cdot}g(X_0)-g)=0$ $\lambda_F$-a.e. on 
$D^n$. With this, $T_1(t)$ converges to zero as $t\downarrow 0$ by 
(ix) of the present proposition.  Let us focus on $T_2(t)$. Taking 
into consideration condition (8) of this section, we get 
\begin{eqnarray*}
(T_2(t))^2&&\hspace{-.5cm}\le\int f^2\left(E_{t,\cdot}g(X_0)-g 
\right)^2\, dx\times\\ 
&&\hspace{.5cm}\times\int\left(E_\cdot\left(\frac{m\left(W_0^{-t} 
\circ u^{-1}\right)-m}{t}\left(\prod_{i=1}^{n\cdot d}\left|{\sf e}+ 
\nabla_{d,W_0}Z_{-t}\circ u^{-1}\vphantom{\dot{f}}\right|_i-\1\right) 
\right)\right)^2 \, dx \\ 
&&\hspace{-.5cm}\le\int f^2\left(E_{t,\cdot}g(X_0)-g\right)^2\, dx 
\times\\ 
&&\hspace{.5cm}\times\int E_\cdot\left(\frac{m\left(W_0^{-t}\circ 
u^{-1}\right)-m}{t}\right)^2\cdot E_\cdot\left(\prod_{i=1}^{n\cdot d} 
\left|{\sf e}+\nabla_{d,W_0}Z_{-t}\circ u^{-1}\vphantom{\dot{f}}\right 
|_i-\1\right)^2\, dx \\ 
&&\hspace{-.5cm}\le 4\|m\|^2\cdot\left(\int f^2\left(E_{t,\cdot}g 
(X_0)-g\right)^2\, dx\right)^2\times \\ 
&&\hspace{.5cm}\times\frac{1}{t^2}\int E_\cdot\left(\prod_{i=1}^{n 
\cdot d}\left|{\sf e}+\nabla_{d,W_0}Z_{-t}\circ u^{-1}\vphantom{\dot 
{f}}\right|-\1\right)^2\, dx\stack{t\downarrow 0}{\lra}0\vphantom{\int}
\end{eqnarray*} 
the limit in the last line is a consequence of (\ref{4.1}), (\ref{4.12}) 
and condition (viii) of Theorem \ref{Theorem4.2}. For $T_3(t)$, we recall 
condition (7) of this section. Relations (\ref{4.1}), (\ref{4.12}) and  
assumption (viii) of Theorem \ref{Theorem4.2} show also here that 
\begin{eqnarray*} 
(T_3(t))^2&&\hspace{-.5cm}\le\int f^2\left(E_{t,\cdot}g(X_0)-g\right)^2 
\, d\bmu\times\\ 
&&\hspace{.5cm}\times\frac{1}{t^2}\int E_\cdot\left(\prod_{i=1}^{n\cdot d} 
\left|{\sf e}+\nabla_{d,W_0}Z_{-t}\circ u^{-1}\vphantom{\dot{f}}\right|-\1 
\right)^2\, d\bmu\stack{t\downarrow 0}{\lra}0\, . 
\end{eqnarray*} 
From (\ref{4.13}) we get now 
\begin{eqnarray*} 
\lim_{t\downarrow 0}\frac1t\int f\left(E_{t,\cdot}g(X_0)-g\right) 
(E_\cdot\rho_{-t}-\1)\, d\bmu=0\, ,\quad f,g\in C_b(D^n). 
\end{eqnarray*} 
Furthermore, it holds that 
\begin{eqnarray*} 
P_{\sbmu}(X_t\in dx)=\left.E_{\sbmu}\left(\rho_{-t}\right|X_0=x 
\right)\cdot\bmu(dx)=E_x\rho_{-t}\cdot\bmu(dx)\, . 
\end{eqnarray*} 
We get with this  
\begin{eqnarray}\label{4.14} 
&&\hspace{-.5cm}-\lim_{t\downarrow 0}\frac1t\int f\left(E_{t,\cdot} 
g(X_0)-g\right)\, d\bmu\nonumber \\ 
&&\hspace{.5cm}=-\lim_{t\downarrow 0}\frac1t\int f\left(E_{t,\cdot} 
g(X_0)-g\right)E_\cdot\rho_{-t}\, d\bmu\nonumber \\ 
&&\hspace{1.0cm}+\lim_{t\downarrow 0}\frac1t\int f\left(E_{t,\cdot} 
g(X_0)-g\right)(E_\cdot\rho_{-t}-\1)\, d\bmu\nonumber \\ 
&&\hspace{.5cm}=-\lim_{t\downarrow 0}\frac1t\int fE_{t,\cdot}\left( 
g(X_0)-g\right)\, dP_{\sbmu}(X_t\in\, \cdot\, )\nonumber \\ 
&&\hspace{.5cm}=-\lim_{t\downarrow 0}\frac1t\int f(X_t)(g(X_0)-g(X_t) 
)\, dP_{\sbmu}\nonumber \\ 
&&\hspace{.5cm}=\lim_{t\downarrow 0}\frac1t\int f(X_0)(g(X_t)-g(X_0) 
)\, dP_{\sbmu}\nonumber \\ 
&&\hspace{1.0cm}+\lim_{t\downarrow 0}\frac1t\int\left(f(X_t)-f(X_0) 
\right)\left(g(X_t)-g(X_0)\right)\, dP_{\sbmu}\nonumber \\ 
&&\hspace{.5cm}=A_wg(f)+R(f,g)\, ,\quad f\in {\cal C},\ g\in D(A_w) 
\vphantom{\int}. 
\end{eqnarray} 
We obtain $g\in D(A^\ast_w,{\cal C})$ as well as $A_w^\ast g(f)=A_w 
g(f)+R(f,g)$. 
\medskip 

\nid 
(b) By (x) and the definition of $D(A_w)$ we have $m\in C_b^1(D^n)$ 
and by (v) it holds that 
\begin{eqnarray}\label{4.15} 
B^\ast_wm(h)&&\hspace{-.5cm}=-\lim_{t\downarrow 0}\frac1t\int\frac{h} 
{m}\left(E_\cdot m(X_{-t})-m\right)\, \bmu(dx)\nonumber \\ 
&&\hspace{-.5cm}=-\lim_{t\downarrow 0}\frac1t\int\frac{h(X_0)}{m(X_0)} 
\left(m(X_{-t})-m(X_0)\right)\, dP_{\sbmu}\nonumber \\ 
&&\hspace{-.5cm}=-\lim_{t\downarrow 0}\frac1t\int\frac{h(X_0)}{m(X_0)} 
\left(m(X_{-t})-m(X_0)\right)\rho_{-t}(X)\, dP_{\sbmu}\nonumber \\ 
&&\hspace{-.0cm}-\lim_{t\downarrow 0}\frac1t\int\frac{h(X_0)}{m(X_0)} 
\left(m(X_{-t})-m(X_0)\right)(\1-\rho_{-t}(X))\, dP_{\sbmu}\nonumber \\ 
&&\hspace{-.5cm}=-\lim_{t\downarrow 0}\frac1t\int\frac{h(X_t)}{m(X_t)} 
\left(m(X_0)-m(X_t)\right)\, dP_{\sbmu}\nonumber \\ 
&&\hspace{-.0cm}-\lim_{t\downarrow 0}\frac1t\int\frac{h(X_0)}{m(X_0)} 
\left(m(X_{-t})-m(X_0)\right)(\1-\rho_{-t}(X))\, dP_{\sbmu} 
\end{eqnarray} 
for all $h\in C_b(D^n)$ such that $h/m\in C_b(D^n)$. The last line of 
(\ref{4.15}) is zero, by 
\begin{eqnarray*} 
&&\hspace{-.5cm}\left(\frac1t\int\frac{h(X_0)}{m(X_0)}\left(m(X_{-t}) 
-m(X_0)\right)(\1-\rho_{-t}(X))\, dP_{\sbmu}\right)^2 \\ 
&&\hspace{.5cm}=\left(\frac1t\int\frac{h(X_0)}{m(X_0)}\left(m(X_{-t}) 
-m(X_0)\right)(m(X_0)-m(X_0)\, \rho_{-t}(X))\, dP_{\lambda_F}\right)^2 
 \\ 
&&\hspace{.5cm}\le\left\|\frac{h}{m}\right\|^2\frac1t\int(\left(m(X_{ 
-t})-m(X_0)\right)^2\, dP_{\lambda_F}\cdot\frac1t\int(m(X_0)-m(X_0)\, 
\rho_{-t}(X))^2\, dP_{\lambda_F} \\ 
&&\hspace{.5cm}\le\left\|\frac{h}{m}\right\|^2\frac1t\int(\left(m(X_{ 
-t})-m(X_0)\right)^2\, dP_{\lambda_F}\times\\ 
&&\hspace{1.0cm}\times\frac1t\int\left(m(X_0)-m\left(X_{-t}\right)\cdot 
\prod_{i=1}^{n\cdot d}\left|{\sf e}+\nabla_{d,W_0}Z_{-t}\circ u^{-1} 
\vphantom{\dot{f}}\right|_i\right)^2\, dP_{\lambda_F}\stack{t 
\downarrow 0}{\lra}0 
\end{eqnarray*} 
as well as condition (8) of this section and hypotheses (vi) and (viii) 
of Theorem \ref{Theorem4.2}. Comparing the second last line of 
(\ref{4.15}) with (\ref{4.14}) it turns out that 
\begin{eqnarray*} 
B^\ast_wm(h)=A_wm(h/m)+R(h/m,m)=A^\ast_wm(h/m)
\end{eqnarray*} 
if $m\in D(A_w)$ and $h/m\in {\cal C}$ which implies $m\in D(A^\ast_w, 
{\cal C})$. 
\qed 
 
\section{Relative Compactness of Particle Systems} 
\setcounter{equation}{0} 

Let us assume that for all $n\in {\Bbb N}$ we are given a system 
consisting of $n$ particles $\left\{X^n_1,\ldots ,X^n_n\right\}$, each 
of which being a stochastic process with time domain ${\Bbb R}$ taking 
values in $D$. We consider the measure valued stochastic process ${\bf 
X}$ given by ${\bf X}^n_t:=\frac1n\sum_{i=1}^n\delta_{(X^n_i)_t}$, $t\in 
{\Bbb R}$.  

In addition, we introduce the $n\cdot d$-dimensional process $X^n\equiv 
(X^n_t)_{t\in {\Bbb R}}=\left((X^n_1)_t,\ldots ,\right.$ $\left.(X^n_n 
)_t\right)_{t\in {\Bbb R}}$ which we assume to be of the form $X^n=W^n+ 
A^n$ satisfying, for every $n\in {\Bbb N}$, the hypotheses of Section 1. 
Furthermore, we suppose that the distribution of $X^n$ is invariant 
under the permutation of the $n$ particles provided that the distribution 
of $W^n_0=\left((W_1)_0,\ldots ,(W_n)_0\right)$ is invariant under the 
permutation of the $n$ $d$-dimensional random variables $\left\{(W_1)_0, 
\ldots ,(W_n)_0\right\}$. 

Let ${\cal M}_1(\overline{D})$ denote the space of all probability 
measures on $(\overline{D},{\cal B}(\overline{D}))$. Moreover, abbreviate 
for $g\in C(\overline{D})$ and $\mu\in {\cal M}_1(\overline{D})$ the 
integral $\int g\, d\mu$ by $(g,\mu)$ and define 
\begin{eqnarray*}
&&\hspace{-.5cm}\t C^1({\cal M}_1(\overline{D})) \\ 
&&\hspace{.5cm}:=\left\{F(\mu)=\vp\left((g_1,\mu),\ldots ,(g_k,\mu)\right), 
\ g_1,\ldots ,g_k\in C^1(\overline{D}),\ \vp\in C_b^1({\Bbb R}^k),\ k\in 
{\Bbb Z}_+\right\}\, . 
\end{eqnarray*}

It is our primary goal to formulate conditions in order to verify relative 
compactness of the family $F({\bf X}^n)\equiv (F({\bf X}^n_t))_{t\ge 0}$, 
$F\in\t C^1({\cal M}_1(\overline{D}))$, $n\in {\Bbb N}$. To derive relative 
compactness of the family ${\bf X}^n\equiv ({\bf X}^n_t)_{t\ge 0}$, $n\in 
{\Bbb N}$, from this is then standard. 

As in Section 4, we use the symbol $Z\equiv Z^n$ to abbreviate $A-Y\equiv 
A^n-Y^n$. Referring to the number of particles, here we also use the 
notation $A^n=(A^1)^n+(A^2)^n$. In contrast to Section 4, where we considered 
$X$ the given process, we look in this section at the particle systems as at 
processes constructed from Brownian motion. Therefore we also refer to the 
measures $Q_{\sbnu}\equiv Q^{(n)}_{\sbnu_n}$ and $Q_{\cdot}\equiv Q^{(n)}_{ 
\cdot}$ rather than to $P_{\sbmu}$ and $P_{\cdot}$. 
\begin{theorem}\label{Theorem5.1} 
For every $n\in {\Bbb N}$, let $X^n$ be a stochastic process on the 
probability space $(\Omega^{(n)},{\cal F}^{(n)},Q^{(n)}_{\sbnu_n})$ 
satisfying the hypotheses of Theorem \ref{Theorem1.11} or, respectively, 
of Theorem \ref{Theorem1.12}. Let $m_n$ be as in  Corollary 
\ref{Corollary1.13}. In addition to condition (2) (i) of Subsection 1.2 
or, respectively, condition (ii) of Theorem \ref{Theorem1.12}, suppose 
that for every $n\in {\Bbb N}$ and $Q^{(n)}_{\sbnu_n}$-a.e. $W\in 
\Omega^{(n)}$, $(A^1)^n(W)$ is equi-continuous on ${\Bbb R}$. Assume the 
following. 
\begin{itemize} 
\item[(i)] For all $n\in {\Bbb N}$, $Q^{(n)}_{\sbnu_n}$-almost never two 
of the particles $\{X^n_1,\ldots ,X^n_n\}$ jump at the same time. 
\item[(ii)] There exists a sequence $c_n>0$, $n\in {\Bbb N}$, such that 
for all $k\in {\Bbb N}$, all $j\in\{1,\ldots ,n\cdot d\}$, all $\delta\in 
[0,1]$, and all $x\in D^n$, 
\begin{eqnarray*} 
{\rm ess}\sup_{\hspace{-.6cm}\Omega}\left(\sup_{t\in [0,\delta]}\left| 
\left\langle e_j,\nabla_{d,W_0}Z^n_{-t}\right\rangle_F\vphantom{\dot{f}} 
\right|^k\right)\le c_n\delta^k 
\end{eqnarray*} 
and $\lim_{\delta\downarrow 0}\limsup_{n\in {\Bbb N}}c_n(1+\delta)^n=0$. 
Here, the {\rm ess \hspace{-.5mm}sup} refers to the measure $Q^{(n)}_{ 
\sbnu_n}$. 
\item[(iii)] Denoting by $E^{(n)}$ the expectation with respect to $Q^{(n) 
}$, we have 
\begin{eqnarray*} 
\lim_{\delta\downarrow 0}\limsup_{n\in {\Bbb N}}\int_{D^n} E^{(n)}_x\left 
(\sup_{t\in [0,\delta]}\left|m_n(X^n_{-t}-Y^n_{-t})-m_n(x)\vphantom{\dot 
{f}}\right|\right)\, dx=0\, . 
\end{eqnarray*} 
\end{itemize} 
(a) Let $F\in\t C^1({\cal M}_1(\overline{D}))$. The family of stochastic 
processes $F({\bf X}^n)\equiv (F({\bf X}^n_t))_{t\ge 0}$, $n\in {\Bbb N}$, 
is relatively compact with respect to the topology of weak convergence of 
probability measures over the Skorokhod space $D_{[-\|F\|,\|F\|]}({\Bbb 
R})$. \\ 
(b) The family of stochastic processes ${\bf X}^n\equiv ({\bf X}^n_t)_{t\ge 
0}$, $n\in {\Bbb N}$, is relatively compact with respect to the topology of 
weak convergence of probability measures over the Skorokhod space $D_{{\cal 
M}_1(\overline{D})}({\Bbb R})$.
\end{theorem}
Proof. {\it Step 1 } Let $f\in C^1\left(\overline{D^n}\right)$ be defined 
by $f(x_1,\ldots ,x_n):=F\left(\frac1n\sum_{i=1}^n\delta_{x_i}\right)$, 
$x_1,\ldots ,x_n\in\overline{D}$. In this step we introduce the objects 
used in Chapter 3 of \cite{EK86}, Theorem 8.6 and Remark 8.7. 

We fix $n\in {\Bbb N}$ and drop the index $n$ from the notation since no 
ambiguity is possible. Let $0<\delta<1$ and introduce 
\begin{eqnarray*}
&&\hspace{-.5cm}\gamma(\delta)\equiv\gamma_n(\delta;f):=\sup_{t\in {\Bbb 
R},\, 0\le u\le\delta}\left(f(X_{t+u})-f(X_t)\right)^2\, . 
\end{eqnarray*}

Furthermore, set ${\cal A}=[-\delta,\delta]$. For $\alpha\in {\cal A}$, 
let $\Omega_\alpha$ be the collection of all $u^{-1}(X)\in\Omega$ 
satisfying the following. 
\begin{itemize}
\item[(j)] There are sequences $s_m\in {\Bbb R}$ and $\alpha_m\in {\Bbb 
R}$, $m\in {\Bbb N}$, with $\lim_{m\to\infty}\alpha_m=\alpha$ and 
\begin{eqnarray*} 
\lim_{m\to\infty}\left(f(X_{s_m+|\alpha_m|})-f(X_{s_m})\right)^2=\sup_{t 
\in {\Bbb R},\, 0\le u\le\delta}\left(f(X_{t+u})-f(X_t)\right)^2=:A\, . 
\end{eqnarray*}
\item[(jj)] 
\begin{eqnarray*} 
\mbox{\rm sign}(\alpha)\cdot\lim_{m\to\infty}\left(f(X_{s_m+|\alpha_m|}) 
-f(X_{s_m})\right)\ge 0\, . 
\end{eqnarray*}
\item[(jjj)] If $\alpha\le 0$ there is no pair of sequences $s_m'\in{\Bbb 
R}$ and $\alpha'_m\in {\Bbb R}$, $m\in {\Bbb N}$, with $\lim_{m\to\infty} 
\alpha'_m=\alpha$ such that 
\begin{eqnarray*} 
\lim_{m\to\infty}\left(f(X_{s_m'+|\alpha_m'|})-f(X_{s_m'})\right)>0\, . 
\end{eqnarray*}
\item[(jv)] 
\begin{eqnarray*} 
\limsup_{m\to\infty}\left(f(X_{t_m+|\beta_m|})-f(X_{t_m})\right)^2<A
\end{eqnarray*}
for all sequences $t_m\in {\Bbb R}$ and $\beta_m\in {\Bbb R}$, $m\in 
{\Bbb N}$, with $\lim_{m\to\infty}\beta_m=\beta$ and $|\beta|<|\alpha|$. 
\end{itemize}

Let $\ve >0$. Set $t(\alpha):=|\alpha|-\ve$, $\alpha\in {\cal A}$, let 
$\tau^\alpha\equiv\tau^\alpha (\ve)$ be the random time defined in 
Corollary \ref{Corollary1.15}, and denote by $\1$ the function which 
is constant one. We have the hypotheses (i) and (ii) of Corollary 
\ref{Corollary1.15} and obtain 
\begin{eqnarray}\label{5.1}
\rho_{-\tau^\alpha (\ve)}(X)&&\hspace{-1cm}\stack{\ve\downarrow 0}{\lra} 
\frac{\D m\left(X_{-|\alpha|}-Y_{-|\alpha|}\circ u^{-1}\vphantom{l^1} 
\right)}{\D m\left(u^{-1}_0\vphantom{l^1}\right)}\cdot\vphantom{\left\{ 
\sum_{\tau_k}\right\}}\prod_{i=1}^{n\cdot d}\left|{\sf e}+\nabla_{d,W_0} 
Z_{-|\alpha|}\circ u^{-1}\vphantom{\dot{f}}\right|_i\nonumber \\ 
&&\hspace{-.5cm}=:\rho_{-|\alpha|}(X)\quad P_{\sbmu}{\rm -a.e.}\vphantom 
{\int}
\end{eqnarray}
In addition, using the notation of the definition of the sets 
$\Omega_\alpha$, we get 
\begin{eqnarray*}
&&\hspace{-.5cm}E_{\sbmu}\left(\gamma(\delta)\right)=E_{\sbmu}\left( 
\sup_{t\in {\Bbb R}\, ,0\le u\le\delta}\left(f(X_{t+u})-f(X_t)\right)^2 
\right) \\ 
&&\hspace{.5cm}\le E_{\sbmu}\left(\lim_{m\to\infty}\left(f(X_{s_m+| 
\alpha_m|})-f(X_{s_m})\right)^2\right) \\ 
&&\hspace{.5cm}=\lim_{m\to\infty}E_{\sbmu}\left(\vphantom{\lim_{m\to 
\infty}}\left(f(X_{s_m+|\alpha_m|})-f(X_{s_m})\right)^2\right) \\ 
&&\hspace{.5cm}=\lim_{m\to\infty}\left(\vphantom{\lim_{m\to\infty}}E_{ 
\sbmu}\left((f(X_{s_m+|\alpha_m|}))^2-(f(X_{s_m}))^2\right)\right. 
\vphantom{\int} \\ 
&&\hspace{4.5cm}\left.-2E_{\sbmu}\left(f(X_{s_m})\cdot\left(f(X_{s_m+ 
|\alpha_m|})-f(X_{s_m})\right)\right)\vphantom{\lim_{m\to\infty}}\right) 
\vphantom{\int} \\ 
&&\hspace{.5cm}\le\liminf_{m\to\infty}\left(\vphantom{\lim_{m\to\infty} 
}\left|E_{\sbmu}\left((f(X_{{s_m}+|\alpha_m|}))^2-(f(X_{s_m}))^2\right) 
\right|\right.\vphantom{\int} \\ 
&&\hspace{4.5cm}\left.+2\|f\|\cdot E_{\sbmu}\left|f(X_{s_m+|\alpha_m|}) 
-f(X_{s_m})\right|\vphantom{\lim_{m\to\infty}}\right)\vphantom{\int} \\ 
&&\hspace{.5cm}=\liminf_{m\to\infty}\left|E_{\sbmu}\left((f(X_{{s_m}+ 
\alpha_m}))^2-(f(X_{s_m}))^2\right)\right|\vphantom{\int} \\ 
&&\hspace{4.5cm}+2\|f\|\cdot\lim_{m\to\infty}E_{\sbmu}\left({\rm sign} 
(\alpha)\left(f(X_{s_m+|\alpha_m|})-f(X_{s_m})\right)\vphantom{l^1_1} 
\right)\vphantom{\int} \\ 
&&\hspace{.5cm}\le\liminf_{m\to\infty}\left(\vphantom{\lim_m}\left|E_{ 
\sbmu}\left((f(X_{{s_m}}))^2\left(\rho_{-\tau^\alpha}-\1\right)\right) 
\right|\right.\vphantom{\int} \\ 
&&\hspace{4.5cm}\left.+E_{\sbmu}\left|(f(X_{{s_m}+|\alpha_m|}))^2-(f( 
X_{s_m+\tau^\alpha}))^2\right|\vphantom{\lim_m}\right)\vphantom{\int}\\ 
&&\hspace{2.5cm}+2\|f\|\cdot\lim_{m\to\infty}\left(\vphantom{\lim_m} 
E_{\sbmu}\left({\rm sign}(\alpha)f(X_{s_m})\left(\rho_{-\tau^\alpha}- 
\1\right)\vphantom{l^1_1}\right)\right.\vphantom{\int} \\ 
&&\hspace{4.5cm}\left.+2\|f\|\cdot E_{\sbmu}\left|f(X_{s_m+|\alpha_m|}) 
-f(X_{s_m+\tau^\alpha})\right|\vphantom{\lim_m}\right)\, .\vphantom 
{\int} 
\end{eqnarray*}
We recall $\tau^\alpha=|\alpha|-\ve$ and (\ref{5.1}). Letting $\ve 
\downarrow 0$ we obtain 
\begin{eqnarray*} 
&&\hspace{-.5cm}E_{\sbmu}\left(\gamma(\delta)\right)\le 3\|f\|^2\cdot 
E_{\sbmu}\left|\rho_{-|\alpha|}(X)-\1\right|+4\|f\|\cdot\liminf_{m\to 
\infty}E_{\sbmu}\left|f(X_{s_m+|\alpha_m|})-f(X_{s_m+|\alpha|-})\right| 
\, .\vphantom{\int}
\end{eqnarray*}
For the second item of the right-hand side, we use now $\lim_{m\to\infty 
}|\alpha_m|=|\alpha|$ and the fact that $W_t+A^1_t$ is equi-continuous 
on $t\in [s_m,s_m+|\alpha|+1]$, $m\in {\Bbb N}$, by Paul L\'evy's modulus 
of continuity of $W$ and by the equi-continuity of $A^1$. Furthermore we 
pay attention to the particular representation of $f$, $f(x_1,\ldots ,x_n 
):=F\left(\frac1n\sum_{i=1}^n\delta_{x_i}\right)$, $x_1,\ldots ,x_n\in 
\overline{D}$ where $F(\mu)=\vp\left((g_1,\mu),\ldots ,(g_k,\mu)\right)$, 
$g_1,\ldots ,g_k\in C^1(\overline{D})$, $\vp\in C_b^1({\Bbb R}^k)$, $k\in 
{\Bbb Z}_+$. Also, we take into consideration that every jump of $X^n$ is 
$P_{\sbmu}$-a.e. caused by only one of the $d$-dimensional processes 
$X^n_1,\ldots ,X^n_n$, cf. condition (i) of this theorem. Denoting by 
diam$(D)$ the diameter of $D$ we get 
\begin{eqnarray}\label{5.2} 
&&\hspace{-.5cm}E_{\sbmu}\left(\gamma(\delta)\right)\le 3\|f\|^2\cdot E_{ 
\sbmu}\left|\rho_{-|\alpha|}(X)-\1\right|+\frac{{\rm diam}(D)}{n}\, 4\|f 
\|\cdot\|\nabla\vp\|\cdot\sum_{j=1}^k\|\nabla g_j\|\, . 
\end{eqnarray}

\nid
{\it Step 2 } We re-add the index $n$ to the notation and note that 
the density function $m$ depends on the number of particles $n$. In 
this step we aim to prove 
\begin{eqnarray}\label{5.3}
\lim_{\delta\downarrow 0}\limsup_{n\to\infty}E^n_{\sbmu_n}\left( 
\gamma_n(\delta)\right)=0\, . 
\end{eqnarray} 
We note that, independent of $n\in {\Bbb N}$, we have $\|f\|=\|F\|$. 
From (\ref{5.1}) and (\ref{5.2}) we obtain 
\begin{eqnarray}\label{5.4}
&&\hspace{-.5cm}\frac{1}{3\|F\|^2}\lim_{\delta\downarrow 0}\limsup_{n 
\to\infty}E^n_{\sbmu_n}\left(\gamma_n(\delta)\right)\le\lim_{\delta 
\downarrow 0}\limsup_{n\to\infty}E^n_{\sbmu_n}\left(\sup_{t\in [0, 
\delta]}\left|\frac{\D m_n\left(X^n_{-t}-Y^n_{-t}\circ u^{-1}\vphantom 
{l^1}\right)}{\D m_n\left(u^{-1}_0\vphantom{l^1}\right)}\vphantom 
{\prod_{i=1}^{n\cdot d}}\, \times\right.\right.\nonumber \\ 
&&\hspace{1.0cm}\left.\left.\times\prod_{i=1}^{n\cdot d}\left|{\sf e}+ 
\nabla_{d,W_0}Z^n_{-t}\circ u^{-1}\vphantom{\dot{f}}\right|_i-\1\right| 
\right)\nonumber \\ 
&&\hspace{.5cm}=\lim_{\delta\downarrow 0}\limsup_{n\to\infty}E^{(n) 
}_{\sbnu_n}\left(\sup_{t\in [0,\delta]}\left|\frac{\D m_n\left(X^n_{-t} 
-Y^n_{-t}\vphantom{l^1}\right)}{\D m_n\left(W_0\vphantom{l^1}\right)} 
\cdot\prod_{i=1}^{n\cdot d}\left|{\sf e}+\nabla_{d,W_0}Z^n_{-t} 
\vphantom{\dot{f}}\right|_i-\1\right|\right)\nonumber \\ 
&&\hspace{.5cm}\le\lim_{\delta\downarrow 0}\limsup_{n\to\infty}E^{ 
(n)}_{\sbnu_n}\left(\sup_{t\in [0,\delta]}\left|\frac{m_n(X^n_{-t}- 
Y^n_{-t})-m_n\left(W_0\vphantom{l^1}\right)}{\D m_n\left(W_0\vphantom 
{l^1}\right)}\right|\right)\nonumber \\ 
&&\hspace{1.0cm}+\lim_{\delta\downarrow 0}\limsup_{n\to\infty}E^{(n 
)}_{\sbnu_n}\left(\sup_{t\in [0,\delta]}\left|\frac{m_n(X^n_{-t}-Y^n_{ 
-t})-m_n\left(W_0\vphantom{l^1}\right)}{\D m_n\left(W_0\vphantom{l^1} 
\right)}\vphantom{\left(\prod_{i=1}^{n\cdot d}\right)}\times\right. 
\right.\nonumber \\ 
&&\hspace{5.0cm}\left.\left.\times\left(\prod_{i=1}^{n\cdot d}\left| 
{\sf e}+\nabla_{d,W_0}Z^n_{-t}\vphantom{\dot{f}}\right|_i-\1\right) 
\right|\right)\nonumber \\ 
&&\hspace{1.0cm}+\lim_{\delta\downarrow 0}\limsup_{n\to\infty}E^{(n 
)}_{\sbnu_n}\left(\sup_{t\in [0,\delta]}\left|\prod_{i=1}^{n\cdot d} 
\left|{\sf e}+\nabla_{d,W_0}Z^n_{-t}\vphantom{\dot{f}}\right|_i-\1 
\right|\right)\nonumber \\ 
&&\hspace{.5cm}=\lim_{\delta\downarrow 0}\limsup_{n\to\infty}T_1'( 
\delta)+\lim_{\delta\downarrow 0}\limsup_{n\to\infty}T_2'(\delta)+ 
\lim_{\delta\downarrow 0}\limsup_{n\to\infty}T_3'(\delta)\, . 
\vphantom{\int}
\end{eqnarray}

Let us take a look at the items of the right-hand side. It holds that 
\begin{eqnarray}\label{5.5}
T_1'(\delta)=\int_{D^n}E^{(n)}_x\left(\sup_{t\in [0,\delta]}\left|m_n 
(X^n_{-t}-Y^n_{-t})-m_n(x)\vphantom{\dot{f}}\right|\right)\, dx\, ,  
\end{eqnarray}
\begin{eqnarray}\label{5.6}
T_2'(\delta)&&\hspace{-.5cm}=\int_{D^n}E^{(n)}_x\left(\sup_{t\in [0, 
\delta]}\left|\left(m_n(X^n_{-t}-Y^n_{-t})-m_n(x)\right)\vphantom{ 
\left(\prod_{i=1}^{n\cdot d}\right)}\times\right.\right.\nonumber \\ 
&&\hspace{3.5cm}\times\left.\left.\left(\prod_{i=1}^{n\cdot d}\left| 
{\sf e}+\nabla_{d,W_0}Z^n_{-t}\vphantom{\dot{f}}\right|_i-\1\right) 
\right|\right)\, dx\nonumber \\ 
&&\hspace{-.5cm}\le\int_{D^n}E^{(n)}_x\left(\sup_{t\in [0,\delta]} 
\left|m_n(X^n_{-t}-Y^n_{-t})-m_n(x)\vphantom{\dot{f}}\right|\right)\, 
dx\times\nonumber \\ 
&&\hspace{3.5cm}\times{\rm ess}\sup_{\hspace{-.6cm}\Omega}\left( 
\sup_{t\in [0,\delta]}\left|\prod_{i=1}^{n\cdot d}\left|{\sf e}+ 
\nabla_{d,W_0}Z^n_{-t}\vphantom{\dot{f}}\right|_i-\1\right|\right)\, 
,\qquad
\end{eqnarray}
and 
\begin{eqnarray}\label{5.7}
T_3'(\delta)\le{\rm ess}\sup_{\hspace{-.6cm}\Omega}\left(\sup_{t\in 
[0,\delta]}\left|\prod_{i=1}^{n\cdot d}\left|{\sf e}+\nabla_{d,W_0} 
Z^n_{-t}\vphantom{\dot{f}}\right|_i-\1\right|\right)\, . 
\end{eqnarray}
We apply the hypotheses (ii) and (iii) of the present proposition to 
(\ref{5.5})-(\ref{5.7}). From (\ref{5.4}) we obtain (\ref{5.3}) in this 
way. 
\medskip 

\nid 
{\it Step 3 } We finish the proof of part (a). Let $-\infty<S<T<\infty$ 
and $0<\delta\le 1$ and let $({\cal F}^n_t)_{t\in {\Bbb R}}$ denote the 
filtration generated by $X^n$ on $(-\infty,t]$. For $0\le u\le\delta$ and 
$S\le t\le T$ we have 
\begin{eqnarray*}
E^n_{\sbmu_n}\left.\left(\left(f(X^n_{t+u})-f(X^n_t)\right)^2\right|{\cal 
F}^n_t\right)\le E^n_{\sbmu_n}\left(\gamma_n(\delta)\left.\vphantom{{\cal 
F}^n_{S,t}}\right|{\cal F}^n_t\right)\, . 
\end{eqnarray*}
We get now the claim from (\ref{5.3}), Chapter 3 of \cite{EK86}, Theorem 
8.6 and Remark 8.7. 
\medskip 

\nid 
{\it Step 4 } Part (b) follows now from (a) and Chapter 3 of \cite{EK86}, 
Theorem 9.1. 
\qed 
\medskip 

The above theorem and the corollary cover particle systems with a certain 
abstract drift, Brownian noise such that the drift may not be adapted, and, 
for given Brownian trajectory, non-random jumps. As discussed in the 
introduction, the jump mechanism is adopted from particle systems 
approximating Boltzmann type equations. 

The conditions formulated here in order to verify weak compactness of the 
particle system are of abstract character. They may give a guideline in 
more concrete situations. Conditions (i) and (ii) refer to the 
construction of the particle system. They have to be verified by taking 
into consideration 	the particular situation. The interesting condition 
of Theorem \ref{Theorem5.1} is (iii). We are interested in the following 
simple but fairly instructive example. 
\medskip

Let $\nabla_d$ and $\Delta_d$ denote the $d$-dimensional gradient and 
Laplace operator. For the sake of simplicity, let us suppose that $D$ 
is of Lebesgue measure one which also means that $\lambda^{n\cdot d} 
(D^n)\equiv\lambda_F(D^n)=1$. Choose $c\in (0,1)$ and let Let $\hat{H} 
\equiv\hat{H}_c$ be the set of all probability measures $d(x)\, dx$ on 
$(\overline{D},{\cal B}(\overline{D}))$ such that 
\begin{eqnarray}\label{5.8}
d\in C^2\left(\overline{D}\right)\ \mbox{\rm with}\ c\le d\le c^{-1}\ 
\mbox {\rm and}\ \|\nabla_d d\|\le c^{-1}\, .
\end{eqnarray} 
Let $\bnu$ be a probability measure on 
$\left({\cal M}_1\left(\overline{D}\right),{\cal B}\left({\cal M}_1\left 
(\overline{D}\right)\right)\right)$ such that $\bnu\left({\cal M}_1\left 
(\overline{D}\right)\setminus\hat{H}\right)=0$. We define the measures 
$\t\bnu_n(dx)=m_n(x)\, dx$ on $(\overline{D^n},{\cal B}(\overline{D^n})) 
$ and $\bnu_n(d\mu^n)$ on the set of all empirical measure $\mu^n=\frac 
1n\sum_{i=1}^n\delta_{x_i}$, $x_1,\ldots,x_n\in\overline{D}$, by $\bnu_n 
(d\mu^n)=\t\bnu_n(dx)$, $x=(x_1,\ldots ,x_n)$ and 
\begin{eqnarray}\label{5.9}
m_n(x)&&\hspace{-.5cm}:=\int_{\mu=d(x)\, dx\in\hat{H}}\frac{1}{\|d^{\frac 
1n}\|_{L^1(D)}^n}\prod_{i=1}^n\left(d(x_i)\vphantom{l^1}\right)^{\frac1n} 
\, \bnu(d\mu)\nonumber \\ 
&&\hspace{-.5cm}=\int_{\mu=d(x)\, dx\in\hat{H}}\frac{1}{\|d^{\frac1n}\|_{ 
L^1(D)}^n}\exp\left\{(\log(d),\mu^n)\vphantom{l^1}\right\}\, \bnu(d\mu)\, 
,\quad n\in {\Bbb N}. 
\end{eqnarray} 
The following two remarks shall discuss the choice of the initial 
distribution relative to $n\in {\Bbb N}$. On the one hand the distribution 
of a typical individual particle has order $d^{1/n}$. On the other hand, 
Remark (1) below shows that the overall initial distribution behaves 
moderately with respect to $n\in {\Bbb N}$. In addition, Remark (2) below 
says that even the gradient of the overall initial density as well as the 
generator of the Brownian motion part, applied to the overall initial 
density, behave moderately with respect to $n\in {\Bbb N}$. These are 
important technical features. 
\medskip 

\nid 
{\bf Remarks. (1)} Since $1+(\log(d))/n\le d^{\frac1n}\le 1+(\log(d))/n+  
\left((\log(d))/n\right)^2$ for sufficiently large $n\in {\Bbb N}$ for 
all $d\in\hat{H}$ by (\ref{5.8}), we have for such $n$ 
\begin{eqnarray}\label{5.10}
\left(1+\frac1n\int\log(d)\, dx\right)^n\le\left\|d^{\frac1n}\right 
\|^n_{L^1(D)}\le\left(1+\frac1n\int\log(d)\, dx+\frac{1}{n^2}\int(\log 
(d))^2\, dx\right)^n\, . 
\end{eqnarray} 
Relation (\ref{5.9}) and the weak law of large numbers show now that 
$\bnu_n$ converges to $\bnu$ as $n\to\infty$ in the weak topology of 
probability measures over ${\cal M}_1(\overline{D})$. 
\medskip 

\nid 
{\bf (2)} Denoting by $\nabla$ and $\Delta$ the $n\cdot d$-dimensional 
gradient and Laplace operator we obtain
\begin{eqnarray*}
\nabla m_n(x)=\left(\frac1n\int_{\mu=d(x)\, dx\in\hat{H}}\frac{1}{\|d^{ 
\frac1n}\|_{L^1(D)}^n}\cdot\frac{\nabla_d d(x_i)}{d(x_i)}\prod_{j=1}^n 
\left(d(x_j)\vphantom{l^1}\right)^{\frac1n}\, \bnu(d\mu)\right)_{i=1, 
\ldots ,n}
\end{eqnarray*} 
and 
\begin{eqnarray*}
\frac12\Delta m_n(x)&&\hspace{-.5cm}=\frac1n\left(-\sum_{i=1}^n\int_{\mu 
=d(x)\, dx\in\hat{H}}\frac{n-1}{\|d^{\frac1n}\|_{L^1(D)}^n}\cdot\frac{{\D 
\left(\nabla_d d(x_i),\nabla_d d(x_i)\vphantom{l^1}\right)_{{\Bbb R}^d}}} 
{2n(d(x_i))^2}\prod_{j=1}^n\left(d(x_j)\vphantom{l^1}\right)^{\frac1n}\, 
\bnu(d\mu)\right. \\ 
&&\hspace{.5cm}+\left.\sum_{i=1}^n\int_{\mu=d(x)\, dx\in\hat{H}}\frac{1} 
{\|d^{\frac1n}\|_{L^1(D)}^n}\cdot\frac{\frac12\Delta_d d(x_i)}{d(x_i)} 
\prod_{j=1}^n\left(d(x_j)\vphantom{l^1}\right)^{\frac1n}\, \bnu(d\mu)\right) 
\, ,\quad n\in {\Bbb N}. 
\end{eqnarray*} 
Keeping in mind that the left-hand side of (\ref{5.10}) is increasing in 
$n$ for sufficiently large $n\in {\Bbb N}$, we verify the existence of 
$c_0>0$ and $n_0>0$ such that $\|d^{\frac1n}\|_{L^1(D)}^n\ge c_0$ for $n 
\ge n_0$ and all $d\in\hat{H}$ by (\ref{5.8}). By assumption (\ref{5.8}) 
it holds also that 
\begin{eqnarray}\label{5.11}
\|\nabla m_n\|\le\frac{n^{-\frac12}}{c_0 c^3}\, ,\quad n\ge n_0,
\end{eqnarray} 
where c is the constant in (\ref{5.8}). It follows that 
\begin{eqnarray}\label{5.12}
\lim_{n\to\infty}\int_{D^n}\left(\nabla m_n,\nabla m_n\vphantom{l^1} 
\right)_{{\Bbb R}^{n\cdot d}}\, dx=0 
\end{eqnarray} 
By (\ref{5.9}), (\ref{5.10}), and the weak law of large numbers we have 
$\lim_{n\to\infty}\int_{D^n}{\T\frac12\Delta m_n}\, dx=\int_{\mu=d(x)\, 
dx\in\hat{H}}\int_D{\T\frac12}\Delta_d\log(d)\, dx\, \bnu(d\mu)$ as well 
as 
\begin{eqnarray}\label{5.13}
&&\hspace{-.5cm}\limsup_{n\to\infty}\int_{D^n}\left|{\T\frac12\Delta 
m_n}\right|\, dx\le\int_{\mu=d(x)\, dx\in\hat{H}}\int_D\left|{\T\frac12} 
\Delta_d\log(d)\right|\, dx\, \bnu(d\mu)<\infty\, . 
\end{eqnarray} 
\medskip 

In order to verify (iii) of Theorem \ref{Theorem5.1} under weak additional 
conditions for the example (\ref{5.9}) it is beneficial to apply the following 
proposition. 
\begin{proposition}\label{Proposition5.2}
For every $n\in {\Bbb N}$, let $X^n$ be a stochastic process  satisfying 
the hypotheses of Theorem \ref{Theorem1.11} or, respectively, of Theorem 
\ref{Theorem1.12}. Let $m_n$ be as in Corollary \ref{Corollary1.13}. In 
addition, suppose that for every $n\in {\Bbb N}$ and $Q^{(n)}_{\sbnu_n} 
$-a.e. $W\in\Omega^{(n)}$, $(A^1)^n(W)$ is equi-continuous on ${\Bbb R}$. 
Let $m_n$, $n\in {\Bbb N}$, be defined by (\ref{5.9}). If 
\begin{itemize} 
\item[(iv)] for every $n\in {\Bbb N}$, there is a random variable $b_n$
with 
\begin{eqnarray*}
\left|A^n_{-t}-Y^n_{-t}\vphantom{l^1}\right|\le n^{\frac12}\cdot b_nt\, , 
\quad t\in [0,1], 
\end{eqnarray*} 
and $\int_{D^n}E^{(n)}_x(b_n)\, dx$ is bounded in $n\in {\Bbb N}$ 
\end{itemize} 
then we have (iii) of Theorem \ref{Theorem5.1}. 
\end{proposition} 
Proof. In condition (iii) of Theorem \ref{Theorem5.1} we look at 
\begin{eqnarray*}
E^{(n)}_{\1_{D^n}\lambda_F}\left(\sup_{t\in [0,\delta]}\left|m_n(X^n_{-t} 
-Y^n_{-t})-m_n(x)\vphantom{\dot{f}}\right|\right)\, . 
\end{eqnarray*} 
Let $\delta\in [0,1]$. We have 
\begin{eqnarray}\label{5.14}
&&\hspace{-.5cm}\sup_{t\in [0,\delta]}\left|m_n(X^n_{-t}-Y^n_{-t})-m_n(x) 
\vphantom{\dot{f}}\right|\nonumber \\ 
&&\hspace{.5cm}\le\sup_{t\in [0,\delta]}\left|m_n(W^n_{-t})-m_n(x)\vphantom 
{\dot{f}}\right|+\sup_{t\in [0,\delta]}\left|m_n\left(W^n_{-t}+(A^n_{-t}- 
Y^n_{-t})\vphantom{l^1}\right)-m_n(W^n_{-t})\vphantom{\dot{f}}\right| 
\vphantom{\prod_{i=1}^n}\nonumber \\ 
&&\hspace{.5cm}=\sup_{t\in [0,\delta]}T_1''(t)+\sup_{t\in [0,\delta]}T_2'' 
(t)\vphantom{\prod_{i=1}^n}\, .  
\end{eqnarray} 

We extend $m_n$ outside of $D^n$ with zero and denote this extension to 
${\Bbb R}^{n\cdot d}\equiv F$ also with $m_n$. The random subset $T_n$ of 
$[0,1]$ in which $W_{-t}\in D^n$ is $Q_{\sbnu_n}$-a.e. an open subset of 
$[0,1]$ in the trace topology of ${\Bbb R}$. Thus $T_n\cap (0,1)$ is a 
countable union of open subintervals of $(0,1)$. According to It\^ o's 
formula it holds that 
\begin{eqnarray*}
T_1''(t)&&\hspace{-.5cm}=\left|m_n(W^n_{-t})-m_n(x)\vphantom{\dot{f}}\right 
| \\ 
&&\hspace{-.5cm}=\left|\int_{s\in [0,t]\cap T_n}\left\langle\nabla m_n(W^n_{ 
-s}),dW_{-s}\right\rangle_F+\int_{s\in [0,t]\cap T_n}{\T\frac12}\Delta m_n( 
W^n_{-s})\, ds\right|\, ,\quad t\in [0,\delta]. 
\end{eqnarray*} 
By Doob's maximal inequality 
\begin{eqnarray*}
&&\hspace{-.5cm}E^{(n)}_{\1_{D^n}\lambda_F}\left(\sup_{t\in [0,\delta]} 
\left|T_1''\right|\right)\le\left(E^{(n)}_{\1_{D^n}\lambda_F}\left(\int_{s 
\in [0,\delta]\cap T_n}\left\langle\nabla m_n(W^n_{-s}),dW_{-s}\right 
\rangle_F\right)^2\right)^{1/2} \\ 
&&\hspace{1.0cm}+E^{(n)}_{\1_{D^n}\lambda_F}\int_{s\in [0,\delta]\cap T_n} 
\left|{\T\frac12}\Delta m_n(W^n_{-s})\right|\, ds \\  
&&\hspace{.5cm}\le \left(\int_{s\in [0,\delta]}E^{(n)}_{\1_{D^n}\lambda_F 
}\left\langle\nabla m_n(W^n_{-s}),\nabla m_n(W^n_{-s})\right\rangle_F\, ds 
\right)^{1/2} \\ 
&&\hspace{1.0cm}+\int_{s\in [0,\delta]}E^{(n)}_{\1_{D^n}\lambda_F}\left| 
{\T\frac12}\Delta m_n(W^n_{-s})\right|\, ds 
\end{eqnarray*}
where we treat $\nabla m_n(W^n_{-s})$ as well as ${\T -\frac12}\Delta m_n 
(W^n_{-s})$ as random elements which are zero outside the random set $T_n$. 
This is equivalent to considering $\nabla m_n$ as well as ${\T -\frac12} 
\Delta m_n$ as functions which are zero outside of $D^n$. In this sense we 
obtain 
\begin{eqnarray*} 
&&\hspace{-.5cm}E^{(n)}_{\1_{D^n}\lambda_F}\left(\sup_{t\in [0,\delta]} 
\left|T_1''\right|\right)\le e^{1/2}\left(\int_{s\in [0,\delta]}e^{-s/ 
\delta}E^{(n)}_{\1_{D^n}\lambda_F}\left\langle\nabla m_n(W^n_{-s}),\nabla 
m_n(W^n_{-s})\right\rangle_F\, ds\right)^{1/2} \\ 
&&\hspace{1.0cm}+e\int_{s\in [0,\delta]}e^{-s/\delta}E^{(n)}_{\1_{D^n} 
\lambda_F}\left|{\T\frac12}\Delta m_n(W^n_{-s})\right|\, ds  \\ 
&&\hspace{.5cm}\le e^{1/2}\left(\int_{s\in {\Bbb R}_+}e^{-s/\delta}E^{ 
(n)}_{\1_{D^n}\lambda_F}\left\langle\nabla m_n(W^n_{-s}),\nabla m_n(W^n_{ 
-s})\right\rangle_F\, ds\right)^{1/2} \\ 
&&\hspace{1.0cm}+e\int_{s\in {\Bbb R}_+}e^{-s/\delta}E^{(n)}_{\1_{D^n} 
\lambda_F}\left|{\T\frac12}\Delta m_n(W^n_{-s})\right|\, ds \, . 
\end{eqnarray*}
Denoting by $G_{\beta}^{(n)}$, $\beta>0$, the resolvent of the $n\cdot d 
$-dimensional Brownian motion on ${\Bbb R}^{n\cdot d}=F$ this yields 
\begin{eqnarray*} 
&&\hspace{-.5cm}E^{(n)}_{\1_{D^n}\lambda_F}\left(\sup_{t\in [0,\delta]} 
\left|T_1''\right|\right)\le e^{1/2}\left(\int_{x\in D^n}G_{1/\delta}^{ 
(n)}\left(\left\langle\nabla m_n,\nabla m_n\right\rangle_F\right)(x)\, d 
x\right)^{1/2} \\ 
&&\hspace{1.0cm}+e\int_{x\in D^n}G_{1/\delta}^{(n)}\left|{\T\frac12} 
\Delta m_n(x)\right|\, dx \\  
&&\hspace{.5cm}\le (e\delta)^{1/2}\left(\int_{D^n}\left\langle\nabla m_n, 
\nabla m_n\right\rangle_F\, dx\right)^{1/2}+e\delta\int_{D^n}\left|{\T 
\frac12}\Delta m_n\right|\, dx 
\end{eqnarray*}
where, for the last line, we have used the symmetry of $G_{1/\delta}^{(n) 
}$ in $L^2(F,\lambda_F)$ restricted to the bounded functions of compact 
support. With (\ref{5.12}) and (\ref{5.13}) we get 
\begin{eqnarray}\label{5.15}
\lim_{\delta\downarrow 0}\limsup_{n\to\infty}E^{(n)}_{\1_{D^n}\lambda_F} 
\left(\sup_{t\in [0,\delta]}\left|T_1''\right|\right)=0\, . 
\end{eqnarray}

For $T_2''$ introduced in (\ref{5.14}) we have 
\begin{eqnarray*}
T_2''(t)&&\hspace{-.5cm}=\left|m_n\left(W^n_{-t}+(A^n_{-t}-Y^n_{-t}) 
\vphantom{l^1}\right)-m_n(W^n_{-t})\vphantom{\dot{f}}\right|\le\|\nabla 
m_n\|\cdot\left|A^n_{-t}-Y^n_{-t}\vphantom{l^1}\right|\, . 
\end{eqnarray*} 
Relation (\ref{5.11}) and condition (iv) of this proposition give directly 
\begin{eqnarray}\label{5.16}
\lim_{\delta\downarrow 0}\limsup_{n\to\infty}E^{(n)}_{\1_{D^n}\lambda_F} 
\left(\sup_{t\in [0,\delta]}\left|T_2''\right|\right)=0\, . 
\end{eqnarray} 
The claim follows now from (\ref{5.14})-(\ref{5.16}). 
\qed 
\section{Appendix: Basic Malliavin Calculus for Brownian Motion with 
Random Initial Data} 
\setcounter{equation}{0} 

In this section, we are going to identify the logarithmic derivative 
relative to a Brownian motion $B$ with random initial data. We follow 
the usual concept which means we establish a related Cameron-Martin 
type formula, a gradient operator, and a stochastic integral as its 
dual. 

It turns out that the whole analysis could be separated into the well 
known Wiener space case, which is the case with initial datum zero, and 
an independent part related to the finite dimensional initial variable. 
This is due to the fact that we assume $B$ at time zero to be 
independent of the process at all other times. We treat both parts 
simultaneously and mention that the pure Wiener space part is presented 
here in a very compressed way. For this case it is recommended to 
compare the following with \cite{U95} and the thorough introduction 
to differentiation and integration on the Wiener space in \cite{UZ00}, 
Appendix B. For the case including the random initial condition, we 
refer to \cite{BM-W99} as a source of stochastic calculus relative to 
abstract non-Gaussian measure spaces. We assume the process $B=\left( 
(B_s)_{s\in {\cal T}},({\cal F}_u^v\times \sigma(B_0))_{u,v\in {\cal 
T},\ u\le v},(P_x)_{x\in F}\right)$, ${\cal T}$ either $[0,1]$ or $[-1 
,1]$, to be endowed with an initial random variable $B_0$ independent 
of $(B_s-B_0)_{s\in {\cal T}\setminus \{0\}}$ whose distribution is 
denoted by $\bnu$. Let us assume that $\bnu$ admits a density 
\begin{eqnarray}\label{6.1} 
0<m\in C^1 (F)\quad\mbox{\rm with }\quad\lim_{\|x\|_F\to\infty}m(x)=0\, . 
\end{eqnarray}
Also, let us introduce the shift by $y\in F$,
\begin{eqnarray*} 
U_y(x)=x+y\, , \quad x\in F.    
\end{eqnarray*} 
Furthermore, let us assume that there exist $q\in (1,\infty)$ such that, 
for all $y\in F$, 
\begin{eqnarray}\label{6.2} 
\frac{m\circ U_y}{m}\in L^q(F,\bnu)\quad\mbox{\rm and}\quad\left.\frac 
{d}{d\lambda}\right|_{\lambda =0}\frac{m\circ U_{\lambda y}}{m}=\left 
\langle\frac{\nabla m}{m},y\right\rangle_F\quad \mbox{\rm exists in }L^q 
(F,\bnu)\, . 
\end{eqnarray} 
Let $\Omega:=C({\cal T};F)$ denote the set of trajectories associated with 
$(B_s)_{s\in {\cal T}}$ starting from $s=0$ to either direction in case 
of ${\cal T}=[-1,1]$. Also introduce $P_{\sbnu}:=\int P_x\, \bnu(dx)$ and 
let $E_\nu$ denote the corresponding mathematical expectation. Furthermore, 
we suppose that the filtration $\{{\cal F}_u^v=\sigma(W_\alpha-W_\beta :u 
\le\alpha,\beta\le v)\times\sigma(W_0):u,v\in {\cal T},\ u<v\}$ is completed 
by the $P_{\sbnu}$-completion of the $\sigma$-algebra ${\cal F}$ of all Borel 
subsets of $\Omega$ relative to uniform convergence on ${\cal T}$. 
\bigskip 

{\bf Cameron-Martin type space, embedding, shift relative to $B_0$. } Let 
\begin{eqnarray*} 
H:=\left\{(f,x):f\in L^2({\cal T};F),\ x\in F\right\}\, , 
\end{eqnarray*} 
be equipped with the inner product 
\begin{eqnarray*} 
\left\langle (f,x),(g,y)\right\rangle_H:=\langle f,g\rangle_{L^2}+\langle 
x,y\rangle_F\, . 
\end{eqnarray*} 
Furthermore, define 
\begin{eqnarray*} 
W:=\left\{(\omega ,x):\omega\in C({\cal T};F),\ \omega(0)=0,\ x\in F\right 
\}\, . 
\end{eqnarray*} 
We embed $H$ into $W\equiv C({\cal T};F)$ by 
\begin{eqnarray*} 
j(f,x):=\left(\int_0^\cdot f(s)\, ds,x\right)\equiv x+\int_0^\cdot f(s)\, ds 
\, , \quad (f,x)\in H\, . 
\end{eqnarray*} 

{\bf Cameron-Martin type formula. } Let ${\cal M}^{f,s}({\cal T};F)$ denote 
the set of all $F$-valued finite signed measures on $({\cal T},{\cal B}({\cal 
T}))$ and 
\begin{eqnarray*} 
W^\ast:=\left\{(\omega^\ast,x^\ast):\omega^\ast\in {\cal M}^{f,s}({\cal T};F) 
\ \mbox{\rm with }\omega^\ast({\cal T})=0,\ x^\ast\in F\right\}\, . 
\end{eqnarray*} 
We have  
\begin{eqnarray*} 
_{W^\ast}\hspace{-1mm}\left\langle (\omega^\ast ,x^\ast),(\omega ,x) 
\right\rangle_W&&\hspace{-.5cm}=\int_{\cal T}\langle\omega,d\omega^\ast 
\rangle_F+\langle x,x^\ast\rangle_F \\ 
&&\hspace{-.5cm}=\int_{\cal T}\langle\omega^\ast((s,1]),d\omega_s 
\rangle_F+\langle x,x^\ast\rangle_F 
\, . 
\end{eqnarray*} 
For $(\omega ,x)\in W$ and $(f,y)\in H$ introduce 
\begin{eqnarray*} 
S_{(f,y)}(\omega ,x):=\left(\omega +\int_0^\cdot f(s)\, ds\, ,\, U_y(x)
\right)
\end{eqnarray*} 
and  
\begin{eqnarray*} 
{\cal E}_{(f,y)}(\omega ,x)&&\hspace{-.5cm}:=\exp\left\{\int_{\cal T} 
\langle f(s),d\omega_s\rangle_F-\frac12\int_{\cal T}\langle f(s),f(s) 
\rangle_F\, ds\right\}\frac{m(U_{-y}(x))}{m(x)}\, . 
\end{eqnarray*} 
\begin{theorem}\label{Theorem6.1} 
We have 
\begin{itemize} 
\item[(i)] $\displaystyle\qquad\frac{dP_{\sbnu}\circ S_{(f,y)}
}{dP_{\sbnu}}((\omega ,x))={\cal E}_{(-f,-y)}(\omega ,x)$.
\end{itemize} 
\end{theorem} 
Proof. Let $(f,y)\in H$. First, verify that the map 
\begin{eqnarray*} 
L^1(W,P_{\sbnu})\ni\Phi\lra E\sbnu\left(\exp\left\{i\, _{W^\ast} 
\hspace{-1mm}\left\langle\, \cdot\, ,S_{(f,y)}\right\rangle_W\right\} 
\Phi\vphantom{\sum}\right)
\end{eqnarray*} 
is an injection. Then focus on 
\begin{eqnarray*} 
&&\hspace{-.5cm}E\sbnu\left(\exp\left\{i\, _{W^\ast}\hspace{-1mm} 
\left\langle (\omega^\ast ,x^\ast),S_{(f,y)}(\omega ,x)\right\rangle_W
\right\}\, {\cal E}_{(-f,-y)}(\omega ,x)\vphantom{\sum}\right)\vphantom{ 
\int_0^0} \\ 
&&\hspace{.5cm}=E\sbnu\left(\exp\left\{i\, \left(\int_{\cal T}\left\langle 
\omega^\ast((s,1]),d\left(\omega_s +\int_0^s f(v)\, dv\right)\right 
\rangle_F+\left\langle x^\ast,U_y(x)\right\rangle_F\right)\right.\right. 
 \\ 
&&\hspace{1.5cm}\left.\left.-\int_{\cal T}\langle f(s),d\omega_s\rangle_F 
-\frac12\int_{\cal T}\langle f(s),f(s)\rangle_F\, ds\right\}\, \frac{m( 
U_y(x))}{m(x)}\right) \\ 
&&\hspace{.5cm}=\exp\left\{-\frac12\int_{\cal T}\left\langle\omega^\ast 
((s,1]),\omega^\ast((s,1])\right\rangle_F\, ds\right\}\times \\ 
&&\hspace{1.5cm}\times E_{\sbnu}\left(\exp\left\{i \left\langle x^\ast, 
U_y(x)\right\rangle_F\right\}\frac{m(U_y(x))}{m(x)}\right) \\ 
&&\hspace{.5cm}=\exp\left\{-\frac12\int_{\cal T}\left\langle\omega^\ast 
((s,1]),\omega^\ast((s,1])\right\rangle_F\, ds\right\}\times \\ 
&&\hspace{1.5cm}\times\int\exp\left\{i \left\langle x^\ast,U_y(z)\right 
\rangle_F\right\}m(U_y(z))\, dz \\ 
&&\hspace{.5cm}=\exp\left\{-\frac12\int_{\cal T}\left\langle\omega^\ast 
((s,1]),\omega^\ast((s,1])\right\rangle_F\, ds\right\}\, \int\exp\left\{i 
\left\langle x^\ast,z\right\rangle_F\right\}m(z)\, dz\, . 
\end{eqnarray*} 
The claim is an immediate consequence of the independence of the 
right-hand side of $(f,y)$.\qed 
\begin{corollary}\label{Corollary6.2}
For all bounded measurable $\vp:W\to {\Bbb R}$, we have 
\begin{itemize} 
\item[(ii)]  $\displaystyle\qquad E_{\sbnu}\left(\vp\left((\omega ,x)+j( 
f,U_y(x)-x)\right)\cdot{\cal E}_{(-f,-y)}(\omega ,x)\right)=E_{\sbnu}\vp( 
\omega,x)$ and 
\item[(iii)]  $\displaystyle\qquad E_{\sbnu}\vp\left((\omega ,x)+j(f,U_y( 
x)-x)\right)=E_{\sbnu}\left(\vp(\omega,x)\cdot{\cal E}_{(f,y)}(\omega ,x) 
\right)$. 
\end{itemize} 
\end{corollary} 
\bigskip 

{\bf Gradient operator. } Let $k\in {\Bbb N}$, $l_1,\ldots ,l_k\in W^\ast$, 
$f\in C_b^1({\Bbb R^k})$, and 
\begin{eqnarray}\label{6.3}
\vp(\omega,x):=f\left(l_1(\omega ,x),\ldots ,l_k(\omega ,x)\right)\, . 
\end{eqnarray} 
In other words, let $\vp$ be a bounded differentiable cylindrical function 
on $W$. As usual, define 
\begin{eqnarray}\label{6.4}
{\Bbb D}\vp(\omega,x)&&\hspace{-.5cm}:=\sum_{i=1}^kf_{x_i}\left(l_1(\omega 
,x),\ldots ,l_k(\omega ,x)\right)\cdot _{W^\ast}\hspace{-1mm}\left\langle 
l_i,\, \cdot\, \right\rangle_W \nonumber \\ 
&&\hspace{-.35cm}\equiv\sum_{i=1}^kf_{x_i}\left(l_1(\omega ,x),\ldots ,l_k 
(\omega ,x)\right)\cdot l_i 
\end{eqnarray} 
where $f_{x_i}$ is the partial derivative of $f$ with respect to the $i$-th 
entry, $i\in \{1,\ldots ,k\}$. Let ${\cal C}$ denote the set of all cylindrical 
functions of the form (\ref{6.3}). As a particular case, let us assume the 
linear functionals to be of the form $l_i(\omega ,x):=\left\langle a_i,\omega 
(t_i)\right\rangle_F+\left\langle b_i,x\right\rangle_F$, $a_i,b_i\in F$, $t_1 
\le\ldots\le t_k$, all in ${\cal T}$. Then 
\begin{eqnarray}\label{6.5}  
&&\hspace{-.5cm}{\Bbb D}\vp(\omega,x)(\rho,y) \nonumber \\ 
&&\hspace{.0cm}=\sum_{i=1}^kf_{x_i}\left(\left\langle a_1,\omega(t_1)\right 
\rangle_F+\left\langle b_1,x\right\rangle_F,\ldots ,\left\langle a_k,\omega 
(t_k)\right\rangle_F+\left\langle b_k,x\right\rangle_F\right)\cdot\left(\left 
\langle a_i,\rho(t_i)\right\rangle_F+\left\langle b_i,y\right\rangle_F\right) 
\, , \nonumber \\ 
&&\hspace{11.0cm}(\rho ,y)\in W,
\end{eqnarray} 
which makes it comparable with the common representations for the pure Wiener 
process case as, for example, presented in \cite{N06}. In either case we have 
the usual relation to directional derivatives 
\begin{eqnarray*}
{\Bbb D}\vp(\omega,x)(\rho,y)\equiv\, _{W^\ast}\hspace{-1mm}\left\langle 
{\Bbb D}\vp(\omega,x),(\rho,y)\right\rangle_W =\frac{\partial\vp(\omega,x)} 
{\partial (\rho,y)}\, , \quad (\rho ,y)\in jH=\{jh:h\in H\}. 
\end{eqnarray*} 
\begin{proposition}\label{Proposition6.3} 
Let $q$ be the number used in (\ref{6.2}) and $1/p+1/q\le 1$. 
\medskip 

\nid
(a) The set of all cylindrical functions 
\begin{eqnarray*} 
\vp(\omega,x)=f\left(\left\langle a_1,\omega(t_1)\right\rangle_F 
+\left\langle b_1,x\right\rangle_F,\ldots ,\left\langle a_k,\omega(t_k) 
\right\rangle_F+\left\langle b_k,x\right\rangle_F\right) 
\end{eqnarray*} 
$f\in C_b^1({\Bbb R^k})$, $a_i,b_i\in F$, $i\in \{1,\ldots ,t\}$, $t_1\le 
\ldots\le t_k$, all in ${\cal T}$, is dense in $L^p(W,P_{\sbnu})$. \\ 
(b) Let $\t {\cal C}\subseteq {\cal C}$ be a set of cylindrical functions 
which is dense in $L^p(W,P_{\sbnu})$. Then the operator 
\begin{eqnarray}\label{6.6} 
(D ,\t {\cal C}):=(j^\ast\circ {\Bbb D},\t {\cal C}) 
\end{eqnarray} 
is closable on $L^p(W,P_{\sbnu};H)$. 
\end{proposition} 
Proof. (a) This is explained in \cite{U95} Preliminaries 5 i), for the pure 
Wiener space case. Anything else is trivial. \\ 
(b) Assume $\t {\cal C}\ni\vp_n\stack{n\to\infty}{\lra}0$ in $L^p(W,P_{ 
\sbnu})$ and $D \vp_n$, $n\in {\Bbb N}$, be Cauchy in $L^p(W,P_{\sbnu};H)$. 
Then $D\vp_n\stack{n\to\infty}{\lra}\xi$ in $L^p(W,P_{\sbnu};H)$ for some 
$\xi\in L^p(W,P_{\sbnu};H)$. We have to show $\xi=0$ $P_{\sbnu}$-a.e. For 
$\vp\in\t {\cal C}$ and $(f,y)\in H$, we have by Corollary \ref{Corollary6.2} 
(iii) 
\begin{eqnarray}\label{6.7} 
E_{\sbnu}\left(\left\langle D\vp_n,(f,y)\right\rangle_H\cdot\vp\right) 
&&\hspace{-.5cm}=E_{\sbnu}\left( _{W^\ast}\hspace{-1mm}\left\langle {\Bbb D} 
\vp_n,j(f,y)\right\rangle_W\cdot\vp\right)\vphantom{\int}\nonumber \\ 
&&\hspace{-4cm}=\left.\frac{d}{d\lambda}\right|_{\lambda =0}E_{\sbnu}\left( 
\vp_n\left((\omega ,x)+\lambda j(f,y)\right)\cdot\vp(\omega ,x)\right) 
\nonumber \\ 
&&\hspace{-4cm}=\left.\frac{d}{d\lambda}\right|_{\lambda =0}E_{\sbnu}\left( 
\vphantom{\int_{\cal T}}\vp_n (\omega ,x)\cdot\vp\left((\omega , x)-\lambda 
j(f,y)\right)\times\right.\nonumber \\ 
&&\hspace{-1.5cm}\left.\times \exp\left\{\lambda\int_{\cal T} f(s)\, d 
\omega_s-\frac{\lambda^2}{2}\int_{\cal T}f^2(s)\, ds\right\}\cdot\frac{m 
(U_{-\lambda y}(x))}{m(x)}\right)\, .  
\end{eqnarray}
The derivative of the last line with respect to $\lambda$ at $\lambda =0$ 
exists in $L^q(W,P_{\sbnu})$ by (\ref{6.2}) and is equal to 
\begin{eqnarray}\label{6.8}
\left.\frac{d}{d\lambda}\right|_{\lambda =0}{\cal E}_{\lambda (f,y)}(\omega,x) 
&&\hspace{-.5cm}=\int_{\cal T}f\, d\omega -\left\langle\frac{\nabla m(x)}{m(x)}, 
y\right\rangle_F\, . 
\end{eqnarray} 
Therefore, 
\begin{eqnarray}\label{6.9} 
&&\hspace{-.5cm}E_{\sbnu}\left(\left\langle D\vp_n,(f,y)\right\rangle_H 
\cdot\vp\right)\nonumber \\ 
&&\hspace{.5cm}=E_{\sbnu}\left(\vp_n\cdot _{W^\ast}\hspace{-1mm}\left\langle 
{\Bbb D}\vp,-j(f,y)\right\rangle_W\vphantom{I^I_I}\right)+E_{\sbnu}\left( 
\vp_n\cdot\vp\cdot\left.\frac{d}{d\lambda}\right|_{\lambda =0}{\cal E}_{ 
\lambda(f,y)}\right)\, . 
\end{eqnarray} 
Letting $n\to\infty$ we obtain $E_{\sbnu}\left(\left\langle\xi,(f,y)\right 
\rangle_H\cdot\vp\right)=0$ for all $\vp\in\t {\cal C}$ and $(f,y)\in H$ 
from (\ref{6.9}). Thus $\xi=0$ $P_{\sbnu}$-a.e. In other words, $D$ is closable 
on $L^p(W,P_{\sbnu};H)$. \qed 
\begin{definition}\label{Definition6.4} 
{\rm Let $q$ be the number used in (\ref{6.2}) and $1/p+1/q\le 1$. We say 
$\vp\in {\rm Dom}_p(D)$ if there is a sequence $\vp_n\in {\cal C}$, $n\in 
{\Bbb N}$, such that ${\cal C}\ni\vp_n\stack{n\to\infty}{\lra}\vp$ in $L^p 
(W,P_{\sbnu})$ and $D\vp_n$, $n\in {\Bbb N}$, is Cauchy in $L^p(W,P_{\sbnu} 
;H)$. In this case 
\begin{eqnarray*} 
D\vp:=\lim_{n\to\infty}D(\vp_n)\, . 
\end{eqnarray*} 
Let $D_{p,1}\equiv D_{p,1}(P_{\sbnu})$ be the space of all $\vp\in {\rm 
Dom}_p(D)$ equipped with the norm $\|\vp\|_{p,1}:=\|\vp\|_{L^p(W,P_{\sbnu})} 
+\|D\vp\|_{L^p(W,P_{\sbnu};H)}$. 
}
\end{definition} 
\bigskip 

{\bf Gradient and directional derivative } 
\begin{proposition}\label{Proposition6.5} 
(a) Let $q$ be the number used in (\ref{6.2}) and $1/p+1/q\le 1$. 
Furthermore, let $\psi\in D_{p,1}$ and $k\in H$. Suppose the existence of 
the limit 
\begin{itemize}
\item[(i)] $\qquad\D\frac{\partial\psi}{\partial jk}=\lim_{\lambda\to 0} 
\frac1\lambda\left(\psi (\cdot+\lambda jk)-\psi\vphantom{\dot{f}}\right)$ 
\end{itemize} 
in $L^p(W,P_{\sbnu})$. Then 
\begin{eqnarray}\label{6.10}
\frac{\partial\psi}{\partial jk}=\left\langle D\psi ,k\right\rangle_H\in 
L^p(W,P_{\sbnu})\, . 
\end{eqnarray} 
(b) Let $1<p<\infty$ and $\psi\in L^p(W,P_{\sbnu})$. Assume the existence of 
the limit (i) for all $k\in H$ in the measure $P_{\sbnu}$. Also suppose the 
existence of $C\in L^p(W,P_{\sbnu};H)$ such that 
\begin{itemize}
\item[(ii)] $\qquad\D\frac{\partial\psi}{\partial jk}= \langle C,k\rangle_H 
\quad$ for all $k\in H$ and  
\item[(iii)] $\qquad\D\psi((\omega,x) +jk)-\psi((\omega,x))=\int_0^1\langle C 
((\omega,x)+tk),k\rangle_H\, dt\quad$ for all $k\in H$ and $P_{\sbnu}$-a.e. 
$(\omega,x)\in W$. 
\end{itemize} 
Then $\psi\in D_{p,1}$, $C=D\psi$, and we have 
\begin{eqnarray}\label{6.11}
\frac{\partial\psi}{\partial jk}=\left\langle D\psi ,k\right\rangle_H\in 
L^p(W,P_{\sbnu})\, . 
\end{eqnarray} 
\end{proposition} 
Proof. (a) Let $\vp\in {\cal C}$ and let $\psi_n\in {\cal C}$, $n\in {\Bbb 
N}$, be a sequence with $\psi_n\stack{n\to\infty}{\lra}\psi$ in $D_{p,1}$. 
Following (\ref{6.7})-(\ref{6.9}) back and forth we obtain for $k=(f,y)\in 
H$
\begin{eqnarray*} 
E_{\sbnu}\left(\left\langle D\psi,k\right\rangle_H\cdot\vp\right)
&&\hspace{-.5cm}=\lim_{n\to\infty}E_{\sbnu}\left(\left\langle D\psi_n,k 
\right\rangle_H\cdot\vp\right)\vphantom{\left.\frac{d}{d\lambda}\right|_{ 
\lambda =0}} \nonumber \\ 
&&\hspace{-.5cm}=\lim_{n\to\infty}\left(-E_{\sbnu}\left(\psi_n\cdot\left 
\langle D\vp,k\right\rangle_H\vphantom{I^I_I}\right)+E_{\sbnu}\left( 
\psi_n\cdot\vp\cdot\left.\frac{d}{d\lambda}\right|_{\lambda =0}{\cal E}_{ 
\lambda(f,y)}\right)\right) \nonumber \\ 
&&\hspace{-.5cm}=-E_{\sbnu}\left(\psi\cdot\left\langle D\vp,k\right 
\rangle_H\vphantom{I^I_I}\right)+E_{\sbnu}\left(\psi\cdot\vp\cdot\left. 
\frac{d}{d\lambda}\right|_{\lambda =0}{\cal E}_{\lambda(f,y)}\right) 
\nonumber \\ 
&&\hspace{-.5cm}=\left.\frac{d}{d\lambda}\right|_{\lambda =0}E_{\sbnu} 
\left(\psi\left(\cdot+\lambda jk\right)\cdot\vp\right) \nonumber \\ 
&&\hspace{-.5cm}=E_{\sbnu}\left(\left.\frac{d}{d\lambda}\right|_{\lambda 
=0}\psi\left(\cdot+\lambda jk\right)\cdot\vp\right)\, ,  
\end{eqnarray*} 
the last line by hypothesis (i). From here, we conclude 
(\ref{6.10}). 
\medskip 

\nid 
(b) Let us first assume that the condition (i) in the measure and the 
conditions (ii), (iii) are satisfied only for $C=(C_1,0)$ with $C_1\in 
L^p\left(W,P_{\delta_0};L^2({\cal T};F)\right)$ and $k=(f,0)$ with $f\in 
L^2({\cal T};F)$. For the conclusion (\ref{6.11}) for all such $k=(f,0)$, 
we refer to \cite{UZ00}, Appendix B.6. Now part (b) of the proposition 
follows from the fact that (i) in the measure and (ii), (iii) together 
with the definition of $\nu(dx)=m(x)\, dx$ in (\ref{6.1}) implies 
(\ref{6.11}) for $k=(0,x)$ with $x\in F$ and $C=(0,C_2)$ where $C_2\in 
L^p(W,\bnu;F)$. For this one may also consult \cite{KO84}, Theorem 1.11. 
\qed
\bigskip

{\bf The stochastic integral } We are going to define an anticipating 
integral which in its specification to the adapted case on the Wiener 
space is the It\^o integral. The following definition takes into 
consideration the Definition \ref{Definition6.4} in which the spaces  
$D_{p,1}$ are meaningfully introduced only for $1/p\le 1-1/q$. 
\begin{definition}\label{Definition6.6} 
{\rm Let $q$ be the number appearing in (\ref{6.2}) and $1/p+1/q=1$. Let 
$q'\ge p$ and $1/{p'}+1/{q'}=1$. Furthermore, let $\xi\in L^{p'}\left( 
W,P_{\sbnu};H\right)$. We say that $\xi\in {\rm Dom}_{p'}(\delta)$ if 
there exists $c_{p'}(\xi)>0$ such that 
\begin{eqnarray*} 
E_{\sbnu}\left(\langle D\vp,\xi\rangle_H\right)\le c_{p'}(\xi)\cdot \|\vp 
\|_{L^{q'}(W,P_{\sbnu})}\, , \quad\vp\in D_{q',1}\, . 
\end{eqnarray*} 
In this case, we define the {\it stochastic integral} $\delta(\xi)$ by 
\begin{eqnarray*} 
E_{\sbnu}\left(\delta(\xi)\cdot\vp\right)=E_{\sbnu}\left(\langle\xi,D 
\vp\rangle_H\right)\, , \quad\vp\in D_{q',1}\, . 
\end{eqnarray*} 
}
\end{definition} 

Replacing $\nu$ by $\delta_0$ and concentrating on the case ${\cal T}= 
[0,1]$, $\delta(\xi)$ is the {\it Skorohod integral}. 
\begin{theorem}\label{Theorem6.7} 
Let $q$ appearing in (\ref{6.2}) and $1/p+1/q=1$. 
\medskip 

\nid
(a) Let $(f,y)\in H$. We have $(f,y)\in {\rm Dom}_q(\delta)$ and 
\begin{eqnarray*} 
\delta(f,y)=\int_{\cal T}\langle f,d\omega\rangle_F-\left\langle y, 
{\displaystyle\frac{\nabla\, m(x)}{m(x)}}\right\rangle_F\, . 
\end{eqnarray*} 
(b) Let $[s_1,s_2)\subseteq {\cal T}$, $i,j\in \{1,\ldots,n\cdot d\}$, 
and $a\in D_{q',1}$ for some $q'>p$. Then $a(\omega,x)\cdot (e_j\1_{[ 
s_1,s_2)},e_i)\in {\rm Dom}_w(\delta)$ with $1/w=1/{q'}+1/q$ and 
\begin{eqnarray}\label{6.12}
&&\hspace{-.5cm}\delta \left(a\cdot (e_j\1_{[s_1,s_2)},e_i)\vphantom{ 
I^I_I}\right)(\omega,x) \nonumber \\ 
&&\hspace{.5cm}=a(\omega ,x)\cdot(\omega_{s_2}-\omega_{s_1})_j-a(\omega 
,x)\cdot\frac{(\nabla m(x))_i}{m(x)}-\left\langle Da,(e_j\1_{[s_1,s_2)} 
,e_i)\right\rangle_H\, . 
\end{eqnarray} 
Assume the existence of the directional derivatives 
\begin{eqnarray*} 
\left.\frac{d}{d\lambda}\right|_{\lambda =0}a\left(\omega +\lambda{ 
\int_0^\cdot}e_j\1_{[s_1,s_2)}(s)\, ds,x\right)\quad\mbox{\rm and} 
\quad\left.\frac{d}{d\lambda}\right|_{\lambda =0}a(\omega,x+\lambda 
\cdot e_i)\, . 
\end{eqnarray*} 
in $L^{q'}(W,P_{\sbnu})$. Then 
\begin{eqnarray}\label{6.13}
\delta \left(a\cdot (e_j\1_{[s_1,s_2)},e_i)\vphantom{I^I_I}\right) 
(\omega,x)&&\hspace{-.5cm}=a(\omega ,x)\cdot(\omega_{s_2}-\omega_{s_1} 
)_j-a(\omega,x)\cdot\frac{(\nabla m(x))_i}{m(x)}\nonumber \\ 
&&\hspace{-3.5cm}-\left.\frac{d}{d\lambda}\right|_{\lambda =0}a\left( 
\omega +\lambda{\int_0^\cdot}e_j\1_{[s_1,s_2)}(s)\, ds,x\right)-\left. 
\frac{d}{d\lambda}\right|_{\lambda =0}a(\omega,x+\lambda\cdot e_i)\, . 
\end{eqnarray} 
\end{theorem}
Proof. (a) Looking at (\ref{6.9}), choosing $\vp=1$ there, and renaming 
$\vp_n$ by $\psi \in {\cal C}$ we get 
\begin{eqnarray*} 
E_{\sbnu}\left(\left\langle D\psi,(f,y)\right\rangle_H\right) 
&&\hspace{-.5cm}=E_{\sbnu}\left(\psi\cdot\left.\frac{d}{d\lambda}\right|_{ 
\lambda =0}{\cal E}_{\lambda(f,y)}\right) \\ 
&&\hspace{-.5cm}\le\|\psi\|_{L^p(W,P_{\sbnu})}\cdot\left\|\left.\frac{d}{d 
\lambda}\right|_{\lambda =0}{\cal E}_{\lambda(f,y)}\right\|_{L^q(W,P_{\sbnu})}
\, . 
\end{eqnarray*} 
Part (a) follows now by taking now into consideration $q\ge 2$ and (\ref{6.8}), 
\begin{eqnarray*}
\left.\frac{d}{d\lambda}\right|_{\lambda =0}{\cal E}_{\lambda (f,y)}(\omega,x) 
=\int_{\cal T}f\, d\omega -\left\langle\frac{\nabla m(x)}{m(x)},y\right\rangle_F 
\, . 
\end{eqnarray*} 
(b) Relation (\ref{6.12}) is a consequence of part (a) and the following. For 
$\xi\in H$, $\alpha\in D_{q,1}$ we have $\alpha\delta (\xi)-\langle D\alpha,\xi 
\rangle_H\in L^w$ where $1/w=1/q+1/{q'}$ and 
\begin{eqnarray*}
\alpha\xi\in{\rm Dom}_w(\delta)\quad\mbox{\rm and}\quad\delta(\alpha\xi)= 
\alpha\delta (\xi)-\langle D\alpha,\xi\rangle_H\, . 
\end{eqnarray*} 
Indeed the latter becomes evident for $\vp\in {\cal C}$ by adapting a calculation 
from \cite{U95}, Section 1.1, 
\begin{eqnarray*}
E_\nu\left(\langle a\xi,D\vp\rangle_H\vphantom{\dot{f}}\right)&&\hspace{-.5cm}= 
E_\nu\left(\langle\xi,aD\vp\rangle_H\vphantom{\dot{f}}\right) \\ 
&&\hspace{-.5cm}=E_\nu\left(\langle\xi,D(a\vp)-\vp Da\rangle_H\vphantom{\dot{f}} 
\right) \\ 
&&\hspace{-.5cm}=E_\nu\left(a\delta(\xi)\cdot\vp\vphantom{\dot{f}}\right) 
-E_\nu\left(\langle Da,\xi\rangle_H\cdot\vp\vphantom{\dot{f}}\right)\, .
\end{eqnarray*} 
Relation (\ref{6.13}) follows now from Proposition \ref{Proposition6.5} (a). 
\qed
\medskip 

Let us continue with an object being, up to its signum, nothing but the 
stochastic integral. 
\begin{definition}\label{Definition6.8} 
{\rm Let $h\in {\rm Dom}_p(\delta)$ for some $p\in (1,\infty)$. We define the 
{\it logarithmic derivative} $\beta_{jh}$ of the measure $P_{\sbnu}$ {\it in 
direction of} $jh$ by $\beta_{jh}=-\delta (h)$. }
\end{definition} 
\medskip 

{\bf The $L^2$-integral } There is another type of anticipating integral 
whose specification to the adapted case on the Wiener space is the 
Stratonovich integral. 
\begin{definition}\label{Definition6.9} 
{\rm (a) Let $L_{2,1}$ denote the space of all $(f,y)\in L^2(W,P_{\sbnu};H)$ 
such that 
\begin{itemize}
\item[(i)] $f_t,y\in D_{2,1}$ for a.e. $t\in {\cal T}$, and $P_{\sbnu}$-a.e. 
it holds that 
\item[(ii)] $\int_{t\in {\cal T}}\|Df_t\|_H^2\, dt+\|Dy\|_H^2<\infty$.  
\end{itemize} 
(b) Let $\vp_i$, $i\in {\Bbb N}$, be a complete orthonormal basis of $L^2( 
{\cal T};F)$ and $(f,y)\in L_{2,1}$. The {\it $L^2$-integral} 
\begin{eqnarray*} 
\int_{\cal T}\langle (f,y)\circ d(\omega,x)\rangle_F\equiv\int_{\cal T} 
\langle f\circ d\omega\rangle_F-\left\langle y,\frac{\nabla\, m(x)}{m(x)} 
\right\rangle_F 
\end{eqnarray*} 
is defined by 
\begin{eqnarray*} 
\sum_{i=1}^\infty\langle f,\vp_i\rangle_{L^2({\cal T};F)}\cdot\int_{\cal 
T}\langle\vp_i,d\omega\rangle_F-\left\langle y,\frac{\nabla\, m(x)}{m(x)} 
\right\rangle_F 
\end{eqnarray*} 
provided that this sum converges in the measure $P_{\sbnu}$ independent of 
the chosen orthonormal basis $\vp_i$, $i\in {\Bbb N}$, of $L^2({\cal T};F)$. 
}
\end{definition} 

Replacing $\nu$ by $\delta_0$ and choosing ${\cal T}=[0,1]$, $\int_{\cal T} 
\langle f\circ d\omega\rangle_F$ is the {\it Ogawa integral}. 
\begin{theorem}\label{Theorem6.10} 
We have the fundamental relation between the stochastic integral and the 
$L^2$-integral. Let $(f,y)\in L_{2,1}$ and suppose that $P_{\sbnu}$-a.e. 
$D(f,y)(\cdot):H\to H$ is a nuclear operator. Then 
\begin{eqnarray*} 
\int_{\cal T} (f,y)\circ d(\omega,x)=\delta (f,y)+{\rm trace}(D(f,y))\, .
\end{eqnarray*} 
\end{theorem} 
Proof. For $\nu$ replaced by $\delta_0$ and ${\cal T}=[0,1]$, this is due to 
\cite{N06}, Section 3.1.2. Anything else is left as an exercise. \qed
\medskip

For expositions beyond the scope of this paper, we refer to \cite{B98}, 
\cite{BM-W99}, \cite{Bu94}, \cite{KP00}, \cite{N06}, \cite{NZ89}, \cite{Ro89}, 
\cite{S00}, \cite{SW99}.

\end{document}